\documentclass[11pt]{book}
\usepackage[latin1]{inputenc}
\usepackage[english]{babel}
\usepackage{amsmath}
\usepackage[ntheorem]{empheq}
\usepackage{amssymb}
\usepackage{graphics}
\usepackage{amsthm}
\usepackage{mathrsfs}
\usepackage{mathdots}
\usepackage{comment}
\usepackage{longtable}
\usepackage[tableposition=below]{caption}
\usepackage[margin=3cm]{geometry}
\usepackage{makeidx}
\makeindex
\numberwithin{equation}{section}

\usepackage{marginnote}
\usepackage[pagebackref,colorlinks,citecolor=red,urlcolor=blue,
bookmarks=false,hypertexnames=true]{hyperref}

\renewcommand{\thefootnote}{}

\usepackage{authblk}

\begin{document}

\title{Conjugacy in finite classical groups}
\author[]{Giovanni De Franceschi} 
\author[]{Martin W. Liebeck} 
\author[]{E.A.\ O'Brien}
\affil[]{\vspace*{1cm} De Franceschi and O'Brien, Department of Mathematics, University of Auckland,
Auckland, New Zealand} 
\affil[]{\vspace*{0.05cm} Liebeck, Imperial College, London SW7 2AZ, UK} 

\newcommand{\Ft}{\mathbb{F}_q[t]_{\lambda}}
\newcommand{\Ftunit}{\mathbb{F}_q[t]_{\lambda}^*}
\newcommand{\rmN}{\mathrm N}
\newcommand{\w}{\widetilde}

\def\P{{\cal P}}
\newcommand{\GL}{\mathrm{GL}}
\newcommand{\SL}{\mathrm{SL}}
\newcommand{\U}{\mathrm{U}}
\newcommand{\Sp}{\mathrm{Sp}}
\newcommand{\Or}{\mathrm{O}}
\newcommand{\SO}{\mathrm{SO}}
\newcommand{\SU}{\mathrm{SU}}
\newcommand{\N}{\mathrm{N}}
\newcommand{\C}{\mathscr{C}}
\newcommand{\im}{\mathrm{im}}
\newcommand{\Spec}{\mathscr{S}}
\newcommand{\tr}{\mathrm{t}}
\newcommand{\rk}{\mathrm{rk}\,}
\newcommand{\ged}{generalized elementary divisor }
\newcommand{\geds}{generalized elementary divisors }
\newcommand\lreqn[2]{\noindent\makebox[\textwidth]{$\displaystyle#1$\hfill(#2)}\vspace{2ex}}

\newtheorem{thm}{Theorem}[section]
\newtheorem{theorem}[thm]{Theorem}
\newtheorem{lemma}[thm]{Lemma}
\newtheorem{lem}[thm]{Lemma}
\newtheorem{prop}[thm]{Proposition}
\newtheorem{cor}[thm]{Corollary}
 \newtheorem{conj}[thm]{Conjecture}
\newtheorem{definition}[thm]{Definition}
\newtheorem{proposition}[thm]{Proposition}
\newtheorem{corollary}[thm]{Corollary}
\theoremstyle{definition}
\newtheorem{remark}[thm]{Remark}
\newtheorem{example}[thm]{Example}
\newtheorem{theor}{Theorem}[chapter]


\def\myX{Z}
\def\myx{z}

\def\onemu{u}
\def\onenu{v}
\def\mymu{i}
\def\mynu{j}
\def\myr{s}
\def\mys{m}
\def\myV{V}
\def\myW{U}
\def\m_gamma{m} 
\def\myexpR{m} 
\def\myN{m} 
\def\myz{D} 
\def\myP{P} 
\def\matalpha{\mathtt{a}}
\def\matbeta{\mathtt{b}}
\def\matgamma{\mathtt{c}}
\def\matH{\mathtt{h}}

\def\a{\alpha}
\def\om{\omega}
 \def\b{\beta}
 \def\e{\epsilon}
 \def\d{\delta}
  \def\D{\Delta}
 \def\c{\chi}
 \def\k{\kappa}
 \def\g{\gamma}
 \def\t{\tau}
\def\ti{\tilde}
 \def\N{\mathbb N}
 \def\Q{\mathbb Q}
 \def\Z{\mathbb Z}
 \def\F{\mathbb F}
 \def\ovF{\overline\F}
 \def\bfN{\mathbf N}
 \def\cG{\mathcal G}
 \def\cT{\mathcal T}
 \def\cX{\mathcal X}
 \def\cY{\mathcal Y}
 \def\C{\mathcal C}

 \def\cD{\mathcal D}
 \def\cZ{\mathcal Z}
 \def\cO{\mathcal O}
 \def\cW{\mathcal W}
 \def\cL{\mathcal L}
 \def\bfC{\mathbf C}
\def\NN{{\rm N}}
 \def\bfZ{\mathbf Z}
 \def\bfO{\mathbf O}
 \def\G{\Gamma}
 \def\bO{\boldsymbol{\Omega}}
 \def\bgo{\boldsymbol{\omega}}
 \def\go{\rightarrow}
 \def\do{\downarrow}
 \def\ra{\rangle}
 \def\la{\langle}
 \def\fix{{\rm fix}}
 \def\ind{{\rm ind}}
 \def\rfix{{\rm rfix}}
 \def\diam{{\rm diam}}
 \def\uni{{\rm uni}}
 \def\diag{{\rm diag}}
 \def\Irr{{\rm Irr}}
 \def\Syl{{\rm Syl}}
\def\det{{\rm det}}
 \def\Gal{{\rm Gal}}
 \def\Tr{{\rm Tr}}
 \def\M{{\cal M}}
 \def\cE{{\mathcal E}}
\def\td{\tilde\delta}
\def\tx{\tilde\xi}
\def\DC{D^\circ}

\def\Ker{{\rm Ker}}
 \def\rank{{\rm rank}}
 \def\soc{{\rm soc}}
 \def\Cl{{\rm Cl}}
 \def\A{{\sf A}}
 \def\sP{{\sf P}}
 \def\sQ{{\sf Q}}
 \def\SSS{{\sf S}}
  \def\SQ{{\SSS^2}}
 \def\St{{\sf {St}}}
 \def\p{\ell}
 \def\ps{\ell^*}
 \def\SC{{\rm sc}}
 \def\supp{{\sf{supp}}}
  \def\cR{{\mathcal R}}
 \newcommand{\tw}[1]{{}^#1}

 \def\Sym{{\rm Sym}}
 \def\PSL{{\rm PSL}}
 \def\SL{{\rm SL}}
 \def\Sp{{\rm Sp}}
 \def\GL{{\rm GL}}
 \def\SU{{\rm SU}}
 \def\GU{{\rm GU}}
 \def\SO{{\rm SO}}
 \def\PO{{\rm P}\Omega}
 \def\Spin{{\rm Spin}}
 \def\PSp{{\rm PSp}}
 \def\PSU{{\rm PSU}}
 \def\PGL{{\rm PGL}}
 \def\PGU{{\rm PGU}}
 \def\Iso{{\rm Iso}}
 \def\Stab{{\rm Stab}}
 \def\GO{{\rm GO}}
 \def\Ext{{\rm Ext}}
 \def\E{{\cal E}}
 \def\l{\lambda}
 \def\ve{\varepsilon}
 \def\Lie{\rm Lie}
 \def\s{\sigma}
 \def\O{\Omega}
 \def\o{\omega}
 \def\ot{\otimes}
 \def\op{\oplus}
 \def\oc{\overline{\chi}}
 \def\pf{\noindent {\bf Proof.$\;$ }}
 \def\Proof{{\it Proof. }$\;\;$}
 \def\no{\noindent}
\def\hal{\unskip\nobreak\hfil\penalty50\hskip10pt\hbox{}\nobreak
 \hfill\vrule height 5pt width 6pt depth 1pt\par\vskip 2mm}
\setcounter{MaxMatrixCols}{20}

 \renewcommand{\thefootnote}{}

\maketitle

\cleardoublepage
\phantomsection

\chapter*{Abstract}
\addcontentsline{toc}{chapter}{Abstract}
Let $G$ be a classical group defined over a finite field. 
We consider the following fundamental problems concerning conjugacy in $G$:
\begin{itemize}
\item List a representative for each conjugacy class of $G$.
\item Given $x \in G$, describe the centralizer of $x$ in $G$, 
by giving its group structure and a generating set.
\item Given $x,y \in G$, establish whether $x$ and $y$ are conjugate in $G$ 
and, if so, then find explicit $z \in G$ such that $z^{-1}xz = y$.
\end{itemize}

We present comprehensive theoretical solutions to all three problems, 
and use our solutions to formulate practical algorithms. 
 In parallel to our theoretical work, 
we have developed in \textsc{Magma} complete implementations of our algorithms.
They form a critical component of various general algorithms in computational group theory -- for example,
computing character tables and solving conjugacy problems in arbitrary finite groups.

\cleardoublepage
\phantomsection
\addcontentsline{toc}{chapter}{Acknowledgments}

\chapter*{Acknowledgments}
Liebeck and O'Brien were supported by the 
Marsden Fund of New Zealand and by the 2020 and 2022 programmes ``Groups, representations and applications: new perspectives" at the Isaac Newton Institute, Cambridge. 

Some of this work formed part of the PhD thesis of De 
Franceschi \cite{thesis}, supported by a University of Auckland PhD scholarship. 

We thank Professor Donald Taylor of the University of Sydney for 
many helpful discussions about the project and related {\sc Magma} code.

\tableofcontents

\chapter{Introduction and background} \label{introchap}
This book is a study of conjugacy in finite classical groups: the general and special linear, symplectic, orthogonal and unitary groups over finite fields. 
In this chapter we introduce the three main conjugacy problems studied. 
In Section \ref{mainprob}
we describe the problems and some of their history. 
In Section \ref{Sect12} we present both our strategy 
to solve them and a summary of our main results. 
Finally, in Sections \ref{Sect13}--\ref{Sect14} we discuss our related algorithms and their 
performance, and the role they play in the solution of 
these problems for arbitrary finite groups.

\section{The main problems} \label{mainprob}
Let $G$ be a classical group defined over a finite field. Our aim is to 
obtain a precise and explicit solution to each of the following closely 
related problems:
\begin{itemize}
\item[(1)] List a representative for each conjugacy class of $G$.
\item[(2)] Given $x \in G$, describe the centralizer $C_G(x)$ of $x$ in $G$, 
by giving its group structure and a generating set.
\item[(3)] Given $x,y \in G$, establish whether $x$ and $y$ are conjugate in $G$ 
and, if so, then find explicit $z \in G$ such that $z^{-1}xz = y$.
\end{itemize}

\noindent We also aim to provide algorithms that implement our solution of each of the problems.

We achieve all of these aims in this work. First we present comprehensive theoretical solutions to all three problems,
and then use our theoretical solutions to formulate practical algorithms 
to solve them. 
We have developed implementations of our algorithms
in \textsc{Magma} \cite{MAGMA}. 
These algorithms form a critical component of various general algorithms in computational group theory: as one example, algorithms for all three 
problems are vital in computing the character table of a classical group; as another, the algorithms are fundamental to the computational solution of conjugacy problems for arbitrary finite groups, as outlined in Section \ref{Sect14}.  

Various parts of the conjugacy problems (1)-(3) for finite classical groups have been much studied in the literature \cite{JBRIT, Carterclass, FS, GLOB, Huppert, LS, Milnor, Murray, GEW, SCCSS}. 
One of the most significant such works is the 1963 paper  of Wall \cite{GEW}. 
In that paper, for a given classical group $G$ on a finite vector space $V$, the author gave criteria for an arbitrary element of the general linear group $\GL(V)$ to lie in $G$; he gave necessary and sufficient conditions for two elements of $G$ to be conjugate; and he determined the total number of conjugacy classes in $G$ and their sizes. The theory developed in \cite{GEW} was set in the much broader context of classical groups over division rings, and does not lend itself well to our explicit and algorithmic requirements. The same comment applies to the approach of Carter \cite{Carterclass}, where the semisimple classes and centralizers in finite  classical groups are determined using the theory of algebraic groups, and also to the work in \cite{LS}, where the same is done for the unipotent classes. In summary, 
while various treatments in the literature give partial solutions to our problems, none is constructive:  
for example, the structure of the centralizer of an arbitrary element of a finite classical group is known, but nowhere is a method given to construct its generating set. 

Nevertheless, we found various parts of the literature very useful for our approach to these problems. In particular, our work is influenced by the article of Milnor \cite{Milnor}, as interpreted and developed by Britnell \cite{JBRIT}, and a short paper of Wall \cite{SCCSS}; also we make full use of \cite{GLOB, LS} for representatives and centralizers of unipotent classes.

\section{Background and summary of results}   \label{Sect12}

We first introduce some notation that we use throughout. 
If $A,A_1,\ldots,A_r$ are square matrices, then 
$$\bigoplus_{i=1}^rA_i = A_1\oplus \cdots \oplus A_r$$ denotes the block 
diagonal matrix with diagonal blocks $A_i$, and $A^{[s]}$ is the block 
diagonal matrix $A\oplus \cdots \oplus A$, where there are $s$ blocks.

We begin in Chapter \ref{glsl} by describing the 
conjugacy classes and centralizers in $\GL_n(q)$. Although the theory for 
this case is well known -- see for example \cite[\S 3.2]{PGPG}, 
\cite[Chap.\ 3]{PFGG}, \cite[\S 4.2]{MD} and \cite{Murray} -- we 
make a new contribution: an explicit description of the generators of 
the centralizer of an 
element. Moreover, the notation and concepts introduced in this 
chapter underpin our approach to the general problems, so we present the 
solution in detail. 

Here is a brief summary. Let $q=p^a$, where $p$ is prime and $a$ is a positive integer, and let $F=\mathbb{F}_q$. Let $V=F^n$ 
be an $n$-dimensional vector space over $F$, and set $G=\GL(V) \cong \GL_n(F)$. For $x \in G$ with minimal polynomial $\prod_1^h f_i(t)^{e_i}$, where $f_1,\ldots,f_h \in F[t]$ are monic irreducible polynomials, there is a decomposition
\begin{eqnarray}     \label{Vdecompositionintro}
V = V_1 \oplus V_2 \oplus \cdots \oplus V_h,
\end{eqnarray}
where $V_i = \hbox{ker}(f_i(x)^{e_i})$. Moreover for each $i$ we can write $V_i = V_{i,1} \oplus \cdots \oplus V_{i,k_i}$, a sum of $x$-invariant subspaces such that $x$ acts cyclically on each $V_{i,j}$ with minimal polynomial $f_i(t)^{e_{ij}}$, where 
$1 \leq e_{i1} \leq \cdots \leq e_{ik_i} = e_i$. The polynomials $f_i^{e_{ij}}$ are the \emph{elementary divisors} of $x$. It is well known that elements of
$G$ are conjugate if and only if they have the same elementary divisors.

The choice of representative for each class is important for the solution of the conjugacy problems (2) and (3) from the previous section. If $x$ and $y$ are conjugate in $G$, then the easiest way to compute a conjugating element is to find matrices $z_x$ and $z_y$ such that $z_x^{-1}xz_x = J = z_y^{-1}yz_y$, where $J$ is the representative for the conjugacy class of $x$ and $y$. The desired conjugating element $z$ is then $z_xz_y^{-1}$. 
Moreover, if $z^{-1}xz=y$, then $C_G(y)=z^{-1}C_G(x)z$; thus, to describe centralizers in $G$ it is sufficient to describe the centralizers of the representatives of conjugacy classes, so it is convenient to choose a suitable canonical form. We use the 
\emph{Jordan form} \index{Jordan form} of $x$: it is a block diagonal matrix $J = \bigoplus_{i,j} J_{f_i,e_{ij}}$, where $J_{f_i,e_{ij}}$ is the matrix of the restriction of $x$ to $V_{i,j}$ and
$$
J_{f_i,e_{ij}} = \begin{pmatrix}
C_i & \mathbb{I} && \\ & \ddots & \ddots & \\ && \ddots & \mathbb{I} \\ &&& C_i
\end{pmatrix},
$$
where $C_i$, the companion matrix of $f_i$, appears $e_{ij}$ times. The matrix $J_{f_i,e_{ij}}$ is the \textit{Jordan block}\index{Jordan block}
 of order $e_{ij}$ relative to $f_i$. This is a generalization of the canonical Jordan form (defined in \cite[\S 6.6]{PFGG} only for matrices whose 
characteristic polynomial splits into linear factors), and is used also by Macdonald \cite[\S 4.2]{MD} and Murray \cite{Murray}.

Recall that $x \in G \cong \GL_n(q)$ is \textit{semisimple}\index{semisimple} 
if its order is coprime to 
$p$, and $x$ is \textit{unipotent}\index{unipotent} 
if its order is a power of $p$. 
If $x$ is semisimple, then $e_{ij}=1$ for all $i,j$; 
if $x$ is unipotent, then $h=1$ and $f_1(t)=t-1$. The 
\textit{Jordan decomposition}\index{Jordan decomposition}
states that every $x \in G$ can be written uniquely as 
\begin{equation}\label{jor}
x=su=us,
\end{equation}
 where $s$ is semisimple and $u$ is unipotent. In particular, $u \in C_G(s)$. Two important consequences are the following:
\begin{itemize}
\item Let $x,y \in G$ and let $x=s_xu_x$ and $y=s_yu_y$ be their Jordan decompositions. Then $x$ and $y$ are conjugate in $G$ if and only if $s_x$ and $s_y$ are conjugate in $G$ and $z^{-1}u_xz$ and $u_y$ are conjugate in $C_G(s_y)$, where $z \in G$ is such that $z^{-1}s_xz = s_y$.
\item If $x \in G$ and $x=su$ is its Jordan decomposition, then $C_G(x) = C_G(s) \cap C_G(u) = C_{C_G(s)}(u)$.
\end{itemize}
Hence, to solve conjugacy problems in $G$, it is usually sufficient to solve them in $G$ for semisimple elements $s$ and in $C_G(s)$ for unipotent elements. One additional simplification is the following. Let $X_i = \bigoplus_{j=1}^{k_i} J_{f_i,e_{ij}}$ be the matrix of the restriction of $x$ to $V_i$. If $y \in C_G(J)$, then $V_iy \subseteq V_i$, so $y$ has block diagonal form $y = \bigoplus_{i=1}^m y_i$, with $y_iX_i = X_iy_i$ for every $i$. Hence, we can solve the centralizer problem separately on each $V_i$. 

Thus we may replace $J$ by a single $X_i$. Having done this, let $J = J_uJ_s$ be the Jordan decomposition of $J$. Then $J_s = C^{[m]} = C\oplus \cdots \oplus C$, with $C$ the companion matrix of an irreducible $f \in F[t]$ of degree $d$, and $C_G(J) \cong \GL_m(q^d)$ (see Section \ref{embedding}). Moreover, $J_u$ is a unipotent element of this centralizer. In \cite[\S 2.3]{Murray}, Murray describes the centralizer in a general linear group of a unipotent element. 
In Section \ref{Cue} we expand this description and obtain an algorithm 
which lists an explicit generating set for such a centralizer. 
For completeness, we discuss in Section \ref{SLconjclasses} 
classes and centralizers in the special linear group $\SL_n(q)$. 

In the rest of the work,  from Chapter \ref{LinearChapter} onwards, 
we consider the other classical groups -- namely, the symplectic, 
orthogonal and unitary groups. We define these groups 
and record some of their basic properties in Section \ref{fincla}.  
We then establish which conjugacy classes 
in $\GL_n(F)$ contain elements in these groups of isometries. This problem 
was considered by Britnell \cite[Chap.\ 5]{JBRIT}, Milnor \cite{Milnor} and 
Wall \cite[\S 2.6,\,\S 3.7]{GEW}. We mainly refer to 
\cite[Chap.\ 5]{JBRIT}, where Britnell solves the membership problem in 
symplectic and orthogonal groups of odd characteristic. 
In Section \ref{membershipChapter} we present his results, and extend them 
to hermitian and quadratic forms (including characteristic 2). The 
main result is Theorem \ref{elementsofCG}, which gives necessary and 
sufficient criteria for $x \in \GL(V)$ to preserve a non-degenerate 
alternating, hermitian or quadratic form on $V$, in terms of the 
elementary divisors of $x$. 

As in the case of the general linear group, the key to analysing conjugacy in 
the classical groups is the Jordan decomposition. As noted already, 
if $G$ is a classical group over a finite field $F$, then to solve the 
conjugacy problems in $G$ it is sufficient to solve it in $G$ for 
semisimple elements $s$ and in $C_G(s)$ for unipotent elements. 
The structure of $C_G(s)$ is described in Chapter \ref{semisimpleChap}: 
as summarised in (\ref{cgxbel}) below, 
when $G$ is a full isometry group, $C_G(s)$ is a direct product of 
classical groups of smaller dimension over extension 
fields of $F$. Hence a key component of our work is the solution of the 
conjugacy problems for unipotent elements in classical groups. 
This is carried out in Chapters \ref{unigoodchap} and \ref{unibadchap}. 
Let $p$ be the characteristic of the field $F$. 
If $G$ is symplectic or orthogonal and $p=2$,
then $p$ is {\it bad} \index{characteristic!bad} for $G$; 
otherwise $p$ is {\it good}\index{characteristic!good} for $G$.
The three conjugacy problems are much more complicated to solve 
for unipotent classes of classical groups in bad characteristic than in good;
they are solved in Chapters \ref{unigoodchap} and \ref{unibadchap} 
for good and bad characteristic respectively.
In both cases, explicit representatives for the unipotent classes in $G$ are 
given by \cite{GLOB}; and the structures and orders of the centralizers can 
be read off from \cite{LS}. However, computing generators for the 
centralizers and solving the other conjugacy problems presents 
new challenges.

Here is a brief description of our approach, first for the case of good 
characteristic. To keep the exposition reasonably brief, we restrict 
the discussion to the symplectic group $G = \Sp(V) \cong \Sp_{2n}(q)$ 
(with $q$ odd). Let $J_i$ denote an $i\times i$ unipotent Jordan block.
It is well known that 
a block diagonal matrix  $\bigoplus_i J_i^{[n_i]}$ is the Jordan form of 
a unipotent element of $G$ if and only if $n_i$ is even for all odd 
values of $i$ -- so all odd-sized Jordan blocks appear with 
even multiplicities. In Section \ref{spgood} we list 
new representatives for the unipotent classes of $G$. These are labelled
\begin{equation}\label{repsp-intro}
u=\bigoplus_{i=1}^r (V_{\b_i}(2k_i) \oplus V_1(2k_i)^{[a_i-1]}) \oplus \bigoplus_{i=1}^s W(2l_i+1)^{[b_i]},
\end{equation}
where each summand $V_1(2k_i)$ or $V_{\b_i}(2k_i)$ is a single Jordan block 
of size $2k_i$, and each summand $W(2l_i+1)$ comprises two Jordan blocks 
of size $2l_i+1$. Each summand in (\ref{repsp-intro}) corresponds to a nilpotent element 
of the symplectic Lie algebra, via the {\it Cayley map}\index{Cayley map}.  
This is a $G$-equivariant bijection from the set of nilpotent elements of 
the Lie algebra to the set of unipotent elements of $G$: a nilpotent element $e$ is mapped to the 
unipotent element $\mbox{$(1-e)(1+e)^{-1}$}$ of $G$. 
For example, a summand $V_1(2k)$ corresponds to the nilpotent element 
\[
\begin{pmatrix} 0&1&&& \\ &0&1&& \\ &&\ddots && \\ &&&0&1 \\ &&&&0
                           \end{pmatrix}
\]
(matrices with respect to a standard symplectic basis). There is a 1-dimensional torus naturally associated with this nilpotent element, namely $\{T(c) : c \in \bar F^*\}$, where $\bar F$ is the algebraic closure of $F$ and $T(c) = {\rm diag}(c^{-(2k-1)},c^{-(2k-3)},\ldots,c^{2k-1})$ (so that $e\,T(c) = c^2e$). We can define a 1-dimensional torus $T$ in this way corresponding to a unipotent element as in (\ref{repsp-intro}), and the stabilizer in $G$ of the flag of $V$ defined by sums of $T$-weight spaces for decreasing weights is a parabolic subgroup $P$. Write $P=QL$, where $Q$ is the unipotent radical and $L = C_G(T)$ is a Levi factor. Then $C_G(u) \le P$ (see Theorem \ref{sympreps}); 
indeed, $C_G(u) = C_Q(u) \rtimes C_L(u)$. We compute generators for both factors in Section \ref{spcent}, solving the centralizer problem. 

Continuing the discussion of symplectic groups in good characteristic, the problem of deciding conjugacy is fairly straightforward: given unipotent $g \in G$, 
we must decide which representative of the form 
(\ref{repsp-intro}) 
it is conjugate to. This is determined by the Jordan form of $g$, apart from computing the values of the parameters $\b_i$, for which an algorithm is given in Section \ref{spconprob}. More challenging is the problem of finding a conjugating element: for example, given $g \in G = \Sp_{2k}(q)$ that is conjugate to a block $u=V_1(2k)$, 
compute $y \in G$ such that $g^y=u$. A key observation is that if $e$ is the associated nilpotent element given above, then $V$ has a basis
\begin{equation}\label{basint}
v,\,ve,\,ve^2,\ldots,ve^{2k-1},
\end{equation}
and this is a standard basis for the symplectic form. 
Consider the nilpotent element $f := (1-g)(1+g)^{-1}$.
We seek $w \in V$ such that $w,\,wf,\,wf^2,\ldots,wf^{2k-1}$ is a basis of $V$ with symplectic form values matching those of  (\ref{basint}). Then the map $y$ sending $wf^i \mapsto ve^i$ for all $i$ is the desired conjugating element. An algorithm to compute such a vector $w$ is given in Section \ref{spconjel}. Carrying this out for each block in an arbitrary unipotent element (\ref{repsp-intro}) leads to the solution of the conjugating element problem in general.

As remarked above, the analysis of unipotent classes in bad characteristic 
is more complicated than in good characteristic. 
Let $G = \Sp(V)$ or $\Or (V)$ in characteristic 2. Representatives for the 
unipotent conjugacy classes were written down in \cite{GLOB}, and take 
the form 
\begin{equation}\label{canonint}
u = \bigoplus_i W(m_i)^{[c_i]} \oplus \bigoplus_j V(2k_j)^{[d_j]}  \oplus  \bigoplus_r W_\a(m_r')  \oplus  \bigoplus_s V_\a(2k_s'),
\end{equation}
where the parameters $m_i,m_r',k_j,k_s'$ satisfy various numerical 
conditions (see (i)--(v) given after (\ref{canon})), 
and $\a$ is a fixed element of $F$ such that the quadratic $x^2+x+\a$ is 
irreducible in $F[x]$. Again, each summand $V(2k_i)$ or $V_{\a}(2k_i)$ is 
a single Jordan block of size $2k_i$, and each summand $W(m_i)$ or 
$W_\a(m_i)$ comprises two Jordan blocks of size $m_i$. 
In characteristic 2 there is no Cayley map.
Nevertheless, we can associate weights 
to the basis vectors for each summand, and use these to define a 
canonical parabolic subgroup $P$ corresponding to the class representative $u$. 
A deep result of Clarke and Premet \cite{CP} implies that $C_G(u) \le P$. 
We do not have a factorization $C_G(u) = C_Q(u) \rtimes C_L(u)$ as in the 
good characteristic case -- indeed, $C_P(u)$ can be a nonsplit extension 
of $C_Q(u)$. However, we can work in $P$ to 
compute generators for $C_G(u)$; this is done 
in Section \ref{centbad}. 

The question of deciding conjugacy is also more complex for 
the bad characteristic case, and is solved in Section \ref{conjprobbad}. 
The problem of computing conjugating elements is 
the most challenging, and is solved 
in Section \ref{badconj}. As an illustration, consider a single 
block $u = V_\a(2k)$: suppose $g \in G = \Sp_{2k}(q)$ (for $q=2^a$) is 
conjugate to $u$; we wish to compute $y \in G$ such 
that $g^y=u$. To do this, we would like to write down a standard 
basis in similar fashion to (\ref{basint}) above, in terms of 
a single vector $v$ and a nilpotent operator $e$. It turns out 
that this can be done if we use three nilpotent operators, 
rather than one: let $f = 1+u$ and $e = f+f^2+\cdots + f^{2k-1}$. 
We show that there is an additional nilpotent operator $h$ such that the 
following sequence of vectors is a standard basis for the symplectic form:
\begin{equation}\label{basbadalint}
v,\,ve,\ldots,ve^{k-2},\,v(e^{k-1}+h),\,v(e^{k-1}+h)f,\ldots,v(e^{k-1}+h)f^k.
\end{equation}
There is a polynomial $p(x)$ of degree $2k-1$ such that $h = p(f)$ -- we cannot write $p(x)$ down explicitly, but it can be computed by machine for any given $k$ and $q$, which is enough for our algorithm.
Now we take a similar approach to that for 
good characteristic. We let 
$f_0 = 1+g$, $e_0=f_0+\cdots +f_0^{2k-1}$, $h_0 = p(f_0)$, and seek $w \in V$ such that 
\[
w,\,we_0,\ldots,we_0^{k-2},\,w(e_0^{k-1}+h_0),\,w(e_0^{k-1}+h_0)f_0,\ldots,w(e_0^{k-1}+h_0)f_0^k
\]
is a basis of $V$ with symplectic form values matching those of  (\ref{basbadalint}). An algorithm to compute such a vector $w$ is given in Section \ref{badconj}. Once again, applying this procedure for each block in an arbitrary unipotent element (\ref{canonint}) leads to the solution of the conjugating element problem in general.

Having dealt with the unipotent classes, 
in Chapter \ref{semisimpleChap} we consider the semisimple classes. Our starting point is a brief paper of Wall \cite{SCCSS} where the semisimple conjugacy classes in symplectic groups of odd characteristic are classified. In Section \ref{semisimpleconjclasses} we extend this work to all sesquilinear and quadratic forms in all positive characteristics.
The structure of the centralizer appears in \cite[Chap.\ 3]{PGPG} and \cite[\S 1]{FS}. Here we briefly discuss our description of the centralizer of a semisimple element for the symplectic and orthogonal groups, given in Theorem \ref{CentrSem}. For a monic polynomial $f(t)=t^d+a_{d-1}t^{d-1}+\cdots +a_0 \in \F_q[t]$ with $a_0\ne 0$, define the {\it dual polynomial} $f^*(t) = a_0^{-1}t^df(t^{-1})$, another monic polynomial of degree $d$. Following \cite[\S 1]{FS}, define
\begin{eqnarray*}
\Phi_1 & := & \{ f: \, f \in F[t] \; | \; f=f^* \mbox{ monic irreducible}, \, \deg{f}=1 \};\\
\Phi_2 & := & \{ f: \, f \in F[t] \; | \; f = gg^*, \, g \neq g^*, \, g \mbox{ monic irreducible} \};\\
\Phi_3 & := & \{ f: \, f \in F[t] \; | \; f=f^* \mbox{ monic irreducible}, \, \deg{f}>1 \}.
\end{eqnarray*}
Define $\Phi := \Phi_1 \cup \Phi_2 \cup \Phi_3$. The elementary divisors of 
every semisimple $s \in G = I_n(q)$ (a symplectic or orthogonal isometry 
group) are either in $\Phi_1  \cup \Phi_3$, or appear in 
pairs $g,g^*$, where $gg^* \in \Phi_2$. 
The polynomials in $\Phi$ are 
the {\it generalized} elementary divisors \index{generalized elementary divisor} of $s$. For each $f \in \Phi$, 
let $m_f$ be its multiplicity as a generalized elementary divisor of $s$, 
and let $d_f = \frac{1}{2}\hbox{deg}(f)$. Then 
\begin{equation}\label{cgxbel}
C_G(s) = \prod_{f \in \Phi_1} I_{m_f}(q) \times \prod_{f \in \Phi_2}\GL_{m_f}(q^{d_f}) \times 
\prod_{f \in \Phi_3}\GU_{m_f}(q^{d_f}).
\end{equation}
An algorithm to compute a generating set for this centralizer is given in 
Section \ref{centgens}.

In Chapter \ref{generalChap} we use the theory developed for semisimple 
and unipotent classes to solve the three conjugacy  problems in the 
general case. The listing of conjugacy classes in $G$ proceeds as follows. 
We first list all semisimple classes; each is identified by a pair $(s,B)$, 
where $s$ and $B$ are the matrices of the representative and of the form 
respectively. If $h_1, \dots, h_r \in \Phi$ are the distinct generalized 
elementary divisors of $s$, then $s = \bigoplus_{i=1}^r s_i$ and 
$B = \bigoplus_{i=1}^r B_i$, where $s_i$ 
and $B_i$ are the restrictions of $s$ and $B$ to the generalized eigenspace 
relative to $h_i$. If $G_i$ is the isometry group of the form $B_i$, 
then $C_{G_i}(s_i)$ is one of the factors 
listed in (\ref{cgxbel}).
Now multiply $s$ by all matrices of the form $u=\bigoplus_{i=1}^r u_i$, 
where $u_i$ runs over the set of representatives for unipotent classes in 
$\C_i := C_{G_i}(s_i)$. 
The set of pairs $(su,B)$ built in this way is a complete set of 
representatives for conjugacy classes in $G$ (see Theorem \ref{CCgen}). 
For every $(su,B)$, the centralizer of $su$ in $G$ is 
$\prod_i C_{\C_i}(u_i)$ (see Theorem \ref{CCgen}). 
An algorithm to compute a generating set for this centralizer is 
given in Section \ref{SectCentrGen}. The conjugacy problem is addressed in Section \ref{conprobgen}.
Finally, given conjugate $x,y \in G$, this theory is used 
in Section \ref{SectConjGen} to compute $z \in G$ such that $x^z=y$.

\section{Providing electronic access to the results} \label{Sect13}
We have developed in \textsc{Magma} 
\cite{MAGMA} implementations of our algorithms to 
solve the conjugacy problems (1)-(3) for classical groups. 
The resulting code is available
publicly at \cite{github} and is also distributed as part 
of {\sc Magma}. 

By default, our conjugacy class representatives are returned as 
elements of the {\it standard copy}\index{standard copy} of the 
classical group $\C = \SL(V), \SU(V), \Sp(V), \SO(V)$, or $\O(V)$, as defined in {\sc Magma}; classes in the corresponding isometry groups are also determined. 
The functions to list conjugacy classes in $\C$ return a sequence of 
triples $ \langle |r|, \, |r^{\C}|, \, r \rangle$,  where $r$ is a 
representative for a conjugacy class of $\C$. 
The semisimple and unipotent classes can be constructed independently. 
We provide a function that, given semisimple $s \in \C$, 
returns representatives for all conjugacy classes of $\C$ having 
$s$ as semisimple part. As discussed in Section \ref{centgens}, 
by applying the algorithm of \cite{GSM}, the data 
is readily translated to any other natural copy of $\C$.

Table \ref{ConjClassTable} records the CPU time in seconds to construct 
the representatives for conjugacy classes of some classical groups $\C$;
it also records the total times taken to construct the centralizer
in $\C$ of a random element from each class, and to construct 
an element of $\C$ which conjugates between two random elements 
from each class.  All calculations were 
carried out using \textsc{Magma} 2.28-13 on a $2.6$ GHz machine.
The relevant commands are listed at \cite{github}.

\begin{table}[thb]\caption{Some calculations in classical groups}\label{ConjClassTable}
\begin{small}
\begin{center}
\begin{tabular}{|r|r|r|r|r|}
\hline
\multicolumn{1}{|c|}{$\C$} & \multicolumn{1}{|c|}{Number of classes} 
& \multicolumn{1}{|r|}{Setup} 
& \multicolumn{1}{|r|}{Centralizer} 
& \multicolumn{1}{|r|}{Conjugation} 
\\
\hline
$\Omega^+_{10}(4)$ & $1543$ & $2$ & 12 & 14 \\
\hline
$\Omega^-_{10}(4)$ & $1593$ & $3$ & 13 & 15 \\
\hline
$\Omega^+_{10}(9)$ & $36177$ & $24$  & 148 & 154\\
\hline
$\Omega^-_{10}(9)$ & $36205$ & $29$  & 160 & 177\\
\hline
$\Omega_{11}(5)$ & $3771$ & $4$  & 38 & 37 \\
\hline
$\Omega_{11}(9)$ & $44933$ & $57$  & 319 & 316 \\
\hline
$\Sp_{10}(4)$ & $2170$ & $2$  & 13 & 11\\
\hline
$\Sp_{10}(9)$ & $107992$ & $46$  & 605 & 519\\
\hline
$\SU_{10}(2)$ & $2340$ & $1$  & 45 & 30\\
\hline
$\SU_{10}(3)$ & $41218$ & $18$  & 537 & 392\\
\hline
$\SL_{10}(4)$ &  349420 &   46 & 644 & 128 \\
\hline
\end{tabular}
\end{center}
\end{small}
\end{table}

\section{The complexity  of our algorithms}
We comment briefly on the theoretical complexity of our solutions.
Of course, many tasks are solved by writing down
explicit solutions.
Most of our remaining tasks for elements of a classical
group $\C$  of degree $n$ defined over a field of size $q$ 
are solved in time $O(n^3 \log q)$.
These include constructing Jordan forms, and constructing and factorising 
minimal polynomials; see for example \cite{vzg, handbook, NP}.
We identify two exceptions. 

\begin{itemize}
\item 
To construct the centralizer in $\C$ of a unipotent element $x$, we construct 
an isomorphic copy $P$ of a unipotent subgroup of $\C$, describe $P$ by 
a power-conjugate presentation, and construct the centralizer in $P$ 
of $x$.  While calculations using such presentations 
are practically very efficient, they rely on ``collection", an 
algorithm not known to run in polynomial time (see \cite[Chap.\ 9]{handbook}).

\item If $x$ and $y$ are unipotent conjugate elements of a classical
group $\C$ of bad characteristic, then we write down 
a system of quadratic equations to construct $z \in \C$ such that
$x^z = y$. This system is solved using a Gr\"obner basis algorithm
(see \cite[Chap.\ 1]{Adams}). 
While our systems are solved readily in practice,
this algorithm is not known to run in polynomial time.
\end{itemize}

\section{Conjugacy in central quotients}

As above, let $\C$ be a standard copy of $\SL(V), \SU(V), \Sp(V), \SO(V)$, or $\O(V)$. Having solved the conjugacy problems for
$\C$, we can solve them for any central quotient $\bar \C = \C/Z$ in an arbitrary representation, 
including the corresponding finite simple classical group. 
For such a group $\bar \C$, there is a constructive recognition algorithm that 
provides an explicit surjective homomorphism $\phi: \C \to \bar \C$;
for details of such algorithms, see for example \cite{black}. 
Using this machinery we can solve the conjugacy problems for $\bar \C$ as follows.

\begin{enumerate}
\item[(1)] Class representatives: Let $(g_i)_{i\in I}$ be a set of class representatives in $\C$. Define
an equivalence relation on these by setting $g_i \sim g_j$ if and only if $g_j$ is $\C$-conjugate to $g_iz$ for some $z \in Z$.
A set of equivalence class representatives $(g_{i_j})_{j\in J}$ can be selected using the solution to the conjugacy 
problem in $\C$. Then $(\phi(g_{i_j}))_{j\in J}$ is a set of class representatives in $\bar \C$.

\item[(2)] Centralizers: Let $\phi(g) \in \bar \C$. For each $z \in Z$ such 
that $g$ is $\C$-conjugate to $gz$, find 
$x_z \in \C$ such that $g^{x_z} = gz$.
Then the centralizer of $\phi(g)$ in ${\bar \C}$
is generated by the image under $\phi$ of $C_{\C}(g)$, together with the images of these conjugating elements $x_z$.

\item[(3)] Conjugacy problem and conjugating element: $\phi(g)$ and $\phi(h)$ 
are conjugate in $\bar \C$ if and only if $g$ is $\C$-conjugate to $hz$ for some $z \in Z$. 
If $\phi(g)$ and $\phi(h)$ are $\bar \C$-conjugate, then 
a conjugating element is $\phi(x)$ where $g^x = hz$.
\end{enumerate}
Our implementations of the resulting algorithms in {\sc Magma}
are available at \cite{github}.

\section{Conjugacy in arbitrary finite groups} \label{Sect14}
We expect that our solution to the conjugacy problems for classical groups 
can assist in their solution for an arbitrary finite group $G$. 
Here is a brief sketch of how this is done in practice.
Existing algorithms follow the ``soluble radical" 
model (see \cite[Chap.\ 10]{handbook}). 
An efficient practical algorithm to construct the necessary 
data structure for this model is described in \cite{BHLO} and 
is available in {\sc Magma}. Its output is a characteristic series for $G$:
$$
1 \leqslant L \leqslant S \leqslant P \leqslant G,
$$
where
\begin{itemize}
\item $L$ is the solvable radical of $G$;
\item $S/L$ is the socle of $G/L$ with $S/L \cong \prod_i T_i^{d_i}$, 
where the $T_i$s are non-abelian, pairwise non-isomorphic simple groups;
\item $P/S \leqslant \prod_i \mathrm{Out}(T_i)^{d_i}$ is solvable;
\item ${G/P \leqslant \prod_i\mathrm{Sym}(d_i)}$.
\end{itemize}
Observe that we have a monomorphism from $G/L$ into the direct product 
$\prod_i W_i$ where $W_i = \mathrm{Aut}(T_i) \wr \mathrm{Sym}(d_i)$. 
\begin{itemize}
\item The problem to extend the solutions from $\mathrm{Aut}(T_i)$ 
to $\mathrm{Aut}(T_i) \wr \mathrm{Sym}(d_i)$, and then to $G/L$, 
was solved by Cannon \& Holt \cite{CannonHolt} and 
Hulpke \cite{Hulpke2000}. 
\item Since $L$ is solvable, there exists a series
$$
L = N_1 \vartriangleright N_2 \vartriangleright \cdots \vartriangleright N_r=1
$$
with $N_i/N_{i+1}$ elementary abelian. For every $i$, the solution of the problem in $G/N_{i+1}$ can be obtained from that in $G/N_i$. This procedure is described in \cite{Hulpke2000, Hulpke2013} for conjugacy classes and \cite[\S 8.8]{handbook} for centralizers.
\end{itemize}
Hence the conjugacy problems in arbitrary finite groups 
reduce to their solutions for finite {\em almost simple} groups. 

\section{Notation}
Our standard references for notation are \cite{KL} and \cite{PG}. 
See also the notation index. 
Symbols used include the following.  \\ 

\noindent 
\makebox[1.7cm][l]{$q=p^a$} for prime number $p$ and  positive integer $a$;\\
\makebox[1.7cm][l]{$\mathbb{F}_q$} finite field of size $q$;\\
\makebox[1.7cm][l]{$F^*, \, \mathbb{F}_q^*$} unit group of $F$, $\mathbb{F}_q$ respectively;\\
\makebox[1.7cm][l]{$F^{*2}, \, \mathbb{F}_q^{*2}$} subgroups of squares in $F^*$, $\mathbb{F}_q^*$ respectively;\\
\makebox[1.7cm][l]{$F[t]$} polynomial ring with coefficients in $F$;\\
\makebox[1.7cm][l]{$M_n(F)$} algebra of $n \times n$ matrices with entries in $F$;\\
\makebox[1.7cm][l]{$\GL_n(F)$} general linear group defined on the vector space $F^n$;\\
\makebox[1.7cm][l]{$V_n(q)$} vector space of dimension $n$ over $\F_q$;\\
\makebox[1.7cm][l]{$\GL_n(q)$} general linear group defined on the vector space $\mathbb{F}_q^n$;\\
\makebox[1.7cm][l]{$\GL(V)$} general linear group defined on the vector space $V$;\\
\makebox[1.7cm][l]{$g^U$} action of matrix $g$ on space $U$; \\
\makebox[1.7cm][l]{$\mathbb{O}, \, \mathbb{O}_n$} zero matrix (dimension not specified), $n \times n$ zero matrix;\\
\makebox[1.7cm][l]{$\mathbb{I}, \, \mathbb{I}_n$} identity matrix (dimension not specified), $n \times n$ identity matrix;\\
\makebox[1.7cm][l]{$A \oplus B$} block diagonal matrix $\left(\begin{array}{cc} A & \mathbb{O}\\ \mathbb{O} & B \end{array}\right)$;\\
\makebox[1.7cm][l]{$\bigoplus M_i$} block diagonal sum of matrices $M_i$; \\
\makebox[1.7cm][l]{$M^{[r]}$}  block diagonal sum $M\oplus \cdots \oplus M$ ($r$ blocks); \\
\makebox[1.7cm][l]{$X^{\tr}$} transpose of matrix $X$;\\
\makebox[1.7cm][l]{$a \mapsto \overline{a}$} automorphism of $\mathbb{F}_q$ of order 1 or 2;\\
\makebox[1.7cm][l]{$\overline{X}$} matrix $(\overline{a_{ij}})$, 
where $X = (a_{ij})$;\\
\makebox[1.7cm][l]{$X^*$} matrix $\overline{X}^{\tr}$;\\
\makebox[1.7cm][l]{$f(t)^*$} polynomial $\overline{a}_0^{-1}\sum_{i=0}^n \overline{a}_{n-i}t^i$, where $f(t) = \sum_{i=0}^n a_it^i$, $a_0 \neq 0$;\\
\makebox[1.7cm][l]{deg $f$} degree of polynomial $f$; \\
\makebox[1.7cm][l]{$\rk{X}$} rank of $X$;\\
\makebox[1.7cm][l]{$\delta_{ij}$} Kronecker delta: $\delta_{ij}=1$ if $i=j$, $0$ otherwise; \\
\makebox[1.7cm][l]{$V \downarrow X$} restriction of a $G$-module $V$ to $X$, an element or subgroup of $G$; \\
\makebox[1.7cm][l] 
{$[n]$}  group of order $n$; \\
\makebox[2.7cm][l] {$\hbox{diag}(\a_1,\ldots,\a_n)$}  diagonal matrix with diagonal entries $\a_1,\ldots,\a_n$. \\

\chapter{General and special linear groups}\label{glsl}

In this chapter we give a complete description of 
conjugacy classes and centralizers in the general
and special linear groups. Although the theory for 
this case is well known, we 
make a new contribution: an explicit 
generating set for
the centralizer of an arbitrary element.

\section{Conjugacy classes in $\GL_n(q)$}
Let $q=p^a$, where $p$ is prime and $a$ is a positive integer, 
and let $F=\mathbb{F}_q$.
Let $V$ be an $n$-dimensional vector space over $F$
and 
let $x \in \GL(V) \cong \GL_n(q)$. Let $f_1(t)^{e_1} \cdots f_h(t)^{e_h}$ be the 
minimal polynomial of $x$, where $f_1,\ldots,f_h \in F[t]$ are distinct 
monic irreducible polynomials. We write
$$
V = V_1 \oplus \cdots \oplus V_h,
$$
where $V_i = \ker{(f_i(x)^{e_i})}$ is the 
\textit{generalized eigenspace}\index{generalized eigenspace} 
corresponding to $f_i(t)$ for every $i$. 
By \cite[4.5.1]{PFGG}, every $V_i$ can be written as a direct sum of 
$x$-invariant subspaces
\begin{eqnarray}\label{decompCSM}
V_i = V_{i,1} \oplus \cdots \oplus V_{i,k_i},
\end{eqnarray}
where $x$ acts cyclically on $V_{i,j}$ with minimal polynomial $f_i(t)^{e_{ij}}$, and $1 \leq e_{i1} \leq \cdots \leq e_{ik_i} = e_i$. The polynomials $f_1(t)^{e_{11}}, \dots, f_1(t)^{e_{1k_1}}, \dots, f_h(t)^{e_{h1}}, \dots, f_h(t)^{e_{hk_h}}$ are the \textit{elementary divisors}\index{elementary divisors} of $x$. 
If $g$ is a power of an irreducible polynomial, then the 
\textit{multiplicity} of $g$ as an elementary divisor of $x$ is 
the number of times $g$ appears in the list of elementary divisors of $x$. 
Since the decomposition in (\ref{decompCSM}) is unique up to rearranging the 
factors, the list of elementary divisors of $x$ is well-defined. 
Elements of $\GL(V)$ are conjugate if and only if they have the 
same elementary divisors; see for example \cite[6.7.3]{PFGG}.

As a representative for each conjugacy class of $\GL(V)$, we choose the Jordan form, defined as follows. 
For every monic irreducible polynomial $f$ of degree $d$ and positive integer $e$, the \textit{Jordan block}\index{Jordan block} 
$J_{f,e}$ 
of order $e$ relative to $f$ is the block matrix
$$
J_{f,e}= \begin{pmatrix}
C & \mathbb{I} && \\ & \ddots & \ddots & \\ && \ddots & \mathbb{I} \\ &&& C
\end{pmatrix},
$$
where $C$, the companion matrix of $f$, appears $e$ times, and
$\mathbb{I}$\index{$\mathbb{I}$, $\mathbb{I}_n$} is the $d \times d$ identity matrix. A Jordan block is \textit{unipotent}\index{Jordan block!unipotent} if $f(t)=t-1$. For every $V_{i,j}$ as in (\ref{decompCSM}) there is a basis such that the matrix of the restriction of $x$ to $V_{i,j}$ is the Jordan block $J_{f_i, e_{ij}}$.
Hence, there exists a basis of $V$ such that the matrix of $x$ is a 
diagonal join of Jordan blocks. 
This matrix is the \textit{Jordan form}\index{Jordan form} of $x$.

We comment briefly on problems (1) and (3) of Section \ref{mainprob}: 
we take as conjugacy class representatives in $\GL_n(q)$ the list of distinct Jordan forms of its elements;
and for $x \in \GL_n(q)$, an element that conjugates $x$ to 
its Jordan form is determined by 
a Jordan basis algorithm (see for example \cite{JF}). 
Hence we focus on the construction of 
the centralizer of $x$ in $\GL_n(q)$.

\section{Centralizers in $\GL_{\MakeLowercase{n}}(\MakeLowercase{q})$}\label{gl-cent}
In this section we record how to construct the centralizer of 
$x \in \GL(V) \cong \GL_n(q)$. The results are well known, except 
for the generation of 
centralizers of unipotent elements; our work is motivated by that of 
Murray \cite{Murray}. 

Let $x \in G := \GL(V)$ be in Jordan form, with
minimal polynomial $f_1(t)^{e_1} \cdots f_h(t)^{e_h}$ 
and associated generalized eigenspaces $V_i$.
Every element of $C_G(x)$ fixes each $V_i$.
Hence, if $x = \bigoplus_{i=1}^h x_i$, where $x_i$ is the restriction of 
$x$ to $V_i$, then 
\begin{eqnarray}\label{SplitCentr}
C_{\GL(V)}(x) = \prod_{i = 1}^h 
C_{\GL(V_i)}(x_i). 
\end{eqnarray}

Recall the {\it Jordan decomposition}\index{Jordan decomposition}:  
if $x \in G$, then 
$x = su=us$ for unique semisimple $s$ and unipotent $u$. 
If $x$ is a diagonal join of Jordan blocks, then $s$ and 
$u$ are the diagonal joins of the semisimple and the unipotent 
parts of each block.
To compute $C_G(x)$, 
we first compute $C_G(s)$ and then, using the fact 
that $u \in C_G(s)$, compute $C_{C_G(s)}(u)$.

\subsection{Centralizer of a semisimple element}\label{embedding}
Let $x \in G := \GL_n(q)$ be semisimple. 
By (\ref{SplitCentr}),
we can assume that $x$ has a unique elementary divisor 
$f(t) \in \mathbb{F}_q[t]$, which is irreducible.

If $\deg{f}=1$, then $x$ is a scalar matrix, so $C_G(x)=\GL_n(q)$.

Now suppose that $f(t) \in \mathbb{F}_q[t]$ is irreducible of degree $r>1$. 
Let $E = \mathbb{F}_q[t]/(f)$ be the splitting field of $f$ 
over $\mathbb{F}_q$ 
and let $\alpha \in E$ be a root of $f$. Every element of $E$ can be written as 
$\phi(\alpha)$ for some polynomial $\phi(t) \in \mathbb{F}_q[t]$ of degree 
smaller than $r$. For every positive integer $m$, 
there is a canonical embedding 
of $\GL_m(E) = \GL_m(q^r)$ into $\GL_{mr}(q)$ sending the 
matrix $(\phi_{ij}(\alpha))$ into the block matrix $(\phi_{ij}(C))$, 
where $C$ is the companion matrix of $f$ (see \cite[2.1.4]{PGPG}). 

Now assume $x$ is semisimple with a unique elementary 
divisor $f(t)$ of degree $r>1$ and multiplicity $m$. We can suppose 
that the matrix of $x$ is a diagonal join of $m$ copies of $C$, 
the companion matrix of $f$. Then $x$ is the embedding 
into $\GL_{mr}(q)$ of the scalar matrix $\alpha \mathbb{I}_m \in \GL_m(E)$. Every matrix of $\GL_m(E)$ commutes with $\alpha \mathbb{I}_m$, so its embedding into $\GL_{mr}(q)$ commutes with $x$. On the other hand, these are the only matrices in $\GL_{mr}(q)$ that commute with $x$ (see \cite[Lemma 3.1.9]{PGPG}). Thus we have proved the following.

\begin{proposition} \label{largerfield}
	Let $x \in \GL_n(q)$ be semisimple with characteristic polynomial $\prod_{i=1}^hf_i(t)^{e_i}$, where $f_1,\ldots,f_h \in \mathbb{F}_q[t]$ are distinct irreducibles and $\deg f_i = r_i$. Then
	$$
	C_{\GL_n(q)}(x) \cong \prod_{i=1}^{h} \GL_{e_i}(q^{r_i}).
	$$
\end{proposition}

\subsection{Centralizer of a unipotent element} \label{Cue}
Let $x \in G := \GL_n(q)$ be unipotent of order $p^b$. 
The structure of $C_G(x)$ is well known;
our new contribution is to describe an explicit generating set in Section \ref{CGgeneration}.
Following the approach 
of \cite{Murray}, it is convenient to work in the matrix algebra 
$M = M_n(\mathbb{F}_q)$ and find the centralizer $C_M(x)$ of $x$ in $M$. 
The centralizer of $x$ in $\GL_n(q)$ is the set of invertible elements 
of $C_M(x)$.

Since the unique eigenvalue of $x$ is $1$, the Jordan form of $x$ is
$$
\left(\begin{array}{cccc}
J_{\lambda_1} & & & \\
& J_{\lambda_2} & & \\
& & \ddots & \\
& & & J_{\lambda_k}
\end{array}\right),
$$
where $J_{\lambda_i}$ is the unipotent Jordan block of 
dimension $\lambda_i$ and $\lambda_1 + \cdots + \lambda_k = n$. 
We suppose $\lambda_1 \leq \lambda_2 \leq \cdots \leq \lambda_k$ and 
denote this partition of $n$ by $\lambda$. 
Let $x$ be this Jordan form.

Take an element $y$ centralizing $x$ and write the matrix of $y$ as
$$
\begin{pmatrix}
B_{11} & \cdots & B_{1k}\\
\vdots & \ddots & \vdots \\
B_{k1} & \cdots & B_{kk}
\end{pmatrix},
$$
where $B_{ij}$ is a block of dimension $\lambda_i \times \lambda_j$ for every $i,j$.

The condition $xy=yx$ implies $J_{\lambda_i}B_{ij} = B_{ij}J_{\lambda_j}$ for every $1 \leq i,j \leq k$. Write $B_{ij} = (b_{\onemu,\onenu})$ for $b_{\onemu,\onenu} \in \mathbb{F}_q$, where $0 \leq \onemu \leq \lambda_i-1$ and $0 \leq \onenu \leq \lambda_j-1$. 
A simple computation shows that $xy=yx$ is equivalent to 
$b_{\onemu,\onenu} = b_{\onemu+1,\onenu+1}$ for every $\onemu,\onenu$; also 
$b_{\onemu,\onenu} = 0$ if $\lambda_i \leq \lambda_j$ and
$\onemu \geq \onenu-(\lambda_j-\lambda_i-1)$, or 
$\lambda_i > \lambda_j $ and $\onemu \geq \onenu + 1$. 
Hence $B_{ij}$ is an upper 
triangular rectangular matrix with constant upper diagonals; it is either 
\begin{eqnarray}   \label{blocksBij}
\left( \mathbb{O} \begin{array}{lllll}
& b_{\lambda_j-\lambda_i} & b_{\lambda_j-\lambda_i+1} & \cdots & b_{\lambda_j-1} \\
& 0 & b_{\lambda_j-\lambda_i} & \ddots & \vdots \\
& \vdots & \ddots & \ddots & b_{\lambda_j-\lambda_i+1}\\
& 0 & \cdots & 0 & b_{\lambda_j-\lambda_i}
\end{array} \right) \mbox{ or } \left(\begin{array}{c}
\begin{array}{llll}
b_0 & b_1 & \cdots & b_{\lambda_j-1} \\
0 & b_0 & \ddots & \vdots \\
\vdots & \ddots & \ddots & b_1\\
0 & \cdots & 0 & b_0 \\
&&&
\end{array}\\
\mathbb{O}
\end{array}\right).
\end{eqnarray}
Denote by $X_{c\times d}^a$ the $c \times d$ matrix whose $(i,j)$-entry is 1 if $j-i = a$ and 0 otherwise. We can write $B_{ij} = \sum_{a=\lambda_j-\lambda_i}^{\lambda_j-1} b_aX_{\lambda_i\times \lambda_j}^a$ (resp.\ $\sum_{a=0}^{\lambda_j-1}b_aX_{\lambda_i\times \lambda_j}^a$) if $\lambda_i\leq \lambda_j$ (resp.\ $\lambda_i > \lambda_j$). It is easy to check that $X_{c\times d}^aX_{d \times e}^b = X_{c \times e}^{a+b}$ and this gives an algebra homomorphism from
$$
\left( \begin{array}{lllll}
\mathbb{F}_q[t] & t^{\lambda_2-\lambda_1}\mathbb{F}_q[t] & t^{\lambda_3-\lambda_1} \mathbb{F}_q[t] & \cdots & t^{\lambda_k-\lambda_1}\mathbb{F}_q[t] \\
\mathbb{F}_q[t] & \mathbb{F}_q[t] & t^{\lambda_3-\lambda_2} \mathbb{F}_q[t] & \cdots & t^{\lambda_k-\lambda_2}\mathbb{F}_q[t]\\
\vdots & \vdots & \vdots & \ddots & \vdots \\
\mathbb{F}_q[t] & \mathbb{F}_q[t] & \mathbb{F}_q[t] & \cdots & \mathbb{F}_q[t]
\end{array} \right)
$$
to $C_M(x)$ sending the element $t^a$ in the $(i,j)$-entry to $X_{\lambda_i\times\lambda_j}^a$ in the $(i,j)$-block and extended by linearity. From our expression for the $B_{ij}$, it is clear that the homomorphism is surjective and the kernel is the set of matrices described by
$$
\left(\begin{array}{llll}
(t^{\lambda_1}) & (t^{\lambda_2}) & \cdots & (t^{\lambda_k})\\
(t^{\lambda_1}) & (t^{\lambda_2}) & \cdots & (t^{\lambda_k})\\
\vdots & \vdots & \ddots & \vdots \\
(t^{\lambda_1}) & (t^{\lambda_2}) & \cdots & (t^{\lambda_k})
\end{array}\right),
$$
where $(t^{\lambda_i})$ is the ideal of $\mathbb{F}_q[t]$ generated 
by $t^{\lambda_i}$. The centralizer of $x$ in the algebra $M$ is 
isomorphic to the quotient algebra
\begin{equation}\label{kbyk}
\mathbb{F}_q[t]_{\lambda} = \left( \begin{array}{llll}
\mathbb{F}_q[t]_{\lambda_1} & t^{\lambda_2-\lambda_1}\mathbb{F}_q[t]_{\lambda_2} & \cdots & t^{\lambda_k-\lambda_1}\mathbb{F}_q[t]_{\lambda_k} \\
\mathbb{F}_q[t]_{\lambda_1} & \mathbb{F}_q[t]_{\lambda_2} & \cdots & t^{\lambda_k-\lambda_2}\mathbb{F}_q[t]_{\lambda_k}\\
\vdots & \vdots & \ddots & \vdots \\
\mathbb{F}_q[t]_{\lambda_1} & \mathbb{F}_q[t]_{\lambda_2} & \cdots & \mathbb{F}_q[t]_{\lambda_k}
\end{array} \right),
\end{equation}
where $\mathbb{F}_q[t]_{\lambda_s} = \mathbb{F}_q[t]/(t^{\lambda_s})$ 
is the truncated polynomial algebra. The multiplication
$$
{\mathbb{F}_q[t]_{\lambda_s} \times \mathbb{F}_q[t]_{\lambda_{s'}} \rightarrow \mathbb{F}_q[t]_{\lambda_{s'}}}
$$
is defined by multiplying the two polynomials and removing all of 
the monomials of degree greater than $\lambda_{s'}$.

Searching for invertible elements of $C_M(x)$ is equivalent to searching for invertible elements of $\mathbb{F}_q[t]_{\lambda}$. The $(\mymu,\mynu)$-entry of an arbitrary element of $\mathbb{F}_q[t]_{\lambda}$ corresponds to a pair $(\lambda_{\mymu},\lambda_{\mynu})$.

Since the $\lambda_{\mymu}$ are not necessarily distinct, 
it is convenient to redefine the notation as follows: 
let $\lambda_1 < \lambda_2 < \cdots < \lambda_h$ be the distinct dimensions 
of the Jordan blocks of $x$ and let $l_i$ be the multiplicity of $\lambda_i$ 
where $1 \leq i \leq h$. 
We assemble the entries sharing the same 
values $(\lambda_i,\lambda_j)$ in a unique 
$l_i\times l_j$ block. With this new notation, the definition of $\Ft$ becomes
$$
\Ft = \begin{pmatrix*}[l]
M_{l_1 \times l_1} (\mathbb{F}_q[t]_{\lambda_1}) & M_{l_1 \times l_2}(t^{\lambda_2-\lambda_1} \mathbb{F}_q[t]_{\lambda_2}) & \cdots & M_{l_1 \times l_h}(t^{\lambda_h-\lambda_1}\mathbb{F}_q[t]_{\lambda_h}) \\
M_{l_2 \times l_1} (\mathbb{F}_q[t]_{\lambda_1}) & M_{l_2 \times l_2}(\mathbb{F}_q[t]_{\lambda_2}) & \cdots & M_{l_2 \times l_h}(t^{\lambda_h-\lambda_2}\mathbb{F}_q[t]_{\lambda_h}) \\
\vdots & \vdots & \ddots & \vdots \\
M_{l_h \times l_1}(\mathbb{F}_q[t]_{\lambda_1}) & M_{l_h\times l_2}(\mathbb{F}_q[t]_{\lambda_2}) & \cdots & M_{l_h \times l_h}(\mathbb{F}_q[t]_{\lambda_h})
\end{pmatrix*},
$$
where $M_{l_i \times l_j}(F)$ is the set of $l_i \times l_j$ matrices with coefficients in $F$. An arbitrary element 
of $\mathbb{F}_q[t]_{\lambda}$ can be written as a block matrix
\begin{eqnarray}\label{LF}
A = \left(\begin{array}{ccc} A_{11} & \cdots & A_{1h}\\ \vdots & \ddots & \vdots \\ A_{h1} & \cdots & A_{hh}\end{array}\right),
\end{eqnarray}
where $A_{ij} \in M_{l_i \times l_j}(t^{\lambda_j-\lambda_i}\mathbb{F}_q[t]_{\lambda_j})$.
Such $A_{ij}$ can be written as
\begin{eqnarray}\label{matrix-polynomials}
A_{ij} = A_{ij}^{(0)} + tA_{ij}^{(1)} + t^2A_{ij}^{(2)} + t^3A_{ij}^{(3)} + \cdots + t^{\lambda_j-1}A_{ij}^{(\lambda_j-1)},
\end{eqnarray}
where $A_{ij}^{(s)}$ is the $l_i \times l_j$ matrix of the coefficients of $t^s$ in the entries of $A_{ij}$.

Now take an arbitrary $A \in \mathbb{F}_q[t]_{\lambda}$ as in (\ref{LF}), 
and write it as 
$$A=A^{(0)}+tA^{(1)}+\cdots + t^{\lambda_h-1}A^{(\lambda_h-1)}$$ 
with
\begin{eqnarray}\label{Bdecomposition}
A^{(s)} = \left(\begin{array}{ccc} A_{11}^{(s)} & \cdots & A_{1h}^{(s)}\\ \vdots & \ddots & \vdots \\ A_{h1}^{(s)} & \cdots & A_{hh}^{(s)}\end{array}\right),
\end{eqnarray}
where the $A_{ij}^{(s)}$ are those defined in (\ref{matrix-polynomials}). It is clear that $A$ is invertible in $\mathbb{F}_q[t]_{\lambda}$ if and only if $A^{(0)}$ is invertible. In $\mathbb{F}_q[t]_{\lambda}$ polynomials in the blocks above the main diagonal have leading term zero -- in other words,
\begin{eqnarray}\label{B0decomposition}
A^{(0)} = \left(\begin{array}{ccc} A_{11}^{(0)} &  & 0\\  \vdots &  \ddots & \\ A_{h1}^{(0)} & \cdots & A_{hh}^{(0)}\end{array}\right).
\end{eqnarray}
Thus $A$ is invertible if and only if $A_{ii}^{(0)}$ is invertible for every $i$, equivalently $A_{ii}^{(0)} \in \GL_{l_i}(q)$.

Define the following subgroups of $\mathbb{F}_q[t]_{\lambda}^*$,
the unit group of $\Ft$:
\begin{itemize}
\item $R$ is the subgroup of $\mathbb{F}_q[t]_{\lambda}^*$ consisting 
of the matrices of the form
$$
\left(\begin{array}{ccc} A_{11}^{(0)} &  & \\  &  \ddots & \\  &  & A_{hh}^{(0)}\end{array}\right)
$$
with $A_{ii}^{(0)} \in \GL_{l_i}(q)$.
\item $U$ is the subgroup of $\mathbb{F}_q[t]_{\lambda}^*$ 
consisting of the matrices of the form (\ref{LF}) where $A_{ii}^{(0)}$ is the identity matrix for every $i$; equivalently, 
$U$ consists of the matrices
$$
\left(\begin{array}{llll}
1+tU_{11} & U_{12} & \cdots & U_{1h}\\
U_{21} & \ddots & \ddots & \vdots \\
\vdots & \ddots & \ddots & U_{h-1,h}\\
U_{h1} & \cdots & U_{h,h-1} & 1+tU_{hh} \end{array}\right),
$$
where $U_{ij}$ is an arbitrary $l_i\times l_j$ block with entries in $t^{\lambda_j-\lambda_i}\mathbb{F}_q[t]_{\lambda_j}$ (with $t^{\lambda_j-\lambda_i}=1$ if $\lambda_i>\lambda_j$). 
\end{itemize}

\begin{proposition}
Let $R$ and $U$ be defined as above. 
\begin{itemize}
\item[{\rm (i)}] $R \cap U$ is the trivial subgroup.
\item[{\rm (ii)}] $\mathbb{F}_q[t]_{\lambda}^* \subseteq UR$.
\item[{\rm (iii)}] $U$ is a normal subgroup of $UR$.
\end{itemize}
Hence $\mathbb{F}_q[t]_{\lambda}^* = U \rtimes R$.
\end{proposition}
\begin{proof}
Part (i) is clear. For (ii), let $A \in \mathbb{F}_q[t]_{\lambda}^*$ and write $A = A^{(0)} + tA^{(1)} + t^2A^{(2)} + \cdots$ as in (\ref{Bdecomposition}). If $D$ is the block diagonal matrix $\bigoplus_{i=1}^h A_{ii}^{(0)}$, where $A_{ii}^{(0)}$ are as in (\ref{B0decomposition}), then $D \in R$ and $AD^{-1} \in U$.

Finally, consider (iii). 
 Let $B \in U$ and $X \in \mathbb{F}_q[t]_{\lambda}^*$ arbitrary. Write $B = B^{(0)} + tB^{(1)} + t^2B^{(2)}+\cdots$ and ${X = X^{(0)} + tX^{(1)} + t^2X^{(2)} + \cdots}$ as we did for $A$. Note that $X^{-1} = {X^{(0)}}^{-1} + tX'$ for some $X' \in \mathbb{F}_q[t]_{\lambda}$, and $X^{-1}BX = {X^{(0)}}^{-1}B^{(0)}X^{(0)} + tB'$ for some $B' \in \mathbb{F}_q[t]_{\lambda}$. As shown in (\ref{B0decomposition}), $X^{(0)}$ and $B^{(0)}$ are lower triangular block matrices; so is ${X^{(0)}}^{-1}B^{(0)}X^{(0)}$. Moreover, the blocks on the main diagonal of ${X^{(0)}}^{-1}B^{(0)}X^{(0)}$ are ${X_{ii}^{(0)}}^{-1}B_{ii}^{(0)}X_{ii}^{(0)} = 1$, since $B_{ii}^{(0)} = 1$. This is exactly the condition $X^{-1}BX \in U$.
\qedhere
\end{proof}

As in the previous chapter, in the next result we use the notation $\bigoplus_{i=1}^h J_{\l_i}^{[l_i]}$ to denote a block diagonal matrix having $l_i$ diagonal blocks $J_{\l_i}$ for $i=1\ldots,h$.

\begin{theorem} \label{Card}
Let $x = \bigoplus_{i=1}^h J_{\l_i}^{[l_i]}$ be a unipotent element of $\GL_n(q)$. Then $C_{\GL_n(q)}(x) =U \rtimes R$, where $R \cong \prod_{i=1}^h \GL_{l_i}(q)$ and $|U|=q^\m_gamma$ with 
\begin{eqnarray} \label{CentStruct}
 \m_gamma = 2\sum_{i<j}\lambda_il_il_j + \sum_i (\lambda_i-1)l_i^2.
\end{eqnarray}
\end{theorem}
\begin{proof}
Let $G = \GL_n(q)$. Working in $\mathbb{F}_q[t]_{\lambda}$ instead of $C_G(x)$, the 
groups $R$ and $U$ are the subgroups described above, and $C_G(x) \cong U \rtimes R$. 
That $R \cong \prod_{i=1}^h \GL_{l_i}(q)$ follows by definition of $R$. 

It remains to compute the cardinality of $U$, for which we need to determine how many choices there are for 
the $U_{ij}$. The term $1+tU_{ii}$ equals 
$1+tA_{ii}^{(1)} + \cdots + t^{\lambda_i-1}A_{ii}^{(\lambda_i-1)}$ and 
$A_{ii}^{(s)}$ is a matrix in $M_{l_i \times l_i}(\mathbb{F}_q)$ which can be chosen 
arbitrarily; thus, for $U_{ii}$ there are $q^{(\lambda_i-1)l_i^2}$ choices and 
these give the second sum in (\ref{CentStruct}).

If $i > j$, then $\lambda_i>\lambda_j$ and $U_{ij} = A_{ij}^{(0)} + \cdots + t^{\lambda_j-1}A_{ij}^{(\lambda_j-1)}$. Every matrix $A_{ij}^{(s)}$ is an arbitrary $l_i \times l_j$ matrix with entries in $\mathbb{F}_q$, so the number of possible choices for $U_{ij}$ is $q^{\lambda_i l_il_j}$.

If $i<j$, then $\lambda_i < \lambda_j$ and $U_{ij} = t^{\lambda_j-\lambda_i}A_{ij}^{(\lambda_j-\lambda_i)} + \cdots + 
t^{\lambda_j-1}A_{ij}^{(\lambda_j-1)}$. Again, $A_{ij}^{(s)}$ is an arbitrary $l_i \times l_j$ matrix with entries in $\mathbb{F}_q$, so the number of possible choices for $U_{ij}$ is $q^{\lambda_j l_il_j}$.

Summing over all $i \neq j$, we get the first sum in  (\ref{CentStruct}).
\end{proof}

\subsection{Generators for the centralizer of a unipotent element}\label{CGgeneration}

As above, let $x = \bigoplus_{i=1}^h J_{\l_i}^{[l_i]}$
be a unipotent element of $\GL_n(q)$. 
We now describe a generating set for 
$C_{\GL_n(q)}(x)$ by  working in $\Ft$. We continue the notation of the previous section.

Observe that every element of $\Ft$ can be written as a block matrix $A = (A_{ij})$ as in (\ref{LF}). For example, if $A$ is the identity matrix, then $A_{ii}$ is the $l_i \times l_i$ identity matrix with coefficients in $\mathbb{F}_q[t]_{\lambda_i}$ and $A_{ij}=0$ for $i \neq j$.

We introduce the following notation:
\begin{itemize}
\item $\matalpha_i(y)$, for $1 \leq i \leq h$ and 
$y \in \GL_{l_i}(\mathbb{F}_q[t]_{\lambda_i})$,
is the matrix obtained by taking the identity matrix in $\Ft$
and replacing the block $A_{ii}$ by $y$;
\item $\matbeta_i(\mu)$, for $1 \leq i \leq h$ and 
$\mu \in \mathbb{F}_q[t]_{\lambda_i}$, 
is the matrix obtained by 
taking the identity matrix in $\Ft$
and replacing the block $A_{ii}$ by the diagonal matrix
$$
\left( \begin{array}{cccc} \mu &&& \\ & 1 && \\ && \ddots & \\ &&& 1 \end{array}\right);
$$
\item $\matgamma_i^+(\mu)$ (resp.\ $\matgamma_i^-(\mu)$), for $1 \leq i \leq h-1$ and $\mu \in \mathbb{F}_q[t]_{\lambda_{i+1}}$ (resp.\ $\mathbb{F}_q[t]_{\lambda_i}$), is the 
matrix obtained by taking the identity matrix in $\Ft$
and replacing the block $A_{i,i+1}$ (resp.\ $A_{i+1,i}$) by the matrix 
having $\mu$ in the bottom-left (resp.\ top-right) corner and 0 elsewhere. We write $\matgamma_i^+$ and $\matgamma_i^-$ for $\matgamma_i^+(1)$ and $\matgamma_i^-(1)$ respectively.
\end{itemize}

Recall, for example from \cite[p.\ 59]{handbook},
that a {\it primitive element}\index{primitive element}
of $\F_q$ is a generator for the cyclic group $\F_q^*$.

\begin{theorem}\label{ThmGensUnipotent}
Let $\omega$ be a primitive element of $\mathbb{F}_q$. 
The following collection of matrices forms a generating set for the unit group $\Ftunit$:
\begin{itemize}
\item the set of $\matalpha_i(y)$, where $1 \leq i \leq h$ and $y$ runs over a generating set for $\GL_{l_i}(q)$ 
$($as subgroup of $\GL_{l_i}(\mathbb{F}_q[t]_{\lambda_i})$; 
\item the set of $\matbeta_i(\mu)$, where $1 \leq i \leq h$ and $\mu$ runs over all elements of the set 
$$
\{ 1+\omega^jt^l: 0 \leq j \leq [\mathbb{F}_q:\mathbb{F}_p], 1 \leq l \leq \lambda_i-1\};
$$
\item the set of all $\matgamma_i^+$ and $\matgamma_i^-$, where $1 \leq i \leq h-1$.
\end{itemize}
\end{theorem}

To prove the theorem we  first establish some preliminary facts.
The following lemma shows that, 
with the generators listed above, we can get all elements of 
the form $\matalpha_i(y)$ for $y \in \GL_{l_i}(q)$, and all 
elements of the form $\matbeta_i(\mu)$ for  
$\mu \in \mathbb{F}_q[t]_{\lambda_i}$.
\begin{lemma}\label{LemmaGensOfFqt}
Let $\omega$ be a primitive element of $\mathbb{F}_q$. A generating set 
for $\mathbb{F}_q[t]^{*}_{\lambda_i}$ is given by
\begin{eqnarray}\label{generatorsOfFqt}
\Delta := \{\omega\} \cup \{ 1+\omega^jt^l: 0 \leq j \leq [\mathbb{F}_q:\mathbb{F}_p], 1 \leq l \leq \lambda_i-1\}.
\end{eqnarray}
\end{lemma}
\begin{proof}
Let $r = [\mathbb{F}_q:\mathbb{F}_p]$. Let $a_0 +a_1t+a_2t^2+\cdots +a_{\lambda_i}t^{\lambda_i}$ be an element of $\mathbb{F}_q[t]_{\lambda_i}^*$, with $a_j \in \mathbb{F}_q$ and $a_0 \neq 0$. 
We claim that for any $\myN \le \l_i$ there is a product of elements of $\Delta$ such that 
the first $\myN+1$ coefficients coincide with $a_0, \dots, a_\myN$. We prove this by induction on $\myN$.

For $\myN=0$, this is obvious since $a_0 \neq 0$, so $a_0$ is a power of $\omega$. Now suppose there exists a product $\myP(t)$ of elements of $\Delta$, say $\myP(t) = a_0 + a_1t + \cdots + a_\myN t^\myN + b_{\myN+1}t^{\myN+1} + \cdots + b_{\lambda_i}t^{\lambda_i}$ for some $b_{\myN+1}, \dots, b_{\lambda_i} \in \mathbb{F}_q$. Since $a_{\myN+1}-b_{\myN+1} \in \mathbb{F}_q$, we can find $y_0, \dots, y_{r-1} \in \mathbb{F}_p$ such that $
a_0^{-1}(a_{\myN+1}-b_{\myN+1}) = y_0 + y_1\omega+y_2\omega^2 + \cdots + y_{r-1}\omega^{r-1}$.
A straightforward computation shows that 
\begin{eqnarray*}
\myP(t) \cdot \prod_{j=0}^{r-1} (1+\omega^j t^{\myN+1})^{y_j} & = & \myP(t) \cdot (1 + t^{\myN+1}(y_0 + y_1 \omega + \cdots + y_{r-1}\omega^{r-1}) +t^{\myN+2}Q(t))\\
& = & (a_0 + a_1t + \cdots + a_\myN t^\myN + b_{\myN+1}t^{\myN+1} + t^{\myN+2}R(t)) \cdot \\
& & (1+a_0^{-1}(a_{\myN+1}-b_{\myN+1})t^{\myN+1} + t^{\myN+2}Q(t))\\
& = & a_0 + a_1t + \cdots + a_{\myN+1}t^{\myN+1} + t^{\myN+2}S(t),
\end{eqnarray*}
where $Q(t)$, $R(t)$ and $S(t)$ are elements of $\mathbb{F}_q[t]$. 
Hence, we have exhibited a product of elements of $\Delta$ whose 
first $\myN+1$ coefficients are exactly $a_0, \dots, a_{\myN+1}$.
\end{proof}

\begin{proposition}\label{PropMuDiagonal}
Every diagonal matrix in $\Ftunit$ can be written as a product of elements 
of the form $\matalpha_i(y)$  for $y \in \GL_{l_i}(q)$ and  $\matbeta_i(\mu)$
for $\mu \in \mathbb{F}_q[t]_{\lambda_i}$.
\end{proposition}
\begin{proof}
It is sufficient to prove that we can build every diagonal matrix where all 
but one element on the main diagonal are 1. Say $\myz$ is a diagonal matrix whose diagonal has $\mu \neq 0,1$ in a certain entry and 1 elsewhere. If $\mu$ belongs to the block $A_{ii}$, then $\myz$ can be obtained by conjugating the matrix $\matbeta_i(\mu)$ by an appropriate permutation matrix in the block $A_{ii}$ (permutation matrices have coefficients in $\mathbb{F}_q$, so they are products of elements of the form $\matalpha_i(y)$).
\end{proof}

For the rest of the section, we regard elements of $\Ft$ as $k\times k$ matrices, as in (\ref{kbyk}). 
Let $\mathbb{I}$ be the identity matrix in $\Ft$ and let $E_{j\ell}$ be the 
matrix in $\Ft$ having 1 in the $(j,\ell)$-entry where $1 \leq j,\ell \leq k$
and 0 elsewhere. 
\begin{proposition}
Every matrix of the form $\mathbb{I} + \nu E_{j\ell}$, with $1 \leq j,\ell \leq k$ and $\nu \in \mathbb{F}_q[t]_{\lambda_s}$ for some $s$, can be obtained as a product of elements of the form $\matalpha_i(y)$, $\matbeta_i(\mu)$, $\matgamma_i^+$ and $\matgamma_i^-$.
\end{proposition}
\begin{proof}
Consider the elements with $j<\ell$ (namely, the entry $\nu$ is above 
the main diagonal). Observe first that every matrix of the form 
$\mathbb{I}+E_{j,j+1}$ is either an element of the form $\matalpha_i(y)$ 
(if the $(j,j+1)$-entry is in one of the blocks $A_{ii}$) or 
$\matgamma_i^+$ (otherwise). If $j < m < \ell$, then computation shows that
$$
(\mathbb{I} + E_{jm})(\mathbb{I} + E_{m\ell})(\mathbb{I}-E_{jm})(\mathbb{I}-E_{m\ell}) = \mathbb{I}+E_{j\ell};
$$
this shows by induction on $|\ell-j|$ that every element of the form $\mathbb{I}+E_{j\ell}$ can be obtained as a product of the given generators. Finally, 
for every invertible $\nu \in \mathbb{F}_q[t]_{\lambda_s}$ and $1 \leq j,\ell \leq k$, observe that
$$
\mathbb{I}+\nu E_{j\ell} = (\mathbb{I} + (\nu-1)E_{jj})(\mathbb{I}+E_{j\ell})(\mathbb{I}+(\nu^{-1}-1)E_{jj}),
$$
where the second factor 
is obtained as above, and the first and 
third factors are obtained as explained in Proposition \ref{PropMuDiagonal}. If $\nu$ is not invertible, then it is the sum of two invertible elements 
$\nu_1,\nu_2$ (for example $1$ and $\nu-1$), so
\[
\mathbb{I}+\nu E_{j\ell} = (\mathbb{I}+\nu_1 E_{j\ell}) (\mathbb{I}+\nu_2 E_{j\ell}),
\]
and both terms in the product can be obtained as above.

The argument is the same for the elements with $j>\ell$.
\end{proof}

\vspace{2mm}
\no {\bf Proof of Theorem \ref{ThmGensUnipotent}. }
Let $Y = (y_{j\ell})$ be an element of the unit group 
$\Ft^*$ (so $Y$ is $k\times k$ as in (\ref{kbyk})). We can suppose without loss of generality that $y_{11}$ is invertible, otherwise we switch columns by multiplying on the right by permutation matrices on the block $A_{11}$. Observe that
$$
\prod_{j=2}^k (\mathbb{I}-(y_{j1}y_{11}^{-1})E_{j1}) \cdot Y \cdot \prod_{j=2}^k (\mathbb{I}-(y_{1j}y_{11}^{-1})E_{1j}) = \left( \begin{array}{c|ccc} y_{11} & \cdots & 0 & \cdots \\ \hline \vdots & \ddots & & \\ 0 & & Y' & \\ \vdots & & & \ddots \end{array}\right).
$$
We can now apply the same argument to the matrix $Y'$ and proceed inductively 
until we get a diagonal matrix. Thus $Y$ is a product of diagonal matrices, 
matrices of the form $\matalpha_i(y)$ where $y$ is a permutation matrix in 
$\GL_{l_i}(q)$, and matrices of the form $\mathbb{I} + \mu E_{j\ell}$ 
where $\mu \in \mathbb{F}_q[t]_{\lambda_s}$ for some $s$. 
As shown in the previous results, each of these matrices
can be written as a product of the generators listed in the 
theorem. $\;\;\;\;\Box$

\subsection{Centralizer of a general element}
Using Theorem \ref{Card} together with   (\ref{SplitCentr}), we deduce the following result giving the structure of the centralizer of an arbitrary element of 
 $\GL_n(q)$. 
 
\begin{theorem}   \label{GLCentralizer}
Let $x \in \GL_n(q)$ have minimal polynomial 
$\prod_{i=1}^hf_i(t)^{e_i}$, where 
$f_1,\ldots,f_h\in \mathbb{F}_q[t]$ are distinct and irreducible, and $\deg f_i=d_i$. 
Let the Jordan form of $x$ be 
$$
\bigoplus_{i=1}^h ( J_{f_i, \lambda_{i1}}^{[l_{i1}]} 
\oplus \cdots \oplus  J_{f_i, \lambda_{ik_i}}^{[l_{ik_i}]}),
$$
where $\lambda_{i1} < \cdots < \lambda_{ik_i} = e_i$ for each $i$. 
Then $C_G(x) = U \rtimes R$, where 
$R \cong \prod_{i=1}^h \left( \prod_{j=1}^{k_i} 
\GL_{l_{ij}}(q^{d_i})\right)$ and $|U| = q^{\m_gamma}$ with 
$$
\m_gamma = \sum_{i=1}^h d_i\left( 2\sum_{a<b} \lambda_{ia}l_{ia}l_{ib} + \sum_j (\lambda_{ij}-1)l_{ij}^2\right).
$$
\end{theorem}

\section{Conjugacy and centralizers in $\SL_{\MakeLowercase{n}}(\MakeLowercase{q})$} \label{SLconjclasses}

For completeness we include the solutions of the main conjugacy 
problems (1)-(3) of Section \ref{mainprob}  for the special linear group. 
In particular, we construct a generating set for the centralizer 
in $\SL_n(q)$ of an arbitrary element, 
since this is not immediate from our solution for $\GL_n(q)$. 

\begin{theorem} \label{MainResult2}
Let $x \in \SL_n(q)$ have minimal polynomial 
$\prod_{i=1}^hf_i(t)^{e_i}$, where 
$f_1,\ldots,f_h\in \mathbb{F}_q[t]$ are distinct and irreducible, and $\deg f_i=d_i$.
Let the Jordan form of $x$ be 
\begin{equation}\label{jf}
\bigoplus_{i=1}^h ( J_{f_i, \lambda_{i1}}^{[l_{i1}]} 
\oplus  \cdots \oplus  J_{f_i, \lambda_{ik_i}}^{[l_{ik_i}]}),
\end{equation}
where 
$\lambda_{i1} < \cdots < \lambda_{ik_i} = e_i$ for each $i$. Set $r = \mathrm{gcd}(\lambda_{11}, \dots, \lambda_{hk_h}, q-1)$.
Then
\[
|C_{\GL_n(q)}(x):C_{\SL_n(q)}(x)| = \frac{q-1}{r}.
\]
The conjugacy class $x^{\GL_n(q)}$ splits into $r$ classes in $\SL_n(q)$, with representatives $x,x^z,\ldots,x^{z^{r-1}}$, where 
$z \in \GL_n(q)$ is a fixed element of 
determinant $\o$ $($a primitive element of $\F_q)$.
\end{theorem}

\begin{proof} 
We first deal with the case where $h=1$. To simplify notation, write $f = f_1$ of degree $d$, and let the distinct elementary divisors of $x$ be  $f^{\l_1},\ldots,f^{\l_k}$, with multiplicities $l_1,\ldots,l_k$, so that $\sum_i \l_i l_i = n/d$. If $x = su=us$ is the Jordan decomposition of $x$ with $s$ semisimple and $u$ unipotent, then by Proposition \ref{largerfield}, $C_{\GL_n(q)}(s) = \GL_e(q^d)$, where $de=n$. Then $u \in \GL_e(q^d)$ is unipotent with Jordan form $\bigoplus_i J_{\l_i}^{[l_i]}$. By Theorem \ref{Card},
\[
C_{\GL_n(q)}(x) = C_{\GL_e(q^d)}(u) = UR,
\]
where $R \cong \prod_{i=1}^k \GL_{l_i}(q^d)$. 
If we view $V$ as $V_e(q^d)$, then the action of $R$ on $V$ is as 
$$\prod_{i=1}^k (\GL_{l_i}(q^d) \otimes I_{\l_i}).$$ 
Let $\nu$ be a primitive element of $\F_{q^d}$ such that $\NN(\nu) = \o$, where $\NN$ is the norm map $\F_{q^d} \to \F_q$. 
If $A \in \GL_{l_i}(q^d)$ has determinant $\nu$, then the 
determinant of $A \otimes I_{\l_i}$ as an element of $\GL_e(q^d)$ is $\nu^{\l_i}$, and so its determinant as an element of $\GL_n(q)$ is $\NN(\nu^{\l_i})$ (see \cite[(4.3.13)]{KL}), which is $\o^{\l_i}$. 

Thus the subgroup of $\F_q^*$ consisting of all determinants of elements of $R$ is generated by $\o^{\l_i}$ for $i=1,\ldots,k$. Hence it equals $\la \o^r\ra$ 
where $r  = \mathrm{gcd}(\lambda_{1}, \ldots, \lambda_{k}, q-1)$. 
Thus the image of the determinant map $C_{\GL_n(q)}(x) \to \F_q^*$ has order $(q-1)/r$. This completes the proof when $h=1$.

\vspace{2mm}
Now consider the general case. We know from (\ref{SplitCentr}) that $C_{\GL(V)}(x) = \prod_{i=1}^h C_{\GL(V_i)}(x_i)$, where $V_i = \ker \,f_i(x)^{e_i}$ are the generalized eigenspaces of $x$. By the $h=1$ case, the image of the determinant map on $C_{\GL(V_i)}(x_i)$ is $\la \o^{r_i}\ra$, where $r_i  = \mathrm{gcd}(\lambda_{i1}, \ldots, \lambda_{ik_i}, q-1)$. It follows that the 
image of the determinant map on $C_{\GL(V)}(x)$ is $\la \o^{r}\ra$, where $r  = \mathrm{gcd}(r_1,\ldots,r_h)$. The conclusion follows.
\end{proof}

We now use this theorem to decide conjugacy and compute conjugating elements in $S:=\SL_n(q)$.
First, as conjugacy class representatives in $S$, for each Jordan form $J \in S$ as in (\ref{jf}), we choose representatives 
\[
J,J^z,\ldots,J^{z^{r-1}}, 
\]
where  $r = \mathrm{gcd}(\lambda_{11}, \dots, \lambda_{hk_h}, q-1)$ and 
$z \in \GL_n(q)$ is a fixed element of determinant $\o$ (a primitive element of $\F_q$). 

Given $x \in S$ with Jordan form $J$ as in (\ref{jf}), we find its
class representative by computing $g \in \GL_n(q)$ such that $J = x^g$; then $x$ is $S$-conjugate to $J^{z^i}$, where $\det(g) \in \o^{-i}\la\o^r\ra$.

To decide whether $x,y \in S$, both with Jordan form (\ref{jf}), are $S$-conjugate, we compute $g \in \GL_n(q)$ such that $y=x^g$; then $x$ and $y$ are $S$-conjugate if and only if $\det (g) \in \la \o^r\ra$. 

Finally, given $x \in S$ that is $S$-conjugate to a representative $J^{z^i}$, we find a conjugating element as follows. First, compute $g \in \GL_n(q)$ such that $x^g = J$. Then $\det (g) \in \o^{-i}\la\o^r\ra$, say $\det (g) = \o^{-i+kr}$. From the proof of 
Theorem \ref{MainResult2}, we can find $c \in C_{\GL_n(q)}(J)$ of determinant $\o^{-kr}$ with $c$ a product of scalars in the factors $\GL_{l_{ij}}(q^{d_i})$ of the subgroup $R$ of $C(J)$. Then $x^{gcz^i} = J^{z^i}$, and $gcz^i \in S$.

\subsubsection{Generators for centralizer in $\SL_n(q)$}
Write $G = \GL_n(q)$ and $S = \SL_n(q)$, and let $x \in S$ be as in (\ref{jf}). 
Using the notation of 
Theorem \ref{GLCentralizer}, $C_G(x) = U \rtimes R$, and so 
$C_S(x) = U \rtimes (R \cap S)$, since every element 
of $U$ has determinant 1. So our task is to compute the elements of $R$ of determinant 1.

Recall Theorem \ref{GLCentralizer}: there is a group isomorphism 
$\varphi$ between $\prod_{i=1}^h \left(\prod_{j=1}^{k_i} \GL_{l_{ij}}(q^{d_i}) \right)$ 
and $R$. We realise this isomorphism as follows.
Map a sequence 
$[\myx_{ij}: \, 1 \leq i \leq h, \, 1 \leq j \leq k_i]$ 
with $\myx_{ij} \in \GL_{l_{ij}}(q^{d_i})$ to the block diagonal matrix
\begin{eqnarray}          
\label{homProdGLintoR}
\begin{pmatrix} \ddots && \\ & \myX_{ij} & \\ && \ddots \end{pmatrix},
\end{eqnarray} 
where $\myX_{ij}$ is defined as follows: if $\myx_{ij} = (\alpha_{\mu\nu})$, 
then $\myX_{ij}$ is the matrix obtained by substituting each 
entry $\alpha_{\mu\nu}$ of $\myx_{ij}$ by the embedding into 
$\GL_{\lambda_{ij}d_i}(q)$ of a $\lambda_{ij}$-dimensional scalar 
matrix with diagonal $\alpha_{\mu\nu}$.
Note that if $\myx \in R$ is such that $\myx = \varphi(\myx_{11}, \dots, \myx_{hk_h})$, 
where $\myx_{ij} \in \GL_{l_{ij}}(q^{d_i})$ corresponds to a Jordan block of 
dimension $\lambda_{ij}$, then 
\begin{eqnarray}  \label{eqprod}
\det(\myx) = \prod_{i,j} \NN_i(\det({\myx_{ij}}))^{\lambda_{ij}},
\end{eqnarray}
where $\NN_i: \mathbb{F}_{q^{d_i}} \rightarrow \mathbb{F}_q$ is the norm map. 
Let $\omega$ be a primitive element of $\mathbb{F}_q$. 
Then $\NN_i(\det ({\myx_{ij}})) = \omega^{a_{ij}}$ for some 
$a_{ij} \in \mathbb{Z}_{q-1}$, the ring of residue classes modulo $q-1$, 
so (\ref{eqprod}) becomes
$
\det (\myx) = \prod_{i,j}\omega^{a_{ij}\lambda_{ij}}.
$
The condition $\det(\myx)=1$ is equivalent to
$
\sum_{i=1}^h \sum_{j=1}^{k_i} a_{ij}\lambda_{ij} = 0.
$

Let $V_i$ be the generalized eigenspace $\ker\,f_i(x)^{e_i}$, so that 
$x = x_1 \oplus \cdots \oplus x_h$, where $x_i$ is the restriction of $x$ to $V_i$.
 Every $y \in C_S(x)$ can also be written as $y = y_1 \oplus \cdots \oplus y_h$,
but the centralizer in $S$ of $x$ is not necessarily the direct product of the centralizers of the $x_i$. 
We list generators for $C_S (x)$ in two steps.

\begin{enumerate}
\item[1.] For each $i$, we generate $C_{S_i}(x_i)$, where $S_i = \SL(V_i) \cong \SL_{n_id_i}(q)$ and $n_i = \sum_j l_{ij}\lambda_{ij}$. 
Putting them together, we generate a subgroup 
\begin{eqnarray} \label{productOfCS}
\prod_{i=1}^h C_{S_i}(x_i)
\end{eqnarray}
of $C_S(x)$.
\item[2.] 
To complete a generating set for $C_S(x)$, we add elements 
$$
\matH_\mys = \bigoplus_{i,j} \matH_{{\mys}ij}, \, 1 \leq \mys \leq \myr,
$$
where each $\matH_{{\mys}ij}$ commutes with $J_{f_i,\lambda_{ij}}^{[l_{ij}]}$, and together with the subgroup (\ref{productOfCS}),
$\matH_1, \dots, \matH_\myr$ generate the kernel of the 
determinant map $R \to \F_q^*$.
\end{enumerate}

We now describe how to write down generators for each of the above steps: call $\mathcal{S}_i$ the set of generators in Step $i$ for $i=1,2$. 

We first discuss how to construct $\mathcal{S}_2$. 
For every $1 \leq i \leq h$ and $1 \leq j \leq k_i$, 
take $\widetilde{\matH}_{ij} \in \GL_{l_{ij}}(q^{d_i})$ having 
determinant of order $q^{d_i}-1$ and such that its image $\matH_{ij}$ 
under application 
of the homomorphism $\varphi$ defined in (\ref{homProdGLintoR}) has 
determinant  $\omega^{\lambda_{ij}}$.
Let $k = \sum_{i=1}^h k_i$ and let $Z(0)$ be the subgroup of $\mathbb{Z}_{q-1}^k$ 
consisting of the solutions $(x_{11}, \dots, x_{hk_h})$ of the equation 
$\sum_{i,j} x_{ij}\lambda_{ij}=0$ in $\mathbb{Z}_{q-1}$. 
Compute a generating set $\{ (a_{{\mys}11}, \dots, a_{{\mys}hk_h}) \, : 
\, 1 \leq \mys \leq \myr\}$ for $Z(0)$. Then $\mathcal{S}_2$ consists of 
all the matrices
$$
\matH_\mys = \bigoplus_{i,j} \matH_{ij}^{a_{{\mys}ij}}, \quad 1 \leq \mys \leq \myr.
$$

We construct $\mathcal{S}_1$ by describing separately 
each $C_{S_i}(x_i)$ for every $i$.
Let $d_i = \deg{f_i}$ and let $n_i = \sum_{j=1}^{k_i} \lambda_{ij}l_{ij}$. 
Let $x_i$ be the embedding into 
$\GL_{n_id_i}(q)$ of some $\w{x}_i \in \GL_{n_i}(q^{d_i})$ and let 
$C_{\GL_{n_id_i}(q)}(x_i)$ be the embedding into $\GL_{n_id_i}(q)$ of 
$C_{\GL_{n_i}(q^{d_i})}(\w{x}_i)$. 
If $\w{x}_i = \w{s}_i\w{u}_i$ is the Jordan decomposition of $\w{x}_i$, 
then $\w{s}_i$ is a scalar matrix, so 
$C_{\GL_{n_i}(q^{d_i})}(\w{x}_i) = C_{\GL_{n_i}(q^{d_i})}(\w{u}_i)$. 
So we reduce our task to computing the centralizer of a unipotent element with 
Jordan form $\bigoplus_{j=1}^{k_i} J_{\lambda_{ij}}^{[l_{ij}]}$. 
For every $y \in C_{\GL_{n_id_i}(q)}(x_i)$, if $\w{y}$ is the 
corresponding element 
in $C_{\GL_{n_i}(q^{d_i})}(\w{x}_i)$, then 
$\det(y) = \NN_i(\det(\w{y}))$. Hence $\det(y)=1$ if and only if 
$\det(\w{y})$ has determinant a power of $\omega_i^{q-1}$, where 
$\omega_i$ is a primitive element of $\mathbb{F}_{q^{d_i}}$.
It follows that the construction of a generating set for $C_{S_i}(x_i)$ reduces to 
listing the generators of the centralizer of a unipotent element 
in a subgroup of a specific index in $\GL_{n_i}(q^{d_i})$. We do this by modifying 
Theorem \ref{ThmGensUnipotent}.

To state this modification, we revert to the notation of Theorem \ref{ThmGensUnipotent}. 
Thus let $x = \bigoplus_{i=1}^h J_{\l_i}^{[l_i]}$ be a unipotent element of $\GL_n(q)$
with corresponding partition $\l$ of $n$, and recall that we have an isomorphism between
the centralizer $C_{\GL_n(q)}(x)$ and the unit group $\mathbb{F}_q[t]_{\lambda}^*$.

\begin{theorem} \label{ThmGensUnipotentSL}
Let $\omega$ be a primitive element of $\mathbb{F}_{q}$ and 
let $d$ be a divisor of $q-1$. 
A generating set for the subgroup of $\mathbb{F}_q[t]_{\lambda}^*$ of matrices 
whose determinant is a power of $\omega^d$ consists of the following:
\begin{enumerate}
\item[\rm (1)] the set of $\matalpha_i(y)$, where $1 \leq i \leq h$ and 
$y$ runs 
over a generating set for the unique subgroup of $\GL_{l_i}(q)$ of index $d$;
\item[\rm (2)] the set of $\matbeta_i(\mu)$, where $1 \leq i \leq h$ and 
$\mu$ runs over all elements of the set 
$$
\{ 1+\omega^jt^l: 0 \leq j \leq [\mathbb{F}_q:\mathbb{F}_p], 1 \leq l \leq \lambda_i-1\};
$$
\item[\rm (3)] the set of all $\matgamma_i^+(\xi)$ and $\matgamma_i^-(\xi)$, 
where $1 \leq i \leq h-1$ and 
$\xi \in \{1, \omega, \dots, \omega^{[\mathbb{F}_q:\mathbb{F}_p]-1} \};$
\item[\rm (4)] 
matrices which, together with the elements in $(1)-(3)$ above,
generate the kernel of the determinant map $R \to \F_q^*$.
\end{enumerate}
\end{theorem}

\begin{proof}
The proof is similar to that of Theorem \ref{ThmGensUnipotent}. 
Using elements of the form $\matalpha_i(y)$ and $\matbeta_i(\mu)$ described above, 
we can get in every block $A_{ii}$ a matrix whose determinant is a 
power of $\omega^d$. 
To get diagonal elements, here we cannot conjugate by permutation 
matrices (their determinant is not a power of $\omega^d$ in general), 
but we can conjugate by a monomial matrix of determinant 1.

The elements $\matgamma_i^+$ and $\matgamma_i^-$ are not sufficient 
to repeat the argument of Theorem \ref{ThmGensUnipotent}, because we do 
not have diagonal matrices to get all elements of the form 
$\mathbb{I}+\nu E_{j\ell}$. But for every $1 \leq j, \ell \leq h$, 
$$
(\mathbb{I} + \xi_1 E_{j\ell})(\mathbb{I} + \xi_2 E_{j\ell}) = 
\mathbb{I} + (\xi_1+\xi_2)E_{j\ell}
$$
for every $\xi_1 \in \mathbb{F}_q[t]_{\lambda_{s_1}}$ and 
$\xi_2 \in \mathbb{F}_q[t]_{\lambda_{s_2}}$ for appropriate $s_1$ and $s_2$. 
So we take 
$\matgamma_i^+(\xi)$ and $\matgamma_i^-(\xi)$ for 
enough $\xi$ to generate the additive group of $\mathbb{F}_q$.

Finally, using the matrices listed at (1)-(3) 
we can generate every matrix where each 
block $A_{ii}$ has determinant a power of $\omega^d$. 
But we need additional matrices satisfying the weaker condition that 
$\prod_i \det(A_{ii})$ is a power of $\omega^d$. 
For these we add elements listed at (4); 
these are obtained using the method described earlier for $\mathcal{S}_2$. 
\end{proof}

\section{Some examples}
We illustrate some of the results of this chapter with the following examples.
Matrices are listed in Jordan form and we use the notation of the earlier 
sections throughout.

\begin{example}
We list the conjugacy classes of $\GL_3(2)$. If $x \in \GL_3(2)$, 
then its characteristic polynomial has degree 3 and it is a 
product of irreducible polynomials of degree at most 3 in $\mathbb{F}_2[t]$. 
These are $t+1$, $t^2+t+1$, $t^3+t+1$ and $t^3+t^2+1$. 
Every elementary divisor is 
a power of one of these polynomials and the sum of the degrees must be 3. 
In Table \ref{ExTable}, for each class representative,
we record its elementary divisors, its Jordan form, 
and its minimal and characteristic polynomials. 

\begin{table}[thbp]\caption{Conjugacy classes in $\GL_3(2)$}\label{ExTable}
\begin{tabular}{|c|l|l|l|}
\hline
List of elementary divisors & Representative & Minimal polynomial & Characteristic polynomial \\
\hline
$\begin{array}{l}
t+1,\\ t+1,\\ t+1 \end{array}$ & $\left( \begin{array}{ccc} 1 & 0 & 0\\ 0 & 1 & 0\\ 0 & 0 & 1 \end{array}\right)$ & $t+1$ & $(t+1)^3$\\
\hline
$\begin{array}{l}
t+1,\\ (t+1)^2 \end{array}$ & $\left( \begin{array}{ccc} 1 & 0 & 0\\ 0 & 1 & 1\\ 0 & 0 & 1 \end{array}\right)$ & $(t+1)^2$ & $(t+1)^3$\\
\hline
$(t+1)^3$ & $\left( \begin{array}{lll} 1 & 1 & 0\\ 0 & 1 & 1\\ 0 & 0 & 1 \end{array}\right)$ & $(t+1)^3$ & $(t+1)^3$\\
\hline
$\begin{array}{l}
t+1,\\ t^2+t+1 \end{array}$ & $\left( \begin{array}{ccc} 1 & 0 & 0\\ 0 & 0 & 1\\ 0 & 1 & 1 \end{array}\right)$ & $(t+1)(t^2+t+1)$ & $(t+1)(t^2+t+1)$\\
\hline
$t^3+t+1$ & $\left( \begin{array}{lll} 0 & 1 & 0\\ 0 & 0 & 1\\ 1 & 1 & 0 \end{array}\right)$ & $t^3+t+1$ & $t^3+t+1$\\
\hline
$t^3+t^2+1$ & $\left( \begin{array}{lll} 0 & 1 & 0\\ 0 & 0 & 1\\ 1 & 0 & 1 \end{array}\right)$ & $t^3+t^2+1$ & $t^3+t^2+1$\\
\hline
\end{tabular}
\end{table}
\end{example}

\begin{example} Consider $x \in G := \GL_{12}(3)$ where 
$$
x = \left(\begin{array}{cccccccccccc}
2 & 0 & 0 & 0 & 0 & 0 & 0 & 0 & 0 & 0 & 0 & 0 \\ 
0 & 0 & 1 & 0 & 0 & 0 & 0 & 0 & 0 & 0 & 0 & 0 \\ 
0 & 1 & 1 & 0 & 0 & 0 & 0 & 0 & 0 & 0 & 0 & 0 \\
0 & 0 & 0 & 0 & 1 & 0 & 0 & 0 & 0 & 0 & 0 & 0 \\
0 & 0 & 0 & 1 & 1 & 0 & 0 & 0 & 0 & 0 & 0 & 0 \\
0 & 0 & 0 & 0 & 0 & 0 & 1 & 1 & 0 & 0 & 0 & 0 \\
0 & 0 & 0 & 0 & 0 & 1 & 1 & 0 & 1 & 0 & 0 & 0 \\
0 & 0 & 0 & 0 & 0 & 0 & 0 & 0 & 1 & 0 & 0 & 0 \\
0 & 0 & 0 & 0 & 0 & 0 & 0 & 1 & 1 & 0 & 0 & 0 \\
0 & 0 & 0 & 0 & 0 & 0 & 0 & 0 & 0 & 0 & 1 & 0 \\
0 & 0 & 0 & 0 & 0 & 0 & 0 & 0 & 0 & 0 & 0 & 1 \\
0 & 0 & 0 & 0 & 0 & 0 & 0 & 0 & 0 & 1 & 2 & 2 
\end{array} \right).
$$

The minimal polynomial of $x$ is 
$(t - 2) (t^2 - t - 1)^2 (t^3 + t^2 + t - 1).$
The generalized eigenspaces are $V_1 = \ker{(t-2)}$, $V_2 = \ker{((t^2-t-1)^2)}$ 
and $V_3=\ker{(t^3+t^2+t-1)}$. 
Recall that $C_G(x) = U \rtimes R$. 

Since $V_1$ has dimension 1, the corresponding factor in $R$ is $\GL_1(3)$.

The matrix of the restriction of $x$ to $V_2$ is the embedding 
in $\GL_8(3)$ of 
$$
y := \left(\begin{array}{cccc}
\omega & 0 & 0 & 0\\
0 & \omega & 0 & 0\\
0 & 0 & \omega & 1\\
0 & 0 & 0 & \omega
\end{array}\right) \in \GL_4(9),
$$
where $\omega$ is a root of the polynomial $t^2-t-1$ in $\mathbb{F}_9$. 
Observe $y$ has two Jordan blocks of order 1 and one block of order 2. 
An element in $\GL_4(9)$ centralizing $y$ has the form
$$
\left(\mbox{\begin{tabular}{c|c}
$A$ & $\begin{array}{cc} 0 & * \\ 0 & * \end{array}$\\ \hline
$\begin{array}{cc} * & * \\ 0 & 0 \end{array}$ & $X$
\end{tabular}}\right),
$$
where $A \in \GL_2(9)$ and $X$ is a matrix of the form 
$b_0X_{2\times 2}^0+b_1X_{2\times 2}^1$ with $b_0,b_1 \in \mathbb{F}_9$ 
and $b_0 \neq 0$ (using the notation after (\ref{blocksBij})). 

We can now use the descriptions in Sections \ref{embedding} and \ref{Cue} to 
write down explicit generators for each of $R$ and $U$.
Now $C_G(x) = U \rtimes R$ where 
$
R \cong \GL_1(3) \times \GL_2(9) \times \GL_1(9) \times \GL_1(27)
$
and $U$ has order $3^{\gamma}$ where
$$
\gamma = \sum_{i=1}^h d_i\left( 2\sum_{j<\ell} \lambda_{i,j}l_{i,j}l_{i,\ell} + \sum_j (\lambda_{i,j}-1)l_{i,j}^2\right) = 2[2(1 \cdot 2 \cdot 1) + (2-1)1^2] = 10.
$$
\end{example}

\begin{example}
Consider $x \in \SL_8(3)$ where 
$$
x = \left( \begin{array}{cccccccc} 
0 & 1 & 1 & 0 & 0 & 0 & 0 & 0 \\
1 & 2 & 0 & 1 & 0 & 0 & 0 & 0 \\
0 & 0 & 0 & 1 & 0 & 0 & 0 & 0 \\
0 & 0 & 1 & 2 & 0 & 0 & 0 & 0 \\
0 & 0 & 0 & 0 & 0 & 1 & 1 & 0 \\
0 & 0 & 0 & 0 & 1 & 2 & 0 & 1 \\
0 & 0 & 0 & 0 & 0 & 0 & 0 & 1 \\
0 & 0 & 0 & 0 & 0 & 0 & 1 & 2  
\end{array}\right).
$$
Its minimal polynomial is $(t^2+t-1)^2$. 
As usual, $C_{\GL_8(3)}(x) = U \rtimes R$ where $R \cong \GL_2(9)$. 
Every element of $R$ can be identified with some $z_{1,1} \in \GL_2(9)$ and 
its determinant is $\NN(\det\,{z_{1,1}})^2$, where $\NN$ is the norm map
$\mathbb{F}_9 \rightarrow \mathbb{F}_3$. For example, if $\omega$ is a 
root of $t^2+t-1$ in $\mathbb{F}_9$ and 
$$z_{1,1} = \left(\begin{array}{cc} \omega & 1\\ \omega+2 & 0 
\end{array}\right),$$ then the corresponding element of $R$ is
$$
\left( \begin{array}{cccccccc} 0 & 1 & 0 & 0 & 1 & 0 & 0 & 0\\
1 & 2 & 0 & 0 & 0 & 1 & 0 & 0\\
0 & 0 & 0 & 1 & 0 & 0 & 1 & 0\\
0 & 0 & 1 & 2 & 0 & 0 & 0 & 1\\
2 & 1 & 0 & 0 & 0 & 0 & 0 & 0\\
1 & 1 & 0 & 0 & 0 & 0 & 0 & 0\\
0 & 0 & 2 & 1 & 0 & 0 & 0 & 0\\
0 & 0 & 1 & 1 & 0 & 0 & 0 & 0 \end{array}\right),
$$
and has determinant $\NN(\det\,{z_{1,1}})^2 = \NN(1-\omega)^2=1$.
\end{example}

\begin{example}
Consider $x \in \SL_6(13)$ where 
$$
x = \left(\begin{array}{cccccc} 1 & 1 & 0 & 0 & 0 & 0\\ 0 & 1 & 1 & 0 & 0 & 0\\ 0 & 0 & 1 & 0 & 0 & 0 \\
0 & 0 & 0 & 1 & 1 & 0\\ 0 & 0 & 0 & 0 & 1 & 1\\ 0 & 0 & 0 & 0 & 0 & 1 \end{array}\right).
$$
Its minimal polynomial is $(t - 1)^3$.
Here $x$ has two Jordan blocks of order 3, and $d = \gcd(3,3,13-1) = 3$. 
Thus, the conjugacy class of $x$ in $\GL_6(13)$ splits into 
three classes in $\SL_6(13)$. Observe 2 is a primitive element 
of $\mathbb{F}_{13}$. For $i = 0, 1, 2$, 
$$
y_i := \left(\begin{array}{cccc} 1 & & & \\ & \ddots & & \\ & & 1 & \\ & & & 2^i \end{array}\right).
$$
Hence, the representatives for the three classes in $\SL_6(13)$ are
$$
x^{y_0}=\left(\begin{array}{cccccc} 1 & 1 & 0 & 0 & 0 & 0\\ 0 & 1 & 1 & 0 & 0 & 0\\ 0 & 0 & 1 & 0 & 0 & 0 \\
0 & 0 & 0 & 1 & 1 & 0\\ 0 & 0 & 0 & 0 & 1 & 1\\ 0 & 0 & 0 & 0 & 0 & 1 \end{array}\right), \, x^{y_1} = \left(\begin{array}{cccccc} 1 & 1 & 0 & 0 & 0 & 0\\ 0 & 1 & 1 & 0 & 0 & 0\\ 0 & 0 & 1 & 0 & 0 & 0 \\
0 & 0 & 0 & 1 & 1 & 0\\ 0 & 0 & 0 & 0 & 1 & 2\\ 0 & 0 & 0 & 0 & 0 & 1 \end{array}\right), \, x^{y_2} = \left(\begin{array}{cccccc} 1 & 1 & 0 & 0 & 0 & 0\\ 0 & 1 & 1 & 0 & 0 & 0\\ 0 & 0 & 1 & 0 & 0 & 0 \\
0 & 0 & 0 & 1 & 1 & 0\\ 0 & 0 & 0 & 0 & 1 & 4\\ 0 & 0 & 0 & 0 & 0 & 1 \end{array}\right).
$$
\end{example}

\chapter{Preliminaries on classical groups}   \label{LinearChapter}

In this chapter we first define the symplectic, orthogonal and unitary groups. 
Building on the work of Britnell \cite{JBRIT} and Milnor \cite{Milnor}, we 
then give in Theorem \ref{elementsofCG} 
 necessary and sufficient criteria for 
$X \in \GL(V)$ to preserve a non-degenerate alternating, 
hermitian or quadratic form on $V$ in terms of the elementary divisors of $X$. 

\section{The finite classical groups}\label{fincla}

Let $F = \F_{q^u}$\index{$\F_{q^u}$}, where $u$ is 1 or 2, and 
let $\bar \l = \l^q$ for $\l \in F$, so that 
$\l \mapsto \bar \l$ is a field automorphism of order $u$. Let $V= F^n$ be an 
$n$-dimensional vector space over $F$. For $X = (x_{ij}) \in \GL_n(F)$, define $\bar X = (\bar x_{ij})$, and set $X^* = \bar X^\tr$\index{$X^*$}, the transpose of $\bar X$. We also regard $X$ as an element of $\GL(V)$, namely the map $v \mapsto Xv$; similarly for $\bar X$ and $X^*$.

We now define the classical 
symplectic, orthogonal and unitary groups as isometry groups of sesquilinear forms\index{form!sesquilinear} $\b: V\times V \to F$ or 
quadratic forms\index{form!quadratic} $Q:V \to F$. 
We first recall a few notions about such forms. 
We consider sesquilinear forms 
$\b: V\times V \to F$ such that 
\begin{itemize}
    \item $\b$ is left-linear: $\b(au_1+bu_2,v) = a\b (u_1,v)+b\b(u_2,v)$ for $a,b \in F$, $u_i,v \in V$, and
    \item $\b(v,u) = \overline{\b(u,v)}$ for $u,v \in V$,
\end{itemize}
and quadratic forms \index{form!quadratic} $Q:V \to F$ such that 
\begin{itemize}
    \item $Q(av) = a^2Q(v)$ for $a \in F$, $v \in V$, and 
    \item the function $\b_Q(u,v) = Q(u+v)-Q(u)-Q(v)$ is a bilinear form on $V$.
\end{itemize}

For a basis $B = \{v_1,\ldots,v_n\}$ of $V$, the matrix of a sesquilinear form $\b$ with respect to $B$ is the $n\times n$ matrix $\b_B = (\b(v_i,v_j))$, and the matrix of a quadratic form $Q$ is $A = (a_{ij})$, where $a_{ii} = Q(v_i)$, $a_{ij}=\b_Q(v_i,v_j)$ for $i<j$, and $a_{ij} = 0$ for $i>j$. Note that $Q(v) = vAv^\tr$ (writing vectors relative to the basis $B$), and $A+A^\tr$ is the matrix of $\b_Q$. 
We say that $g \in \GL(V)$ 
is an  {\it isometry}\index{isometry} of $\b$ if $\b(ug,vg) = \b(u,v)$ for all $u,v \in V$. 
Writing matrices with respect to a basis $B$, if $A = \b_B$ is the matrix 
of $\b$, then $X \in \GL_n(F)$ is an isometry if and only 
if $XAX^* = A$. Similarly, $g \in \GL(V)$ is an isometry of a quadratic from $Q$ if $Q(vg) = Q(v)$ for all $v \in V$.

Two sesquilinear forms $\beta_1$, $\beta_2$ on $V$ are \textit{congruent}\index{forms!congruent} if there exists $T \in \GL(V)$ such that $\beta_1(uT,vT) = \beta_2(u,v)$ for all $u,v \in V$. 
If $A_1$ and $A_2$ are the matrices of $\beta_1$ and $\beta_2$ respectively, 
then this condition is equivalent to the existence of 
$S \in \GL_n(F)$ such that $SA_1\bar S^\tr = A_2$. Similarly, quadratic forms $Q_1$, $Q_2$ are \textit{congruent} if there exists $T \in \GL(V)$ such that $Q_1(vT) = Q_2(v)$ for all $v \in V$. 

A sesquilinear form $\b$ is {\it non-degenerate}\index{non-degenerate} 
if its matrix $\b_B$ with respect to every basis $B$ is invertible; and a quadratic form $Q$ is non-degenerate if the associated bilinear form $\b_Q$ is invertible. A subspace $U$ of $V$ is {\it totally singular}\index{totally singular}
with respect to $\b$ (resp.\ $Q$) if $\b$ (resp.\ $Q$) vanishes on 
restriction to $U$. 
The {\it orthogonal complement}\index{orthogonal complement} $U^\perp$ of  
$U$ with respect to $\b$ (resp.\ $Q$) is 
$\{v \in V : \b(u,v) = 0 \;\forall u \in U\}$
(resp.\ $\{v \in V : \b_Q(u,v) = 0 \;\forall u \in U\}$).

We now define the classical groups; all assertions can be found 
in \cite{PG}, for example.

\begin{itemize}
\item {\it Symplectic group}\index{group!symplectic} $\Sp(V)$\index{$\Sp(V)$}: 
here $u=1$, $n$ is even  and $\Sp(V)$ is the isometry group 
of a non-degenerate alternating bilinear form $\b$ on $V$. Up to 
congruence there is a unique such form. The corresponding matrix 
group is denoted $\Sp_n(q)$.
\item {\it Unitary group} \index{group!unitary} $\GU(V)$\index{$\GU(V)$}: 
here $u=2$,  
and $\GU(V)$ is the isometry group of a non-degenerate hermitian 
form $\b$ on $V$. Up to congruence there is a unique such form. 
The corresponding matrix group is denoted $\GU_n(q)$.
\item {\it Orthogonal group}\index{group!orthogonal} 
$\Or(V)$\index{$\Or(V)$}: here $u=1$, and $\Or (V)$ is the isometry group of 
a non-degenerate quadratic form $Q$ on $V$. 
The associated symmetric bilinear form is $\b_Q$, as defined above.
If $q$ is odd, then $\Or (V)$ is also the isometry group of $\b_Q$. Up to congruence there are two non-degenerate quadratic forms on $V$. If $n$ is odd, then $q$ is odd and the forms are $Q$ and $\l Q$, where $\l$ is a non-square in $F$. 
When $n=2k$, the congruence classes are distinguished by the dimension, 
$k$ or $k - 1$, of a maximal totally singular subspace: 
if the dimension is $k$, 
then the orthogonal group is denoted $\Or^+(V)$ \index{$\Or^+(V)$} 
or $\Or^+_{2k}(q)$;  if the 
dimension is $k-1$, then it is denoted $\Or^-(V)$ \index{$\Or^-(V)$} 
or $\Or^-_{2k}(q)$. We denote the sign \index{sign of quadratic form} $\pm$ by $\hbox{sgn}(Q)$, and say that $Q$ is of {\it plus} or {\it minus type} accordingly.\index{plus type}\index{minus type} For $n=2k$ and $q$ odd, $\hbox{sgn}(Q)$ is also determined by the {\it discriminant} \index{discriminant} $D(Q)$, which is the determinant of the matrix of $\b_Q$ modulo the squares $(F^*)^2$: namely, if $D(Q)$ is a square (resp., a non-square), then  $\hbox{sgn}(Q) = (-1)^{k(q-1)/2}$ (resp., $(-1)^{k(q+1)/2}$) -- see \cite[Prop.\ 2.5.10]{KL}.
\end{itemize}

We often denote the isometry group of the sesquilinear 
form $\b$ by the symbol $\C(\b)$\index{$\C(\b)$}, and the isometry group of 
the quadratic form $Q$ by $\C(Q)$\index{$\C(Q)$}, sometimes 
replacing the form by its associated matrix.

The {\it special} \index{special group} unitary and orthogonal groups are the intersections of 
the isometry groups with $\SL(V)$, and denoted by 
$\SU(V) = \GU(V) \cap \SL(V)$\index{$\SU(V)$} 
and $\SO(V)=\Or(V) \cap \SL(V)$\index{$\SO(V)$}. 
Finally, $\SO(V)$ has a subgroup 
$\O(V)$ \index{$\O(V)$} 
of index 2 which is the kernel of the 
{\it spinor norm}\index{spinor norm}
 map $\SO(V) \to C_2$ (see \cite[\S 2.5]{KL}). 

Generating sets of size 2 are recorded in \cite{DET} for 
the linear, symplectic, and unitary groups; 
generating sets of size at most 4 are recorded in \cite{RT} 
for the orthogonal groups.  
These {\it standard generators}\index{standard generators}, 
written with respect to a fixed form, 
are used to define (the standard copies of) the groups in {\sc Magma}.  

\section{Membership of classical groups}   \label{membershipChapter}

Let $\C = \C(\b)$ or $\C(Q)$ be a classical group on $V$, as defined in the previous section. The first step in classifying the  conjugacy classes of $\C$ is to determine which classes of $\GL(V)$ have elements in $\C$, and in this section we do this. The main result is Theorem \ref{elementsofCG}. Our approach is modelled on those of Britnell \cite[\S 5.1]{JBRIT} and Milnor 
\cite[\S 3]{Milnor}.

Recall that $F = \F_{q^u}$, where $u$ is 1 or 2, and 
$\l \mapsto \bar \l$ is a field automorphism of order $u$ 
(where $\bar \l = \l^q$). 
For a monic polynomial 
$f(t) = t^d+a_{d-1}t^{d-1}+\cdots +a_0 \in F[t]$ with $a_0\ne 0$, 
define $\bar f(t) \index{$\bar f(t)$} = 
 t^d+\bar a_{d-1}t^{d-1}+\cdots +\bar a_0$, and define 
the {\it dual polynomial}\index{dual polynomial} $f^*$ \index{$f^*(t)$}
of $f$ by 
\begin{equation}\label{dualdef}
f^*(t) = \bar a_0^{-1}t^d\bar f(t^{-1}).
\end{equation}
Note that $(f^*)^* = f$ and $(fg)^* = f^*g^*$ for monic 
polynomials $f$ and $g$. 
In particular, $f$ is irreducible if and only if $f^*$ is.

The next two results are essentially \cite[Lemma 5.1 and Thm.\ 5.2]{JBRIT}. 
They give information about the elementary divisors and the 
$F[t]$-module structure of $V$ induced by an element of 
a classical group $\C(\b)$ or $\C(Q)$. 

\begin{proposition} \label{necessarycondition}
Let $\beta$ be a non-degenerate alternating, symmetric, or hermitian form 
on $V$, and let $Q$ be a non-degenerate quadratic form. 
Let $X \in \C(\beta)$ or $\C(Q)$, and let $U$ be an $X$-invariant subspace 
of $V$ with $U \cap U^{\bot}=\{0\}$. If $X$ has minimal 
polynomial $f$ on $U$, then $f=f^*$.
\end{proposition}

\begin{proof}
First assume $X \in \C(\b)$. Observe that $\beta(uX,v)=\beta(u,vX^{-1})$ for all $u,v \in U$. 
Write $f(t) = \sum_{i=0}^d a_it^i$, where $a_0 \neq 0$ and $a_d=1$. Now
	\begin{eqnarray*}
		0 & = & \beta(Uf(X),\,U)\\
		& = & \beta(U(\sum_i a_iX^i),\,U) \\
		& = & \beta(U,\,U\sum_i \overline{a_i}X^{-i}) \\
		& = & \beta(U,\,\overline{a}_0^{-1}UX^d\sum_i \overline{a_i}X^{-i}) \\
		& = & \beta(U,\,Uf^*(X)).
	\end{eqnarray*}
Since $U$ is non-degenerate by hypothesis, the identity $\beta(U,\,Uf^*(X))=0$ implies $Uf^*(X)=0$, so $f$ divides $f^*$. But both $f$ and $f^*$ are monic 
of degree $d$, so $f=f^*$.

For $X \in C(Q)$, the above proof goes through verbatim, with $\b_Q$ replacing $\b$.
\end{proof}

\begin{proposition} \label{orthogonaldecomposition}
Let $\beta$, $Q$ be as in Proposition $\ref{necessarycondition}$, and let 
$X \in \C(\beta)$ or $\C(Q)$. There exists an orthogonal  decomposition 
$V = \bigoplus_i U_i$ such that, for each $i$, one of the following holds:
\begin{itemize}
\item[{\rm (i)}]  $X$ acts cyclically on $U_i$ with minimal polynomial $f^e$ for some $e\ge 1$, where $f$ is irreducible and $f=f^*$;
\item[{\rm (ii)}]  $U_i = W \oplus W^*$ and $X$ acts cyclically on $W$ 
$($resp.\ $W^*)$ with minimal polynomial $f^e$ $($resp.\ $f^{*e})$ for some $e\ge 1$, where $f$ is irreducible.
\end{itemize}
\end{proposition}

\begin{proof}
For notational convenience, in the case where $X \in \C(Q)$ we write 
$\b$ for the associated bilinear form $\b_Q$.

Let $U$ be one of the summands in the decomposition of $V$ into 
cyclic $X$-submodules given in (\ref{decompCSM}). 
Then $X$ acts cyclically on $U$ with minimal polynomial $f^e$, 
where $f$ is irreducible and $e\ge 1$. 
Let $W$ be an $X$-invariant subspace of $V$ such that $V = U \oplus W$, 
and let $U^* = W^{\bot}$. Since $\beta$ is non-degenerate, 
\begin{equation}\label{uustar}
U^* \cap U^{\bot}= U^{\bot} \cap W^{\bot} = (U \oplus W)^{\bot}=V^{\bot} = 0,
\end{equation}
so for each non-zero $v \in U^*$ there exists 
$u \in U$ such that $\beta(u,v) \neq 0$. Moreover, 
for every non-zero $u \in U$, there exists $v \in U^*$ such 
that $\beta(u,v) \neq 0$ (otherwise $u \in U^{*\bot} = W$, a contradiction). 
We distinguish three cases.

\vspace{2mm} \no {\bf Case 1:} Assume $U\cap U^\perp = 0$. Then $U$ is non-degenerate, so $f=f^*$ by Proposition \ref{necessarycondition}; thus
$U$ is as in conclusion (i) of the proposition.

\vspace{2mm} \no {\bf Case 2:} Assume $U\cap U^\perp \ne 0$ and $U\cap U^* \ne 0$. 
 Since $U \cap U^*$ is $X$-invariant, $Uf(X)^{e-1} \subseteq U \cap U^*$ (because $X$ acts cyclically on $U$ with minimal polynomial $f^e$). For the same reason $Uf(X)^{e-1} \subseteq U \cap U^{\bot}$. Hence $Uf(X)^{e-1} \subseteq U^* \cap U^{\bot}$, contradicting (\ref{uustar}). 

\vspace{2mm} \no {\bf Case 3:} 
Assume $U\cap U^\perp \ne 0$ and $U\cap U^* = 0$. 
We show that $\beta$ is non-degenerate on $U \oplus U^*$: namely, 
for every $u \in U$, $v \in U^*$ with $u+v \neq 0$, there exists 
$z \in U \oplus U^*$ with $\beta(u+v,\,z) \neq 0$. If $v=0$, 
then take $z \in U^*$ such that $\beta(u,z) \neq 0$; 
if $v \neq 0$ and $u \in U^{\bot}$, then take $z \in U$ such that 
$\beta(z,v) \neq 0$ (note that $z$ exists in both cases 
since $U^{\bot} \cap U^* = 0$). Now suppose $v \neq 0$ and $u \notin U^{\bot}$.
Let $a = \min\{n \in \mathbb{N} \, | \, uf(X)^n \in U^{\bot}\}$. 
By our choice of $u$, observe that $1 \leq a \leq e$. 
Consider $$ (u+v)f(X)^a = uf(X)^a+vf(X)^a.  $$
If $uf(X)^a=0$, then $uf(X)^{a-1} \in \ker{f(X)} = Uf(X)^{e-1} \subseteq U \cap U^{\bot}$, contradicting the minimality of $a$. So $uf(X)^a\ne 0$. If $vf(X)^a=0$, then there exists $w \in U^*$ such that $\beta(uf(X)^a,w) \neq 0$, so $\beta((u+v)f(X)^a,w) \neq 0$, and we can choose $z = w\overline{f}(X^{-1})^a$. 
So we can assume that $vf(X)^a\ne 0$.
Hence there exists $w \in U$ such that $\beta(w,v\overline{f}(X^{-1})^a) \neq 0$. Since $uf(X)^a \in U^{\bot}$, it follows that $\beta((u+v)f(X)^a,w) \neq 0$, so we can again choose $z = w{\overline{f}}(X^{-1})^a$.
Thus $\beta$ is non-degenerate on $U \oplus U^*$. 

Let $U_1 = U \oplus U^*$. For every $n$, 
\[
\beta(Uf(X)^n,\,U^*)=0 \;\; \Leftrightarrow \;\; \beta(U,\,U^*f^*(X)^n) = 0.
\]
Since $U^*\cap U^\perp = 0$, 
\[
Uf(X)^n = 0 \;\; \Leftrightarrow \;\; U^*f^*(X)^n = 0.
\]
This proves that the minimal polynomial of $X$ on $U^*$ is $f^{*e}$, 
and the action of $X$ on $U^*$ is cyclic because $\dim{U} = \dim{U^*}$. 
Hence conclusion (ii) holds for $U$.
	
Since $V = U \oplus U^{\bot}$ and $U^{\bot}$ is $X$-invariant, we can repeat the argument for $U^{\bot}$, and the proposition follows by induction.
\end{proof}

Proposition \ref{orthogonaldecomposition} has the following consequence:
if $f^e$ is an elementary divisor of $X \in \C(\b)$ or $\C(Q)$ of multiplicity $m$, 
then $f^{*e}$ is also an elementary divisor of multiplicity $m$. 
Recall that $X, Y \in \GL(V)$ 
are \textit{similar}\index{similar} if they are conjugate, in which case 
we write $X \sim Y$. \label{similarSymbol}\index{$x \sim y$}

\begin{cor}\label{simi} 
If $X$ is an element of $\C(\b)$ or $\C(Q)$, 
then $X$ is similar to $\bar X^{-1}$.
\end{cor}

\begin{proof} By the definition of the dual polynomial, $f^e$ is an elementary divisor of $X$ if and only if $f^{*e}$ is an elementary divisor of $\bar X^{-1}$. Hence, by the preceding remarks, $X$ and $\bar X^{-1}$ have the same elementary divisors, and so they are similar. 
\end{proof}

We now introduce the following notation, following \cite{FS}. 

\begin{definition}  \label{ThreeCases}\index{$\Phi$}\index{$\Phi_1, \Phi_2, \Phi_3$}
{\rm Let $F = \mathbb{F}_{q^u}$ with $u = 1$ or 2, let $\l \mapsto \bar \l$ be a field automorphism of order $u$, and let $f^*$ be as in (\ref{dualdef}). Define}
	\begin{eqnarray*}
		\Phi_1 & := & \{ f: \, f \in F[t] \; | \; f=f^* \mbox{ monic irreducible}, \, \deg{f}=1 \};\\
		\Phi_2 & := & \{ f: \, f \in F[t] \; | \; f = gg^*, \, g \neq g^*, \, g \mbox{ monic irreducible} \};\\
		\Phi_3 & := & \{ f: \, f \in F[t] \; | \; f=f^* \mbox{ monic irreducible}, \, \deg{f}>1 \}.
	\end{eqnarray*}
{\rm Let $\Phi := \Phi_1 \cup \Phi_2 \cup \Phi_3$ and let $f \in \Phi$. For $X \in \C(\beta)$ or $\C(Q)$ as above, and $m$ a positive integer, $f^m$ is a 
\emph{generalized elementary divisor} 
\index{generalized elementary divisor} 
of $X$ if one of the following holds:}
\begin{itemize}
\item  $f \in \Phi_1 \cup \Phi_3$ {\rm and $f^m$ is an elementary divisor of} $X$; 
\item $f \in \Phi_2$, $f=gg^*$ {\rm and $g^m$ is an elementary divisor of $X$ (and so $g^{*m}$ is as well)}.
\end{itemize}
 \index{generalized elementary divisor}
\end{definition}

\begin{proposition}\label{elem} 
The following hold:
\begin{itemize}
\item[{\rm (i)}] If $u=1$, then $\Phi_1 = \{t+1,\,t-1\}$.
\item[{\rm (ii)}] If $u=2$, then $\Phi_1 = \{t-\l : \l\bar \l=1  \}$.
\item[{\rm (iii)}] If $u=1$ and $f \in \Phi_3$, then $\deg{f}$ is even. 
\item[{\rm (iv)}] If $u=2$ and $f \in \Phi_3$, then $\deg{f}$ is odd. 
\end{itemize}
\end{proposition}

\begin{proof}
Parts (i) and (ii) are clear. 

Consider (iii). If $u=1$ and $f \in \Phi_3$, then, 
for each root $\l$ of $f$ (in the splitting field), 
$\l^{-1}$ is also a root, and $\l \ne \l^{-1}$ as $f \ne t\pm 1$. 
Hence $\deg{f}$ is even. 

Finally consider (iv). Let $u=2$ and $f \in \Phi_3$, and let $R(f)$ be the set of roots of $f$.
The map $\alpha \mapsto \alpha^{-q}$ acts as a permutation on $R(f)$. Consider an orbit $\mathcal{O} = \{ \lambda, \lambda^{-q}, \lambda^{q^2}, \dots, \lambda^{(-q)^r} \}$ under this action, and subset $\mathcal{O'} = \{\lambda, \lambda^{q^2}, \lambda^{q^4}, \dots \}$. The polynomial
	$$
	g(t) = \prod_{\mu \in \mathcal{O'}} (t-\mu) = (t-\lambda)(t-\lambda^{q^2})(t-\lambda^{q^4}) \cdots 
	$$
	is a divisor of $f$ and belongs to $\mathbb{F}_{q^2}[t]$ because its coefficients are fixed by the field automorphism $\alpha \mapsto \alpha^{q^2}$. Since $f$ is irreducible, the only possibility is that $g=f$ and $\mathcal{O'}=R(f)$. It follows that $\mathcal{O}$ is the unique orbit in $R(f)$, and also that $\mathcal{O}=\mathcal{O'}$, which implies that $|\mathcal{O}|$ is odd.
\end{proof}

Now we state the main result of this section. 

\begin{theorem} \label{elementsofCG}
Let $F=\F_{q^u}$ with $u=1$ or $2$, 
and let $\l\mapsto \bar \l$ be an automorphism of $F$ of order $u$. 
Let $V$ be a vector space over $F$ and let $X \in \GL(V)$. 
\begin{itemize}
\item[{\rm (i)}] Suppose $u=2$. There exists a non-degenerate hermitian form 
$\beta$ on $V$ such that $X \in \C(\beta)$ if and only if $X \sim \bar X^{-1}$.
\item[{\rm (ii)}]  Suppose $u=1$. There exists a non-degenerate 
alternating form $\beta$ on $V$ such that $X \in \C(\beta)$ if and only 
if $X \sim X^{-1}$ and every elementary divisor $(t\pm 1)^{2k+1}$ of $X$ with $k\in \N$ has even multiplicity.
\item[{\rm (iii)}]  Suppose $q$ is odd and $u=1$. 
\begin{itemize}
\item[{\rm (a)}]  There exists a non-degenerate quadratic form $Q$ on $V$ 
such that $X \in \C(Q)$ if and only if $X \sim X^{-1}$ and every 
elementary divisor $(t\pm 1)^{2k}$ with $k \in \N^+$ 
has even multiplicity. 
\item[{\rm (b)}] Suppose {\rm (a)} holds for $X$ and 
$\dim V$ is even. Then $X \in \C(Q)$ for $Q$ of both plus and minus 
types 
if and only if $X$ has at least one elementary divisor 
$(t \pm 1)^{2k+1}$ for some $k \in \N$. 
 If this is not the case, then $X \in \C(Q)$ for $Q$ of plus type $($resp.\ minus type$)$ if and only if $\sum_{f,e} e\cdot m(f^e)$ is even $($resp.\ odd$)$, where the sum runs over all $f\in \Phi_3$, and $m(f^e)$ is the multiplicity of $f^e$ as an elementary divisor of $X$.
\end{itemize}
\item[{\rm (iv)}]  Suppose $q$ is even and $u=1$. 
\begin{itemize}
\item[{\rm (a)}]  There exists a non-degenerate quadratic form $Q$ 
on $V$ such that $X \in \C(Q)$ if and only if $X \sim X^{-1}$ and 
every elementary divisor $(t+ 1)^{2k+1}$ with $k\in \N$ has even multiplicity. 
\item[{\rm (b)}] Suppose {\rm (a)} holds for $X$. 
Then $X \in \C(Q)$ for $Q$ of both plus and minus types if and only 
if $X$ has at least one elementary divisor $(t + 1)^{k}$ for some 
$k \in \N^+$. 
 If this is not the case, then $X \in \C(Q)$ for $Q$ of plus type 
$($resp.\ minus type$)$ if and only if $\sum_{f,e} e\cdot m(f^e)$ is even 
$($resp.\ odd$)$, where the sum runs over all $f\in \Phi_3$, and 
$m(f^e)$ is the multiplicity of $f^e$ as an elementary divisor of $X$.
\end{itemize}
\end{itemize}
\end{theorem}

To prove the theorem we require several preliminary lemmas. 
The first is \cite[2.2]{Huppert}; we include a proof for completeness.

\begin{lemma}\label{Prop414}
Suppose $q$ is odd and $X$ acts cyclically on $V$ with minimal polynomial $(t-\e)^{2k}$, where $\e = 1$ or $-1$. Then there is no non-degenerate symmetric bilinear form $\beta$ on $V$ such that $X \in \C(\beta)$.
\end{lemma}

\begin{proof}
We prove this for $\e = 1$; the proof for $\e=-1$ is similar. 
Put $m=2k$ and choose a basis $v_1, \dots, v_m$ for $V$ such that $v_1X=v_1$ and $v_iX = v_{i-1}+v_i$ for all $i=2, \dots, m$. 
Suppose for a contradiction that $\beta$ is a non-degenerate symmetric form 
on $V$ such that $X \in \C(\b)$. 
Since $\langle v_1 \rangle^{\bot}$ is an invariant $F[t]$-submodule of $V$ 
of dimension $m-1$, it must equal the unique submodule of this 
dimension, so 
$\langle v_1 \rangle^{\bot} = \langle v_1, \dots, v_{m-1} \rangle$. Hence 
$\beta(v_1,v_m) \neq 0$. Now,
	\begin{eqnarray*}
		0 \neq \beta(v_1,v_m) & = & \beta(v_1,v_1(X-1)^{m-1}) \\
        & = & \beta(v_1(X^{-1}-1)^{m-1},v_1)\\
		& = & \beta((-1)^{m-1}v_1X^{-m+1}(X-1)^{m-1},v_1)\\
		& = & (-1)^{m-1} \beta(v_mX^{-m+1},v_1)\\
		& = & (-1)^{m-1}\beta(v_m,v_1)\\
		& = & -\beta(v_1,v_m).
	\end{eqnarray*}
	This is impossible in odd characteristic. 
\end{proof}

\noindent 
The next lemma also follows from  \cite{Huppert}. 

\begin{lemma}\label{HuppertResults}
Let $q$ be odd. Let $\beta$ be a non-degenerate symmetric form on $V$, 
and let $X \in \C(\beta)$. Let $U$ be a non-degenerate  direct summand of 
$V \downarrow X$  that satisfies conclusion {\rm (ii)} of Proposition 
$\ref{orthogonaldecomposition}$.
Then one of the following holds: 
\begin{itemize}
\item[{\rm (i)}] there are non-degenerate $X$-submodules $U_1,U_2$ of $U$ 
of dimension $\frac{1}{2} \dim U$ such that $U = U_1+U_2$, an orthogonal sum;
\item[{\rm (ii)}] there exist totally isotropic $X$-submodules $W,W^*$ of $U$ such that $U = W \oplus W^*$, and either $f \ne f^*$ or $f^e = (t\pm 1)^{2k}$.
\end{itemize}
\end{lemma}

\begin{proof}
	If $f \neq f^*$, then it follows directly from Lemma \ref{L26w} (whose proof is independent of this chapter) that (ii) holds. Suppose that $f=f^*$ and (i) does not hold. Then, in the terminology of \cite[1.8]{Huppert}, $U$ is of ``Type 1".
 It follows from \cite[2.1]{Huppert} that $f=t\pm 1$ and $e$ is even; moreover, \cite[2.4]{Huppert} gives the existence of
 totally isotropic $X$-submodules $W,W^*$ such that $U=W \oplus W^*$. Thus (ii) holds.
 \end{proof}

\begin{lemma}  \label{typePhi3}
Suppose that $u=1$, so $F = \F_q$. Let $f \in F[t]$ be an irreducible polynomial with 
$f= f^*$ and $d = \deg{f}>1$, and let $X$ be its companion matrix over $F$. 
If $Q$ is a non-degenerate quadratic form on $V = F^d$ such 
that $X \in \C(Q)$, then $Q$ has minus type.
\end{lemma}

\begin{proof} Note that $d$ is even by Proposition \ref{elem}. 
We show in Theorem \ref{CentrSem} 
(proved independently of this chapter) that 
$C_{\C(Q)}(X) \cong \GU_1(q^{d/2})$, and so it is a subgroup 
of $\C(Q) \cong \Or_d^{\e}(q)$, where $\e=\pm$. 
But $|\GU_1(q^{d/2})| = 1+q^{d/2}$ divides 
$|\Or_d^+(q)|$ only if $(d,q) = (2,3)$ or $(6,2)$. 
In these cases the only possible polynomials $f$ are $t^2+1 \in \mathbb{F}_3[t]$ and $t^6+t^3+1 \in \mathbb{F}_2[t]$; 
we check directly that their companion matrices do not 
preserve non-degenerate quadratic forms of plus type. 
Hence $\e=-$. 
\end{proof}

\vspace{4mm}
\no {\bf Proof of Theorem \ref{elementsofCG}}

\vspace{2mm}
\no {\bf (1) } We begin by proving the left to right implications 
(the ``only if" assertions) in parts (i), (ii), (iii)(a) and (iv)(a) of the theorem. 
Suppose that $X \in \C(\b)$ or $\C(Q)$. Then $X \sim \bar X^{-1}$ by Corollary \ref{simi}. Moreover, if the multiplicity of $(t\pm 1)^{2k+1}$ as an elementary divisor of $X$ is odd for some $k$, then, by Proposition \ref{orthogonaldecomposition}, there must be a non-degenerate subspace $U_1$ on which $X$ acts cyclically with minimal polynomial $(t\pm 1)^{2k+1}$; in particular, $\dim U_1 = 2k+1$ is odd, and so $\b$ cannot be alternating. This completes the proof of the left to right implication in (ii) and (iv)(a). Similarly, if $q$ is odd and the multiplicity of $(t\pm 1)^{2k}$ is odd for some $k$, then $\b$ cannot be symmetric, by Lemma \ref{Prop414}, completing the left to right part of (iii)(a). 

\vspace*{2mm}
\no {\bf (2) } Next we prove the right to left implications (the ``if" assertions) in parts (i), (ii), (iii)(a) and (iv)(a) of the theorem. Suppose that $X \sim \bar X^{-1}$ and $X$ satisfies the condition on elementary divisors $(t\pm 1)^{2k+1}$ in parts (ii) and (iv)(a), and on  elementary divisors $(t\pm 1)^{2k}$ in part (iii)(a). We can assume that $X$ has a unique generalized elementary divisor $f^e$ with $f \in \Phi$, and that this has multiplicity $m = 1$ or 2, where $m=2$ only in the cases where even multiplicity is assumed in  (ii), (iii)(a) and (iv)(a); once we exhibit forms $\b$ or $Q$ fixed by such elements, we can just take direct sums to exhibit forms fixed by $X$ in the general case.

\vspace*{1mm}
\noindent 
\textbf{Case a:} $f \in \Phi_1$. Here $X$ has a unique elementary divisor $(t-\l)^e$, where $\l \bar \l = 1	$. Hence $X = \l u$, where $u$ is unipotent. 
Under the assumed conditions on the multiplicity $m$, unipotent elements are 
defined  in $\C(\b)$ and $\C(Q)$ in the following sections: 
Section \ref{sugood} for unitary groups; 
Sections \ref{spgood} and \ref{sogood} for symplectic and orthogonal groups in 
odd characteristic; 
and Section \ref{badsec} for symplectic and orthogonal groups in characteristic 2.

\vspace*{1mm}
\noindent 
\textbf{Case b:} $f \in \Phi_2$. Here $f=gg^*$, where $g$ is irreducible and $g \ne g^*$, and $f^e$ has multiplicity 1 as a generalized elementary divisor of $X$. Let $d= \deg{g}$. In an appropriate basis, $X$ has block diagonal matrix $X = Y \oplus Y^{*-1}$, with $Y$ a Jordan block relative to $g^e$. Then $X$ is an isometry for the form with matrix $B$ or, in the orthogonal case, for the quadratic form with matrix $A$, where
$$
B = \left(\begin{array}{cc} \mathbb{O} & \mathbb{I} \\ \varepsilon \mathbb{I} & \mathbb{O} \end{array}\right), \quad A = \left(\begin{array}{cc} \mathbb{O} & \mathbb{I} \\ \mathbb{O} & \mathbb{O} \end{array}\right),
$$
$\mathbb{O}$\index{$\mathbb{O}$, $\mathbb{O}_n$} is the $de \times de$ zero matrix, $\mathbb{I}$ is the $de \times de$ identity matrix, and $\varepsilon = -1$ in the symplectic case and $1$ otherwise. 

This completes the argument for Case {\bf b}. Note that $V= U \oplus W$, where $X$ acts cyclically on $U$ (resp.\ $W$) with minimal polynomial $g^e$ (resp.\ $g^{*e}$), and $U$ and $W$ are totally singular.

\vspace*{1mm}
\noindent 
\textbf{Case c:} $f \in \Phi_3$. We extend \cite[Thm.\ 5.4]{JBRIT} to unitary groups and to symplectic and orthogonal groups in even characteristic. In this case $X$ has a unique elementary divisor $f^e$ with multiplicity 1, where $f=f^*$ is irreducible of degree $d>1$. Then $X$ acts cyclically on $V$ with minimal polynomial $f^e$. Let $U_1$ and $U_2$ be isomorphic to $V$ (as $X$-modules) and 
write $U = U_1 \oplus U_2$. Using the same argument as in Case {\bf b}, we see that there exists a non-degenerate form $\beta$ (or quadratic form $Q$) on $U$ such that $X \in \C(\beta)$ and $U_1, U_2$ are totally singular. 

Now suppose there exists $v \in U$ and $i \in \mathbb{Z}$ such that $\beta(v,\,vX^if(X)^{e-1}) \neq 0$. Let $W$ be the cyclic $F[t]$-submodule of $U$ generated by $v$. The minimal polynomial of $X$ on $W$ is $f^c$ for some $c \leq e$. If $c < e$, then $vf(X)^{e-1}=0$, so $\beta(v,vX^if(X)^{e-1})=\beta(v,0)=0$, a contradiction. Hence $c=e$, and the only submodules of $W$ are of the form $Wf(X)^m$ for some $m\le e$. Moreover, $W$ is non-degenerate: if not, then $W\cap W^\perp$ is a nonzero submodule of $W$, hence it equals $Wf(X)^m$ for some $m< e$; but this implies that $vf(X)^{e-1} \in W^\perp$, a contradiction. 
Thus $W$ is non-degenerate and $X$ is an element of 
$\C(\beta_{W})$ or $\C(Q_W)$, where $\b_W$, $Q_W$ are the restrictions of $\b, Q$ to $W$.
Since $W$ is isomorphic to $V$ as $F[t]$-module, such a form must exist on $V$, as required.

By the previous paragraph, the proof for the case $f \in \Phi_3$ is complete 
except when $\beta(v\,,vX^if^{e-1}(X))=0$ for every $v \in V$ and 
$i \in \mathbb{Z}$. Suppose for a contradiction that this is the case. 
For $v = v_1+v_2$ with $v_i \in U_i$, this condition is equivalent to
\begin{equation}\label{bineq}
\beta(v_1,\,v_2X^if(X)^{e-1})+\beta(v_2,\,v_1X^if(X)^{e-1})=0
\end{equation}
for all $v_1 \in U_1$, $v_2 \in U_2$ and $i \in \mathbb{Z}$. In the unitary case, choosing $\l \in \F_{q^2}$ such that $\l \ne \bar \l$, and replacing $v_1$ by $\l v_1$ in (\ref{bineq}) gives
\[
\l\beta(v_1,\,v_2X^if(X)^{e-1})+\bar \l\beta(v_2,\,v_1X^if(X)^{e-1})=0.
\]
Combined with (\ref{bineq}), and taking $i=0$, this implies that $\beta(v_1,\,v_2f(X)^{e-1})=0$ for all $v_1\in U_1$, $v_2 \in U_2$, which implies that the restriction $f(X)^{e-1}\downarrow U_2 = 0$, a contradiction.

Now consider the symplectic and orthogonal cases. By a straightforward computation, (\ref{bineq}) implies the following sequence of identities:
\begin{eqnarray*}
	\beta(v_1,\,v_2[X^if(X)^{e-1}+\varepsilon X^{-i}f(X^{-1})^{e-1}]) & = & 0,\\
	\beta(v_1,\,v_2[X^if(X)^{e-1}+\varepsilon f(0)^{e-1}X^{-i-d(e-1)}f(X)^{e-1}]) & = & 0, \quad (\mbox{using } f=f^*)\\
\beta(v_1,\,v_2f(X)^{e-1}[X^i+\varepsilon f(0)^{e-1}X^{-i-d(e-1)}]) & = & 0,\\
	\beta(v_1,\,v_2f(X)^{e-1}[X^{d(e-1)+2i}+\varepsilon f(0)^{e-1}]) & = & 0.
\end{eqnarray*}
Since $\beta$ is non-degenerate on $U_1 \oplus U_2$, for every $i \in \Z$, the restriction 
\[
f(X)^{e-1}\left(X^{d(e-1)+2i}+\varepsilon f(0)^{e-1}\right) \downarrow U_2 =0.
\]
This implies that $f(t)$ divides 
$t^{d(e-1)+2i}+\varepsilon f(0)^{e-1}$ for every $i$ such 
that $d(e-1)+2i$ is non-negative. Since we are in the symplectic or 
orthogonal case, $d$ is even by Proposition \ref{elem}. Hence,
choosing $i = 1-d(e-1)/2$, we see that $f(t)$ divides 
$t+\varepsilon f(0)^{e-1}$ of degree 1, which is impossible. 

This completes the argument for Case {\bf c}. 
Parts (i), (ii), (iii)(a) and (iv)(a) of the theorem are now established. 

\vspace{2mm}
\no {\bf (3) } It remains to prove parts (iii)(b) and (iv)(b) of the theorem. 
Consider (iii)(b). We suppose (as in (iii)(a)) that $\dim V$ is even, 
$X \sim X^{-1}$, and every elementary divisor 
$(t\pm 1)^{2k}$ with $k\in \N^+$ has even multiplicity. 
Combining Propositions \ref{orthogonaldecomposition} and 
\ref{HuppertResults}, we see that $V\downarrow X$ has an orthogonal 
decomposition $V = \bigoplus_iU_i$ such that, for each $i$, one 
of the following holds:
\begin{itemize}
\item[{\bf (a)}]  $X$ acts cyclically on $U_i$ with minimal polynomial $f^e$ for some $e \geq 1$, where $f$ is irreducible and $f=f^*$;
\item[{\bf (b)}]  $U_i = W \oplus W^*$ and $X$ acts cyclically on $W$ (resp.\ $W^*$) with minimal polynomial $f^e$ (resp.\ $f^{*e}$) for some $e \geq 1$, where $f$ is irreducible; moreover, $W$ and $W^*$ are totally isotropic, and either $f \ne f^*$ or $f^e = (t\pm 1)^{2k}$.
\end{itemize}

Assume first that $X$ has an elementary divisor $(t\pm 1)^{2k+1}$ for 
some $k$. Then there must be a summand $U_i$ on which $X$ acts cyclically 
with minimal polynomial $(t\pm 1)^{2k+1}$. From Section \ref{sogood}, 
we see that there are quadratic forms on $U_i$ of both square and 
non-square discriminants preserved by the action of $X$ (corresponding to 
the unipotent elements $V_1(2k+1)$ and \mbox{$V_\a(2k+1)$)}. 
Hence $X$ preserves 
quadratic forms of both square and non-square discriminants on $V$, 
showing that $X \in \C(Q)$ for $Q$ of both plus and minus types, as required.

Now assume that $X$ has no elementary divisor $(t\pm 1)^{2k+1}$. Consider a summand $U_i$. If $U_i$ is as in {\bf (b)} above, then clearly $X^{U_i} \in \Or^+(U_i)$. Now suppose $U_i$ is as in {\bf (a)}. Here $f=f^*$ by Proposition \ref{necessarycondition}. If $f = t\pm 1$, then $f^e = (t\pm 1)^{2k}$ for some $k$, by our assumption at the beginning of this paragraph. But this contradicts Lemma \ref{Prop414}. Hence $\deg{f} \ge 2$ and $f \in \Phi_3$. Now $X$ acts cyclically on $U_i$ with minimal polynomial $f^e$. Let $X = SU$ be the Jordan decomposition of $X$, where $S$ is semisimple and $U$ unipotent. An application of Lemma \ref{HuppertResults} now shows that $V$ is an orthogonal sum of cyclic $S$-submodules $U_1',\ldots,U_e'$, on which $S$ acts cyclically with minimal polynomial $f$. By Lemma \ref{typePhi3}, each $U_i'$ has minus type, and hence $U_i$ has type $(-1)^e$. It follows that $V$ has type $(-1)^l$, where $l = \sum_{f,e} e \cdot m(f^e)$ and the sum runs over $f \in \Phi_3$. 

This concludes the proof of part (iii)(b) of the theorem. The proof of (iv)(b) is very similar, and in fact easier, as we can work as above with the decomposition $V = \bigoplus U_i$ of Proposition \ref{orthogonaldecomposition}. 

The proof of Theorem \ref{elementsofCG} is now complete.

\chapter{Unipotent classes in good characteristic}\label{unigoodchap}
As described in Chapter \ref{introchap}, our approach to the conjugacy problems (1)-(3) of Section \ref{mainprob} splits
naturally into the analysis of the unipotent and semisimple classes. 
In this chapter and the next we solve the conjugacy problems for unipotent classes in classical groups.
Since we deal separately with the different families of classical groups, we drop the $\C(\b)$ and $\C(Q)$ notation 
of the previous chapter, and let $G$ be $\Sp_{2n}(q)$, $\Or_n^\epsilon(q)$, or $\GU_n(q)$ 
where $q = p^a$ with $p$ prime. 
Recall that the characteristic $p$ is {\it bad} if $G$ is symplectic or 
orthogonal and $p=2$; otherwise $p$ is {\it good}.
We address the classical groups in good characteristic in this chapter
and the bad characteristic cases in Chapter \ref{unibadchap}. 

\section{Unipotent class representatives}\label{goodrep}

The unipotent class representatives of the finite classical groups are given 
(in all positive characteristics) in \cite{GLOB}. For effective computation 
with these classes in good characteristic, we choose 
representatives that differ from those of \cite{GLOB}. 
In this section we present these representatives
and structural information about their centralizers.

As outlined in Section \ref{Sect12}, the choice of representatives and the construction of their centralizers is
based on theory developed in \cite{LS}. We now present a more detailed summary of this theory. Let $K$ be an 
algebraically closed field of odd characteristic $p$, and let $\bar V$ be an $n$-dimensional vector space over $K$.
Let $(\,,\,)$ be a non-degenerate symplectic or orthogonal bilinear form on $\bar V$, with isometry group $\bar G = \Sp(\bar V)$
or $\Or(\bar V)$. The Lie algebra\index{Lie algebra} of $\bar G$ is 
\[
L(\bar G) = \{A \in \hbox{End}(\bar V) : (vA,w)+(v,wA)=0 \;\forall v,w \in \bar V\}
\]
(see \cite[Lemma 2.7]{LS}). By \cite[Lemma 2.15]{LS}, if $e \in L(\bar G)$ is a nilpotent element, then 
$u:=(1-e)(1+e)^{-1}$ is a unipotent element of $\bar G$, and the map
\[
e \mapsto (1-e)(1+e)^{-1}
\]
is a bijection, known as the {\it Cayley map}, from the set of 
nilpotent elements of $L(\bar G)$ to the set of unipotent elements of $\bar G$. 
The Cayley map $\chi$ is {\it $\bar G$-equivariant}\index{equivariant}: 
namely, $\chi(e^g) = \chi(e)^g$ for $g \in \bar G$. 
In particular, $C_{\bar G}(u) = C_{\bar G}(e)$.

In \cite[\S 3.3.2]{LS}, various nilpotent elements $e$ of $L(\bar G)$ are defined, and the corresponding unipotent
elements $u=(1-e)(1+e)^{-1}$ are shown to form a complete set of unipotent class representatives for $\bar G$. 
For each such nilpotent element $e$, a 1-dimensional {\it torus}\index{torus} $T = \{T(c) : c \in K^*\}\le \bar G$ is constructed, with the property that 
$eT(c) = c^2e$ for all $c \in K^*$. The space $\bar V$ is a direct sum of 
$T$-weight spaces $\bar \myV_{(i)}$ (for $i \in \Z$), where
\[
\bar \myV_{(i)} = \{v \in \bar V : vT(c) = c^iv\}.
\]
If we set $\bar \myW_{(i)} = \la \bar \myV_{(j)} : j \ge i\ra$, then the subspaces 
$\bar \myW_{(i)}$ form a flag
$$\cdots \subseteq \bar \myW_{(i+1)} \subseteq \bar \myW_{(i)} \subseteq \bar \myW_{(i-1)} \subseteq \cdots.$$ 
The stabilizer in $\bar G$ of this flag is a parabolic subgroup $P$, and the Levi decomposition of this 
parabolic is $P = QL$, with unipotent radical $Q$ and Levi factor $L = C_{\bar G}(T)$. 
By \cite[Prop.\ 3.10 and Lemma 3.13]{LS}, $C_{\bar G}(u) = C_{\bar G}(e) = UR$, where $U$ is 
the unipotent radical and $R = C_{\bar G}(T,e) = C_L(e)$ is reductive. 
Moreover, by \cite[Lemma 2.29]{LS}, $U = C_Q(e)$, and hence
\begin{equation}\label{cbarg}
C_{\bar G}(u) = C_{\bar G}(e) = C_Q(e)\,C_L(e).
\end{equation}
The structure of the reductive group $R = C_L(e)$ and the dimension of the 
unipotent group $C_Q(e)$ are given by \cite[Prop.\ 3.7, (3.4), (3.5)]{LS}.

This summarizes the theory of unipotent class representatives and centralizers in symplectic and orthogonal groups over algebraically closed fields of odd characteristic. Over finite fields, the theory is similar, but requires various 
tweaks, as we now discuss.

\subsection{Symplectic groups}\label{spgood}
Let $G = \Sp_{2n}(q)$ with $q = p^a$ and $p$ an odd prime, 
and let $(\,,\,)$ be a symplectic form on $V = \F_q^n$ preserved by $G$. 
Let $J_i$ be a unipotent $i\times i$ Jordan block matrix. 
By Theorem \ref{elementsofCG}, the Jordan form of a unipotent element 
of $G$ has an even number of blocks of each odd size, so it is a block diagonal sum
\begin{equation}\label{jsp}
\bigoplus_{i=1}^r J_{2k_i}^{[a_i]} \oplus \bigoplus_{i=1}^s J_{2l_i+1}^{[2b_i]},
\end{equation}
where $k_1,\ldots,k_r$ are distinct, as are $l_1,\ldots,l_s$. 
By \cite[Prop.\ 2.3]{GLOB}, the number of conjugacy classes in $G$ of 
elements with such a Jordan form is $2^r$. 

We now present 
representatives for these classes.
While our representatives differ from those in \cite{GLOB}, 
they are labelled in the same way. 
The finite symplectic Lie algebra is
\[
L(G) = sp_{2n}(q) = \{A \in {\rm End}(V) : (vA,w)+(v,wA) = 0 \;\forall v,w \in V\}.
\]
As discussed above, our representatives are images of nilpotent elements in 
$L(G)$ under the Cayley map. 

\subsubsection{Even blocks $V_\b(2m)$} 
Let $\a$ be a fixed non-square in $\F_q^*$, and let $\b \in \{1,\a\}$. Let $V_{2m}$ be a $2m$-dimensional vector space over $\F_q$ with basis $v_{-(2m-1)},v_{-(2m-3)},\ldots, v_{2m-1}$. Define a symplectic form on $V_{2m}$ by setting $(v_i,v_{-i}) = (-1)^{\lfloor i/2\rfloor}$ for all $i$, and all other values $(v_i,v_j) = 0$. Define $e_\b \in {\rm End}(V_{2m})$ by
\[
\begin{array}{llll}
e_\b: & v_i & \mapsto&  v_{i+2}\;\;(i<2m-1,\,i\ne -1), \\
          & v_{-1}& \mapsto & \b v_1,\\
          & v_{2m-1} & \mapsto & 0.
\end{array}
\]
Now $e_\b$ is a nilpotent element of $sp_{2m}(q)$. With respect to the given basis, $e_\b$ has matrix
\[
\begin{pmatrix} 0&1&&&&&&& \\ &0&1&&&&&& \\ &&\ddots &&&&&& \\ &&&0&1&&&& \\ &&&&0&\b&&& \\
                          &&&&&0&1&& \\ &&&&&&\ddots &&\\ &&&&&&&0&1 \\ &&&&&&&&0 \end{pmatrix}.
\]
We use the Cayley map to define a unipotent element
\[
V_\b(2m) = (1-e_\b)(1+e_\b)^{-1} \in G = \Sp(V_{2m}).
\]
Observe that $V_\b(2m)$ is a single Jordan block. Its matrix relative to the given basis is
\[
\begin{pmatrix*}[r]
 1&-2&2&\cdots & 2\e & -2\e \b & 2\e \b&& \cdots & -2\b \\
                           & 1 & -2 & \cdots & -2\e & 2\e\b & -2\e\b && \cdots & 2\b \\
                           && 1 &&&&&&& \\
                           &&& \ddots &&&&&& \\
                           &&&& 1 & -2\b & 2\b && \cdots & 2\e \b \\
                           &&&&& 1 & -2 &2 & \cdots & 2\e \\
                           &&&&&& 1 & -2 & \cdots & -2\e \\
                           &&&&&&&& \ddots & \\
                            &&&&&&&&& 1 \\
\end{pmatrix*}
\]
where $\e = (-1)^{m+1}$ and the first $\b$ appears in the $(m+1)^{th}$ column. 
Writing $\bar \F_q$ for the algebraic closure of $\F_q$, define a 1-dimensional torus 
$T_1= \{T_1(c) : c \in \bar \F_q^*\} \le \Sp(\bar V_{2m})$, where $\bar V_{2m} = V_{2m} \otimes \bar \F_q$,  
by the action 
\[
T_1(c): v_i \mapsto c^iv_i \hbox{ for all }i.
\]
Note that $e_{\b} T_1(c) = c^2e_\b$. Also one checks that $V_\b(2m)$ is conjugate to the element labelled $V_{2\b}(2m)$ in \cite[\S 2.4]{GLOB}.
                           
\subsubsection{Odd blocks $W(2m+1)$}  
Let $V_{4m+2}$ be a $(4m+2)$-dimensional space over $\F_q$ with basis 
\[
w_{-2m},x_{-2m}, w_{-(2m-2)}, x_{-(2m-2)},\ldots, w_{2m}, x_{2m},
\]
and symplectic form defined by $(w_i,x_{-i}) = (-1)^{\lfloor i/2\rfloor}$ for all $i$ and all other values $(w_i,w_j)$, $(w_i,x_j)$, $(x_i,x_j)$ equal 0. 
Define $e \in {\rm End}(V_{4m+2})$ by
\[
\begin{array}{llll}
e: & w_i & \mapsto & w_{i+2}\;\;(i<2m), \\
   & x_i & \mapsto & x_{i+2}\;\;(i<2m), \\
          & w_{2m}& \mapsto & 0,\\
          & x_{2m} & \mapsto & 0.
\end{array}
\]
Now $e$ is a nilpotent element of $sp_{4m+2}(q)$. Define a unipotent element
\[
W(2m+1)  = (1-e)(1+e)^{-1} \in G = \Sp(V_{4m+2}).
\]
Observe that $W(2m+1)$ has Jordan form $J_{2m+1}^{[2]}$. Its matrix relative to the given basis is
\[
\begin{pmatrix*}[r] I&-2I&2I& \cdots &(-1)^{m+1}2I \\
                          & I & -2I & \cdots & (-1)^m2I \\
                          &&& \ddots & \\
                          &&&&I
\end{pmatrix*}
\]
where $I$ denotes the $2\times 2$ identity matrix. 
Define a 1-dimensional torus $T_1 \le \Sp(\bar V_{4m+2})$ by the action
\[
T_1(c): w_i \mapsto c^iw_i,\,x_i \mapsto c^ix_i \;\;\hbox{ for all }i.
\]

\subsubsection{Class representatives} 
By \cite[Prop. 4.3]{GLOB}, for a given Jordan form (\ref{jsp}), there are $2^r$ conjugacy class representatives in $\Sp(V)$, namely
\begin{equation}\label{vwsp}
\bigoplus_{i=1}^r (V_{\b_i}(2k_i) \oplus V_1(2k_i)^{[a_i-1]}) \oplus \bigoplus_{i=1}^s W(2l_i+1)^{[b_i]},
\end{equation}
where $V = \bigoplus_{i=1}^r V_{2k_i}^{[a_i]} \oplus \bigoplus_{i=1}^s V_{4l_i+2}^{[b_i]}$, and each $\b_i \in \{1,\a\}$ (and the $k_i$ are distinct, as are the $l_i$). As defined in the above notation, each $V_{2k_i}$ (resp.\ $V_{4l_i+2}$) is a subspace of $V$ on which the unipotent element 
$V_{\b_i}(2k_i)$ or $V_1(2k_i)$ (resp.\ $W(2l_i+1)$) acts.
Define $T$ to be the 1-dimensional torus acting as $T_1(c)$ on each subspace $\bar V_{2k_i}, \bar V_{4l_i+2}$.

To define corresponding matrix representatives, choose an ordered basis of $V$ according to the $T$-weights. Namely, list the spaces $V_{2k_j}, V_{4l_j+2}$ in decreasing order of Jordan block sizes. For the $i^{th}$ space, if it is $V_{\b}(2k)$, then label its basis as
\[
v_{-(2k-1)}^{(i)},v_{-(2k-3)}^{(i)},\ldots, v_{2k-1}^{(i)},
\]
and if it is $W(2l+1)$, then label its basis as 
\[
w_{-2l}^{(i)},x_{-2l}^{(i)}, w_{-(2l-2)}^{(i)}, x_{-(2l-2)}^{(i)},\ldots, w_{2l}^{(i)}, x_{2l}^{(i)}.
\]
Now take a basis of $V$ consisting of these vectors $v_j^{(i)}, w_j^{(i)}, x_j^{(i)}$ in increasing order lexicographically with respect to the pairs $(j,i)$. Choose the matrix representative corresponding to the linear map (\ref{vwsp}) to be its matrix with respect to this ordered basis.

\begin{example} \label{sp-example}
Consider the unipotent element 
\[
V_\a(4) \oplus W(3) \oplus V_1(2) \in \Sp(V),
\]
where $\dim V = 12$. The ordered basis is
\[
v_{-3}^{(1)},w_{-2}^{(2)},x_{-2}^{(2)},v_{-1}^{(1)},v_{-1}^{(3)},w_{0}^{(2)},x_{0}^{(2)},
v_{1}^{(1)},v_{1}^{(3)},w_{2}^{(2)},x_{2}^{(2)},v_{3}^{(1)},
\]
and the corresponding matrix representative is
\[
\begin{pmatrix*}[r]  1&&&-2&&&&2\a &&&&-2\a \\
                           &1&&&&-2&&&&2&& \\
                           &&1&&&&-2&&&&2& \\
                           &&&1&&&&-2\a &&&&2\a \\
                           &&&&1&&&&-2&&& \\
                           &&&&&1&&&&-2&& \\
                           &&&&&&1&&&&-2& \\
                           &&&&&&&1&&&&-2 \\
                           &&&&&&&&1&&& \\
                           &&&&&&&&&1&& \\
                           &&&&&&&&&&1& \\
                           &&&&&&&&&&&1
\end{pmatrix*}.
\]
Let $\myV_{(i)}$ denote the {\it $T$-weight space}\index{weight space} for weight $i$: the span of all vectors in the basis with subscript $i$. 
Set $\myW_{(i)} = \la \myV_{(j)} : j \ge i\ra$ -- so 
$\myW_{(3)} = \la v_{3}^{(1)}\ra$, 
$\myW_{(2)} = \la w_{2}^{(2)},x_{2}^{(2)}, v_{3}^{(1)}\ra$ and so on. 
The subspaces $\myW_{(i)}$ form a flag
\[
\myW_{(3)} \subset \myW_{(2)} \subset \myW_{(1)} \subset \cdots \subset \myW_{(-3)} = V,
\]
and the stabilizer in $G = \Sp(V)$ of this flag is a parabolic subgroup $P = QL$ with unipotent radical $Q$ and Levi subgroup $L = C_G(T) = \GL_1(q) \times \GL_2(q) \times \GL_2(q) \times \Sp_2(q)$.
\end{example}

Returning to the general case, let $P$ be the parabolic subgroup stabilizing the flag of $V$ consisting of subspaces 
$\myW_{(i)} = \la \myV_{(j)} : j \ge i\ra$, where $\myV_{(i)}$ is the $T$-weight space 
for weight $i$. As in the example, $P=QL$, with unipotent radical $Q$ and 
Levi subgroup $L= C_G(T)$. Now the discussion around (\ref{cbarg}) (the finite version also uses \cite[Thm.\ 7.1]{LS}) gives the following result.

\begin{thm}\label{sympreps} 
Let $G = \Sp_{2n}(q)$ with $q$ odd, and let $u\in G$ be the unipotent matrix representative defined above for the linear map
\[
\bigoplus_{i=1}^r (V_{\b_i}(2k_i) \oplus V_1(2k_i)^{[a_i-1]}) + \bigoplus_{i=1}^s W(2l_i+1)^{[b_i]}.
 \]
Let $P = QL$ be the corresponding parabolic subgroup. 
\begin{itemize}
\item[{\rm (i)}] $C_G(u) = C_Q(u)C_L(u)$;
\item[{\rm (ii)}] $C_L(u) \cong \prod_{i=1}^s \Sp_{2b_i}(q) \times \prod_{i=1}^r \Or_{a_i}^{\e_i}(q)$;
\item[{\rm (iii)}] rewrite the Jordan form of $u$ as $\bigoplus J_i^{r_i}$; 
then $|C_Q(u)| = q^\myexpR$ where 
\[
\myexpR = \frac{1}{2}\sum_i (i-1)r_i^2 + \sum_{i<j}ir_ir_j + \frac{1}{2}\sum_{i\;even}r_i.
\]
\end{itemize}
\end{thm}

In part (ii), the sign $\e_i=\pm$ is the sign of the orthogonal form with 
Gram matrix ${\rm diag}(\b_i,1,\ldots,1)$. 
In Section \ref{spcent} we justify this, and  show how to construct 
the subgroups $C_Q(u)$ and $C_L(u)$ of the centralizer. 

\subsection{Orthogonal groups}\label{sogood}

Let $V$ be a vector space of dimension $n$ over $\F_q$, where $q=p^a$ is odd, and let $(\,,\,)$ be a non-degenerate symmetric bilinear form on $V$. Let $G = \Or(V)$ be the corresponding orthogonal group. If $n$ is even, then $G \cong \Or_n^\e(q)$ with $\e = \pm$; if $n$ is odd, then there is only one type of orthogonal group $\Or_n(q)$, but for ease of notation we sometimes denote this by $\Or_n^\e(q)$. 
The Jordan form of a unipotent element of $G$ has an even number 
of blocks of each even size, so it is 
\begin{equation}\label{jso}
\bigoplus_{i=1}^r J_{2k_i+1}^{[a_i]} \oplus \bigoplus_{i=1}^s J_{2l_i}^{[2b_i]},
\end{equation}
where $k_1,\ldots,k_r$ are distinct, as are $l_1,\ldots,l_s$. By \cite[Prop.\ 2.4]{GLOB}, the number of conjugacy classes in $G$ of elements with such a Jordan form is as follows:
\begin{itemize}
\item $n$ odd: $2^{r-1}$ classes;
\item $n$ even: $2^{r-1}$ classes in each of $\Or^+_n(q)$ and $\Or^-_n(q)$, with the exception that if $r=0$, then there is just one class in $\Or^+_n(q)$ and none in $\Or^-_n(q)$.
\end{itemize}

To produce class representatives, 
we use the finite orthogonal Lie algebra
\[
L(G) = so_{n}(q) = \{A \in {\rm End}(V) : (vA,w)+(v,wA) = 0 \;\forall v,w \in V\}.
\]

\subsubsection{Odd blocks $V_\b(2m+1)$} 
Let $\a$ be a fixed non-square in $\F_q^*$, and let $\b \in \{1,\a\}$. 
Let $V_{2m+1}$ be a $(2m+1)$-dimensional vector space over $\F_q$, with basis $v_{-2m},v_{-(2m-2)},\ldots, v_{2m}$. Define an orthogonal  form on $V_{2m+1}$ by setting $(v_i,v_{-i}) = (-1)^{ i/2}$ for all $i \ne 0$, $(v_0,v_0)=\b$, and all other values $(v_i,v_j) = 0$. 
Define $e_\b \in {\rm End}(V_{2m+1})$ by
\[
\begin{array}{llll}
e_\b: & v_i & \mapsto & v_{i+2}\;\;(i\ne 0,2m), \\
          & v_{0}& \mapsto & \b v_2,\\
          & v_{2m} & \mapsto & 0.
\end{array}
\]
Now $e_\b$ is a nilpotent element of $so_{2m+1}(q)$. 
As in Section \ref{spgood}, we apply the Cayley map and set
\[
V_\b(2m+1) = (1-e_\b)(1+e_\b)^{-1} \in G = \Or(V_{2m+1}).
\]
Observe that  $V_\b(2m+1)$ is a single Jordan block. 
Its matrix relative to the given basis is 
the $(2m+1)$-dimensional version of the matrix given 
for $V_\b(2m)$ in the symplectic case.  
It is conjugate to the element  
labelled $V_{2\b}(2m+1)$ in  \cite[\S 2.5]{GLOB}.

Define a 1-dimensional torus $T_1= \{T_1(c) : c \in \bar \F_q^*\} \leq 
\Or(\bar V_{2m+1})$, where $\bar V_{2m+1} = V_{2m+1} \otimes \bar \F_q$,  by the action 
\[
T_1(c): v_i \mapsto c^iv_i \hbox{ for all }i.
\]
Note that $e_\b T_1(c) = c^2e_\b$.
                           
\subsubsection{Even blocks $W(2m)$}  
Let $V_{4m}$ be a $4m$-dimensional space over $\F_q$ with basis 
\[
w_{-(2m-1)},x_{-(2m-1)}, w_{-(2m-3)}, x_{-(2m-3)},\ldots, w_{2m-1}, x_{2m-1},
\]
and orthogonal form defined by $(w_i,x_{-i}) = (-1)^{\lfloor i/2\rfloor}$ for all $i$ and all other values $(w_i,w_j)$, $(w_i,x_j)$, $(x_i,x_j)$ equal 0. 
Define $e \in {\rm End}(V_{4m+2})$ by
\[
\begin{array}{llll}
e: & w_i & \mapsto & w_{i+2}\;\;(i<2m-1), \\
   & x_i & \mapsto & x_{i+2}\;\;(i<2m-1), \\
          & w_{2m-1}& \mapsto & 0,\\
          & x_{2m-1} & \mapsto & 0.
\end{array}
\]
Now $e$ is a nilpotent element of $so_{4m}(q)$. Set 
\[
W(2m)  = (1-e)(1+e)^{-1} \in G = \Or(V_{4m}).
\]
Observe that $W(2m)$ has Jordan form $J_{2m}^{[2]}$. 
Its matrix relative to the given basis is 
\[
\begin{pmatrix*}[r] I&-2I&2I& \cdots &-2I \\
                          & I & -2I & \cdots &2I \\
                          &&& \ddots & \\
                          &&&&I
\end{pmatrix*}
\]
where $I$ denotes the $2\times 2$ identity matrix. 
Define a 1-dimensional torus $T_1 \le \Or(\bar V_{4m})$ by the action
\[
T_1(c): w_i \mapsto c^iw_i,\,x_i \mapsto c^ix_i \;\;\hbox{ for all }i.
\]

\subsubsection{Class representatives}
For a given Jordan form (\ref{jso}), the conjugacy class representatives in $\Or(V)$ with this Jordan form are 
\begin{equation}\label{vwso}
\bigoplus_{i=1}^r (V_{\b_i}(2k_i+1) \oplus V_1(2k_i+1)^{[a_i-1]}) \oplus \bigoplus_{i=1}^s W(2l_i)^{[b_i]},
\end{equation}
where $V = \bigoplus_{i=1}^r V_{2k_i+1}^{[a_i]} \oplus \bigoplus_{i=1}^s V_{4l_i}^{[b_i]}$, and each $\b_i \in \{1,\a\}$ (and the $k_i$ are distinct, as are the $l_i$). 

Let $u$ be the element defined in (\ref{vwso}). Then $u$ fixes an orthogonal form of discriminant $D = \prod_{i=1}^r \b_i(-1)^{k_ia_i}$ (set $D=1$ if $r=0$). If $n=\dim V =2m$ is even, then $u \in \Or^\e_{2m}(q)$, where $\e = +$ if and only if $(-1)^mD$ is a square in $\F_q^*$ (see \cite[Prop.\ 2.5.10]{KL}). For the various choices of $\b_i$, these representatives fall into $2^{r-1}$ conjugacy classes in each of $\Or^+_{2m}(q)$ and $\Or^-_{2m}(q)$, except that for $r=0$, there is one class in $\Or^+_{2m}(q)$ and none in 
$\Or^-_{2m}(q)$. If $n$ is odd, then half of the representatives ($2^{r-1}$ of them) lie in a given $\Or_n(q)$ fixing a form 
of square discriminant.

Define $T$ to be the 1-dimensional torus acting as $T_1(c)$ on each subspace $\bar V_{2k_i}$ and  $\bar V_{4l_i+2}$.
To define corresponding matrix representatives, 
choose an ordered basis of $V$ according to the $T$-weights. 
Namely, list the spaces $V_{2k_j+1}, V_{4l_j}$ in decreasing order of 
Jordan block sizes. For the $i^{th}$ space, if it is $V_{\b}(2k+1)$, 
then label its basis as
\[
v_{-(2k)}^{(i)},v_{-(2k-2)}^{(i)},\ldots, v_{2k}^{(i)},
\]
and if it is $W(2l)$, then label its basis as
\[
w_{-(2l-1)}^{(i)},x_{-(2l-1)}^{(i)}, w_{-(2l-3)}^{(i)}, x_{-(2l-3)}^{(i)},\ldots, w_{2l-1}^{(i)}, x_{2l-1}^{(i)}.
\]
Now take a basis of $V$ consisting of these 
vectors $v_j^{(i)}, w_j^{(i)}, x_j^{(i)}$ in increasing order 
lexicographically with respect to the pairs $(j,i)$. 
Choose the matrix representative corresponding to the linear 
map (\ref{vwso}) to be its matrix with respect to this ordered basis.

Let $P$  be the parabolic subgroup stabilizing the flag of $V$ defined by sums of $T$-weight spaces for decreasing weights; so $P=QL$, with unipotent radical $Q$ and Levi subgroup $L= C_G(T)$. As in the 
previous section, we have the following result.

\begin{thm}\label{orthreps} 
Let $G = \Or^\e_n(q)$ with $q$ odd, and let $u\in G$ be the unipotent matrix representative defined above for the linear map
\[
\bigoplus_{i=1}^r (V_{\b_i}(2k_i+1) \oplus V_1(2k_i+1)^{[a_i-1]}) \oplus \bigoplus_{i=1}^s W(2l_i)^{[b_i]}.
 \]
Let $P = QL$ be the corresponding parabolic subgroup. 
\begin{itemize}
\item[{\rm (i)}] $C_G(u) = C_Q(u)C_L(u)$;
\item[{\rm (ii)}] $C_L(u) \cong \prod_{i=1}^s \Sp_{2b_i}(q) \times \prod_{i=1}^r \Or_{a_i}^{\e_i}(q)$;
\item[{\rm (iii)}] rewrite the Jordan form of $u$ as $\sum J_i^{r_i}$; 
then $|C_Q(u)| = q^\myexpR$ where 
\[
\myexpR = \frac{1}{2}\sum_i (i-1)r_i^2 + \sum_{i<j}ir_ir_j - \frac{1}{2}\sum_{i\;even}r_i.
\]
\end{itemize}
\end{thm}

In part (ii), the sign $\e_i=\pm$ is the sign of the orthogonal form with Gram matrix ${\rm diag}(\b_i,1,\ldots,1)$. In Section \ref{ocent} we justify this, and  show how to construct the subgroups $C_Q(u)$ and $C_L(u)$ of the centralizer. 

Theorem \ref{orthreps} gives the unipotent class representatives $u$  in the orthogonal group $G=\Or^\e_n(q)$. We conclude by 
describing how these classes split within the special orthogonal group $\SO^\e_n(q)$ and its subgroup $\O^\e_n(q)$ (see \cite[Prop.\ 2.4]{GLOB}). First, a class $u^G$ splits into two $\SO^\e_n(q)$-classes if and only if $\e=+$ and $r=0$; 
in this case, representatives are $u$ and $u^t$ for $t$ a 
reflection in $G$. 
Secondly, an $\SO^\e_n(q)$-class with representative $u$ (as in (\ref{vwso})) splits into two $\O_n^\e(q)$-classes if and only if either $r=0$, or $r\ge 1$ and the following hold:
\begin{itemize}
\item[(a)] $a_i=1$ for all $i$, and
\item[(b)] $\b_i \equiv (-1)^{k_1+k_i}\b_1 \hbox{ mod } (\F_q^*)^2$ for all $i$.
\end{itemize}
In case of splitting, representatives of the $\O_n^\e(q)$-classes are $u$ and $u^s$, where $s \in \SO^\e_n(q)\setminus \O^\e_n(q)$. 

\subsection{Unitary groups}\label{sugood}
Let $\GU_n(q)$ fix a unitary form on an $n$-dimensional
vector space $V$ over $\F_{q^2}$, 
and let $\a \mapsto \bar \a$ denote the involutory automorphism of $\F_{q^2}$. 
Here we take the unipotent class representatives directly from 
\cite[\S 2.3]{GLOB}. 
To define these, we need to describe unitary Jordan blocks of even and odd sizes.  
Fix $\b,\g \in \F_{q^2}^*$ with $\b+\bar \b = 0$, $\g+\bar \g = -1$.

\subsubsection{Even blocks $V(2m)$}  
Let $V_{2m}$ be a $2m$-dimensional space over $\F_{q^2}$ with basis 
$v_{-(2m-1)}$,
$v_{-(2m-3)}$,$\ldots,v_{2m-1}$, and unitary form defined by $(v_{-i},v_i) = 1$ for all $i$,  and all other $(v_i,v_j)=0$.  Define a unipotent element of the unitary group $\SU(V_{2m})$ as follows:
\[
\begin{array}{llll}
V(2m): & v_{-i} & \mapsto & v_{-i}+v_{-i+2}+\cdots +v_{-1}+\b v_1\;\;(i=1,3,\ldots, 2m-1) \\
           & v_j & \mapsto & v_j-v_{j+2} \;\;(j=1,3,\ldots,2m-1).
\end{array}
\]

\subsubsection{Odd blocks $V(2m+1)$} 
Let $V_{2m+1}$ be a $(2m+1)$-dimensional space over $\F_{q^2}$ with basis $v_{-2m},
v_{-(2m-2)},\ldots,v_{2m}$, and unitary form defined by $(v_{-i},v_i) = 1$ for all $i$, and all other $(v_i,v_j)=0$.  Define a unipotent element of the unitary group $\SU(V_{2m+1})$ as follows:
\[
\begin{array}{llll}
V(2m+1): & v_{-i} & \mapsto & v_{-i} + v_{-i+2}+\cdots +v_{0}+\g v_2\;\;(i=2,4,\ldots, 2m) \\
           & v_j & \mapsto & v_j-v_{j+2} \;\;(j=0,2,\ldots,2m).
\end{array}
\]

\subsubsection{Class representatives}  
Denote by $\bigoplus_{i=1}^s V(m_i)^{[r_i]}$ the orthogonal direct sum of copies of the linear maps defined above acting on $V = \bigoplus_{i=1}^s V_{m_i}^{[r_i]}$, where $m_1>m_2>\cdots >m_s$. This has Jordan form $\bigoplus_{i=1}^s J_{m_i}^{[r_i]}$, and 
there is a unique conjugacy class in $\GU(V)$ having this Jordan form. The number of conjugacy classes in $\SU(V)$ with this Jordan form is ${\gcd}(t,q+1)$, where $t = {\gcd}(m_1,\ldots,m_s)$ (see \cite[Prop.\ 2.2]{GLOB}). Representatives can be found by taking conjugates of $\bigoplus_{i=1}^s V(m_i)^{[r_i]}$ by matrices $d(\a) \in \GU(V)$ of determinant $\a$, where $Z = \{\mu \in \F_{q^2}: \mu \bar \mu =1\}$ and $\a$ ranges over representatives of $Z/Z^t$. For example, we can take $d(\a)$ to act on the basis of one of the blocks as ${\rm diag}(\l,1,\ldots,1,\bar \l^{-1})$, where $\a = \l\bar \l^{-1}$.

To define matrix representatives, list the spaces $V_{m_j}$ in decreasing order of sizes, and for the $i^{th}$ space, 
if it is $V(m)$, label its basis as $v_{-(m-1)}^{(i)},\ldots, v_{m-1}^{(i)}$. 
Now take a basis of $V$ consisting of these vectors 
$v_j^{(i)}$ for all $i,j$, in increasing order lexicographically with respect 
to  the pairs $(j,i)$. Choose the matrix representative $u$ 
corresponding to $\bigoplus_{i=1}^s V(m_i)^{[r_i]}$ to be the matrix of 
the linear map with respect to this ordered basis. 

For each $i$, define $\myV_{(i)}$ to be the span of the basis vectors 
with subscript $i$, and call this the weight space for $i$.
As in previous sections, let $P$ be the parabolic subgroup stabilizing the flag of $V$ 
defined by sums of weight spaces for decreasing weights; 
so $P=QL$, where $Q$ is the unipotent radical and $L$ a Levi subgroup.
The next result follows from \cite[Lemma 2.29 and Props.\ 3.7, 3.8]{LS}.

\begin{thm}\label{sureps} 
Let $G = \GU_n(q)$, let $u\in G$ be the 
unipotent matrix defined above with Jordan form 
$\bigoplus_{i=1}^s J_{m_i}^{[r_i]}$, and let $P = QL$ be 
the corresponding parabolic subgroup. 
\begin{itemize}
\item[{\rm (i)}] $C_G(u) = C_Q(u)C_L(u)$;
\item[{\rm (ii)}] $C_L(u) \cong \prod_{i=1}^s \GU_{r_i}(q)$;
\item[{\rm (iii)}] $|C_Q(u)| = q^\myexpR$, where $\myexpR = \sum_{i} (i-1)r_i^2 + 2\sum_{i<j}ir_ir_j$.
\end{itemize}
\end{thm}

\section{Centralizers of class representatives}\label{gencentgood}
Having described the unipotent class representatives, 
we now show how to construct their centralizers, whose structures 
are given in Theorems 
\ref{sympreps}, \ref{orthreps} and \ref{sureps}. 
In these results, for each class representative $u$ 
of an isometry
group $G$, we defined a parabolic 
subgroup $P = QL$ such that $C_G(u) = C_Q(u)C_L(u)$ and gave
the structure of $C_L(u)$ and the order of $C_Q(u)$. 

We construct $C_Q(u)$ and $C_L(u)$ in the following steps:
\begin{enumerate}
\item[(a)] construct a generating set for $Q$;
\item[(b)] construct a {\it power-conjugate presentation}
\index{power-conjugate presentation}
for $Q$ and use it to compute generators for $C_Q(u)$;
\item[(c)] construct a generating set for $C_L(u)$.
\end{enumerate}

We comment briefly on Step (b). 
Recall that $Q$ is a unipotent matrix group.
A group of prime-power order can be described by a power-conjugate
presentation which reflects a central 
series; various problems -- including deciding conjugation and 
constructing centralizers -- can be solved easily in practice  
by exploiting this presentation.  It is a routine computational exercise to 
construct an isomorphic copy $\bar{Q}$ of $Q$ described by 
such a presentation; the centralizer of an element
of $Q$ can then be readily constructed in $\bar{Q}$,
and so in $Q$.
For related details, see \cite[Chap.\ 9]{handbook}.
 
In the rest of the section, we describe Steps (a) and (c) for 
each type of classical group.

\subsection{Symplectic groups} \label{spcent}
Let $G = \Sp(V) \cong \Sp_{2n}(q)$, and let $u$ be the matrix 
representative of the element 
\begin{equation}\label{urep}
\bigoplus_{i=1}^r (V_{\b_i}(2k_i) \oplus V_1(2k_i)^{[a_i-1]}) \oplus \bigoplus_{i=1}^s W(2l_i+1)^{[b_i]},
\end{equation}
of $\Sp(V)$ as in (\ref{vwsp}), relative to the ordered basis of vectors $v_j^{(i)}, w_j^{(i)}, x_j^{(i)}$ described in Section \ref{spgood}. Let $P = QL$ be the parabolic subgroup defined there.

\subsubsection{Generators of $Q$}  
We define two collections of generators for $Q$.

\vspace*{2mm} \no (1)  
For all $i,j$ with $i<0$, $i\le j<-i$, all $k,l$, and every $\l \in \F_q$ 
define the following generators 
which fix all basis vectors other than those listed.
Set  $\e = (-1)^{\lfloor i/2\rfloor+{\lfloor j/2\rfloor }}$.
\[
\begin{array}{ll}
x_{1ijkl}(\l): & v_i^{(k)} \mapsto v_i^{(k)}+\l v_{-j}^{(l)}, \\
                     & v_j^{(l)}\mapsto v_j^{(l)}+ \e\l v_{-i}^{(k)}\\
\\
x_{2ijkl}(\l): & v_i^{(k)} \mapsto v_i^{(k)}+\l x_{-j}^{(l)}, \\
                     & w_j^{(l)}\mapsto w_j^{(l)}+ \e\l v_{-i}^{(k)}\\
\\
x_{3ijkl}(\l): & v_j^{(k)} \mapsto v_j^{(k)}+\l x_{-i}^{(l)}, \\
                     & w_i^{(l)}\mapsto w_i^{(l)}+ \e\l v_{-j}^{(k)}\\
\\
x_{4ijkl}(\l): & v_i^{(k)} \mapsto v_i^{(k)}+\l w_{-j}^{(l)}, \\
                     & x_j^{(l)}\mapsto x_j^{(l)}- \e\l v_{-i}^{(k)}\\
\\
x_{5ijkl}(\l): & v_j^{(k)} \mapsto v_j^{(k)}+\l w_{-i}^{(l)}, \\
                     & x_i^{(l)}\mapsto x_i^{(l)}- \e\l v_{-j}^{(k)}\\
\\
x_{6ijkl}(\l): & w_i^{(k)} \mapsto w_i^{(k)}+\l x_{-j}^{(l)}, \\
                     & w_j^{(l)}\mapsto w_j^{(l)}+ \e\l x_{-i}^{(k)}\\
\\
x_{7ijkl}(\l): & x_i^{(k)} \mapsto x_i^{(k)}+\l w_{-j}^{(l)}, \\
                     & x_j^{(l)}\mapsto x_j^{(l)}+ \e\l w_{-i}^{(k)}\\
\\
x_{8ijkl}(\l): & w_i^{(k)} \mapsto w_i^{(k)}+\l w_{-j}^{(l)}, \\
                     & x_j^{(l)}\mapsto x_j^{(l)}- \e\l x_{-i}^{(k)}\\
\\
x_{9ijkl}(\l): & x_i^{(k)} \mapsto x_i^{(k)}+\l x_{-j}^{(l)}, \\
                     & w_j^{(l)}\mapsto w_j^{(l)}- \e\l w_{-i}^{(k)}.
\end{array}
\]

\vspace*{2mm} \no (2) For all $i>0$, all $k$, and every $\l \in \F_q$:
\[
\begin{array}{ll}
y_{1ik}(\l): & v_{-i}^{(k)} \mapsto v_{-i}^{(k)}+\l v_{i}^{(k)}, \\
\\
y_{2ik}(\l): & w_{-i}^{(k)} \mapsto w_{-i}^{(k)}+\l x_{i}^{(k)}, \\
\\
y_{3ik}(\l): & x_{-i}^{(k)} \mapsto x_{-i}^{(k)}+\l w_{i}^{(k)}.
\end{array}
\]

Define $X_{cijkl} = \{x_{cijkl}(\l): \l \in \F_q\} \cong \F_q^+$ for each of the generators listed; similarly define  $Y_{mik}$.
Of course, to construct each subgroup 
we need only write down generators 
for those $\l$ in a 
basis of $\F_q$ over $\F_p$. 

To check that these subgroups generate $Q$, 
observe that each lies in $Q$, and the total number of subgroups 
is $m$, where $|Q|=q^m$. 

\begin{example}
We return to Example \ref{sp-example}. Here 
\[
u = V_\a(4) \oplus W(3) \oplus V_1(2) \in \Sp_{12}(q).
\]
The parabolic subgroup $P=QL$ has $L = \GL_1(q)\times \GL_2(q) \times \GL_2(q) \times \Sp_2(q)$, and $|Q| = q^{33}$. We record the number of subgroups $X_{1ijkl}$, $\ldots$, $X_{9ijkl}$, $Y_{1ik}$, $Y_{2ik}$, and $Y_{3ik}$:
\[
\begin{array}{ll}
\hline 
\hbox{subgroups} & \hbox{number} \\
\hline
X_{1ijkl} & 4 \\
X_{2ijkl} & 5 \\
X_{3ijkl} & 4 \\
X_{4ijkl} & 5 \\
X_{5ijkl} & 4 \\
X_{6ijkl} & 1 \\
X_{7ijkl} & 1 \\
X_{8ijkl} & 2 \\
X_{9ijkl} & 2 \\
Y_{1ik}  & 3 \\
Y_{2ik}  & 1 \\
Y_{3ik} & 1 \\
\hline
\end{array}
\]
The total number of subgroups is 33, as claimed.
\end{example}

\subsubsection{Generators of $C_L(u)$} 
For $u\in \Sp_{2n}(q)$ as in (\ref{urep}), Theorem \ref{sympreps} gives 
\[
C_L(u) \cong \prod_{i=1}^s \Sp_{2b_i}(q) \times \prod_{i=1}^r \Or_{a_i}^{\e_i}(q).
\]
Each factor $\Sp_{2b_i}(q)$ acts on the summand $W(2l_i+1)^{[b_i]}$, fixing the orthogonal complement; and each factor 
$ \Or_{a_i}^{\e_i}(q)$ acts on the summand $V_{\b_i}(2k_i) \oplus V_1(2k_i)^{[a_i-1]}$. 
So it suffices to focus on these summands. 
For notational convenience, 
we denote them by  $W(2l+1)^{[b]}$ and $V_{\b}(2k) \oplus V_1(2k)^{[a-1]}$.

\vspace*{4mm}\no {\it Summand $W(2l+1)^{[b]}$.}  
Let $u = W(2l+1)^{[b]}$. We need to find generators for $C_L(u) = \Sp_{2b}(q)$. The ordered basis of $W(2l+1)^{[b]}$ is
\[
w_{-2l}^{(1)},x_{-2l}^{(1)},\ldots, w_{-2l}^{(b)},x_{-2l}^{(b)},\ldots ,w_{2l}^{(1)},x_{2l}^{(1)},\ldots, 
w_{2l}^{(b)},x_{2l}^{(b)}.
\]
With respect to this basis, the symplectic form has matrix
\[
 \begin{pmatrix} &&&K \\&&-K& \\& \iddots && \\ K&&& \end{pmatrix}
\]
where $K = J^{[b]}$ (block diagonal sum of $b$ copies of $J = \begin{pmatrix} &-1 \\ 1& \end{pmatrix}$). With respect to the above basis, 
the factor $L_0 \cong \Sp_{2b}(q)$ of $C_L(u)$ is
\[
L_0 = \left\{ \begin{pmatrix} A&&& \\ &A&& \\ &&\ddots & \\ &&&A \end{pmatrix} : A \in \Sp_{2b}(q) \right\},
\]
where $\Sp_{2b}(q)$ stabilizes the symplectic form with matrix $K$. We choose standard generators for this group (as mentioned at the end of Section \ref{fincla}), and hence obtain generators for $L_0$.

\vspace*{4mm} \no {\it Summand $V_{\b}(2k) \oplus V_1(2k)^{[a-1]}$.} 
Here the basis is   
\[
v_{-(2k-1)}^{(1)},\ldots, v_{-(2k-1)}^{(a)}, \ldots, v_{2k-1}^{(1)},\ldots, v_{2k-1}^{(a)}.
\]
Change the basis as follows: $v_1^{(1)} \mapsto \b v_1^{(1)}, \ldots, v_{2k-1}^{(1)} \mapsto \b v_{2k-1}^{(1)}$. With respect to this basis, 
the factor $L_0 \cong \Or_a^\e(q)$ is 
\[
L_0 = \left\{ \begin{pmatrix} A&&& \\ &A&& \\ &&\ddots & \\ &&&A \end{pmatrix} : A \in \Or_a^\e(q) \right\},
\]
and $\Or_a^\e(q)$ stabilizes the orthogonal form with matrix ${\rm diag}(\b,1,\ldots,1)$. Now choose standard generators for this group, and then change back to the original basis to obtain the required factor $\Or_a^\e(q)$ of 
$C_L(u)$.

\subsection{Orthogonal groups} \label{ocent}
This is similar to the previous section. Let $G = \Or(V) \cong \Or^\e_{2n}(q)$, and let $u$ be the matrix representative of the element 
\begin{equation}\label{orep}
\bigoplus_{i=1}^r (V_{\b_i}(2k_i+1) \oplus V_1(2k_i+1)^{[a_i-1]}) \oplus \bigoplus_{i=1}^s W(2l_i)^{[b_i]},
\end{equation}
of $\Or(V)$ as in (\ref{vwso}), relative to the ordered basis of vectors $v_j^{(i)}, w_j^{(i)}, x_j^{(i)}$ described in Section \ref{sogood}. Let $P = QL$ be the parabolic subgroup defined there.

\subsubsection{Generators of $Q$}  
We define two collections of generators for $Q$.

\vspace*{2mm} \no (1)  For all $i,j$ with $i<0$, $i\le j<-i$ and $j\ne 0$, 
all $k,l$, and every $\l \in \F_q$ define the following generators
which fix all basis vectors other than those listed.
Set  $\e_1 = (-1)^{\lfloor i/2\rfloor+{\lfloor j/2\rfloor }}$, 
$\e_2 = (-1)^{\lfloor i/2\rfloor+{\lfloor -j/2\rfloor }}$, 
and $\e_3 = (-1)^{\lfloor -i/2\rfloor+{\lfloor -j/2\rfloor }}$.
\[
\begin{array}{ll}
x_{1ijkl}(\l): & v_i^{(k)} \mapsto v_i^{(k)}+\l v_{-j}^{(l)}, \\
                     & v_j^{(l)}\mapsto v_j^{(l)}- \e_1\l v_{-i}^{(k)}\\
\\
x_{2ijkl}(\l): & v_i^{(k)} \mapsto v_i^{(k)}+\l x_{-j}^{(l)}, \\
                     & w_j^{(l)}\mapsto w_j^{(l)}- \e_1\l v_{-i}^{(k)}\\
\\
x_{3ijkl}(\l): & v_j^{(k)} \mapsto v_j^{(k)}+\l x_{-i}^{(l)}, \\
                    & w_i^{(l)}\mapsto w_i^{(l)}- \e_1 \l v_{-j}^{(k)}\\
\\
x_{4ijkl}(\l): & v_i^{(k)} \mapsto v_i^{(k)}+\l w_{-j}^{(l)}, \\
                     & x_j^{(l)}\mapsto x_j^{(l)}- \e_1\l v_{-i}^{(k)}\\
\\
x_{5ijkl}(\l): & v_j^{(k)} \mapsto v_j^{(k)}+\l w_{-i}^{(l)}, \\
                     & x_i^{(l)}\mapsto x_i^{(l)}- \e_1 \l v_{-j}^{(k)}\\
\\
x_{6ijkl}(\l): & w_i^{(k)} \mapsto w_i^{(k)}+\l x_{-j}^{(l)}, \\
                     & w_j^{(l)}\mapsto w_j^{(l)}- \e_1\l x_{-i}^{(k)}\\
\\
x_{7ijkl}(\l): & x_i^{(k)} \mapsto x_i^{(k)}+\l w_{-j}^{(l)}, \\
                     & x_j^{(l)}\mapsto x_j^{(l)}- \e_3\l w_{-i}^{(k)}\\
\\
x_{8ijkl}(\l): & w_i^{(k)} \mapsto w_i^{(k)}+\l w_{-j}^{(l)}, \\
                     & x_j^{(l)}\mapsto x_j^{(l)}- \e_2\l x_{-i}^{(k)}\\
\\
x_{9ijkl}(\l): & x_i^{(k)} \mapsto x_i^{(k)}+\l x_{-j}^{(l)}, \\
                     & w_j^{(l)}\mapsto w_j^{(l)}- \e_2 \l w_{-i}^{(k)}.
\end{array}
\]

\vspace*{2mm} \no (2) 
For all $i>0$, all $k,l$, and every $\l \in \F_q$, 
define the following generators
which fix all basis vectors other than those listed.
Set $\mu_\e= -(-1)^{\lfloor \e i/2 \rfloor}\b \l$, 
and $\a_\e = -\frac{1}{2}(-1)^{\lfloor \e i/2 \rfloor }\b\l^2$, 
where $\e = \pm$.
\[
\begin{array}{rlll}
y_{1ikl}(\l): & v_{-i}^{(k)} & \mapsto & v_{-i}^{(k)}+\l v_0^{(l)} + \a_+ v_{i}^{(k)}, \\
                  & v_0^{(l)} & \mapsto & v_0^{(l)}+\mu_+v_i^{(k)} 
\end{array}
\]

\[
\begin{array}{rlll}
y_{2ikl}(\l): & w_{-i}^{(k)} & \mapsto & w_{-i}^{(k)}+\l v_0^{(l)}+ \a_- x_{i}^{(k)}, \\
                  &  v_0^{(l)} & \mapsto & v_0^{(l)}+\mu_-x_i^{(k)} \\
\\
y_{3ikl}(\l): & x_{-i}^{(k)} & \mapsto & x_{-i}^{(k)}+\l v_0^{(l)}+ \a_+ w_{i}^{(k)}, \\
                  &  v_0^{(l)} & \mapsto & v_0^{(l)}+\mu_+w_i^{(k)}.
\end{array}
\]

\vspace*{2mm} As in the previous section, 
define $X_{cijkl} = \{x_{cijkl}(\l): \l \in \F_q\} \cong \F_q^+$ for each of the generators listed; similarly define  $Y_{mikl}$.
Each of these subgroups lies in $Q$, and the total number 
is $m$, where $|Q|=q^m$. 

\subsubsection{Generators of $C_L(u)$} 

For $u\in \Or_n^\e(q)$ as in (\ref{orep}), Theorem \ref{orthreps} gives 
\[
C_L(u) \cong \prod_{i=1}^s \Sp_{2b_i}(q) \times \prod_{i=1}^r \Or_{a_i}^{\e_i}(q).
\]
Each factor $\Sp_{2b_i}(q)$ acts on the summand $W(2l_i)^{[b_i]}$, fixing the orthogonal complement; and each factor 
$ \Or_{a_i}^{\e_i}(q)$ acts on the 
summand $V_{\b_i}(2k_i+1) \oplus V_1(2k_i+1)^{[a_i-1]}$. 
So it suffices to focus on these summands. 
For notational convenience, we denote them by  $W(2l)^{[b]}$ and $V_{\b}(2k+1) \oplus V_1(2k+1)^{[a-1]}$.

\vspace*{4mm}\no {\it Summand $W(2l)^{[b]}$.}  
Let $u = W(2l)^{[b]}$. As in the previous section, with respect to the basis $w_i^{(j)},x_i^{(j)}$, the factor $L_0 \cong \Sp_{2b}(q)$ of $C_L(u)$ is
\[
L_0 = \left\{ \begin{pmatrix} A&&& \\ &A&& \\ &&\ddots & \\ &&&A \end{pmatrix} : A \in \Sp_{2b}(q) \right\},
\]
where $\Sp_{2b}(q)$ stabilizes the symplectic form with matrix $K = J^{[b]}$, as before. We choose standard generators for this group, and 
hence obtain generators for $L_0$.

\vspace*{4mm} \no {\it Summand $V_{\b}(2k+1) \oplus V_1(2k+1)^{[a-1]}$.} 
Here the basis is   
\[
v_{-(2k)}^{(1)},\ldots, v_{-(2k)}^{(a)}, \ldots, v_{2k}^{(1)},\ldots, v_{2k}^{(a)}.
\]
Change the basis as follows: $v_2^{(1)} \mapsto \b v_2^{(1)}, \ldots, v_{2k}^{(1)} \mapsto \b v_{2k}^{(1)}$. 
Let $\Or_a^\d (q)$ be the orthogonal group stabilizing the orthogonal form with 
matrix ${\rm diag}(\b,1,\ldots,1)$.
With respect to the above basis, 
the factor $L_0 \cong \Or_a^\d (q)$ is 
\[
L_0 = \left\{ \begin{pmatrix} A&&& \\ &A&& \\ &&\ddots & \\ &&&A \end{pmatrix} : A \in \Or_a^\d(q) \right\}.
\]
Now choose standard generators for 
this group, and then change back to the original basis to obtain 
the required factor $\Or_a^\d (q)$ of $C_L(u)$.

We can use Schreier's Lemma and the consequent
algorithm\index{Schreier algorithm}
\cite[Chap.\ 2]{handbook} to write down the kernel
of an explicit homomorphism from a group into a cyclic group.
Let $\Spec=\SO^\epsilon(V)$ and $\O=\Omega^\epsilon(V)$. 
Both $C_{G}(u)/C_{\Spec}(u)$ and $C_{\Spec}(u)/C_{\O}(u)$ are cyclic,
and are the images of the determinant and spinor norm
homomorphisms respectively.  We apply this algorithm to construct
generating sets for $C_\Spec(u)$ and $C_\O(u)$.

\subsection{Unitary groups} \label{ucent}
Let $G = \GU(V) \cong \GU_n(q)$, and let $u$ be the matrix representative of the element $\bigoplus_{i=1}^s V(m_i)^{[r_i]}$
relative to the ordered basis of vectors $v_j^{(i)}$ described in Section \ref{sugood}. Let $P = QL$ be the parabolic subgroup defined there.

\subsubsection{Generators of $Q$}  

\vspace*{2mm} \no (1)  For all $i,j>0$ with $-i\le j<i$ and $j\ne 0$, all $k,l$, and every $\l \in \F_{q^2}$ define 
\[
\begin{array}{ll}
x_{1ijkl}(\l): & v_{-i}^{(k)} \mapsto  v_{-i}^{(k)} + \l  v_{-j}^{(l)} ,\\
                   & v_{j}^{(l)} \mapsto  v_{j}^{(l)} - \bar \l  v_{i}^{(k)}.
\end{array}
\]

\no (2) If there are no weight 0 vectors in the basis 
(no vectors $v_0^{(i)}$), then for all $i<0$, all $k$, and 
every $\a \in \F_{q^2}$ such that $\a+\bar \a = 0$, define
\[
\begin{array}{ll}
x_{2ikl}(\l): & v_{-i}^{(k)} \mapsto  v_{-i}^{(k)} + \a  v_{i}^{(k)}.
\end{array}
\]
  
\no (3)  If there are weight 0 vectors in the basis, then for all $i<0$, all $k,l$, and every $\l,\a \in \F_{q^2}$ such that $\l\bar \l+\a+\bar \a = 0$, define
\[
\begin{array}{ll}
x_{3ikl}(\l,\a): & v_{-i}^{(k)} \mapsto  v_{-i}^{(k)} + \l v_0^{(l)}+\a  v_{i}^{(k)} ,\\
                   & v_{0}^{(l)} \mapsto  v_{0}^{(l)} - \bar \l  v_{i}^{(k)}.
\end{array}
\]

\no 
As in previous sections, these elements generate $Q$.

\subsubsection{Generators of $C_L(u)$} 
Let $u=\bigoplus_{i=1}^s V(m_i)^{[r_i]}\in \GU_n(q)$. 
Theorem \ref{sureps} gives 
\[
C_L(u) \cong  \prod_{i=1}^s \GU_{r_i}(q).
\]
We work in a summand $V(m_i)^{[r_i]}$, which we write as $V(m)^{[r]}$ for simplicity, with basis
\[
v_{-(m-1)}^{(1)},\ldots, v_{-(m-1)}^{(r)}, \ldots, v_{m-1}^{(1)},\ldots, v_{m-1}^{(r)},
\]
where the unitary form is as in Section \ref{sugood}. 
Now $C_L(u)$ is 
\[
\left\{ \begin{pmatrix} A&&& \\ &A&& \\ &&\ddots & \\ &&&A \end{pmatrix} : A \in \GU_r(q) \right\},
\]
where $\GU_r(q)$ stabilizes the unitary form with matrix $I_r$. We use 
standard generators for $\GU_r(q)$ to construct $C_L(u)$.

To obtain $C_{\SU_r(q)}(u)$ we must construct $C := C_L(u) \cap\SU_r(q)$.
Using Schreier's algorithm, 
we construct $C$ as the kernel of the determinant homomorphism from $C_L(u)$ into $\F_{q^2}^*$.

\section{The conjugacy problem}\label{conjtestgood}
In this section we solve the conjugacy problem: for a unipotent element
$g$ of a classical group $G$, identify the class representative defined 
in Section \ref{goodrep} that is conjugate to $g$. 
By two applications of our solution, we then solve the general problem: 
given unipotent $g,h \in G$, is $g$ conjugate to $h$?

\subsection{Symplectic groups} \label{spconprob}
Let $G = \Sp(V) \cong \Sp_{2n}(q)$ with $q$ odd, 
preserving a symplectic form $(\,,\,)$, and 
let $g \in G$ be unipotent. We want to compute the representative 
\begin{equation}\label{uvrep}
 \bigoplus_{i=1}^r (V_{\b_i}(2k_i) \oplus V_1(2k_i)^{[a_i-1]}) \oplus \bigoplus_{i=1}^s W(2l_i+1)^{[b_i]},
\end{equation}
as defined in Section \ref{spgood}, such that $g$ is $G$-conjugate to 
this representative.  The values $k_1,\ldots, k_r$ and 
$l_1,\ldots, l_s$ are determined by the
Jordan form of $g$. It remains to determine the 
values of the parameters $\b_i \in \{1,\a\}$ in (\ref{uvrep}).

Our algorithm to do this uses the following two lemmas. 
The {\it radical} \index{radical}
of $F \leq V$ is $\hbox{Rad}(F) = F \cap F^\perp$.

\begin{lem}\label{spalem1} Let $u = V_\b(2k) \in \Sp(V)$ with $V = V_{2k}$, as defined in Section $\ref{spgood}$. Let $0\le r < k$, and $F = V(1-u)^r$. 
\begin{itemize}
\item[{\rm (i)}] $u^{F/{\rm Rad}(F)}$, 
the element induced by the action of $u$ on the symplectic space 
$F/{\rm Rad}(F)$, is $V_\b(2k-2r)$.
\item[{\rm (ii)}] If $r=2$, then $F^\perp$ is totally isotropic.
\end{itemize}
\end{lem}

\begin{proof}
(i)  From the definition,  $V_{2k}$ has basis $v_{-(2k-1)},v_{-(2k-3)},\ldots, v_{2k-1}$, and  symplectic form given by $(v_i,v_{-i}) = (-1)^{\lfloor i/2\rfloor}$ for all $i$, and $u = (1-e_\b(2k))(1+e_\b(2k))^{-1}$ where $e_\b(2k)$ is the nilpotent linear map defined by
\[
\begin{array}{llll}
e_\b(2k): & v_i & \mapsto & v_{i+2}\;\;(i<2m-1,\,i\ne -1), \\
          & v_{-1}& \mapsto & \b v_1,\\
          & v_{2m-1}&  \mapsto & 0.
\end{array}
\]
Since $e_\b(2k) = (1+u)^{-1}(1-u)$, 
\[
F = V(1-u)^r = Ve_\b(2k)^r = \la v_{-(2k-2r-1)},\ldots, v_{2k-1}\ra.
\]
Then ${\rm Rad}(F) = \la v_{2k-2r+1},\ldots, v_{2k-1}\ra$ and so $e_\b(2k)$ acts on $F/{\rm Rad}(F)$ as $e_\b(2k-2r)$. Hence $u$ acts on $F/{\rm Rad}(F)$ as $(1-e_\b(2k-2r))(1+e_\b(2k-2r))^{-1} = V_\b(2k-2r)$. 

(ii) Let $r=2$, so $F = V(1-u)^2 = \la v_{-(2k-5)},\ldots, v_{2k-1}\ra$. Then $F^\perp = \la v_{2k-3},v_{2k-1} \ra$, which is totally isotropic.  \end{proof}

\begin{lem}\label{spalem2} Let $u = W(2l+1) \in \Sp(V)$ with $V = V_{4l+2}$, as defined in Section $\ref{spgood}$. Let $0\le r < 2l+1$, and $F = V(1-u)^r$. 
\begin{itemize}
\item[{\rm (i)}] If $r\le l$, then $u^{F/{\rm Rad}(F)} = W(2l-2r+1)$.
\item[{\rm (ii)}] If $r>l$, then $F$ is totally isotropic.
\item[{\rm (iii)}] If $r= 2$ and $l\ge 2$, then $F^\perp$ is totally isotropic.
\item[{\rm (iv)}] If $r=2$ and $l=1$, then $u^{F^\perp/{\rm Rad}(F^\perp)} = W(1)$.
\item[{\rm (v)}] If $r=2$ and $l=0$, then $F=0$.
\end{itemize}
\end{lem}

\begin{proof} As in the previous lemma, the proof is based on 
simple computations using the basis for $W(2l+1)$ and  the nilpotent element $e$ defined in Section \ref{spgood}. 
\end{proof}

\subsubsection{The algorithm}
Assume that $g \in G = \Sp(V)$ is unipotent  and conjugate to 
\[
u =  \bigoplus_{i=1}^r (V_{\b_i}(2k_i) \oplus V_1(2k_i)^{[a_i-1]}) \oplus \bigoplus_{i=1}^s W(2l_i+1)^{[b_i]},
\]
for some $\b_i \in \{1,\a\}$. We now give an algorithm to compute the parameters $\b_i$. To simplify notation, write 
\[
V_{\b,a}(2k) = V_{\b}(2k) \oplus V_1(2k)^{[a-1]}.
\]

\vspace*{2mm}
\no {\it Subcase} $I$. Assume that $g$ is conjugate 
to $u = V_{\b,a}(2) \oplus W(1)^{[b]}$. We compute $\b$ as follows. 
First, find $e \in V$ such that 
\[
(e,\,e(1-g)) \ne 0.
\]
(The number of such vectors is $1-O(\frac{1}{q})$.) 
Write $f = e(1-g)$ and $\g_1 = (e,f)$, and define 
$F = \la e,f\ra$. Then $F$ is a non-degenerate $g$-invariant subspace, 
and with respect to the basis $[e,-\g_1^{-1}f]$ the matrix of the 
form is $\begin{pmatrix} &-1 \\1& \end{pmatrix}$, and 
\[
g^{F} = \begin{pmatrix*}[r] 1&-\g_1 \\ 0&1 \end{pmatrix*}.
\]
Hence $g$ acts on $F$ as $V_{\frac{1}{2}\g_1}(2)$. 
Now consider $F^\perp$ and repeat. Continuing, 
we compute values $\g_1,\ldots ,\g_a$ such that $g$ is conjugate to 
\[
W(1)^b \oplus \bigoplus_{i=1}^a V_{\frac{1}{2}\g_i}(2).
\]
We conclude that $\b \equiv 2^a\g_1\cdots \g_a \hbox{ mod }\F_q^2$. 

\vspace*{2mm}
\no {\it General case.} Now assume that $g$ is conjugate to 
\[
u =  \bigoplus_{i=1}^r V_{\b_i,a_i}(2k_i) \oplus U,
\]
where $k_1<k_2<\cdots <k_r$, and $U$ is a direct sum of blocks 
of the form $W(2l_i+1)$. We give an algorithm to compute the parameters $\b_i$.

\begin{enumerate}
\item[(1)] Compute $F = V(1-g)^{k_1-1}$ and the action 
of $g$ on $F_0 = F/{\rm Rad}(F)$. By Lemmas \ref{spalem1} and \ref{spalem2}, 
\begin{equation}\label{vf0}
g^{F_0} = V_{\b_1,a_1}(2) \oplus \bigoplus_{i=2}^r V_{\b_i,a_i}(2k_i') \oplus U',
\end{equation}
where $U'$ is a direct sum of blocks of the form $W(2l_i + 1)$, 
and $k_i' = k_i-k_1+1$ (so $1<k_2'< \cdots < k_r'$).

\item[(2)] Compute $F_1 = F_0(1-g)^2$ and the action 
of $g$ on $F_2 = F_1^\perp/{\rm Rad}(F_1^\perp)$. 
By the lemmas, 
\[
g^{F_2} = W(1)^{[b]} \oplus V_{\b_1,a_1}(2).
\]
Apply the algorithm of Subcase $I$ to compute $\b_1$.

\item[(3)]
Consider $g^{F_0}$ in (\ref{vf0}). Compute $F_3 = F_0(1-g)$ and the action of $g$ on $F_4 = F_3/{\rm Rad}(F_3)$. By the lemmas, 
\[
g^{F_4} = \bigoplus_{i=2}^r V_{\b_i,a_i}(2k_i'') \oplus U'',
\]
where $U''$ is a direct sum of blocks of the form $W(2l_i + 1)$, 
and $k_i''=k_i'-1$.
Now repeat Steps (1) and (2) to compute $\b_2$.
\end{enumerate}
Iterating, we  compute all of the $\b_i$. 

\subsection{Orthogonal groups} \label{oconprob}
Let $G = \Or(V) \cong \Or_{n}^\e(q)$ with $q$ odd, preserving an orthogonal form $(\,,\,)$, and let $g \in G$ be unipotent. We want to compute the representative 
\begin{equation}\label{uvorep}
 \bigoplus_{i=1}^r (V_{\b_i}(2k_i+1) \oplus V_1(2k_i+1)^{[a_i-1]}) \oplus \bigoplus_{i=1}^s W(2l_i)^{[b_i]},
\end{equation}
as defined in Section \ref{sogood}, such that $g$ is $G$-conjugate to this representative. (Refining to conjugacy in $\SO(V )$ and $\O(V)$ is discussed at the end of this section.) The values $k_1,\ldots, k_r$ and $l_1,\ldots, l_s$ are determined by the Jordan form of $g$. It remains to determine the values 
of the parameters $\b_i \in \{1,\a\}$ in (\ref{uvorep}).

Our algorithm to do this uses the following lemma.

\begin{lem}\label{olem} 
\mbox{}
\begin{enumerate}
\item[{\rm (a)}]  Let $u = V_\b(2k+1) \in \Or(V)$ with $V = V_{2k+1}$, as defined in Section $\ref{sogood}$. Let $0\le r < k$, and $F = V(1-u)^r$. 
\begin{itemize}
\item[{\rm (i)}] $u^{F/{\rm Rad}(F)} = V_\b(2k-2r+1)$.
\item[{\rm (ii)}] If $r=1$, then $F^\perp$ is totally isotropic.
\end{itemize}
\item[
{\rm (b)}]  Let $u = W(2l) \in \Or(V)$ with $V = V_{4l}$, as defined in Section $\ref{sogood}$. Let $0\le r < 2l$, and $F = V(1-u)^r$. 
\begin{itemize}
\item[{\rm (i)}] If $r\le l$, then $u^{F/{\rm Rad}(F)} = W(2l-2r)$.
\item[{\rm (ii)}] If $r>l$, then $F$ is totally isotropic.
\item[{\rm (iii)}] If $r= 1$, then $F^\perp$ is totally isotropic.
\end{itemize}
\end{enumerate}
\end{lem}

\begin{proof}
As in the proof of Lemma \ref{spalem1}, this is based on simple computations using the bases and nilpotent elements defining  $V_\b(2k+1)$ and $W(2l)$  in Section \ref{sogood}. 
\end{proof}

\subsubsection{The algorithm}
Assume that $g \in G = \Sp(V)$ is unipotent  and conjugate to 
\[
u =  \bigoplus_{i=1}^r (V_{\b_i}(2k_i+1) \oplus V_1(2k_i+1)^{[a_i-1]}) \oplus U
\]
where $k_1<k_2<\cdots <k_r$, and $U$ is a direct sum of blocks 
of the form $W(2l_i)$. 
We now give an algorithm to compute the parameters $\b_i$.
To simplify notation, write $V_{\b,a}(2k+1) = V_\b(2k+1) \oplus V_1(2k+1)^{[a-1]}$. 

\begin{enumerate}
\item[(1)] 
Compute $F = V(1-g)^{k_1}$ and the action of $g$ on $F_0 = F/{\rm Rad}(F)$. By Lemma \ref{olem},
\begin{equation}\label{fou}
g^{F_0} = V_{\b_1,a_1}(1) \oplus  \bigoplus_{i=2}^r V_{\b_i,a_i}(2k_i'+1) \oplus U',
\end{equation}
where $k_i' = k_i-k_1$. 

\item[(2)] 
Compute $F_1 = F_0(1-g)$ and the action of $g$ on $F_2 = F_1^\perp / {\rm Rad}(F_1^\perp)$. By Lemma \ref{olem}, 
\[
g^{F_2} = V_{\b_1,a_1}(1).
\]
The determinant of the Gram matrix of this $a_1$-dimensional space is congruent to $\b_1$ modulo $(\F_q)^2$. So this determinant determines $\b_1$.

\item[(3)] 
Replace $V$ by $F_0$, and $g$ by $g^{F_0}$ as in (\ref{fou}). Compute $F_3=V(1-g)$  and the action of $g$ on $F_3/{\rm Rad}(F_3)$, which is 
\[
g^{F_3/{\rm Rad}(F_3)} = \bigoplus_{i=2}^r V_{\b_i,a_i}(2k_i''+1) \oplus U'',
\]
where $k_i'' = k_i-1$. Now repeat Steps (1) and (2) to compute $\b_2$.
\end{enumerate}
Iterating, we compute all of the $\b_i$. 

\subsubsection{Conjugacy in $\SO(V)$ and $\O(V)$}
As discussed in Section \ref{sogood}, a conjugacy class 
$u^G$ in $G = \Or(V)$ may split into two classes in $G_0=\SO(V)$, 
with representatives $u$ and $u^t$ for a reflection $t$. 
For unipotent $g \in u^G$, 
we determine its $G_0$-class as follows: 
using the method of Section \ref{oconjel}, compute $y \in G$ 
such that $g^y = u$. If $y \in G_0$, then $g$ is $G_0$-conjugate to $u$; 
if $y \in G_0t$, then $g$ is $G_0$-conjugate to $u^t$. 

Similarly, an $\SO(V)$-class $u^{G_0}$ may split into two classes 
in $G_1=\O(V)$, with representatives $u$ and $u^s$ where 
$s \in \SO(V)\setminus \O(V)$. 
For unipotent $g \in u^{G_0}$, 
we determine its $G_1$-class as follows: 
compute $y \in G_0$ such that $g^y = u$.
Now determine which of $G_1$ or $G_1s$ contains $y$.

\subsection{Unitary groups} \label{slconprob}
Let $G = \SU_n(q)$. Let $u \in G$ be the matrix 
representative for $\bigoplus_{i=1}^s V(m_i)^{[r_i]}$ defined in Section \ref{sugood}. 
As described there, representatives of the $G$-conjugacy classes 
of elements with this Jordan form are $u^{d(\a)}$, 
where $d(\a) \in \GU_n(q)$ has determinant $\a$, and $\a$ ranges over 
representatives of $Z_{q+1}/Z_{q+1}^t$ where 
$t = {\gcd}(m_1,\ldots,m_s)$. 

Our algorithm to construct the class representative in $G$ of 
a unipotent element $g$ is the following.
Let $g$ have Jordan form $\bigoplus_{i=1}^s J_{m_i}^{[r_i]}$, and let 
$u$ be the matrix representative for $\bigoplus_{i=1}^s V(m_i)^{[r_i]}$. 
Using the work of Section \ref{uconjel}, find $y \in \GU_n(q)$ 
such that $g^y = u$, and let $\det(g) = \b \in Z_{q+1}$. 
Now $g$ is $G$-conjugate to $u^{d(\a)}$ where $\a$ is congruent 
to $\b^{-1}$ modulo $Z_{q+1}^t$.

\section{Constructing a conjugating element}\label{goodconjel}
In this section we complete our work for classical groups $G$ in 
good characteristic by solving the following problem: 
given unipotent $g \in G$ 
that is conjugate to a class representative $u$,
find $y \in G$ such that $g^y = u$.

\subsection{Symplectic groups}  \label{spconjel}
Let $G = \Sp(V) \cong \Sp_{2n}(q)$ with $q$ odd. Recall from Section \ref{spgood} that the unipotent class representatives in $G$ are the elements
\begin{equation}\label{repsp}
u=\bigoplus_{i=1}^r (V_{\b_i}(2k_i) \oplus V_1(2k_i)^{[a_i-1]}) \oplus \bigoplus_{i=1}^s W(2l_i+1)^{[b_i]},
\end{equation}
where $V = \bigoplus_{i=1}^r V_{2k_i}^{[a_i]} \oplus \bigoplus_{i=1}^s V_{4l_i+2}^{[b_i]}$, and each $\b_i \in \{1,\a\}$ (and the $k_i$ are distinct, as are the $l_i$). 
Given $g \in G$ that is conjugate to $u$, we aim to compute $y \in G$ such that $g^y=u$. We can compute $y$ ``block-by-block", so the main task is to solve the problem when $u$ is a single block $V_{\b}(2k) $ or $W(2l+1)$. 
We handle these cases separately.

\subsubsection{Case 1: $u = V_{\b}(2k)$}
Let $u = V_{\b}(2k) \in G = \Sp(V) \cong \Sp_{2k}(q)$. Recall from Section \ref{spgood} that $u = (1-e_\b)(1+e_\b)^{-1}$, where $e_\b$ is a nilpotent element of $sp(V)$; moreover, if $v = v_{-(2k-1)}$, the first vector in the basis defining $e_\b$, then $V$ has a basis
\[
v,\,ve_\b,\,ve_\b^2,\ldots,ve_\b^{k-1},\b^{-1}ve_\b^k,\ldots,\b^{-1}ve_\b^{2k-1},
\]
with respect to which the symplectic form defining $G$ has matrix $\e B$, where $\e = (-1)^k$ and
\[
B = \begin{pmatrix} &&&&1 \\&&&-1& \\ &&1&& \\ &\iddots &&& \\ -1&&&& \end{pmatrix}.
\]
In terms of the basis, 
\[
(ve_\b^i,\,ve_\b^j) = \left \{ 
\begin{array}{l} 
\e \b (-1)^i, \hbox{ if }i+j=2k-1 \\ 
0, \hbox{ otherwise.}
\end{array}
\right.
\]

Now suppose $g \in G$ is conjugate to $u$. 
We aim to compute $y \in G$ such that $g^y = u$. 
Let 
\[
f = (1-g)(1+g)^{-1},
\]
a nilpotent element of $sp(V)$. We seek $w\in V$ with the property that 
\begin{equation}\label{wbas}
w,\,wf,\ldots,wf^{k-1},\,\b^{-1}wf^k,\ldots, \b^{-1}wf^{2k-1}
\end{equation}
is a basis of $V$ satisfying
\begin{equation}\label{wprod}
(wf^i,\,wf^j) = \left \{ 
\begin{array}{l} 
\e \b (-1)^i, \hbox{ if }i+j=2k-1 \\ 
0, \hbox{ otherwise.}
\end{array}
\right.
\end{equation}
Then the map sending $wf^i \mapsto ve_\b^i$ for all $i$ will lie in $G$ and conjugate $g$ to $u$.

Here is our algorithm to find such a vector $w$. First, we find $z \in V$ such that $(z,\,zf^{2k-1}) = \e\b$; this can be done by random selection.
Now we aim to solve the following equation for $w \in V$ and $a_i\in \F_q$:
\begin{equation}\label{eqz}
z = w\sum_{i=0}^{2k-1} a_if^i
\end{equation}
such that $w$ satisfies (\ref{wbas}) and (\ref{wprod}).

Note that the fact that $(x,\,yf)+(xf,\,y) = 0$ for all $x,y \in V$ implies that for all $i,j$,
\begin{equation}
\begin{array}{lll}
(zf^i,zf^j) & = & (-1)^i(z,zf^{i+j}), \text{ and} \\
(z,zf^{2i}) & = & 0.
\end{array}
\end{equation}
So, to ensure property (\ref{wprod}) for $w$, we let $\a_j = (z,zf^{2j-1})$ for $1\le j\le k$, and solve for $a_i$ the equations
\begin{equation}\label{aeqs}
\begin{array}{lll}
\a_j & = & \left(\sum_r a_rwf^r,\,\sum_s a_swf^{s+2j-1}\right ) \\
      & = & \sum_{r+s+2j=2k} a_ra_s\e \b (-1)^r
\end{array}
\end{equation}
for $j=1,\ldots,k$.
These equations can be solved easily if we take $a_0=1$ and $a_i=0$ for all odd $i$. So
\begin{equation}\label{zw}
z = w(1+a_2f^2+a_4f^4+\cdots +a_{2k-2}f^{2k-2}),
\end{equation}
and the equations (\ref{aeqs}) simplify to 
\[
\a_{k-i} = \e \b \sum_{r=0}^i a_{2i-2r}a_{2r}\;\;\;(i=0,\ldots,k-1).
\]
For $i=0$, the equation is $\a_k = \e\b a_0^2 = \e\b$, which is true by choice of $z$. For $i=1$ the equation is $\a_{k-1} = 2\e\b a_2$, which we solve for $a_2$. Similarly we solve the $i=2$ equation for $a_4$, and continuing, 
we determine $a_{2i}$ for all $i$. Now let 
\[
m = 1+a_2f^2+\cdots + a_{2k-2}f^{2k-2}.
\]
Then $m$ is invertible, and the vector $w=zm^{-1}$ has the required properties (\ref{wbas}) and (\ref{wprod}). This complete the algorithm. 

\vspace*{4mm}
\no {\it Summary of algorithm.} 
Suppose $g \in G = \Sp_{2k}(q)$ is conjugate to $u = V_\b(2k)$. 

\begin{enumerate}
\item[(1)]
Let $f = (1-g)(1+g)^{-1}$. Find $z \in V$ such that $(z,\,zf^{2k-1}) = \e\b$.

\item[(2)]
Compute the values $\a_j = (z,zf^{2j-1})$ for $1\le j\le k$. Solve for $a_i$ the equations (\ref{aeqs}), assuming that $a_0=1$ and $a_i=0$ for all odd $i$.

\item[(3)]
Let $m = 1+a_2f^2+\cdots + a_{2k-2}f^{2k-2}$ and $w = zm^{-1}$. 
Define the linear map
$y:wf^i \mapsto ve_\b^i$ ($0\le i\le 2k-1$). 
Then $y \in G$ and $g^y = u$. 
\end{enumerate}

\subsubsection{Case 2:  $u = W(2l+1)$}

Let $u = W(2l+1) \in G = \Sp(V) \cong \Sp_{4l+2}(q)$. 
Recall from Section \ref{spgood} that $u = (1-e)(1+e)^{-1}$, 
where $e$ is a nilpotent element of $sp(V)$. If 
$w =w_{-2l}$ and $x = x_{-2l}$ are the first two vectors in 
the basis defining $e$, then $V$ has a basis
\[
w,\,we,\ldots,we^{2l},x,\,xe,\ldots, xe^{2l},
\]
where
\[
(we^i,\,xe^j) = \left \{ 
\begin{array}{l} 
\e (-1)^i, \hbox{ if }i+j=2l \\ 
0, \hbox{ otherwise}
\end{array}
\right.
\]
and $(we^i,\,we^j) =(xe^i,\,xe^j) =0$ for all $i,j$.  Also $\e = (-1)^l$.

Now suppose $g \in G$ is conjugate to $u$. We aim to compute $y \in G$ such that $g^y = u$, and adopt a similar approach to the previous case. Let $f = (1-g)(1+g)^{-1}$, a nilpotent element of $sp(V)$. We seek $w',x'\in V$ with the property that 
\begin{equation}\label{wxbas}
w',\,w'f,\ldots,w'f^{2l},\,x',\,x'f,\ldots, x'f^{2l}
\end{equation}
is a basis of $V$ satisfying
\begin{equation}\label{wxprod}
(w'f^i,\,x'f^j) = \left \{ 
\begin{array}{l} 
\e (-1)^i, \hbox{ if }i+j=2l \\ 
0, \hbox{ otherwise.}
\end{array}
\right.
\end{equation}
and $(w'f^i,\,w'f^j) =(x'f^i,\,x'f^j) =0$ for all $i,j$. 
Then the map sending $w'f^i \mapsto we^i$, $x'f^i \mapsto xe^i$ for all $i$ will lie in $G$ and conjugate $g$ to $u$.

Here is our algorithm to find such vectors $w',x'$. 
First, we find $z,t \in V$ such that $\la Vf,\,z,\,t\ra = V$. 
Now we aim to solve the following equations for 
$w',x' \in V$ and $a_i,b_i,c_i,d_i \in \F_q$:
\[
\begin{array}{l}
z = w'\sum_0^{2l} a_if^i+  x'\sum_0^{2l} b_if^i, \\
t = w'\sum_0^{2l} c_if^i+  x'\sum_0^{2l} d_if^i
\end{array}
\]
such that $w',x'$ satisfy (\ref{wxbas}) and (\ref{wxprod}). 
In practice, we can find a solution in which many of the variables are assumed to
be 0, namely
a solution of the form 
\begin{equation}\label{wxform}
\begin{array}{l}
z = w'\sum_{i=0}^{2l} a_if^i+  x'\sum_{j=1}^l b_{2j-1}f^{2j-1}, \\
t = w'\sum_{i=1}^l c_{2i-1}f^{2i-1}+  x'.
\end{array}
\end{equation}
To find this solution, we compute
\begin{equation}\label{abc}
\begin{array}{llll}
\a_i &=& (z,\,zf^{2i-1}) &  (1\le i\le l) \\
\b_i &=& (t,\,tf^{2i-1}) &  (1\le i\le l) \\
\g_j &=& (z,\,tf^j)      &  (0\le j\le 2l). 
\end{array}
\end{equation}
Then for $a_i,b_i,c_i$ as in (\ref{wxform}), the equations (\ref{wxprod}) are
\begin{equation}\label{simabc}
\begin{array}{llll}
\a_{l-i} & = & 2\e \sum_{j=0}^i a_{2j}b_{2i-2j+1} & (i=0,\ldots,l-1) \\
\b_{l-i} & = & -2\e c_{2i+1} & (i=0,\ldots,l-1) \\
\g_{2l-2i} & = & \e\left( a_{2i} + \sum_{j=1}^i b_{2j-1}c_{2i-2j+1} \right) & 
(i=0,\ldots,l) \\
\g_{2l-2i+1} & = & -\e a_{2i-1} & (i=1,\ldots,l). 
\end{array}
\end{equation}
These can be solved easily. First, we observe that 
\[
\begin{array}{lll}
a_0 & =  & \e\g_{2l}, \\
 a_{i} & = & -\e \g_{2l-i} \hbox{ for }i \hbox{ odd},\\
c_{2i-1}  & =  & -\frac{1}{2}\e \b_{l-i+1} \hbox{ for } i=1,\ldots,l.
\end{array}
\]
Now we can solve successively for $b_1,a_2,b_3,a_4,\ldots, b_{2l-1},a_{2l}$. Note that $(z,tf^{2l})\ne 0$ by choice of $z,t$, so $\g_{2l}\ne 0$ and hence $a_0\ne 0$. 

Given this solution for $a_i,b_i,c_i$, let
\[
A = \sum_{i=0}^{2l} a_if^i,\;\;B = \sum_{i=1}^l b_{2i-1}f^{2i-1}, \;\;C = \sum_{i=1}^l c_{2i-1}f^{2i-1}.
\]
Then (\ref{wxform}) is 
\[
\begin{array}{l}
z = w'A+x'B, \\
t = w'C+x'
\end{array}
\]
Solving these for $w',x'$ (noting that $a_0\ne 0$, so $A-CB$ is invertible), we obtain:
\[
w' = (z-tB)(A-CB)^{-1},\;\;x' = t-w'C.
\]
These vectors $w'$, $x'$ satisfy (\ref{wxbas}) and (\ref{wxprod}), as required.

\vspace*{4mm}
\no {\it Summary of algorithm.} 
Suppose $g \in G = \Sp_{4l+2}(q)$ is conjugate to $u = W(2l+1)$. 

\begin{enumerate}
\item[(1)] Let $f = (1-g)(1+g)^{-1}$. Find $z,t \in V$ such that $\la Vf,\,z,\,t\ra = V$.

\item[(2)] Compute the values $\a_i,\b_i,\g_i$ as in (\ref{abc}). 
Solve for $a_i,b_i,c_i$ the equations (\ref{simabc}).

\item[(3)] Let $A = \sum_{i=0}^{2l} a_if^i$, $B = \sum_{i=1}^l b_{2i-1}f^{2i-1}$, $C = \sum_{i=1}^l c_{2i-1}f^{2i-1}$, and let 
\[
w' = (z-tB)(A-CB)^{-1},\;\;x' = t-w'C.
\]
 Define the linear map 
\[
y: \;w'f^i \mapsto we^i,\;\; x'f^i \mapsto e^i \;\;\hbox{ for all }i.
\]
Then $y \in G$ and $g^y = u$.
\end{enumerate}

\subsubsection{General case}
Suppose $g \in \Sp(V)$ is conjugate to 
a unipotent class representative
\begin{equation}\label{ubks}
u=\bigoplus_{i=1}^r (V_{\b_i}(2k_i) \oplus V_1(2k_i)^{[a_i-1]}) \oplus \bigoplus_{i=1}^s W(2l_i+1)^{[b_i]}.
\end{equation}
We find $y \in G$ such that $g^y=u$ as follows:
\begin{enumerate}
\item[(a)] Find a non-degenerate $g$-invariant subspace $X$ on which the action of $g$ is conjugate to $u_1={V_\b(m)}$ or $W(m)$, where $m$ is the largest Jordan block size in (\ref{ubks}).
\item[(b)] By Case 1 or Case 2 above, compute $y_1 \in \Sp(X)$ that conjugates $g^X$ to $u_1$.
\item[(c)] Now work in $X^\perp$ and repeat Steps (a) and (b). 
\item[(d)] Continuing like this, obtain $y = \bigoplus y_i \in \Sp(V)$ conjugating $g$ to $u$.
\end{enumerate}
Note that in this process, when $2k_i$ becomes the largest block size, 
we must perform Step (a) $a_i$ times, and we must
ensure that the corresponding sequence of blocks $V_\b(2k_i)$ obtained 
is $V_1(2k_i),\ldots,V_1(2k_i),V_{\b_i}(2k_i)$.

It remains to give an algorithm to perform Step (a). 
As usual, write $f=(1-g)(1+g)^{-1}$, a nilpotent element of $sp(V)$.
Suppose first that the largest block size in (\ref{ubks}) is $2l+1$. By random choice, we can with high probability find $w,x \in V$ such that $(w,xf^{2l}) \ne 0$ and $(x,wf^{2l}) \ne 0$. Now set 
\[
X = \la wf^i,\,xf^i : 0\le i\le 2l \ra.
\]
Then $X$ is $g$-invariant and non-degenerate, and $g^X$
is conjugate to $W(2l+1)$, as required.
Now suppose the largest block size is $2k$. We find by random choice $v \in V$ such that $(v,vf^{2k-1}) \ne 0$. Now set
\[
X = \la vf^i : 0\le i\le 2k-1 \ra.
\]
Then $X$ is non-degenerate and $g^X$ is conjugate to $V_\b(2k)$ for 
some $\b \in \{1,\a\}$; we can compute $\b$ as in Section \ref{spconprob}. 
As noted above, we require that the sequence of values of $\b$ for 
the blocks of size $2k_i + 1$ should be $1,\ldots,1,\b_i$; so, if necessary, 
then we re-choose the initial vector $v$ in the above process until the 
correct value of $\b$ in this sequence is returned.

This completes our analysis of constructing conjugating elements in $\Sp_{2n}(q)$, $q$ odd. 

\subsection{Orthogonal groups}  \label{oconjel}

Let $G = \Or(V) \cong \Or_n^\e(q)$ with $q$ odd. Recall from Section \ref{sogood} that the unipotent class representatives in $G$ are the elements
\begin{equation}\label{repso}
u=\bigoplus_{i=1}^r (V_{\b_i}(2k_i+1) \oplus V_1(2k_i+1)^{[a_i-1]}) \oplus \bigoplus_{i=1}^s W(2l_i)^{[b_i]},
\end{equation}
where $V = \bigoplus_{i=1}^r V_{2k_i+1}^{[a_i]} \oplus \bigoplus_{i=1}^s V_{4l_i}^{[b_i]}$, and each $\b_i \in \{1,\a\}$ (and the $k_i$ are distinct, as are the $l_i$). 
Given $g \in G$ that is conjugate to $u$, we aim to compute $y \in G$ such that $g^y=u$. We then refine the method to handle conjugacy in $\SO(V)$ and $\O(V)$. 

The approach used is similar to that of Section \ref{spconjel}.
As before, we can compute $y$ ``block-by-block", so 
the main task is to solve the problem when $u$ is a 
single block $V_{\b}(2k+1) $ or $W(2l)$. 
We handle these cases separately.

\subsubsection{Case 1:  $u = V_{\b}(2k+1)$}

Let $u = V_{\b}(2k+1) \in G = \Or(V) \cong \Or_{2k+1}(q)$, and 
let $e_\b = (1-u)(1+u)^{-1}$, a nilpotent element of $so(V)$. 
If $v = v_{-2k}$, the first vector in the basis 
defining $e_\b$, then $V$ has a basis
\[
v,\,ve_\b,\,ve_\b^2,\ldots,ve_\b^{k},\b^{-1}ve_\b^{k+1},\ldots,\b^{-1}ve_\b^{2k},
\]
with respect to which the orthogonal form defining $G$ has values
\[
(ve_\b^i,\,ve_\b^j) = \left \{ 
\begin{array}{l} 
\e \b (-1)^i, \hbox{ if }i+j=2k \\ 
0, \hbox{ otherwise}
\end{array}
\right.
\]
where $\e=(-1)^k$.

Now suppose $g \in G$ is conjugate to $u$. We aim to compute $y \in G$ such 
that $g^y = u$. 
Let $f = (1-g)(1+g)^{-1}$, a nilpotent element of $so(V)$. 
We seek $w\in V$ with the property that 
$w,\,wf,\ldots,wf^{k},\,\b^{-1}wf^{k+1},\ldots, \b^{-1}wf^{2k}$ is 
a basis of $V$ satisfying
\begin{equation}\label{woprod}
(wf^i,\,wf^j) = \left \{ 
\begin{array}{l} 
\e \b (-1)^i, \hbox{ if }i+j=2k \\ 
0, \hbox{ otherwise.}
\end{array}
\right.
\end{equation}
Then the map sending $wf^i \mapsto ve_\b^i$ for all $i$ will lie in 
$G$ and conjugate $g$ to $u$.

Here is our algorithm to find such a vector $w$. First, we find 
$z \in V$ such that $(z,\,zf^{2k}) = \e\b$. Now we aim to solve 
the following equation for $w \in V$ and $a_i\in \F_q$:
\begin{equation}\label{eqoz}
z = w\,(1+a_2f^2+a_4f^4+\cdots +a_{2k}f^{2k})
\end{equation}
such that $w$ satisfies (\ref{woprod}).

Note that in the orthogonal case $(zf^i,zf^j) = (-1)^i(z,zf^{i+j})$, 
and $(z,zf^{2i+1}) = 0$ for all $i,j$.
So to ensure property (\ref{woprod}) for $w$, 
we let $\a_j = (z,zf^{2j})$ for $0\le j\le k$, 
and solve for $a_i$ the equations
\begin{equation}\label{aoeqs}
\begin{array}{ll}
\a_j & = \left(\sum_r a_{2r}wf^{2r},\,\sum_s a_{2s}wf^{2s+2j}\right ) \\
      & = \sum_{r+s=k-j} a_{2r}a_{2s}\e \b
\end{array}
\end{equation}
for $j=0,\ldots,k$.
These equations can easily be solved uniquely for $a_2,a_4,\ldots ,a_{2k}$. For this solution, let 
\[
m = 1+a_2f^2+\cdots + a_{2k}f^{2k}.
\]
Then $m$ is invertible, and the vector $w=zm^{-1}$ has the 
required properties (\ref{woprod}). 

\subsubsection{Case 2:  $u = W(2l)$}

Let $u = W(2l) \in G = \Or(V) \cong \Or^+_{4l}(q)$, and 
let $e = (1-u)(1+u)^{-1}$. If 
$w =w_{-(2l-1)}$ and $x = x_{-(2l-1)}$ are the first two vectors in 
the basis defining $e$, then $V$ has a basis 
$$w,\,we,\ldots,we^{2l-1},x,\,xe,\ldots, xe^{2l-1},$$ where
\[
(we^i,\,xe^j) = \left \{ 
\begin{array}{l} 
\e (-1)^i, \hbox{ if }i+j=2l-1 \\ 
0, \hbox{ otherwise}
\end{array}
\right.
\]
and $(we^i,\,we^j) =(xe^i,\,xe^j) =0$ for all $i,j$. 

Now suppose $g \in G$ is conjugate to $u$. 
We aim to compute $y \in G$ such that $g^y = u$.
We let $f = (1-g)(1+g)^{-1}$, and seek 
$w',x'\in V$ with the property that 
$w',\,w'f,\ldots,w'f^{2l-1},$
$x',\,x'f,\ldots, x'f^{2l-1}$ is a 
basis of $V$ satisfying
\begin{equation}\label{wxoprod}
(w'f^i,\,x'f^j) = \left \{ 
\begin{array}{l} 
\e (-1)^i, \hbox{ if }i+j=2l-1 \\ 
0, \hbox{ otherwise}
\end{array}
\right.
\end{equation}
and $(w'f^i,\,w'f^j) =(x'f^i,\,x'f^j) =0$ for all $i,j$. 
Then the map sending $w'f^i \mapsto we^i$, $x'f^i \mapsto xe^i$ 
for all $i$ will lie in $G$ and conjugate $g$ to $u$.

Here is our algorithm to find such vectors $w',x'$. First, 
we find $z,t \in V$ such that $\la Vf,\,z,\,t\ra = V$. 
In similar fashion to (\ref{wxform}), we aim to solve the following equations for 
$w',x' \in V$ and $a_i,b_i,c_i \in \F_q$:
\begin{equation}\label{wxoform}
\begin{array}{l}
z = w'\sum_{i=0}^{2l-1} a_if^i+  x'\sum_{j=1}^l b_{2j-1}f^{2j-1}, \\
t = w'\sum_{i=1}^l c_{2i-1}f^{2i-1}+  x'
\end{array}
\end{equation}
such that $w',x'$ satisfy (\ref{wxoprod}). 
To find this solution, we compute
\begin{equation}\label{abc2}
\begin{array}{llll}
\a_i & = & (z,\,zf^{2i}) & (0\le i\le l-1) \\
\b_i & = & (t,\,tf^{2i}) & (0\le i\le l-1) \\
\g_j & =&  (z,\,tf^j) & (0\le j\le 2l-1). 
\end{array}
\end{equation}
For $a_i,b_i,c_i$ as in (\ref{wxoform}), the equations (\ref{wxoprod}) are
\begin{equation}\label{simoabc}
\begin{array}{llll}
\a_{l-i} &=& 2\e \sum_{j=0}^{i-1} a_{2j}b_{2i-2j+1} & (i=1,\ldots,l) \\
\b_{l-i} &=& -2\e c_{2i-1} & (i=1,\ldots,l) \\
\g_{2l-2i-1} &=& \e\left( a_{2i} + \sum_{j=1}^i b_{2j-1}c_{2i-2j+1} \right)
& (i=0,\ldots,l-1) \\
\g_{2l-2i} &=& -\e a_{2i-1} & (i=1,\ldots,l). 
\end{array}
\end{equation}
These equations can be solved easily for $a_i,b_i,c_i$. 
For this solution, let
\[
A = \sum_{i=0}^{2l-1} a_if^i,\;\;B = \sum_{i=1}^l b_{2i-1}f^{2i-1}, \;\;C = \sum_{i=1}^l c_{2i-1}f^{2i-1}.
\]
Then the vectors
\[
w' = (z-tB)(A-CB)^{-1},\;\;x' = t-w'C.
\]
satisfy (\ref{wxoprod}), as required.

\subsubsection{General case}
Suppose $g \in \Or(V)$ is conjugate to a unipotent class representative
\begin{equation}\label{uobks}
u=\bigoplus_{i=1}^r (V_{\b_i}(2k_i+1) \oplus V_1(2k_i+1)^{[a_i-1]}) \oplus \sum_{i=1}^s W(2l_i)^{[b_i]}.
\end{equation}
We find $y \in G$ such that $g^y=u$ as follows:
\begin{itemize}
\item[(a)] Find a non-degenerate $g$-invariant subspace $X$ on which the action of $g$ is conjugate to $u_1=V_\b(m)$ or $W(m)$, where $m$ is the largest Jordan block size in (\ref{uobks}).
\item[(b)] By Case 1 or Case 2 above, compute $y_1 \in \Or(X)$ that conjugates $g^X$ to $u_1$.
\item[(c)] Now work in $X^\perp$ and repeat Steps (a) and (b). 
\item[(d)] Continuing like this, obtain $y = \bigoplus y_i \in \Or(V)$ conjugating $g$ to $u$.
\end{itemize}

To perform Step (a), write $f=(1-g)(1+g)^{-1}$, a nilpotent element 
of $so(V)$.
If the largest block size in (\ref{uobks}) is $2l$, then find 
$w,x \in V$ such that $(w,xf^{2l-1}) \ne 0$ and 
$(x,wf^{2l-1}) \ne 0$, and set $X = \la wf^i,\,xf^i : 0\le i\le 2l-1 \ra$; 
now $X$ is $g$-invariant and non-degenerate, and $g^X$ 
is conjugate to $W(2l)$. 
If the largest block size is $2k+1$, then  find 
$v \in V$ such that $(v,vf^{2k}) \ne 0$ and set
$X = \la vf^i : 0\le i\le 2k \ra$; now $X$ is $g$-invariant and 
non-degenerate, and 
$g^X$ is conjugate to $V_\b(2k+1)$ for some $\b \in \{1,\a\}$. 
We require that the sequence of values of $\b$ for the blocks of size 
$2k_i + 1$ should be $1,\ldots,1,\b_i$; so, if necessary, then we re-choose the 
initial vector $v$ in the above process until the correct 
value of $\b$ in this sequence is returned.

This completes our analysis of constructing conjugating elements in $\Or(V)$ 
(in odd characteristic). 

\subsubsection{Conjugation in $\SO(V)$ and $\O(V)$}
We now address finding conjugating elements in $\SO(V)$ and $\O(V)$. 
Let $u$ be as in (\ref{uobks}). 

First we consider $\SO(V)$. Assume that $u^{\Or(V)}$ splits into two $\SO(V)$-classes, with 
representatives $u,u^t$ as described in Section \ref{sogood}. 
Suppose $g$ is $\SO(V)$-conjugate to one of these. 
To construct the conjugating element of $\SO(V)$, first compute
$y \in \Or(V)$ such that $g^y=u$. Now
find $z \in \{1,t\}$ such that $yz \in \SO(V)$. 
Then $g^{yz} = u^z$, which is the required representative.

Now assume that $u^{\Or(V)}$ is a single $\SO(V)$-class, so that $r \ge 1$ in (\ref{uobks}).
We can find $v \in C_{\Or(V)}(u)$ of determinant $-1$: for example, $v = -I_{2k_i+1}$ acting on 
one of the $V(2k_i+1)$ summands in (\ref{uobks}).
Suppose $g$ is $\SO(V)$-conjugate to $u$, and compute $y \in \Or(V)$ such that $g^y=u$. Then $g^y=g^{yv}=u$, 
and either $y$ or $yv$ is in $\SO(V)$.

Now consider $\O(V)$. Let $u \in \O(V)$, and assume first that $u^{\SO(V)}$ splits into two $\O(V)$-classes, with representatives 
$u,u^s$ as described in Section \ref{sogood}. Let $g$ be $\SO(V)$-conjugate to $u$, and compute $y \in \SO(V)$
such that $g^y = u$. Now find $z \in \{1,s\}$ such that $yz \in \O(V)$. 
Then $g^{yz} = u^z$, which is the required representative.

Finally, suppose $u^{\SO(V)}$ is a single $\O(V)$-class. 
By \cite[Prop.\ 2.4]{GLOB}, $r \geq 1$ in (\ref{uobks}) 
and there exists $i$ such that one of the following holds:
\begin{itemize}
\item[(i)] $a_i\ge 2$;
\item[(ii)] $\b_i \not \equiv (-1)^{k_1+k_i}\b_1 \hbox{ mod }(\F_q)^2$.
\end{itemize}
In case (i), let $X$ be a sum of two blocks of dimension $2k_i+1$, and let 
$u_X$ be the restriction of $u$ to $X$. Now $C_{\SO(X)}(u_X)\,\O(X) = \SO(X)$. 
Compute representatives $\{1,s\}$ in $C_{\SO(X)}(u_X)$ for 
$\SO(X)/\O(X)$. Extend these to act trivially on $X^\perp$. 
Proceed as follows: suppose $g$ is $\O(V)$-conjugate to $u$; compute $y \in \SO(V)$ such that $g^y=u$, and then find 
$z \in \{1,s\}$ such that $yz \in \O(V)$. Now $g^{yz} = u$, as required. 

Case (ii) is similar. Let $X$ be the subspace corresponding to the 
direct sum of the blocks $V_{\b_1}(2k_1+1)$ and $V_{\b_i}(2k_i+1)$. 
The condition on $\b_i$ implies that 
$-1_X \not \in \O(X)$. So if we let $s=-1_X$, then $\{1,s\}$ are representatives in 
$C_{\SO(X)}(u_X)$ for $\SO(X)/\O(X)$. Now proceed as in case (i).

\subsection{Unitary groups}  \label{uconjel}
Let $G = \GU(V) \cong \GU_n(q)$, and let 
\[
u = \bigoplus_{i=1}^s V(m_i)^{[r_i]}
\]
be a unipotent class representative in $G$, as defined 
in Section \ref{sugood}. 
Given $g \in G$ that is conjugate to $u$, we aim to compute $y \in G$ such that $g^y=u$. 
Our approach is similar to that for the symplectic and orthogonal cases. 
Again we compute $y$ ``block-by-block". 
The analysis differs for blocks of even and odd size.

\subsubsection{Case 1:  $u = V(2k)$}
Let $u = V(2k) \in G = \GU(V) \cong \GU_{2k}(q)$, and let $v = v_{-(2k-1)}$, 
the first vector in the basis defining $u$ in Section \ref{sugood}. 
For the symplectic and orthogonal cases, 
we described the basis in terms of $v$ and a single nilpotent element of the corresponding Lie algebra. Here, we use two nilpotent elements, namely $e_1$ and $e_2$, where
\[
e_2 = u-1,\;\;e_1=1-u^{-1} = e_2(1+e_2)^{-1}.
\]
Recall that $\b \in \F_{q^2}^*$ satisfies $\b+\bar \b=0$. 
The basis defining $u$ as in Section \ref{sugood} is 
\[
v,\,ve_1,\ldots,ve_1^{k-1},
\b^{-1}ve_1^{k-1}e_2, -\b^{-1}ve_1^{k-1}e_2^2,
\ldots ,(-1)^{k-1}\b^{-1}ve_1^{k-1}e_2^k,
\]
with respect to which the unitary form defining $G$ has matrix
\begin{equation}\label{mat}
\begin{pmatrix} &&&1 \\&&1& \\ &\iddots&& \\ 1&&& \end{pmatrix}.
\end{equation}
Note that $e_1 = e_2(1+e_2)^{-1} = e_2(1-e_2+e_2^2-\cdots)$, so $e_1^k = e_1^{k-1}e_2(1-e_2+e_2^2-\cdots)$.

Now suppose $g \in G$ is conjugate to $u$. We aim to compute $y \in G$ 
such that $g^y = u$. 
Let $f_2 = g-1$ and $f_1 = f_2(1+f_2)^{-1}$. We seek a vector $w$ such that 
\begin{equation}\label{ubase}
w,\,wf_1,\ldots,wf_1^{k-1},\b^{-1}wf_1^{k-1}f_2,-\b^{-1}wf_1^{k-1}f_2^2,\ldots ,(-1)^{k-1}\b^{-1}wf_1^{k-1}f_2^k
\end{equation}
is a basis of $V$ with respect to which the unitary form 
has matrix (\ref{mat}). Then the map sending $wf_1^i\mapsto ve_1^i$,  $wf_1^{k-1}f_2^i \mapsto ve_1^{k-1}e_2^i$ for all $i$ will lie in $G$ and conjugate $g$ to $u$.

Our algorithm to find such a vector $w$ is the following.
First, we find $z \in V\setminus Vf_1$ such that 
$(z,\,\b^{-1}zf_1^{k-1}f_2^k) = (-1)^{k+1}$. Now we aim to solve 
the following equation for $w \in V$ and $a_i,b_i\in \F_{q^2}$:
\begin{equation}\label{equz-one}
z = w\left(1+\sum_{i=1}^{k-1} a_if_1^i + \sum_{i=1}^k b_i\b^{-1}(-1)^{i-1}f_1^{k-1}f_2^i\right),
\end{equation}
such that $w$ gives a basis  (\ref{ubase}) with respect to which the unitary form has matrix (\ref{mat}). 
Set  
\[
\begin{array}{lll}
\a_i = & (z,\,zf_1^i)&  (i=0,\ldots,k-1) \\
\b_i = & (z,\, \b^{-1}(-1)^{i-1}zf_1^{k-1}f_2^i) & (i=1,\ldots,k-1).
\end{array}
\]
To ensure that the matrix of the unitary form is (\ref{mat}), we need to solve $2k-1$ equations for $a_c,b_d$; in each equation, the left hand side is $\a_i$ or $\b_i$ and the right hand side is a function of $a_c,b_d$. The simplest equation is 
\[
\b_{k-1} = a_1-\bar a_1.
\]
We solve for $a_1$ using an additive Hilbert 90 algorithm. 
The next equation is 
\[
\b_{k-2} = a_2+\bar a_2 - \bar a_1 -a_1\bar a_1
\]
This is solved for $a_2$. Now we proceed to the $\b_{k-3}$ equation, solving for $a_3$, and so on. In this way we can find a solution for $a_c,b_d$  to all of the $2k-1$ equations. For this solution, define
\[
m = 1+\sum_{i=1}^{k-1} a_if_1^i + \sum_{i=1}^k b_i\b^{-1}(-1)^{i-1}f_1^{k-1}f_2^i.
\]
Then $w = zm^{-1}$ is the required vector.

\subsubsection{Case 2:  $u = V(2k+1)$}
Let $u = V(2k+1) \in G = \GU(V) \cong \GU_{2k+1}(q)$, and let $v = v_{-2k}$, 
the first vector in the basis defining $u$ in Section \ref{sugood}. 
We use three nilpotent elements to describe the basis in terms 
of $v$: recall that $\g \in \F_{q^2}$ satisfies $\g+\bar \g = -1$, and define
\[
e_2 = u-1,\;\;e_1=1-u^{-1} = e_2(1+e_2)^{-1},\;\;e_0 = e_2(1-\g e_2)^{-1}.
\]
The basis defining $u$ as in Section \ref{sugood} is 
\[
v,\,ve_1,\ldots,ve_1^{k-1},\,ve_1^{k-1}e_0,\,
-ve_1^{k-1}e_0e_2,
\, ve_1^{k-1}e_0e_2^2,
\ldots, (-1)^{k}ve_1^{k-1}e_0e_2^k,
\]
with respect to which the unitary form defining $G$ has matrix (\ref{mat}).

Now suppose $g \in G$ is conjugate to $u$. We aim to compute $y \in G$ 
such that $g^y = u$. 
Let $f_2 = g-1$, $f_1 = f_2(1+f_2)^{-1}$ and 
$f_0 = f_2(1-\g f_2)^{-1}$. We seek a vector $w$ such that 
\begin{equation}\label{ubase2}
w,\,wf_1,\ldots,wf_1^{k-1},\,wf_1^{k-1}f_0,\,-wf_1^{k-1}f_0f_2, \, wf_1^{k-1}f_0f_2^2,
\ldots, (-1)^{k}wf_1^{k-1}f_0f_2^k,
\end{equation}
is a basis of $V$ with respect to which the unitary form 
has matrix (\ref{mat}). Then the map 
sending $wf_1^i\mapsto ve_1^i$,  $wf_1^{k-1}f_0f_2^i \mapsto ve_1^{k-1}e_0e_2^i$ for all $i$ will lie in $G$ and conjugate $g$ to $u$.

Our algorithm to find such a vector $w$ is similar to that of 
Case 1. First, we find $z \in V$ such that 
$(z,\,zf_1^{k-1}f_0f_2^k) = (-1)^{k}$. Now we aim to solve the following equation for $w \in V$ and  $a_i,b_i\in \F_{q^2}$:
\begin{equation}\label{equz-two}
z = w\left(1+\sum_{i=1}^{k-1} a_if_1^i + \sum_{i=0}^k b_i(-1)^{i}f_1^{k-1}f_0f_2^i\right),
\end{equation}
such that $w$ gives a basis  (\ref{ubase2}) with respect to which the unitary form has matrix (\ref{mat}). 
To do this, set 
\[
\begin{array}{lll}
\a_i = & (z,\,zf_1^i) &  (i=0,\ldots,k-1) \\
\b_i = & (z,\,zf_1^{k-1}f_0f_2^i) & (i=0,\ldots,k-1).
\end{array}
\]
This leads to $2k$ equations for $a_i,b_i$, the simplest of which is the $\b_{k-1}$-equation, which is $\b_{k-1} = 
(-1)^k(\bar a_1- a_1)$. We solve for $a_1$, 
then proceed to the $\b_{k-2}$-equation, and so on. 
Having solved for all $a_i,b_i$, we define
\[
m = 1+\sum_{i=1}^{k-1} a_if_1^i + \sum_{i=0}^k b_i(-1)^{i}f_1^{k-1}f_0f_2^i.
\]
Then $w = zm^{-1}$ is the required vector.

\subsubsection{General case}
Suppose $g \in \GU(V)$ is conjugate to a unipotent class representative
\begin{equation}\label{uubks}
u=\bigoplus_{i=1}^s  V(m_i)^{[r_i]}, 
\end{equation}
where $m_1>m_2>\cdots >m_s$.
As before, we find $y \in G$ such that $g^y=u$ proceeding 
``block-by-block". 
The main step is to compute a largest non-degenerate block for $g$, starting with $V(m_1)$. 
To do this, we find $v \in V$ such that $(v,v(g-1)^{m_1-1}) \ne 0$ and set
$$X = \la v(g-1)^i : 0\le i\le m_1-1 \ra.$$
Then $X$ is non-degenerate and 
$g^X$ is conjugate to $V(m_1)$. Now, as in Case 1 
or 2 above, we compute $y_1 \in \GU(X)$ conjugating $g^X$ to the chosen 
representative for $V(m_1)$, and repeat the process in $X^\perp$, ending up 
with $y = \bigoplus y_i\in \GU(V)$ that conjugates $g$ to $u$.

\subsubsection{Conjugation in $\SU(V)$}
We now address the issue of finding conjugating elements in $\SU(V)$. 
Let $G = \SU(V) \cong \SU_n(q)$.
Recall from Section \ref{sugood} that our representatives for the unipotent classes of $G$ are $u^{d(\a)}$, where $u$ is as in (\ref{uubks}), $d(\a)\in \GU(V)$ has determinant $\a$, and $\a$ ranges over representatives of $Z/Z^t$, where $Z = \{\mu \in \F_{q^2}: \mu\bar \mu=1\}$ and $t = {\gcd}(m_1,\ldots,m_s)$. 

Suppose $g \in G$ is $G$-conjugate to $u^{d(\a)}$. As above, 
we can find $y \in \GU(V)$ such that $g^y = u$. Let $\b= \det(y)$. 
If $\b^{-1}\a^{-1} \in Z^t$, then we can find 
$c \in C_{\GU(V)}(u)$ of determinant $\b^{-1}\a^{-1}$ (acting as a scalar 
on each block). Now $g^{ycd(\a)} = u^{d(\a)}$, and 
the conjugating element $ycd(\a)$ has determinant 1; hence it lies in $G$ as required. 

\section{Some examples}
We close this chapter by listing the unipotent 
class representatives and centralizer orders
for the 7-dimensional orthogonal groups and 
the 8-dimensional orthogonal and symplectic groups 
over fields of odd size $q$.
The structures of the centralizers are given by 
Theorems \ref{sympreps} and \ref{orthreps}. 

Notation is as in earlier sections. 
Recall that $\a \in \F_q$ is a fixed non-square, and both $\b$ 
and $\b_i$ are in $\{1, \a \}$. 
The discriminant of the orthogonal form fixed by $G = \SO_8^\e(q)$ is a 
square in $\F_q$ if and only if $\e = +$ (see \cite[Prop.\ 2.5.13]{KL}).
In Table \ref{so8} we use the notation $\D$ for the value of the discriminant 
modulo $(\F_q^*)^2$; column 3 lists the conditions which 
determine the sign $\e = \pm$ of the group $\SO_8^\e(q)$ containing $u$.
Tables \ref{so7odd} and \ref{so8} have a 
column indicating which unipotent classes split in the corresponding group 
$\O_7(q)$ or $\O_8^\e(q)$; this information is recorded in 
Section \ref{sogood}.

\begin{table}[htbp] 
\caption{Unipotent class representatives in $G = \SO_7(q)$, $q$ odd} 
\label{so7odd}
\[
\begin{array}{|l|l|l|}
\hline 
\hbox{Representative } u & |C_G(u)| & \hbox{class splits in }\O_7(q)? \\
\hline
V_1(1)^{[7]} & |G| &  \hbox{no}\\
W(2) \oplus V_1(1)^{[3]} & q^{7}|\Sp_2(q)|\,|\SO_3(q)| &  \hbox{no}\\
W(2) \oplus V_1(3) & q^{6}|\Sp_2(q)| &  \hbox{yes} \\

V_\b(3) \oplus V_\b(1) \oplus V_1(1)^{[3]} \,(\b \in \{1,\a\}) & q^{5}|\Or_4^{\pm}(q)| &  \hbox{no} \\
V_\b(3) \oplus V_1(3) \oplus V_\b(1) \,(\b \in \{1,\a\})&  2q^{6}(q\pm 1) &  \hbox{no} \\
                                                                        
V_\b(5) \oplus V_\b(1) \oplus V_1(1)\,(\b\in\{1,\a\})  & 2q^4(q\pm 1) &  \hbox{no} \\

V_1(7)  &  q^3  &  \hbox{yes}\\
\hline
\end{array}
\]
\end{table}

\begin{table}[htbp] 
\caption{Unipotent class representatives in $G = \Sp_8(q)$, $q$ odd} \label{sp8odd}
\[
\begin{array}{|l|l|}
\hline 
\hbox{Representative } u & |C_G(u)| \\
\hline
W(1)^{[4]} & |G| \\

W(3) \oplus W(1) & q^{8}|\Sp_2(q)|^2 \\
W(3) \oplus V_\b(2)\,(\b \in \{1,\a\}) & 2q^{9}|\Sp_2(q)| \\

W(1)^{[3]} \oplus V_\b(2) \,(\b \in \{1,\a\}) &  2q^7|\Sp_6(q)| \\

W(1)^{[2]} \oplus V_\b(2) \oplus V_1(2) \,(\b \in \{1,\a\}) &  2q^{11}|\Sp_4(q)|(q\pm 1) \\
W(1)^{[2]} \oplus V_\b(4)\,(\b\in\{1,\a\})  &  2q^6|\Sp_4(q)| \\

W(1) \oplus V_\b(2) \oplus V_1(2)^{[2]} \,(\b \in \{1,\a\}) &  q^{12}|\Sp_2(q)|\,|\Or_3(q)| \\
W(1) \oplus V_\b(4) \oplus V_\g(2)\,(\b,\g\in\{1,\a\})  &  4q^9|\Sp_2(q)| \\
W(1) \oplus V_\b(6) \,(\b \in \{1,\a\}) &  2q^5|\Sp_2(q)| \\

V_\b(2) \oplus V_1(2)^{[3]} \,(\b \in \{1,\a\}) & q^{10}|\Or_4^{\pm}(q)| \\
V_\b(4) \oplus V_\g(2) \oplus V_1(2) \,(\b,\g \in \{1,\a\})&  4q^{9}(q\pm 1) \\
                                                                         & (2 \hbox{ of each size}) \\
V_\b(4) \oplus V_1(4)\,(\b\in\{1,\a\})  & 2q^7(q\pm 1) \\
V_\b(6) \oplus V_\g(2)\,(\b,\g\in\{1,\a\})  &  4q^6 \\
V_\b(8)\,(\b\in\{1,\a\})  &  2q^4 \\
\hline
\end{array}
\]
\end{table}

\begin{table}[htbp] 
\begin{center}
\caption{Unipotent class representatives in $G = \SO_8^\e(q)$, $q$ odd}
\label{so8} 
\[
\begin{array}{|l|l|l|l|}
\hline 
\hbox{Representative } u & |C_G(u)| & \hbox{condition} &  \hbox{class splits} \\

 &&& \hbox{ in }\O_8^\e(q)? \\
\hline
V_\b(1) \oplus V_1(1)^{[7]} & |G| & \b \equiv \D \hbox{ mod }(\F_q^*)^2 & \hbox{no} \\

W(4),\,W(4)^t \;(t \hbox{ reflection}) & q^5|\Sp_2(q)| & \e=+ & \hbox{yes} \\

W(2)^{[2]},\,W(2)^t \oplus W(2) & q^6|\Sp_4(q)| & \e=+ & \hbox{yes} \\

W(2) \oplus V_\b(1) \oplus V_1(1)^{[3]} & q^{9}|\Sp_2(q)|\,|\SO_4^\e(q)| &  \b \equiv \D & \hbox{no} \\

W(2) \oplus V_{\b_1}(3) \oplus  V_{\b_2}(1)   & 2q^{9}|\Sp_2(q)| & -\b_1\b_2 \equiv \D & \hbox{yes, }\e=+ \\
&&&  \hbox{no, }\e=- \\

V_{\b_1}(3) \oplus V_{\b_2}(1)  \oplus  V_1(1)^{[4]}  & 2q^{6}|\SO_5(q)| & -\b_1\b_2 \equiv \D & \hbox{no} \\

V_{\b_1}(3) \oplus V_1(3) \oplus V_{\b_2}(1) \oplus V_1(1) &  2q^{8}(q\pm 1)^2,\,\e=+ & -\b_1\b_2 \equiv \D & \hbox{no} \\
                                                                      & 2q^{8}(q^2-1),\,\e=- && \\
  
V_{\b_1}(5) \oplus V_{\b_2}(1) \oplus V_1(1)^{[2]}  & 2q^5|\SO_3(q)| & \b_1\b_2 \equiv \D & \hbox{no} \\

V_{\b_1}(5) \oplus V_{\b_2}(3)   & 2q^6 & -\b_1\b_2 \equiv \D & \hbox{yes, }\e=+ \\
&&&  \hbox{no, }\e=- \\

V_{\b_1}(7) \oplus V_{\b_2}(1)    &  2q^4 &  -\b_1\b_2 \equiv \D & \hbox{yes, }\e=+ \\
&&&  \hbox{no, }\e=- \\
\hline
\end{array}
\]
\end{center}
\end{table}

\chapter{Unipotent classes in bad characteristic}\label{unibadchap}

In this chapter we solve the conjugacy problems (1)-(3) of Section 
\ref{mainprob} for unipotent elements 
of symplectic and orthogonal groups in characteristic 2. 
Many of the methods of 
Chapter \ref{unigoodchap} do not apply: there is no Cayley 
map between unipotent and nilpotent elements, nor is there 
a factorization of the centralizer of a unipotent element 
in a canonical parabolic subgroup. Hence, although we 
can exploit some features of the good characteristic case, 
our algorithms are more complicated.

\section{Unipotent class representatives}\label{badsec}

We begin by describing the conjugacy class representatives. 
Let $G = \Sp(V)$ or $\Or(V)$, where $V = V_{2n}(q)$, a vector space 
of dimension $2n$ over $\F_q$ and $q=2^a$. 
Denote by $(\,,\,)$ a symplectic form on $V$ preserved by $G$, 
and by $Q$ a quadratic form in the orthogonal case. 

The representatives are defined as orthogonal direct sums of 
various indecomposable blocks, as defined in \cite{GLOB} 
(with a correction for the block $W_\b(2l+1)$, see Section 
\ref{correctionforW}).
These  blocks are denoted $V_\b(2k)$, $W(k)$ and 
$W_\b(2l+1)$ for $\b\in \{0,\a\}$, where $\a$ is a fixed element 
of $\F_q$ such 
that the quadratic $x^2+x+\a$ is irreducible over $\F_q$. We define 
these blocks, using bases labelled by vectors 
$v_i^{(j)}$, $w_i^{(j)}$, $x_i^{(j)}$; these will be used to define a 
certain parabolic subgroup that contains the centralizer of the 
unipotent element, as in Section \ref{goodrep}. 

\subsection{Blocks $V_\b(2k)$} 
Let $k \geq 1$ and let $V_{2k}$ be a $2k$-dimensional vector space over $\F_q$, with basis 
\[
v_{-(2k-1)},v_{-(2k-3)},\ldots, v_{2k-1}.
\]
Define a symplectic form on $V_{2k}$ by setting $(v_i,v_{-i}) =1$ for all $i$, 
and all other values $(v_i,v_j) = 0$. For $k \ge 2$, define $V_\b(2k)$ to be the  linear map acting as follows:
\[
\begin{array}{llll}
v_i &\mapsto &v_i+v_{i+2}+\cdots +v_1+ \b v_3 & \hbox{ for } -(2k-1) \le i\le -3, \\
v_i &\mapsto &v_i+v_{i+2} & \hbox{ for }-1 \le i<2k-1, \\
v_{2k-1} &\mapsto &v_{2k-1}. &
\end{array}
\]
For $k=1$, define $V_\b(2)$ as follows:
$v_{-1} \mapsto v_{-1}+v_1$ and $v_1\mapsto v_1$. 

Observe that 
$V_\b(2k)$ lies in $\Sp(V_{2k})$ and acts as a single Jordan block. If we define a quadratic form $Q_\b$ on $V$ associated to $(\,,\,)$ by setting $Q_\b(v_{-1}) = \b$, $Q_\b(v_1) = 1$ and  $Q$-values 0 on all other basis vectors, then $V_\b(2k)$ is in the orthogonal group preserving $Q_\b$. Specifically, $V_0(2k) \in \Or^+_{2k}(q)$ and $V_\a(2k) \in \Or^-_{2k}(q)$. For $\b=0$ we 
write just $V(2k)$ instead of $V_0(2k)$. For all these assertions, see \cite[\S 3]{GLOB}.

\subsection{Blocks $W(k)$} \label{wkcase}

Let $k\ge 1$ and let $V_{2k}$ be a $2k$-dimensional vector space over $\F_q$, with basis 
\[
w_{-(k-1)},x_{-(k-1)},w_{-(k-3)},x_{-(k-3)},\ldots, w_{k-1},x_{k-1}.
\]
Define a symplectic form on $V_{2k}$ by setting $(w_i,x_{-i}) =1$ for all $i$, 
and all other values $(w_i,w_j)$, $(w_i,x_j)$, $(x_i,x_j)$ equal 0; 
and let $Q$ be the quadratic form associated to $(\,,\,)$ such that 
$Q(w_i)=Q(x_i)=0$ for all $i$. Now define $W(k)$ to be the linear map 
acting as follows:
\[
\begin{array}{llll}
 w_i & \mapsto & w_i+w_{i+2} & \hbox{ for } -(k-1) \le i<k-1, \\
w_{2k-1} & \mapsto & w_{2k-1}, & \\
x_i & \mapsto & x_i+x_{i+2}+\cdots + x_{k-1} & \hbox{ for } -(k-1) \le i\le k-1.
\end{array}
\]
Then $W(k)$ lies in $\Sp(V_{2k})$, and also in $\O^+(V_{2k})$, the orthogonal group of plus type preserving $Q$. It acts as a sum of two Jordan blocks of size $k$.

\subsection{Blocks $W_\b(2l+1)$ }  \label{correctionforW}
Let $l\ge 1$ and let $V_{4l+2}$ be a $(4l+2)$-dimensional vector space over $\F_q$, with basis 
\[
w_{-2l},x_{-2l},w_{-(2l-2)},x_{-(2l-2)},\ldots, w_{2l},x_{2l}.
\]
Define a symplectic form on $V_{4l+2}$ as in Section \ref{wkcase},
and let $Q_\b'$ be the quadratic form associated to $(\,,\,)$ such that $Q_\b'(w_{-2})=Q_\b'(w_0)=Q_\b'(x_0) =\b$, and $Q_\b'$ takes the value 0 on all other basis vectors. Now define $W_\b(2l+1)$ to be the linear map acting as follows:
\[
\begin{array}{llll}
w_i & \mapsto & w_i+w_{i+2} & \hbox{ for }i \ne -4,-2,2l, \\
w_{-4} & \mapsto & w_{-4}+w_{-2}+\b x_2, \\
w_{-2} & \mapsto &  w_{-2}+w_{0}+\b x_2, \\
w_{2l} & \mapsto & w_{2l}, \\
x_i & \mapsto & x_i+x_{i+2}+\cdots + x_{2l}+\b w_2 & \hbox{ for } -2l \le i\le -2, \\
x_i & \mapsto &  x_i+x_{i+2}+\cdots + x_{2l} & \hbox{ for } 0 \le i\le 2l-2, \\
x_{2l} & \mapsto &  x_{2l}.
\end{array}
\]
This is a correction to the block $W_\b(2l+1)$ as defined in \cite{GLOB}. 
Note that $W_0(2l+1)$ is just the element $W(2l+1)$ of $\O_{4l+2}^+(q)$ 
defined in the previous case, and $W_\a(2l+1) \in \O^-_{4l+2}(q)$.

\subsection{The general case}
By \cite[Thm.\ 3.1]{GLOB}, 
every unipotent element of $G = \Sp(V)$ or $\Or(V)$
is conjugate to exactly one element of the form
\begin{equation}\label{canon}
\bigoplus_i W(m_i)^{[c_i]}  \oplus  \bigoplus_j V(2k_j)^{[d_j]}  \oplus  \bigoplus_r W_\a(m_r') \oplus \bigoplus_s V_\a(2k_s')
\end{equation}
satisfying the following conditions:
\begin{itemize}
\item[(i)] the $m_r'$ are odd and distinct, and the $k_s'$ are distinct;
\item[(ii)] $d_j\le 2$, and $d_j\le 1$ if there exist $j,s$ such that $k_j=k_s'$;
\item[(iii)] there exist no $j,s$ such that $k_s'-k_j=1$ or $k_s'-k_j'=1$;
\item[(iv)] there exist no $j,r$ such that $m_r' = 2k_j\pm 1$ or $m_r'=2k_j'\pm 1$;
\item[(v)] for $G = \Sp(V)$, each $m_r'\ge 3$ and each $k_s' \ge 2$.
\end{itemize}
In the orthogonal case, the element is in
$\O(V)$ if and only if the total number of $V$- and $V_\a$-blocks is even; moreover, the only $G$-classes which split into two $\O(V)$-classes are those of the form $\bigoplus W(m_i)^{[a_i]}$ with all $m_i$ even, and for these a second class representative can be obtained by replacing one summand $W(m_i)$ by $W(m_i)^t$ with $t$ a reflection.

\subsection{Class representatives} 
We choose a matrix representative for the element (\ref{canon}) in 
a manner similar to Section \ref{goodrep}. List the $V$-blocks 
and $W$-blocks in (\ref{canon}) in decreasing order of Jordan block sizes.  
For the $i^{th}$ space, if it is $V_{\b}(2k)$, then label its basis as
\[
v_{-(2k-1)}^{(i)},v_{-(2k-3)}^{(i)},\ldots, v_{2k-1}^{(i)};
\]
if it is $W(k)$ or $W_\b(k)$, then label its basis as
\[
w_{-(k-1)}^{(i)},x_{-(k-1)}^{(i)}, \ldots, w_{k-1}^{(i)}, x_{k-1}^{(i)}.
\]
Take a basis of $V$ consisting of these vectors 
$v_j^{(i)}, w_j^{(i)}, x_j^{(i)}$ in increasing order lexicographically with 
respect to the pairs $(j,i)$. Choose the matrix representative corresponding 
to the linear map (\ref{canon}) to be its matrix with respect to 
this ordered basis.

\section{Centralizers of class representatives}\label{centbad}
We now discuss how to construct the centralizers of these 
unipotent class representatives. 
We begin by describing the structure of the centralizer of 
a representative (\ref{canon}); this is given by \cite[Thm.\ 7.3]{LS}. 

To make the description easier, we rewrite (\ref{canon}) slightly differently. 
For $m_i$ odd, let $a_i$ be the total number of $W$-blocks of size $m_i$ (including $W_\a(m_i)$ if it is present), and write
\[
W_{\b,a_i}(m_i) = \left\{ \begin{array}{l} W(m_i)^{[a_i]}, \hbox{ if }\b = 0 \\
W(m_i)^{[a_i-1]} \oplus W_\a(m_i), \hbox{ if } \b = \a. \end{array} \right.
\]
Similarly, for the $V$-blocks define 
\[
V_{\b,b_i}(2k_i) = \left\{ \begin{array}{l} V(2k_i)^{[b_i]}, \hbox{ if }\b = 0 \\
V(2k_i)^{[b_i-1]} \oplus V_\a(2k_i), \hbox{ if } \b = \a. \end{array} \right.
\]
Now rewrite (\ref{canon}) as
\begin{equation}\label{newcan}
\bigoplus_{m_i\;\text{even}} W(m_i)^{[a_i]}  \oplus  \bigoplus_{m_i\;\text{odd}}W_{\b_i,a_i}(m_i)  \oplus  \bigoplus_j V_{\b_j,b_j}(2k_j).
\end{equation}
Thus $a_i = c_i$ or $c_i+1$, where $c_i$ is as in (\ref{canon}); and $a_i=c_i+1$ if and only if there is a summand $W_\a(m_i)$. 
Similarly, $b_j = d_j$ or $d_j+1$.

Let $u$ be the unipotent element (\ref{newcan}), and let ${\mathcal J}_u$ be the 
set of Jordan block sizes in the Jordan form of $u$. 
Define a function $\c_u: {\mathcal J}_u \to \N$ as follows:
\[
\begin{array}{llll}
G = \Sp(V): & \c_u(n) & = & \left\{ \begin{array}{l} \frac{1}{2}n, \hbox{ if }\exists\, V(n) \hbox{ or }V_\a(n) \hbox{ in } (\ref{canon}) \\
\lfloor \frac{1}{2}(n-1)\rfloor, \hbox{otherwise.} \end{array} \right.   \\ 
& & \\
G = \Or(V):  &\c_u(n) & = & \left\{ \begin{array}{l} \frac{1}{2}(n+2), \hbox{ if }\exists\, V(n) \hbox{ or }V_\a(n) \hbox{ in } (\ref{canon}) \\
\lfloor \frac{1}{2}(n+1)\rfloor, \hbox{otherwise.} \end{array} \right.
\end{array}
\]
Let $f_1,\ldots,f_s$ be the list of Jordan block sizes in (\ref{newcan}), 
including multiplicities, where $f_1\ge f_2\ge \cdots \ge f_s$, and define
\begin{equation}\label{dimc}
N_u = \sum_{i=1}^s (if_i - \c_u(f_i)).
\end{equation}
Then $N_u$ is the dimension of the centralizer of $u$ in the algebraic group 
$\bar G = \Sp(\bar V)$ or $\SO(\bar V)$ 
over $\bar \F_q$ corresponding to $G$, where $\bar V = V \otimes \bar \F_q$ (see \cite[(4.1)]{LS}). 

The following is \cite[Thm.\ 7.3]{LS}.
\begin{thm}\label{uni} 
Let $G = \Sp(V)$ or $\Or(V)$, and let $u$ be as in $(\ref{newcan})$; in the orthogonal case, assume $u \in \O(V)$.
Let $C = C_G(u)$.
\begin{itemize}
\item[{\rm (i)}] $C$ has a normal subgroup $R$ of order $q^{M_u}$, 
where $M_u$ is given in part $(iii)$.
\item[{\rm (ii)}] 
\begin{equation}\label{cov}
C/R \cong \prod_{m_i\;\text{even}} \Sp_{2a_i}(q) \times \prod_{m_i\;\text{odd}}  I_{2a_i}(q) \times Z_2^{t+\d},
\end{equation}
where 
\[
I_{2a_i}(q) = \left\{ \begin{array}{ll}
               \Sp_{2a_i}(q), & \hbox{ if }m_i\pm 1 = 2k_j \hbox{ or }2k_j' \hbox{ for some }j, \\
               \Sp_{2a_i}(q), & \hbox{ if }m_i = 1 \hbox{ and }G = \Sp(V), \\
               \Or^{\e_i}_{2a_i}(q), & \hbox{ otherwise;}
\end{array}
\right.
\]
and  $t$ is the number of values of $j$ such that 
$k_j-k_{j+1}\ge 2$ where we list the $k_j$ in $(\ref{newcan})$ in decreasing order; and $\d \in \{0,1\}$ with $\d = 0$ if and only if either there are no $V$-blocks, or  $G=\Sp(V)$ and there is a block $V(2)$.
\item[{\rm (iii)}]  $|R|=q^{M_u}$, where 
\[
M_u = N_u - \sum_{m_i\;\text{even}} \dim \Sp_{2a_i} - \sum_{m_i\;\text{odd}}  \dim I_{2a_i}.
\]
\end{itemize}
\end{thm}

For the remainder of this section 
$G$ denotes either $\Sp(V)$ or $\O(V)$ (rather than $\Or(V)$). 
Let $F$ be a Frobenius endomorphism of $\bar G$ such that $G = \bar G^F$ (see \cite[\S 2.4]{LS}). Note that 
the normal subgroup $R$ in part (i) is the fixed point group under $F$ of the unipotent radical of $C_{\bar G}(u)$.

As in good characteristic, the key to 
computing $C_G(u)$ is to construct a canonical parabolic subgroup 
that contains this centralizer. We define this parabolic subgroup to be 
the stabilizer of a flag defined by the ``weights" of the basis vectors, 
as we now describe. Note that we do not define a 1-dimensional torus (as we did in the previous chapter), so these
weights are just a notational device, and are not the weights of a torus.

First, we tweak the bases for the blocks $W_\a(2l+1)$ for the orthogonal group. For such a block, with basis $w_{-2l},x_{-2l},\ldots, w_{2l},x_{2l}$, we change the basis vectors $w_{-2}$ and $x_0$ as follows:
\[
w_{-2} \mapsto w_{-2}+w_0,\;\;x_0 \mapsto x_0+x_2.
\]
With this change, the values of the quadratic form $Q_\b'$ are nonzero only for the basis vectors $w_0$ and $x_0$.

Now we define a canonical parabolic subgroup associated to $u$, 
the representative 
(\ref{newcan}) in $G = \Sp(V)$ or $\O(V)$. For a basis vector $v_j^{(i)}$, $w_j^{(i)}$ or  $x_j^{(i)}$, the subscript $j$ is its {\it weight}; the span of all basis vectors of weight $j$ is the {\it weight space} \index{weight space} for $j$. Define $P$ to be the parabolic subgroup stabilizing the flag of $V$ defined by sums of weight spaces for decreasing weights. Write $P=QL$, with unipotent radical $Q$ and Levi subgroup $L$. 

The next result follows from \cite[Thm.\ 5.3]{CP}. 

\begin{thm}\label{par2} The centralizer in $G$ of $u$ is a subgroup of $P$.
\end{thm}

Unlike the good characteristic case, it is not true in general 
that $C_P(u) = C_Q(u)C_L(u)$. 
From Theorem \ref{uni}(i) we know that $C_G(u)$ 
has a normal subgroup $R$ of order $q^{M_u}$. 
We conjecture the order of $RC_Q(u)/C_Q(u)$. 

\begin{conj}\label{rcon} Let $u$ be as in $(\ref{newcan})$, and 
let $J = \{j: b_j=2\}$. For $j \in J$, let $a_j$ be the 
multiplicity of $W(2k_j)$ in $(\ref{newcan})$. 
Then $|RC_Q(u)/C_Q(u)| = q^t$, where $t = \sum_{j\in J}(2a_j+1)$.
\end{conj}

In our algorithm given below to construct the centralizer, we assume the conjecture. In any given example, it can be verified 
by checking that the subgroup $R$ of $C_G(u)$ constructed in (e) 
has order $q^{M_u}$, as given in Theorem \ref{uni}(i).

We construct $C_G(u)$ in the following steps:
\begin{itemize}
\item[(a)] construct $C_Q(u)$;
\item[(b)] construct subgroups of $C_G(u)$ that cover the factors $\Sp_{2a_i}(q)$ in (\ref{cov});
\item[(c)] construct subgroups of $C_G(u)$ that cover the factors $I_{2a_i}(q)$ in (\ref{cov});
\item[(d)] construct the subgroup of $C_G(u)$ that covers the factor
$Z_2^{t+\d}$ in (\ref{cov});
\item[(e)] construct the subgroup $R$ of $C_G(u)$.
\end{itemize}
Each step typically requires more effort than the corresponding 
step in good characteristic.

\subsection{Constructing $C_Q(u)$}\label{cqubad}
As in good characteristic, the main task is to write down 
a generating set for the unipotent group $Q$. 
As discussed in Section \ref{gencentgood}, we then construct 
a power-conjugate presentation for $Q$ and so construct $C_Q(u)$.

\subsubsection{Generators for $Q$} 
For $G = \Sp(V)$, we can use exactly the generators given in 
Section \ref{spcent} (the characteristic is 2, so the sign $\e$ can be ignored).

Now consider $G = \O(V)$. The generators for $Q$ 
differ from those in good characteristic. 
We define ten collections of generators for $Q$.

\vspace{2mm} \no (1)  For all $i>0$, all $k,l$, and every $\l \in \F_q$ 
define the following generators 
which fix all basis vectors other than those listed.
\[
\begin{array}{ll}
x_{1ikl}(\l): & w_{-i}^{(k)} \mapsto w_{-i}^{(k)}+\l x_0^{(l)} + \l^2 Q(x_0^{(l)})x_i^{(k)}, \\
                     & w_0^{(l)}\mapsto w_0^{(l)}+ \l x_i^{(k)} \\
\\
x_{2ikl}(\l): & w_{-i}^{(k)} \mapsto w_{-i}^{(k)}+\l w_0^{(l)} + \l^2 Q(w_0^{(l)})x_i^{(k)}, \\
                     & x_0^{(l)}\mapsto x_0^{(l)}+ \l x_i^{(k)} \\
\\
x_{3ikl}(\l): & x_{-i}^{(k)} \mapsto x_{-i}^{(k)}+\l x_0^{(l)} + \l^2 Q(x_0^{(l)})w_i^{(k)}, \\
                     & w_0^{(l)}\mapsto w_0^{(l)}+ \l w_i^{(k)} \\
\\
x_{4ikl}(\l): & x_{-i}^{(k)} \mapsto w_{-i}^{(k)}+\l w_0^{(l)} + \l^2 Q(w_0^{(l)})w_i^{(k)}. \\
                     & x_0^{(l)}\mapsto x_0^{(l)}+ \l w_i^{(k)} \\
\end{array}
\]
\vspace{2mm} \no (2)  For all $i,j\ne 0$ with $j<i$, all $k,l$, and every $\l \in \F_q$:
\[
\begin{array}{ll}
x_{5ijkl}(\l): & w_{-i}^{(k)} \mapsto w_{-i}^{(k)}+\l w_{-j}^{(l)}, \\
                     & x_j^{(l)}\mapsto x_j^{(l)}+ \l x_{i}^{(k)}\\
\\
x_{6ijkl}(\l): & w_{-i}^{(k)} \mapsto w_{-i}^{(k)}+\l x_{-j}^{(l)}, \\
                     & w_j^{(l)}\mapsto w_j^{(l)}+ \l x_{i}^{(k)}\\
\\
x_{7ijkl}(\l): & x_{-i}^{(k)} \mapsto x_{-i}^{(k)}+\l w_{-j}^{(l)}. \\
                     & x_j^{(l)}\mapsto x_j^{(l)}+ \l w_{i}^{(k)}
\end{array}
\]
\vspace{2mm} \no (3)  For all $i>1$, $j<i$ with $j\ne \pm 1$, all $k,l$, and every $\l \in \F_q$:
\[
\begin{array}{ll}
x_{8ijkl}(\l): & v_{-i}^{(k)} \mapsto v_{-i}^{(k)}+\l v_{-j}^{(l)}. \\
                     & v_j^{(l)}\mapsto v_j^{(l)}+ \l v_{i}^{(k)}
\end{array}
\]
\vspace{2mm} \no (4)  For all $j<i$ with $j\ne \pm 1$, all $k,l$, and every $\l \in \F_q$:
\[
\begin{array}{ll}
x_{9ijkl}(\l): & w_{-i}^{(k)} \mapsto w_{-i}^{(k)}+\l v_{-j}^{(l)}, \\
                     & v_j^{(l)}\mapsto v_j^{(l)}+ \l x_{i}^{(k)}
\\
x_{10,ijkl}(\l): & x_{-i}^{(k)} \mapsto x_{-i}^{(k)}+\l v_{-j}^{(l)}. \\
                     & v_j^{(l)}\mapsto v_j^{(l)}+ \l w_{i}^{(k)}
\end{array}
\]
\vspace{2mm} \no (5)  For all $i\ge 1$, all $k,l$, and every $\l \in \F_q$:
\[
\begin{array}{ll}
x_{11,ikl}(\l): & w_{-i}^{(k)} \mapsto w_{-i}^{(k)}+\l v_1^{(l)} + \l^2 x_{i}^{(k)}, \\
                     & v_{-1}^{(l)}\mapsto v_{-1}^{(l)}+ \l x_{i}^{(k)} \\
\\
x_{12,ikl}(\l): & x_{-i}^{(k)} \mapsto x_{-i}^{(k)}+\l v_1^{(l)} + \l^2 w_{i}^{(k)}.  \\
                     & v_{-1}^{(l)}\mapsto v_{-1}^{(l)}+ \l w_{i}^{(k)}
\end{array}
\]
\vspace{2mm} \no (6)  For all $i> 1$, all $k,l$, and every $\l,\mu \in \F_q$:
\[
\begin{array}{ll}
x_{13,ikl}(\l,\mu): & v_{-i}^{(k)} \mapsto v_{-i}^{(k)}+\l w_0^{(l)}+\mu x_0^{(l)} + (\l\mu +(\l^2+\mu^2)Q(w_0^{(l)})) v_{i}^{(k)}, \\
                     & w_0^{(l)}\mapsto w_0^{(l)}+ \mu v_{i}^{(k)} \\
                      & x_0^{(l)}\mapsto x_0^{(l)}+ \l v_{i}^{(k)} \\
\\
x_{14,ikl}(\l): & v_{-i}^{(k)} \mapsto v_{-i}^{(k)}+\l v_i^{(l)}.  \\
                     & v_{-i}^{(l)}\mapsto v_{-i}^{(l)}+ \l v_{i}^{(k)}
\end{array}
\]
\vspace{2mm} \no (7)  For all $i> 1$, all $k,l$, and every $\l,\mu \in \F_q$:
\[
\begin{array}{ll}
x_{15,ikl}(\l,\mu): & v_{-i}^{(k)} \mapsto v_{-i}^{(k)}+\l v_{-1}^{(l)}+\mu v_1^{(l)} + (\l\mu +\mu^2+\l^2Q(v_{-1}^{(l)})) v_{i}^{(k)}. \\
                     & v_{-1}^{(l)}\mapsto v_{-1}^{(l)}+ \mu v_{i}^{(k)} \\
                      & v_1^{(l)}\mapsto v_1^{(l)}+ \l v_{i}^{(k)} 
\end{array}
\]
\vspace{2mm} \no (8)  For all $k,l$, and every $\l,\mu\in \F_q$ satisfying $\mu^2 = \l^2+\l$:
\[
\begin{array}{ll}
x_{16,kl}(\l): & v_{-1}^{(k)} \mapsto v_{-1}^{(k)}+\l v_1^{(k)}+\mu v_1^{(l)}. \\
                     & v_{-1}^{(l)} \mapsto v_{-1}^{(l)}+\mu v_1^{(k)}+\l v_1^{(l)}
\end{array}
\]
\vspace{2mm} \no (9)  For all $k,l,m$, and every $\l,\mu,\g,\d\in \F_q$ satisfying $\d^2 = \l\mu+\g^2+\g+(\l^2+\mu^2)Q(w_0^{(l)})$:
\[
\begin{array}{ll}
x_{17,klm}(\l,\mu,\g): & v_{-1}^{(k)} \mapsto v_{-1}^{(k)}+\l w_0^{(l)}+\mu x_0^{(l)} + \g v_1^{(k)} +
\d v_{1}^{(m)}. \\
                     & v_{-1}^{(m)} \mapsto v_{-1}^{(m)}+\l w_0^{(l)}+\mu x_0^{(l)} + \d v_1^{(k)} +
\g v_{1}^{(m)} \\
                     & w_0^{(l)}\mapsto w_0^{(l)}+ \mu v_{1}^{(k)} +  \mu v_{1}^{(m)} \\
                      & x_0^{(l)}\mapsto x_0^{(l)}+ \l v_{1}^{(k)}+ \l v_{1}^{(m)} \\
\end{array}
\]
\vspace{2mm} \no (10)  For all $k,l,m$:
\[
\begin{array}{ll}
x_{18,klm}: & v_{-1}^{(k)} \mapsto v_{-1}^{(k)}+v_1^{(l)}+v_1^{(m)}. \\
             & v_{-1}^{(l)} \mapsto v_{-1}^{(l)}+v_1^{(k)}+v_1^{(m)} \\
             & v_{-1}^{(m)} \mapsto v_{-1}^{(m)}+v_1^{(k)}+v_1^{(l)}
\end{array}
\]

\subsection{Factors $\Sp_{2a_i}(q)$}\label{fsp}

These factors arise from a summand $W(m_i)^{[a_i]}$ in (\ref{newcan}) with $m_i$ even. Let  $u = W(2l)^{[a]}$, with basis labelled 
\[
w_{-(2l-1)}^{(1)},x_{-(2l-1)}^{(1)},\ldots, w_{-(2l-1)}^{(a)},x_{-(2l-1)}^{(a)},\ldots ,w_{2l-1}^{(1)},x_{2l-1}^{(1)},\ldots, 
w_{2l-1}^{(a)},x_{2l-1}^{(a)}.
\]
Relative to this basis, the symplectic form has matrix
\[
F = \begin{pmatrix} &&&K \\&&K& \\& \iddots && \\ K&&& \end{pmatrix}
\]
where $K = J^{[a]}$ (the block diagonal sum of $a$ copies of $J = \begin{pmatrix} &1 \\ 1& \end{pmatrix}$). 

Let $P = QL$ be the associated parabolic subgroup of $G = \Sp_{4al}(q)$ or $\O^+_{4al}(q)$. Elements of $C_P(u)$ project to elements of $L$ 
having shape $x^{[l]} = x\oplus \cdots \oplus x$, where $x \in \Sp_{2a}(q)$ 
preserves the symplectic form with matrix $K$. 
Choose standard generators $x_i\,(1\le i\le r)$ for this $\Sp_{2a}(q)$, and write $X_i = x_i^{[l]}$. For each $i$, we compute 
$v_i \in Q$ such that $u^{X_i} = u^{v_i}$. 
(This computation is carried out in the isomorphic copy of $Q$ described 
by the power-conjugate presentation.) 
Then $X_iv_i^{-1} \in C_P(u)$, and 
\[
\la X_iv_i^{-1} : 1\le i\le r \ra \le C_P(u)
\]
maps onto a subgroup $\Sp_{2a}(q)$ of the Levi subgroup $L$. 
This gives us the factor $\Sp_{2a}(q)$ of $C_G(u)$.

\subsection{Factors $I_{2a_i}(q)$ for $\Sp(V)$}\label{fi2a}

Let $G = \Sp(V)$ and let $u$ be as in (\ref{newcan}). Here we construct the factors $I_{2a_i}(q) = \Sp_{2a_i}(q)$ or 
$\Or_{2a_i}^{\e_i}(q)$ of $C_G(u)$, given by (\ref{cov}).

First, if  $m_i=1$ for some $i$, then there is a summand $W(1)^{[a_i]}$ in (\ref{newcan}), and clearly $C_G(u)$ has a subgroup $\Sp_{2a_i}(q)$ acting on this summand (and trivially on its orthogonal complement). We choose standard generators to generate this subgroup.

The other case where there is a factor $I_{2a_i}(q)=\Sp_{2a_i}(q)$ occurs when (\ref{newcan}) has a summand 
$W(m_i)^{[a_i]} \oplus V_{\b}(2k_j)$, where $m_i$ is odd and $m_i = 2k_j\pm 1$. 
To simplify notation, denote the summand by $W(2l+1)^{[a]} \oplus V_{\b}(2k)$, where $2l+1 = 2k\pm 1$. This is handled as in Section \ref{fsp}: we choose standard generators $x_i$ ($1\le i \le r$) of $\Sp_{2a}(q)$ preserving the symplectic form with matrix $K$, and let $X_i = x_i^{[2l+1]}$ act on $W(2l+1)^{[a]}$ and 
trivially on $V_\b(2k)$. 
For each $i$, compute $v_i \in Q$ such that $u^{X_i} = u^{v_i}$. Then $\la X_iv_i^{-1} : 1\le i\le r \ra$ is a subgroup of  $C_P(u)$ that maps onto $\Sp_{2a}(q)$, as required.

Finally, we deal with the case where the factor $I_{2a_i}(q)$ is $\Or^{\e_i}_{2a_i}(q)$. This occurs when (\ref{newcan}) has a summand $W_{\b_i,a_i}(m_i)$ with $m_i$ odd, but no summand $V_{\b}(2k_j)$ such that $m_i = 2k_j\pm 1$. Write the summand as $u = W_\b(2l+1) \oplus W(2l+1)^{[a-1]}$. In this case,  elements of $C_P(u)$ project to elements of $L$ having shape 
$x^{[2l+1]}$, where $x \in \Or^{\e}_{2a}(q)$ preserves the quadratic form with matrix 
\begin{equation}\label{squad}
S = \begin{pmatrix} \b&1&&&&& \\ 0&\b&&&&& \\&&0&1&&&\\ &&0&0&&& \\&&&& \ddots && \\ &&&&&0&1 \\ &&&&&0&0 \end{pmatrix}.
\end{equation}
 So we choose standard generators $x_i$ ($1\le i \le r$) of $\Or^{\e}_{2a}(q)$, and let $X_i = x_i^{[2l+1]}$ act on $W_\b(2l+1) \oplus W(2l+1)^{[a-1]}$. 
For each $i$, compute $v_i \in Q$ such that $u^{X_i} = u^{v_i}$. Then $\la X_iv_i^{-1} : 1\le i\le r \ra$ is a subgroup of  $C_P(u)$ that maps onto $\Or^\e_{2a}(q)$, as required.

\subsection{Factors $I_{2a_i}(q)$ for $\O(V)$}\label{fi2ao}
Let $G = \O(V)$ and let $u$ be as in (\ref{newcan}). 
Here we construct the factors $I_{2a_i}(q) = \Sp_{2a_i}(q)$ or 
$\Or_{2a_i}^{\e_i}(q)$ of $C_G(u)$, given by (\ref{cov}).

The orthogonal factors 
occur when (\ref{newcan}) has a summand $W_{\b_i,a_i}(m_i)$ with $m_i$ odd, but no summand $V_{\b}(2k_j)$ such that $m_i = 2k_j\pm 1$. 
Write the summand as $u = W_\b(2l+1) \oplus W(2l+1)^{[a-1]}$. We produce these 
factors as in the previous section, with a small tweak. 
Let $\Or^{\e}_{2a}(q)$ be the orthogonal group preserving the quadratic form with matrix $S$ as in (\ref{squad}). Choose standard generators $x_i$ ($1\le i \le s$) of $\O^{\e}_{2a}(q)$, and a reflection $r \in \Or^{\e}_{2a}(q)\setminus \O^{\e}_{2a}(q)$. 
Let $X_i = x_i^{[2l+1]}$ and $Y = r^{[2l+1]}$ act 
on $W_\b(2l+1) \oplus W(2l+1)^{[a-1]}$. 
Now compute $v_i,w \in Q$ such that $u^{X_i} = u^{v_i}$ and $u^Y=u^w$. Then $\la X_iv_i^{-1} : 1\le i\le s \ra$ is a subgroup of  $C_P(u)$ that maps onto $\O^\e_{2a}(q)$, 
and $Yw^{-1}$ is an element of $C_{\Or(V)}(u)\setminus C_{\O(V)}(u)$ that maps to $r$. 
Write $r'=Yw^{-1}$, and let $r_1',\ldots,r_t'$ be the collection of such elements, one for each orthogonal factor in (\ref{cov}). 
(We use these elements in Section \ref{z2fac}.)

The symplectic factors $I_{2a_i}(q) = \Sp_{2a_i}(q)$
occur when there is a summand $W(m_i)^{[a_i]} \oplus V_{\b}(2k_j)$, where $m_i = 2k_j\pm 1$. Since $G = \O(V)$, there must be another summand $V_\g(2k')$. So let $u = W(2l+1)^{[a]} \oplus V_{\b}(2k) \oplus V_\g(2k')$, where $2l+1 = 2k\pm 1$. First, as in the previous 
paragraph, we produce a subgroup $X$ of $C_P(u)$ acting trivially on $V_{\b}(2k) \oplus V_\g(2k')$ such that the action of $X$ on $W(2l+1)^{[a]}$  maps onto $\O^\e_{2a}(q)$. Secondly, we produce two additional elements of $C_P(u)$ as follows. 
With the usual labelling of basis vectors, 
let $w_i^{(j)},x_i^{(j)}$ ($-2l\le i\le 2l$, $1\le j\le a$) be 
the basis of $W(2l+1)^{[a]}$, and label the bases of the two $V$-blocks as
\[
\begin{array}{ll}
v_{-(2k-1)},\ldots,v_{2k-1} & \hbox{for }V_{\b}(2k), \\
v'_{-(2k'-1)},\ldots,v'_{2k'-1} & \hbox{for }V_{\g}(2k').
\end{array}
\]
Now define the following linear maps 
which fix all basis vectors other than those listed:
\[
\begin{array}{llll}
s_1:&  w_0^{(1)} & \mapsto & w_0^{(1)}+x_0^{(1)}+v_1', \\
    &    v_{-1}' & \mapsto & v_{-1}'+x_0^{(1)}, \\
    &    w_i^{(1)} & \mapsto & w_i^{(1)}+ x_i^{(1)} \hbox{ for }i\ge 1 \\
\\
s_2:&  x_0^{(1)} &\mapsto & x_0^{(1)}+w_0^{(1)}+v_1, \\
    &    v_{-1} &\mapsto & v_{-1}+w_0^{(1)}, \\
    &    x_i^{(1)} &\mapsto & x_i^{(1)}+ w_i^{(1)} \hbox{ for }i\ge 1.
\end{array}
\]
Observe that $s_1,s_2 \in \O(V)$ and $\la X,s_1,s_2\ra$ maps onto 
$\Sp_{2a}(q)$ in the Levi subgroup $L$. So if we compute $y_i\in Q$ 
such that $u^{s_i}=u^{y_i}$ for $i=1,2$, 
then $\la X,\,s_1y_1^{-1}, \,s_2y_2^{-1}\ra$ is 
a subgroup of $C_P(u)$ mapping onto $\Sp_{2a}(q)$, as required. 

\subsection{The factor $Z_2^{t+\d}$} \label{z2fac}
Let $u$ be as in (\ref{newcan}). For $G = \Sp(V)$, we cover the 
factor $Z_2^{t+\d}$ simply by adding as generators of $C_G(u)$ each 
$V$-block in (\ref{newcan}), acting trivially on the orthogonal complement. 
Let $v_1,\ldots,v_s$ be these generators. 
For $G = \O(V)$ these elements lie in $\Or(V) \setminus G$;  so we add as 
generators of $C_G(u)$ the pairwise products among 
$v_1,\ldots,v_s,r_1',\ldots,r_t'$, where the $r_i'$ are the generators 
defined in Section \ref{fi2ao}. 
Let $S$ be the subgroup generated by these generators. Note that to construct the 
centralizer of $u$ in $\Or(V)$, we include all of the generators 
$v_1,\ldots,v_s,r_1',\ldots,r_t'$.

\subsection{Unipotent radical $R$}
Let $u$ be as in (\ref{newcan}). 
We now compute the unipotent radical $R$ of $C_G(u)$; it has order $q^{M_u}$ as 
in Theorem \ref{uni}(i). We already computed $C_Q(u)$ in Section \ref{cqubad}. 
By Conjecture \ref{rcon}, which we assume, $|RC_Q(u)/C_Q(u)| = q^t$, 
where $J = \{j: b_j=2\}$ and $t = \sum_{j\in J}(2a_j+1)$.
Each summand of (\ref{newcan}) of the form $W(2k_j)^{[a_j]} \oplus V_{\b_j,2}(2k_j)$ 
contributes $q^{2a_j+1}$ to this order. We let 
\begin{equation}\label{2a1}
u = W(2k)^{[a]} \oplus V(2k)  \oplus V_{\b}(2k),
\end{equation}
and aim to construct a subgroup $R$ of $C_G(u)$ such that 
$|RC_Q(u)/C_Q(u)| = q^{2a+1}$.

\subsubsection{Symplectic case } 
Let $G = \Sp(V)$, and let $u$ be as in (\ref{2a1}). We compute $R$ in four steps.

\begin{enumerate}
\item  In Section \ref{fsp} we constructed a subgroup, say $X$,
of $C_P(u)$ acting on the summand $W(2k)^a$  that maps onto $\Sp_{2a}(q)$. 

\item 
Consider the summand $U=V(2k)  \oplus  V_{\b}(2k)$ of (\ref{2a1}), and let 
$u_0 = u\downarrow U$. We work in $\Sp(U)$. Our basis of $U$ is 
\begin{equation}\label{vvbas}
v_{-(2k-1)}^{(1)}, v_{-(2k-1)}^{(2)}, \ldots, v_{2k-1}^{(1)}, v_{2k-1}^{(2)}, 
\end{equation}
with respect to which the symplectic form has matrix 
\begin{equation}\label{symkk}
\begin{pmatrix} && I_2 \\ & \iddots & \\ I_2 && \end{pmatrix}.
\end{equation}
If $A^{[2k]} = A \oplus \cdots \oplus A$ preserves the form (where $A$ is $2\times 2$), then $A$ lies in $Y_0$, where 
\[
Y_0 = \{A \in \GL_2(q) : AA^\tr=I\} = \left\{ \begin{pmatrix} a+1&a \\ a&a+1 \end{pmatrix} : a \in \F_q \right\}
\cong \F_q^+.
\]
Let $$M_a = \begin{pmatrix} a+1&a \\ a&a+1 \end{pmatrix}.$$ 
For $a \in \F_q$, let $l_a = M_a^{[2k]}$, and compute $y_a \in Q$ such that $u_0^{l_a} = u_0^{y_a}$. Now define
\[
Y = \la l_ay_a^{-1} : a \in \F_q\ra,
\]
a subgroup of $C_P(u_0)$ that maps onto a subgroup of $L$ isomorphic to $\F_q^+$.

\item Suppose $a\ge 1$, and consider a summand $U_1=W(2k) \oplus V(2k)  \oplus V_{\b}(2k)$ of (\ref{2a1}); let 
$u_1 = u\downarrow U_1$. 
We work in $\Sp(U_1)$. Our basis of $U_1$ is 
\begin{equation}\label{wxbass}
w_{-(2k-1)}^{(1)}, x_{-(2k-1)}^{(1)}v_{-(2k-1)}^{(2)}, v_{-(2k-1)}^{(3)}, \ldots, w_{2k-1}^{(1)}, x_{2k-1}^{(1)}v_{2k-1}^{(2)}, v_{2k-1}^{(3)}, 
\end{equation}
with respect to which the symplectic form has matrix 
\[
\begin{pmatrix} && M \\ & \iddots & \\ M && \end{pmatrix}, \hbox{ where }M = \begin{pmatrix} 0&1&& \\1&0&& \\ &&1&0 \\ &&0&1 \end{pmatrix}.
\]
If $A^{[2k]}$ preserves the form (where $A$ is $4\times 4$), then $AMA^\tr=M$. One such matrix is 
\[
A = \begin{pmatrix} 1&0&1&1 \\ 0&1&1&1 \\ 1&1&1&0 \\ 1&1&0&1 \end{pmatrix}.
\]
Let $l =  A^{[2k]}$, and compute $v \in Q$ such that $u_1^l=u_1^v$. Let $l_0 = lv^{-1} \in C_P(u_1)$.
Let $X_0$ be the subgroup of $C_P(u_1)$ acting on $W(2k)$ that maps onto $\Sp_2(q)$ (constructed in Step 1). Let  
$Y_0$ be the subgroup of $C_P(u_1)$ acting on $V(2k) + V_{\b}(2k)$ that maps 
onto $\F_q^+$ (constructed in Step 2). 
Now $Z=\la X_0,Y_0,l_0\ra$ is a subgroup of $C_P(u_1)$ that maps onto the following subgroup $Z_0$ of $\GL_4(q)$ in the Levi subgroup $L$:
\[
Z_0 = \left \la \begin{pmatrix} A& \\  I \end{pmatrix},\, \begin{pmatrix} I& \\ & M_a \end{pmatrix}, \,l_0 : A \in \Sp_2(q),\,a \in \F_q \right \ra.
\]
Observe that $Z_0$ has the form $[q^3].\Sp_2(q)$ (where $[q^3]$ denotes a group of order $q^3$). The subgroup of $Z$ mapping to the normal subgroup of order $q^3$ gives the subgroup $R$ for the case $a=1$ (in which case $|RC_Q(u)/C_Q(u)| = q^3$).

\item 
Let $u$ be as in (\ref{2a1}). If $a=0$, 
then we constructed $R$ in Step 2. Assume that $a \ge 1$. 
Let $Z$ be the subgroup constructed in Step 3, acting on a summand 
$W(2k) \oplus V(2k)  \oplus  V_{\b}(2k)$ and trivially on the orthogonal complement, and let $X$ be the subgroup in Step 1 that maps onto 
$\Sp_{2a}(q)$. Define
\[
C = \la X,\,Z \ra \le C_P(u).
\]
Now $C$ maps onto a subgroup $[q^{2a+1}].\Sp_{2a}(q)$ of $L$. 
Hence, adding the subgroup $S$ of Section \ref{z2fac},  
\[
\la C_Q(u),\,C,\,S \ra = C_P(u).
\]
\end{enumerate}
This completes the construction of $C_G(u)$ when $G = \Sp(V)$.

\subsubsection{Orthogonal case} 
Let $G = \O(V)$, and let $u$ be as in (\ref{2a1}). We again compute 
$R$ in four steps, but some details are different. 

\begin{enumerate}
\item 
This step is identical to that for the symplectic case: 
in Section \ref{fsp} we constructed a subgroup, say $X$, of $C_P(u)$ acting on 
the summand $W(2k)^a$  that maps onto $\Sp_{2a}(q)$. 

\item 
We consider the summand $U=V(2k)  \oplus  V_{\b}(2k)$ of (\ref{2a1}), with basis as in (\ref{vvbas}), and let 
$u_0 = u\downarrow U$. We work in $\O(U)$, with quadratic form $Q_\b$ as given in Section \ref{badsec}: namely, $Q_\b$ takes values $1,1,\b,\b$ on 
$v_{1}^{(1)}$, $v_1^{(2)}$, $v_{-1}^{(1)}$, $v_{-1}^{(2)}$ 
respectively and 0 on the other basis vectors, and the associated symplectic form is (\ref{symkk}). When $\b \ne 0$, the matrix $l_a = M_a^{[2k]}$ defined in Step 2 of the symplectic case
does not preserve the quadratic form $Q_\b$. 
We adjust it as follows: for $a \in \F_q$, define $l_a' \in \O(U)$ by
\[
\begin{array}{llll}
l_a': & v_i^{(1)} & \mapsto & (a+1)v_i^{(1)}+av_i^{(2)} \hbox{ for all }i \ne -1 \\
       & v_i^{(2)} & \mapsto & av_i^{(1)}+(a+1)v_i^{(2)} \hbox{ for all }i \ne -1 \\
       & v_{-1}^{(1)} & \mapsto&  (a+1)v_{-1}^{(1)}+av_{-1}^{(2)} + a^2\b v_1^{(1)}+a^2\b v_1^{(2)} \\
       & v_{-1}^{(2)} & \mapsto & av_{-1}^{(1)}+(a+1)v_{-1}^{(2)} + a^2\b v_1^{(1)}+a^2\b v_1^{(2)},
\end{array}
\]
and it
fixes all other basis vectors. Compute $y_a \in Q$ such that $u_0^{l_a'} = u_0^{y_a}$. Now $Y = \la l_a'y_a^{-1} : a \in \F_q\ra$ is a subgroup of $C_P(u_0)$ that maps onto a subgroup of $L$ isomorphic to $\F_q^+$.

\item 
We work in $\O(U_1)$, where $U_1=W(2k) \oplus V(2k)  \oplus  V_{\b}(2k)$ is a summand 
of (\ref{2a1}), with basis as in (\ref{wxbass}). Let $u_1 = u\downarrow U_1$. 
Again, we need to tweak the element $l =  A^{[2k]}$ defined 
in Step 3 of the symplectic case, 
as it does not preserve the quadratic form. Define $l' \in \O(U_1)$ by
\[
l' = A^{[k-1]} \oplus B \oplus A^{[k-1]},
\]
where $A$ is the $4\times 4$ matrix in the symplectic case, 
and $B$ is the $8\times 8$ matrix corresponding to the 
basis vectors of weights $-1$ and 1,
$$w_{-1}^{(1)}, x_{-1}^{(1)}, v_{-1}^{(2)}, v_{-1}^{(3)}, 
w_{1}^{(1)}, x_{1}^{(1)}, v_{1}^{(2)}, v_{1}^{(3)};$$
so 
\[
B = \begin{pmatrix} 1&0&1&1&0&\b&& \\ 0&1&1&1&\b&0&& \\ 1&1&1&0&\b&\b&& \\ 
1&1&0&1&\b&\b&& \\  &&&&1&0&1&1 \\  &&&&0&1&1&1 \\  &&&&1&1&1&0 \\  &&&&1&1&0&1 \end{pmatrix}.
\]
Now compute $v \in Q$ such that $u_1^{l'}=u_1^v$, and let $l_0 = l'v^{-1} \in C_P(u_1)$. Using this element we construct the subgroup $R$ for the case $a=1$ in similar fashion to the corresponding step for the symplectic case.

\item 
This step is identical to that for the symplectic case. 
\end{enumerate}
This completes the construction of $C_G(u)$ when $u \in G = \O(V)$. To construct $C_{\Or(V)}(u)$, we include the extra 
generators indicated in Section \ref{z2fac}.

Let $G = \Or(V)$.  If $u \in G \setminus \O(V)$, then we cannot write down 
directly $C_{G}(u)$.  
But we use our approach to reduce its construction 
to a smaller problem:  
since $u^2 \in \O(V)$, we construct $H := C_{G}(u^2)$; 
if $H$ is sufficiently ``small", then 
we can use a standard centralizer algorithm \cite[Chapter~4]{handbook} to 
construct $C_H(u)$.

\section{The conjugacy problem}\label{conjprobbad}
Let $G = \Sp(V)$ or $\Or(V)$ in characteristic 2. In this section 
we address the conjugacy problem: given unipotent $g \in G$, 
find the canonical class representative as in (\ref{canon}) that 
is conjugate to $g$. Before presenting the algorithm to do this, 
we prove five necessary lemmas.

\begin{lem}\label{reduc} Let $u = V_\b(2k)$ or $W_\b(2l+1)$, as defined in Section $\ref{badsec}$, where $\b \in \{0,\a\}$. For $r\ge 1$, define $F_r = V(1-u)^r$. The actions of $u$ on the spaces $F_r/{\rm Rad}(F_r)$ and 
$F_r^\perp/{\rm Rad}(F_r^\perp)$ are as in Table $\ref{actions}$.
\end{lem}

\begin{table}[ht]
\caption{Actions of $u$ on spaces}\label{actions} 
\[
\begin{array}{|l|l|l|l|}
\hline
u & r & u^{F_r/{\rm Rad}(F_r)} & u^{F_r^\perp/{\rm Rad}(F_r^\perp)} \\
\hline
V_\b(2k) & r\le k & V_\b(2k-2r) & 0 \\
             & k<r\le 2k & 0 & V_\b(2r-2k) \\
W_\b(2l+1) & r\le l & W_\b(2l+1-2r) & 0 \\
                 & l < r \le 2l+1 & 0 & W_\b(2r-2l-1) \\
\hline
\end{array}
\]
\end{table}
\no In the table, the convention is that for $\Sp(V)$, the elements $W_\b(1)$ and $V_\b(2)$ are just $W(1)$ and $V(2)$, whatever the value of $\b$.

\vspace{2mm}
\begin{proof} Let $u = V_\b(2k)$, defined as in Section \ref{badsec} relative to the basis $v_{-(2k-1)},\ldots, v_{2k-1}$. If $r\le k$, then $F_r = V(1-u)^r$ has basis $v_{-(2k-2r-1)},\ldots, v_{2k-1}$, and ${\rm Rad}(F_r) = \la v_{2k-2r+1},\ldots, v_{2k-1}\ra$. Hence  $u$ acts on $F_r/{\rm Rad}(F_r)$ as $V_\b(2k-2r)$, as claimed. The other cases are similar.  \end{proof}

\begin{lem}\label{j2test} Let $u = W(2)^{[a]} \oplus V(2)^{[b]} \oplus V_\b(2)^{[c]} \in \Sp(V)$ or $\Or(V)$, where $b+c \le 2$ and $\dim V = 4a+2b+2c$. Let $w_1,\ldots,w_r$ be a basis of $V(u-1)$, and, for each $i$, let $v_i \in V$ be such that $v_i(u-1)=w_i$. 
Then $b+c =0$ if and only if $(v_i,w_i)=0$ for all $i$.
\end{lem}

\begin{proof} If $b+c=0$, then $u=W(2)^{[a]}$, so $(v,\,v(u-1))=0$ for all $v \in V$.

For the converse, suppose $(v_i,w_i)=0$ for all $i$. Assume for a contradiction that $b+c\ge 1$. We assume that $b+c=2$, the case $b+c=1$ being similar and easier. Let $X$ be a non-degenerate 4-space on which $u$ acts as $V(2)^{[b]} \oplus V_\b(2)^{[c]}$. Then $u$ acts as $W(2)^{[a]}$ on $X^\perp$, so $(v,\,v(u-1))=0$ for all $v \in X^\perp$.
Choose a basis $e_1,e_2,f_1,f_2$ of $X$ such that $u-1$ maps $e_1\mapsto f_1 \mapsto 0$, $e_2\mapsto f_2 \mapsto 0$. 
There are two vectors, say $w_1$ and $w_2$, in the basis of $V(u-1)$
whose projections to $X$ are linearly independent. So we may write
\[
w_1=\a_1f_1+\a_2f_2+x,\;\;w_2=\b_1f_1+\b_2f_2+x',
\]
where $x,x' \in X^\perp$ and $\a_1f_1+\a_2f_2,\,\b_1f_1+\b_2f_2$ are linearly independent. Since $v_i(u-1)=w_i$, 
\[
\begin{array}{l}
v_1=\a_1e_1+\a_2e_2+\g_1f_1+\g_2f_2+y, \\
v_2=\b_1e_1+\b_2e_2+\d_1f_1+\d_2f_2+y',
\end{array}
\]
where $y,y' \in X^\perp$ and $y(u-1)=x$, $y'(u-1)=x'$. Then $(y,x)=(y',x')=0$, and hence
\[
\begin{array}{l}
(v_1,w_1)=0 \Rightarrow \a_1^2+\a_2^2 = 0 \Rightarrow \a_1=\a_2, \\
(v_2,w_2)=0 \Rightarrow \b_1^2+\b_2^2 = 0 \Rightarrow \b_1=\b_2.
\end{array}
\]
This contradicts the linear independence 
of $\a_1f_1+\a_2f_2$ and $\b_1f_1+\b_2f_2$. 
\end{proof}

\begin{lem}\label{j4test} Let $u = V_\b(4) \in \Sp(V)$ or $\Or(V)$, where $\dim V=4$. Suppose $v \in V$ is such that $v(u-1)^3 \ne 0$ and $(v(u-1),\,v(u-1)^2) \ne 0$. Write $v_i = v(u-1)^{i-1}$ for $1\le i\le 4$, and define
\[
\l = (v_2,v_3),\;\mu=(v_1,v_2).
\]
Then $\b=\a$ if and only if $x^2+x+\l^{-1}\mu$ is irreducible over $\F_q$.
\end{lem}

\begin{proof} 
Relative to the basis $v_1,v_2,v_3,v_4$, the matrix of $u$ is $J_4$ and the matrix of the form $(\,,\,)$ is
\[
\begin{pmatrix} 0&\mu&\l&\l \\ \mu & 0&\l&0 \\ \l&\l&0&0 \\ \l&0&0&0 \end{pmatrix}.
\]
We change the basis as follows: let $\d \in \F_q$ with $\d^2=\l^{-1}$, and set
\[
\begin{array}{l}
e_1=\d v_1+\d v_2+\d^3\mu v_3, \\
e_2 = \d v_2, \\
f_2 = \d v_3, \\
f_1 = \d v_4.
\end{array}
\]
Now, relative to the basis $B = [e_1,e_2,f_2,f_1]$, 
\[
[u]_B = \begin{pmatrix} 1&1&1&\l^{-1}\mu \\ &1&1&0 \\&&1&1 \\&&&1 \end{pmatrix},\;\;\;
(\,,\,)_B = \begin{pmatrix} &&&1 \\&&1& \\&1&& \\1&&& \end{pmatrix}.
\]
The matrix of $u$ is that of $V_{\l^{-1}\mu}(4)$, and the conclusion follows. 
\end{proof}

\begin{lem}\label{w3test} Let $u = W_\b(3)\in \Sp(V)$ or $\Or(V)$, where $\dim V=6$. Suppose $v_1,w_1 \in V$ are such that \[
v_1,v_1(u-1),v_1(u-1)^2,\,w_1,w_1(u-1),w_1(u-1)^2
\]
 is a basis of $V$. Let $v_2=v_1(u-1)$, $w_2=w_1(u-1)$, and set 
\[
\g_1 = (v_2,w_2),\;\g_2=(v_1,v_2),\;\g_3=(w_1,w_2).
\]
Then $\b=\a$ if and only if the quadratic $x^2+\g_1x+\g_2\g_3$ is irreducible over $\F_q$.
\end{lem}

\begin{proof}
Recall the basis defining $W_\b(3)$ given in Section \ref{badsec}, namely
\[
w_{-2},x_{-2},w_0,x_0,w_2,x_2.
\]
Rewrite this basis as $e_1,f_3,e_2,f_2,e_3,f_1$ respectively, so that $(e_i,f_j) = \d_{ij}$. We claim first that 
\begin{equation}\label{nonz}
\b=0 \Leftrightarrow \exists\, v \in V \hbox{ such that }(v,\,v(u-1)) = 0 \hbox{ and }v(u-1)^2 \ne 0.
\end{equation}
The left to right implication is clear: if $\b=0$, then $v=e_1$ has the required property. Suppose $\b=\a$. Let $v\in V$ be such that $v(u-1)^2 \ne 0$. Writing $v = \sum \l_ie_i+\sum \mu_if_i$, this means that $\l_1\ne 0$ or $\mu_3\ne 0$. We now compute that 
\[
(v,\,v(u-1)) = \l_1^2\a +\l_1\mu_3+\mu_3^2\a.
\]
Taking $\l_1\ne 0$, this gives 
\[
(v,\,v(u-1))  = \l_1^2\a^{-1}(\a^2 + \mu_3\l_1^{-1}\a + ( \mu_3\l_1^{-1}\a)^2),
\]
which is nonzero since the quadratic $x^2+x+\a$ (hence also $x^2+x+\a^2$) is irreducible over $\F_q$. If $\l_1=0$, then also 
$(v,\,v(u-1)) \ne 0$. This proves (\ref{nonz}).

Next we establish one of the implications in the lemma: namely
\begin{equation}\label{nonz2}
x^2+\g_1x+\g_2\g_3 \hbox{ irreducible over } \F_q \Rightarrow \b=\a.
\end{equation}
We prove the contrapositive. 
Assume $\b=0$, and let $v_1,w_1$ be as in the hypothesis of the lemma. Write 
\[
v_1 = \l_1e_2+\l_3f_3+w,\;w_1 = \mu_1e_2+\mu_3f_3+w',
\]
where $w,w' \in \la e_2,e_3,f_1,f_2\ra$. As $\b=0$, 
observe that $u=W(3)$ and we compute that 
\[
\g_1 = \l_1\mu_3+\l_3\mu_1,\;\g_2=\l_1\l_3,\;\g_3=\mu_1\mu_3.
\]
This gives $(\l_3\mu_1)^2+\g_1(\l_3\mu_1)+\g_2\g_3=0$, and hence 
the quadratic $x^2+\g_1x+\g_2\g_3$ is reducible 
over $\F_q$, proving (\ref{nonz2}). 

Now we complete the proof of the lemma by showing
\begin{equation}\label{nonz3}
x^2+\g_1x+\g_2\g_3 \hbox{ reducible over } \F_q \Rightarrow \b=0.
\end{equation}
Suppose that $x^2+\g_1x+\g_2\g_3$ is reducible. Then there exists $\l \in \F_q$ such that $(\l\g_3)^2+\g_1(\l\g_3)+\g_2\g_3=0$. Hence 
\[
(v_1+\l w_1,\; (v_1+\l w_1)(u-1)) = \l^2\g_3+\l\g_1+\g_2 = 0,
\]
and therefore $\b=0$ by (\ref{nonz}). 
\end{proof}

\begin{lem}\label{conpa}
Let $G = \Sp(V)$ or $\Or(V)$. If $u$ and $v$ are unipotent elements of $G$ 
as in Table $\ref{Asp-conj}$, then $u$ and $v$ are conjugate by an element of $\Sp(V)$ or $\O(V)$, respectively.
\end{lem}

\begin{table}[ht]
\caption{Conjugate unipotent elements}\label{Asp-conj}
\[
\begin{array}{|l|l|}
\hline
u & v  \\
\hline
V(2k) \oplus V(2k) & V_\a(2k) \oplus V_\a(2k)  \\
V(2k+2) \oplus V(2k) & V_\a(2k+2) \oplus V_\a(2k)  \\
V_\a(2k+2) \oplus V(2k) & V(2k+2) \oplus V_\a(2k)  \\
V(2k)^{[3]} & W(2k) \oplus V(2k) \\
V(2k)^{[2]} \oplus V_\a(2k) & W(2k) \oplus V_\a(2k)   \\
W(2l+1) \oplus W(2l+1) & W_\a(2l+1) \oplus W_\a(2l+1)  \\
W(2l+1) \oplus V(2l+2) & W_\a(2l+1) \oplus V_\a(2l+2)   \\
W(2l+1) \oplus V_\a(2l+2) & W_\a(2l+1) \oplus V(2l+2) \\
W(2l+1) \oplus V(2l) & W_\a(2l+1) \oplus V_\a(2l)   \\
W(2l+1) \oplus V_\a(2l) & W_\a(2l+1) \oplus V(2l)   \\
\hline
\end{array}
\]
\end{table}

In the table $k,l\ge 1$, and, as before, our convention is that for $G = \Sp(V)$, the elements $W_\b(1)$ and $V_\b(2)$ are just $W(1)$ and $V(2)$, whatever the value of $\b$.

\begin{proof} For $G =\Sp(V)$, this is \cite[Lemma 3.4]{GLOB}. For $G = \Or(V)$ it follows from the same result, combined with \cite[Thm.\ 3.1]{GLOB}. 
\end{proof}

\subsubsection{The algorithm} 
We now describe our algorithm for conjugacy testing in 
$G = \Sp(V)$ or $\Or(V)$. Let $g \in G$ be unipotent. 
We need to compute a representative $u$ of the form (\ref{canon}) 
(satisfying conditions (i)--(v) listed after (\ref{canon})) such 
that $g$ is conjugate to $u$. 
In the orthogonal case, the main step is to handle conjugacy in $\Or(V)$; 
we address conjugacy in $\O(V)$ below.

\vspace{4mm} \no {\it Step $1$.} 
Let $r$ be the largest Jordan block size for $g$. 
From the rational canonical form \cite[\S 6.7]{PFGG} of $g$, 
we obtain a homocyclic component, say $U$, of the restriction $V\downarrow g$ corresponding 
to the block size $r$. (A {\it homocyclic component}
\index{homocyclic component} is the 
sum of all cyclic summands isomorphic to a given one.)
Then $U$ contains all the blocks of size $r$ for $g$. 
Moreover, $U$ is non-degenerate: suppose for a contradiction that $U\cap U^\perp \ne 0$. Since $U\cap U^\perp$  is $g$-invariant, there exists $0 \ne w \in U\cap U^\perp$ such that $w(g-1)=0$. Then $w = x(g-1)^{r-1}$ for some $x \in U$. Hence for every $v \in V$ such that $v(g-1)^{r-1}=0$, 
\[
(w,v) = (x(g-1)^{r-1},v) = (x,\,v(g^{-1}-1)^{r-1}) = 0.
\]
It follows that $w \in V^\perp$, a contradiction.

\vspace{4mm} \no {\it Step $2$.} Suppose $r$ is odd, so 
$U \downarrow g = \bigoplus_i W_{\b_i}(r)$, where each $\b_i \in \{0,\a\}$. We determine the $\b_i$ as follows. 

Compute $v,w \in U$ such that $(v,\,w(g-1)^{r-1}) \ne 0$. 
Now $U_0 = \la vg^i,\,wg^i : 0\le i\le r-1\ra$ is  non-degenerate and $U_0\downarrow g = W_\b(r)$ for some $\b \in \{0,\a\}$.  If $r=1$ and $G = \Or(V)$, then $\b = 0$ or $\a$ according as the sign of $\Or(U_0)$ is $+$ or $-$, respectively; and if $r=1$ and $G=\Sp(V)$, then $\b=0$. Now assume $r\ge 3$. Let $F = U_0(g-1)^{(r-3)/2}$. By Lemma \ref{reduc}, $g$ acts on $F/{\rm Rad}(F)$ as $W_\b(3)$. We can compute $\b$ using Lemma \ref{w3test}.

Next we work with $U\cap U_0^\perp$ and repeat this process. 
The outcome is a sequence of values $\b_i$ such that 
\begin{eqnarray}\label{rodd}
U \downarrow g = \bigoplus_{i} W_{\b_i}(r),
\end{eqnarray}
where, as always, for $G = \Sp(V)$ the element $W_\b(1)$ is just $W(1)$, 
whatever the value of $\b$.

\vspace{4mm} \no {\it Step $3$.} Suppose $r$ is even, say $r=2k$, so 
\[
U\downarrow g = W(2k)^{[a]}  \oplus   \bigoplus_{i=1}^b V_{\g_i}(2k),
\]
where $b\le 2$. We first use Lemma \ref{j2test} to determine whether $b=0$ or not, as follows. Let $F = U(g-1)^{k-1}$. By Lemma \ref{reduc}, $g$ acts on $F/{\rm Rad}(F)$ as $W(2)^{[a]} \oplus  \bigoplus_{i=1}^b V_{\g_i}(2)$. We can use Lemma \ref{j2test} to determine whether $b=0$.

If $b=0$, then $U\downarrow g = W(2k)^{[a]}$, and we proceed to Step 4. 
If $b>0$, then we determine the $\g_i$ as follows. We find $v \in U$ 
such that $(v,\,v(g-1)^{2k-1})\ne 0$. Now 
$U_0 = \la vg^i : 0\le i\le 2k-1\ra$ is  non-degenerate, 
and $U_0\downarrow g = V_\g(2k)$ for some $\g \in \{0,\a\}$. If $k=1$ and $G = \Or(V)$, then $\b = 0$ or $\a$ according as the sign of $\Or(U_0)$ is $+$ or $-$, respectively. Now assume $k\ge 2$. Let $F = U_0(g-1)^{k-2}$. By Lemma \ref{reduc}, $g$ acts on $F/{\rm Rad}(F)$ as $V_\g(4)$. We can compute $\g$ using Lemma \ref{j4test}.

Next we work with $U\cap U_0^\perp$ and repeat this process. 
The outcome is a sequence of values $\gamma_1, \ldots, \gamma_b$ such that 
\begin{equation}\label{larg}
U\downarrow g = W(2k)^{[a]}  \oplus   \bigoplus_{i=1}^b V_{\g_i}(2k),
\end{equation}
where, as always, for $G = \Sp(V)$ the element $V_\g(2)$ is just $V(2)$, whatever the value of $\g$.

\vspace{4mm} \no {\it Step $4$.}  Having determined (\ref{rodd}) or (\ref{larg}) for the largest homocyclic component, we now work in $U^\perp$, and repeat Steps 1--3 for the next largest component. Repeating this for all the components, we end up with values of $\b_i,\g_i \in \{0,\a\}$ such that $g$ is conjugate to 
\begin{equation}\label{approx}
\bigoplus_{r\;\text{odd}}\bigoplus_i W_{\b_i}(r)  \oplus  \bigoplus_{r\;\text{even}} \left(W(r)^{[a_r]}  \oplus  \bigoplus_i V_{\g_i}(r)\right).
\end{equation}

Our final task is to use Lemma \ref{conpa} to convert this description
to a canonical representative as in (\ref{canon}). We do this in the next step.

\vspace{4mm} \no {\it Step $5$.} Let the Jordan form of $g$ be $\bigoplus_r J_r^{[n_r]}$, and define the following sets:
\[
\begin{array}{l}
T = \{r : n_r>0\}, \\
S = \{r \in T : r \hbox{ odd}\} \cup \{r \in T : r \hbox{ even and }\exists\, 
V_{\g_i}(r) \hbox{ in } (\ref{approx})\}.
\end{array}
\]
For $r \in S$ and $\b \in \{0,\a\}$, define $R_\b(r)$ as follows:
\begin{itemize}
\item[] $r$ odd: $R_\b(r) = W_\b(r) \oplus W(r)^{[\frac{1}{2}n_r-1]}$;
\item[] $r$ even: $R_\b(r) = V_\b(r) \oplus V(r)^{[b]} \oplus W(r)^{[c]}$, where $b\le 1$, $b+2c=n_r-1$.
\end{itemize}
(Thus the Jordan form of $R_\b(r)$ is $J_r^{[n_r]}$ in both cases.)

For $r,s \in S$, write $r \sim s$ if either $|r-s|=1$, or $|r-s|=2$ and 
both $r$ and $s$ are even. Define $\sim$ to be the minimal equivalence relation on $T$ extending this (so each element of $T\setminus S$ 
is its own equivalence class). 

Let $S_0$ be an equivalence class in $T$, and let $n(S_0)$ be the total number of subscripts $\b_i$ or $\g_i$ in (\ref{approx}) that equal $\a$ for block sizes in $S_0$. Define $\a(S_0) \in \{0,\a\}$ as follows:
\[
\begin{array}{ll}
G=\Or(V): & \a(S_0) = \left\{ \begin{array}{l} 0, \hbox{ if }n(S_0) \hbox{ is even} \\
                               \a, \hbox{ if }n(S_0) \hbox{ is odd} \end{array} \right.
\\ 
& \\
G=\Sp(V): & \a(S_0) = \left\{ \begin{array}{l} 0, \hbox{ if }n(S_0) \hbox{ is even, or if }2 \in S_0 \\
                               \a, \hbox{ if }n(S_0) \hbox{ is odd and }2 \not \in S_0. \end{array} \right.
\end{array}
\]
We now identify the canonical representative $R(S_0)$ for 
each equivalence class $S_0$:
\begin{itemize}
\item[(a)] if $S_0 = \{t\}$ with $t \in T\setminus S$, then let $R(S_0) = W(t)^{[n_t/2]}$;
\item[(b)] if $S_0 = \{t\}$ with $t$ odd, then let $R(S_0) = R_{\a(S_0)}(t)$;
\item[(c)] otherwise $S_0 \subseteq S$ and $\exists$ even $t \in S_0$: let $r_0 = \hbox{min}\{r \in S_0: r \hbox{ even}\}$ and set
\[
R(S_0) = R_{\a(S_0)}(r_0)  \oplus  \bigoplus_{r \in S_0, r\ne r_0} R_0(r).
\]
\end{itemize}

\begin{thm}\label{canrep} The canonical representative for the element $g$ as in $(\ref{approx})$ is
\[
\bigoplus_{S_0} R(S_0),
\]
where the sum is over all the equivalence classes $S_0$. 
\end{thm}

\begin{proof}
By definition of the equivalence relation, $\bigoplus_{S_0} R(S_0)$ 
satisfies conditions (i)--(v) listed after (\ref{canon}). By Lemma \ref{conpa}, the sum of the blocks in (\ref{approx}) with block sizes in $S_0$ is conjugate to $R(S_0)$. Hence $g$ is conjugate to $\bigoplus_{S_0} R(S_0)$, as required.
\end{proof}

To complete this section, we address 
conjugacy in $\O(V)$ versus $\Or(V)$. As mentioned in 
Section \ref{badsec}, the only unipotent classes in $\Or(V)$ that 
split in $\O(V)$ are those of the 
form $\bigoplus W(m_i)^{[a_i]}$ with all $m_i$ even; in such a case, a 
second class representative is obtained by replacing one 
summand $W(m_i)$ by $W(m_i)^t$ where $t$ is a reflection. 
Denote these representatives by $u$ and $u^t$ respectively. 
Given $g \in \O(V)$ that is $\Or(V)$-conjugate to $u$, we determine 
whether $g$ is $\O(V)$-conjugate to $u$ or $u^t$, as follows. 
Using the results of Section \ref{badconj}, 
compute $y \in \Or(V)$ such that $g^y = u$. 
Then $g$ is $\O(V)$-conjugate to $u$ if and only if $y \in \O(V)$.

\section{Constructing a conjugating element}\label{badconj}
In this section we solve the conjugation problem for classical groups 
$G$ in bad characteristic: given unipotent $g \in G$ that is 
conjugate to a class representative $u$, find $y \in G$ such that 
$g^y = u$. Recall that the class representatives are given in (\ref{canon}). 
As in Section \ref{goodconjel}, we compute the conjugating 
element $y$ ``block-by-block", so the main task is to solve the problem 
when $u$ is a single block $V_{\b}(2k) $ or $W_\b(m)$. 
We address the symplectic and 
orthogonal groups separately, although both methods are similar. 

\subsection{Symplectic groups}  \label{spbadconjel}
Let $G = \Sp(V)$ in characteristic 2.  The problem for each
of the four types of blocks is challenging.

\subsubsection{Case 1:  $u = V(2k)$}
Let $u = V(2k) \in G = \Sp(V) \cong \Sp_{2k}(q)$. 
For odd characteristic, in Section \ref{spconjel} we used 
the nilpotent element 
of the Lie algebra $sp(V)$ 
corresponding to $u$ (via the Cayley map) 
to express the basis $v_{-(2k-1)},\ldots, v_{2k-1}$. 
In characteristic 2, there is no Cayley map, but we mimic this method. 
Define the following nilpotent elements of ${\rm End}(V)$:
\[
f = 1+u,\;\;e=f+f^2+\cdots +f^{2k-1}.
\]
The basis $v_{-(2k-1)},\ldots, v_{2k-1}$ with respect to which $V(2k)$ is defined in Section \ref{badsec} can be expressed as follows, where $v = v_{-(2k-1)}$:
\begin{equation}\label{basbad1}
v,\,ve,\,ve^2,\ldots,ve^{k-1},\,ve^{k-1}f,\,ve^{k-1}f^2,\ldots,ve^{k-1}f^k.
\end{equation}
With respect to this basis, the symplectic form defining $G$ has matrix
\begin{equation}\label{formm}
\begin{pmatrix} &&&1 \\&&1& \\  &\iddots && \\ 1&&& \end{pmatrix}.
\end{equation}

Now suppose $g \in G$ is conjugate to $u$. 
We aim to compute $y \in G$ such that $g^y = u$. 
Let 
\[
f_0 = 1+g,\;\;e_0=f_0+f_0^2+\cdots +f_0^{2k-1}.
\]
We seek $w\in V$ with the property that 
\begin{equation}\label{basbad2}
w,\,we_0,\,we_0^2,\ldots,we_0^{k-1},\,we_0^{k-1}f_0,\,we_0^{k-1}f_0^2,\ldots,we_0^{k-1}f_0^k
\end{equation}
is a basis of $V$, with respect to which the symplectic form also has matrix (\ref{formm}). Then the element $y$ sending each vector in (\ref{basbad1}) to the corresponding vector in (\ref{basbad2}) will lie in $G$ and conjugate $g$ to $u$, as required. 

To compute such a vector $w$, we adopt a similar strategy to that in Section \ref{spconjel}. To simplify notation, denote the sequence of operators 
\[
1,\,e_0,\,e_0^2,\ldots,e_0^{k-1},\,e_0^{k-1}f_0,\,e_0^{k-1}f_0^2,\ldots,e_0^{k-1}f_0^k
\]
by the symbols 
\[
m_0,\,m_1,\ldots , m_{2k-1}.
\]
Then the vectors in the sequence (\ref{basbad2}) are $wm_0,\ldots,wm_{2k-1}$, and the requirement that the form has matrix (\ref{formm}) can be expressed as
\begin{equation}\label{req}
(wm_i,\,wm_j) = \left \{ 
\begin{array}{l} 
1, \hbox{ if }i+j=2k-1 \\ 
0, \hbox{ otherwise}.
\end{array}
\right.
\end{equation}

Choose $z \in V\setminus V(g-1)$. We aim to solve the following 
equation for $w \in V$ and $a_i \in \F_q$:
\[
z = w\left(a_0+ \sum_{i\;\text{odd},\, 1\le i \le 2k-3} a_im_i \right),
\]
such that $w$ satisfies (\ref{req}). To do this, we compute the values
\[
\a_i = (z,\,zm_i)\;\;(i\hbox{ odd}, 1\le i\le 2k-1).
\]
The requirement (\ref{req}) gives $k$ quadratic equations in the $a_i$. The 
two simplest equations are
\[
a_0^2 = \a_{2k-1},\;a_1^2+a_0a_1=\a_{2k-3}.
\]
These arise from evaluating $\a_{2k-1}$ and $\a_{2k-3}$ imposing (\ref{req}), as follows:
\[
\begin{array}{ll}
\a_{2k-1} &= (z,\,ze_0^{k-1}f_0^k) \\
               & = (a_0w+a_1we_0+\cdots, a_0we_0^{k-1}f_0^k+\cdots ) \\
               &= a_0^2\\
\\
\a_{2k-3} & = (z,\,ze_0^{k-1}f_0^{k-2}) \\
               & = (a_0w+a_1we_0+\cdots, \, a_0we_0^{k-1}f_0^{k-2}+
a_1we_0^{k-1}(f_0^{k-1}+f_0^k)+\cdots ) \\
              & = a_0a_1+a_1^2.
\end{array}
\]
Given a solution for $a_0,a_1$ to these equations, we can successively solve for $a_3$, $a_5$, $\ldots$, $a_{2k-3}$  the ensuing equations with left hand sides 
$\a_{2k-5},\,a_{2k-7},\ldots,\a_1$. For example, when $k=4$, the equations are
\[
\begin{array}{l}
\a_7 = a_0^2, \\
\a_5 = a_0a_1+a_1^2, \\
\a_3 = a_3a_0+a_1^2+a_0a_1,\\
\a_1 = a_5a_0+a_3^2+a_1a_3+a_0a_3.
\end{array}
\]
Having solved for the $a_i$, we let 
\[
M = a_0I+ \sum_{i\;\text{odd},\,1\le i \le 2k-3} a_im_i, 
\]
and set $w = zM^{-1}$. Now $w$ is a vector satisfying (\ref{req}), as required. 

\subsubsection{Case 2:  $u = V_\a(2k)$}
Let $u = V_\a(2k) \in G = \Sp(V) \cong \Sp_{2k}(q)$. This is trickier than the $V(2k)$ case because we cannot list the basis $v_{-(2k-1)},\ldots,v_{2k-1}$ in terms of two nilpotent elements, as in (\ref{basbad1}). We adjust the procedure 
for $V(2k)$ as follows. As before, define
\[
f = 1+u,\;\;e=f+f^2+\cdots +f^{2k-1}.
\]
Writing $v = v_{-(2k-1)}$, the first part of the basis $v_{-(2k-1)},\ldots,v_{-3}$ is $v,ve,\ldots,ve^{k-2}$. 
We compute an endomorphism $h'$, a polynomial in $f$, such that 
\[
v_{-1} = v_{-3}(e+h').
\]
Thus $v_{-1} = ve^{k-2}(e+h')$. Set $h = e^{k-2}h'$. 
Then the basis $v_{-(2k-1)},\ldots,v_{2k-1}$ can be expressed as follows:
\begin{equation}\label{basbadal1}
v,\,ve,\ldots,ve^{k-2},\,v(e^{k-1}+h),\,v(e^{k-1}+h)f,\ldots,v(e^{k-1}+h)f^k.
\end{equation}
Note that $h = p(f)$ is a polynomial in $f$ of degree $2k-1$. 
Here are the polynomials $p(f)$ for the first few values of $k$:
\[
\begin{array}{l|l}
\hline
k & p(f) \\
\hline
2 & \a f^3 \\
3 & \a f^4+\a f^5 \\
4 & \a f^5 + \a^2 f^7 \\
5 & \a f^6 + \a f^7 + (\a^2+\a)f^8 + \a f^9 \\
6 & \a f^7 + (\a^2+\a)f^9+\a^2f^{10}+ \a^3 f^{11} \\
\hline
\end{array}
\]

The rest of the algorithm is similar to the previous case. 
Suppose $g \in G$ is conjugate to $u$. 
We aim to compute $y \in G$ such that $g^y = u$. 
Let 
\[
f_0 = 1+g,\;\;e_0=f_0+f_0^2+\cdots +f_0^{2k-1},\;\;h_0 = p(f_0),
\]
where $p(x)$ is the polynomial such that $h=p(f)$. 
We seek $w\in V$ with the property that 
\begin{equation}\label{basbadal2}
w,\,we_0,\ldots,we_0^{k-2},\,w(e_0^{k-1}+h_0),\,w(e_0^{k-1}+h_0)f_0,\ldots,w(e_0^{k-1}+h_0)f_0^k
\end{equation}
is a basis of $V$, with respect to which the 
symplectic form has matrix (\ref{formm}). 
Then the element $y$ sending each vector in (\ref{basbadal1}) to 
the corresponding vector in (\ref{basbadal2}) will lie in $G$ and 
conjugate $g$ to $u$, as required. 

Denote the sequence of operators $1,e_0,\ldots ,(e_0^{k-1}+h_0)f_0^k$ in (\ref{basbadal2}) by the symbols 
$m_0$, $\ldots$, $m_{2k-1}$, so that 
the requirement that the form has matrix (\ref{formm}) can again be expressed as (\ref{req}).
Choose $z \in V\setminus V(g-1)$. We aim to solve the following equation for $w \in V$ and $a_i \in \F_q$:
\begin{equation}\label{zwa0}
z = w\left(a_0+ \sum_{i\;\text{odd},\, 1\le i \le 2k-3} a_im_i \right),
\end{equation}
such that $w$ satisfies (\ref{req}). To do this, we compute the values
\[
\a_i = (z,\,zm_i)\;\;(i\hbox{ odd},\, 1\le i\le 2k-1).
\]
The requirement (\ref{req}) gives $k$ quadratic equations in the $a_i$. 
As before, these can be solved for the $a_i$. We let 
\[
M = a_0I+ \sum_{i\;\text{odd},\,1\le i \le 2k-3} a_im_i, 
\]
and set $w = zM^{-1}$. Now $w$ is a vector satisfying (\ref{req}), as required. 

\subsubsection{Case 3:  $u = W(m)$}
Let $u = W(m) \in G = \Sp(V) \cong \Sp_{2m}(q)$, defined as in Section \ref{badsec} with respect to a basis $w_{-(m-1)},x_{-(m-1)},\ldots, w_{m-1},x_{m-1}$. In a 
similar manner to Case 1, we express this basis in terms of the 
vectors $w=w_{-(m-1)}$, $x = x_{-(m-1)}$ and the following nilpotent operators:
\[
e=1+u,\;\;f=e+e^2+\cdots+e^{m-1}.
\]
The basis $w_{-(m-1)},x_{-(m-1)},\ldots, w_{m-1},x_{m-1}$ is
\begin{equation}\label{basbad3}
w,\,we,\,\ldots,we^{m-1},\,x,\,xf,\ldots ,xf^{m-1},
\end{equation}
and the symplectic form takes values
\begin{equation}\label{syval}
(we^i,\,xf^j) = \left \{ 
\begin{array}{l} 
1, \hbox{ if }i+j=m-1 \\ 
0, \hbox{ otherwise}
\end{array}
\right.
\end{equation}
and $(we^i,\,we^j) = (xf^i,\,xf^j) = 0$ for all $i,j$. 

Now suppose $g \in G$ is conjugate to $u$. 
We aim to compute $y \in G$ such that $g^y = u$. Let 
\[
e_0 = 1+g,\;\;f_0=e_0+e_0^2+\cdots +e_0^{m-1}.
\]
We seek $w',x'\in V$ with the property that 
\begin{equation}\label{basbad4}
w',\,w'e_0,\,\ldots,w'e_0^{m-1},\,x',\,x'f_0,\ldots ,x'f_0^{m-1}
\end{equation}
is a basis of $V$ with symplectic form values
\begin{equation}\label{valbad}
(w'e_0^i,\,x'f_0^j) = \left \{ 
\begin{array}{l} 
1, \hbox{ if }i+j=m-1 \\ 
0, \hbox{ otherwise}
\end{array}
\right.
\end{equation}
and $(w'e_0^i,\,w'e_0^j) = (x'f_0^i,\,x'f_0^j) = 0$ for all $i,j$. 
Then the element $y$ sending each vector in (\ref{basbad3}) to the 
corresponding vector in (\ref{basbad4}) will lie in $G$ and 
conjugate $g$ to $u$, as required. 

Here is our algorithm to find such vectors $w',x'$. 
First, find $z,t \in V$ such that $\la Vf_0,\,z,\,t\ra = V$ 
and $(z,\,tf_0^{m-1}) \ne 0$. Now compute 
\begin{equation}\label{abcbad}
\begin{array}{llll}
\a_i & = & (z,\,ze_0^i) & (1\le i\le \lfloor (m-1)/2\rfloor), \\
\b_i & = & (t,\,tf_0^i) & (1\le i\le \lfloor (m-1)/2\rfloor), \\
\g_i & = & (z,\,tf_0^j) & (0\le j\le m-1). 
\end{array}
\end{equation}
We consider separately the cases where $m$ is even  and odd.

\vspace*{4mm} \no {\bf (3a)}. Let $m=2k$.
We aim to solve the following equations for $w',x' \in V$ and $b_i,c_i,d_i \in \F_q$:
\begin{equation}\label{wxeqbad}
\begin{array}{l}
z = w' + x'\sum_{i=0}^{k-2}b_if_0^i, \\
t = w'\sum_{i=1}^{k-1} c_{2i-1}e_0^{2i-1} + x'\sum_{i=0}^{2k-1}d_if_0^i,
\end{array}
\end{equation}
such that $w',x'$ satisfy (\ref{valbad}). The values $\a_i,\b_i,\g_i$, combined with (\ref{valbad}), give $4k-2$ quadratic equations in $b_i,c_i,d_i$, and it turns out that these have a unique solution for $b_i,c_i,d_i$. For example, when $k=3$ these equations are:
\[
\begin{array}{l}
\a_2 = b_1,\\
\a_1=b_0+b_1,\\
\b_2=c_1d_0,\\
\b_1 = c_1d_0+c_3d_0+c_1d_1+c_1d_2,\\
\g_5=d_0,\\
\g_4= d_1+b_0c_1,\\
\g_3 = d_2+b_0c_1+d_1c_1,\\
\g_2 = d_3+b_0c_1+b_0c_3,\\
\g_1 = d_4+b_0c_1+b_0c_3+b_1c_1+b_1c_3,\\
\g_0 = d_5.
\end{array}
\]
Note that (in general) $d_0 = \g_{2k-1}$, which is nonzero by the 
choice of $z$ and $t$. For the case $k=3$ we can solve the above equations 
successively for $d_0$, $b_1$, $b_0$, $c_1$, $d_1$, $d_2$, $c_3$, $d_3$, $d_4$ 
and finally $d_5$. A similar procedure applies in the general case.

For this solution for $b_i,c_i,d_i$, define
\[
B = \sum_{i=0}^{k-2}b_if_0^i,\;C = \sum_{i=1}^{k-1} c_{2i-1}e_0^{2i-1},\; D=\sum_{i=0}^{2k-1}d_if_0^i.
\]
By (\ref{wxeqbad}), $z=w'+x'B$ and $t=w'C+x'D$. Hence, letting $M = BC+D$, and noting that $D$ is invertible since $d_0\ne 0$, we deduce that 
\[
x' = (zC+t)M^{-1},\;\;w' = z+x'B.
\]
These vectors $w',x'$ satisfy (\ref{valbad}), as required.

\vspace*{4mm} \no {\bf (3b)}. Let $m=2k+1$.
We aim to solve the following slightly different equations for 
$w',x' \in V$ and $b_i,c_i,d_i \in \F_q$:
\begin{equation}\label{wxeqbad2}
\begin{array}{l}
z = w' + x'\sum_{i=0}^{k-1}b_{2i}f_0^{2i}, \\
t = w'\sum_{i=0}^{k-1} c_{2i}e_0^{2i} + x'\sum_{i=0}^{2k}d_if_0^i,
\end{array}
\end{equation}
such that $w',x'$ satisfy (\ref{valbad}). This time the values $\a_i,\b_i,\g_i$, combined with (\ref{valbad}), give $4k+1$ quadratic equations in the $4k+1$ variables $b_i,c_i,d_i$. However, these are not as straightforward to solve as in the $m$ even case. We solve the equations when $k=2$ -- a case which shows the necessary features of the general case. For $k=2$ the equations are:
\[
\begin{array}{l}
\a_2 = b_0,\\
\a_1=b_0+b_2,\\
\b_2=c_0d_0,\\
\b_1 = c_0d_0+c_2d_0+c_0d_1+c_0d_2,\\
\g_4= d_0+b_0c_0,\\
\g_3 = d_1+b_0c_0,\\
\g_2 = d_2+b_0c_0+b_0c_2+b_2c_0,\\
\g_1 = d_3+b_0c_0+b_0c_2+b_2c_0,\\
\g_0 = d_4+b_2c_2.
\end{array}
\]
Note first that the $1^{st}$, $3^{rd}$ and $5^{th}$ 
equations in this list show that $d_0$ is a root of the quadratic $x^2+\g_4x+\a_2\b_2$. A tedious exercise shows that this quadratic is indeed reducible. 
The choice of $z$ and $t$ implies that $\g_4 \ne 0$, 
so we can take $d_0$ to be a nonzero root. Given this, we can solve the equations uniquely for the rest of the $b_i,c_i,d_i$. 

In the general case, the same features persist for solving the $4k+1$ quadratic equations in the $b_i,c_i,d_i$, which is carried out computationally using a 
Gr\"obner basis algorithm \cite[Chap.\ 1]{Adams}. 
A certain quadratic must be reducible and implies that 
the equations have a unique solution. For this solution, define 
\[
B = \sum_{i=0}^{k-1}b_{2i}f_0^{2i},\;C = \sum_{i=0}^{k-1} c_{2i}e_0^{2i},\; D=\sum_{i=0}^{2k}d_if_0^i.
\]
Letting $M = BC + D$, we deduce that 
\[
x' = (zC+t)M^{-1},\;\;w' = z+x'B.
\]
These vectors $w',x'$ satisfy (\ref{valbad}), as required.

\subsubsection{Case 4:  $u = W_\a(2l+1)$}
Let $u = W_\a(2l+1) \in G = \Sp(V) \cong \Sp_{4l+2}(q)$, defined 
as in Section \ref{badsec} with respect to the basis 
$w_{-2l},x_{-2l},\ldots, w_{2l},x_{2l}$. Let $e = 1+u$. 
Note that $e^{2l+1}=0$. It is 
particularly difficult to describe the basis in terms of the 
nilpotent operator $e$ (as we did in (\ref{basbad3}) in the previous case). 
We let $w=w_{-2l}$, $x = x_{-2l}$, and compute 
polynomials $f_i(t)$, $g_i(t)$, $h_i(t)$, $k_i(t) \in \F_q[t]$ for 
$0\le i\le 2l$, each of degree at most $2l$, such that 
\begin{equation}\label{gb}
w_{-2l+2i} = wf_i(e) + xg_i(e),\;\;x_{-2l+2i} = wh_i(e) + xk_i(e)\;\;\hbox{ for }0\le i \le 2l.
\end{equation}
Note that $f_0(t)=k_0(t)=1$, $g_0(t)=h_0(t)=0$. The symplectic form takes values 
\begin{equation}\label{valbadal}
(wf_i(e) + xg_i(e),\, wh_j(e) + xk_j(e)) = \left \{ 
\begin{array}{l} 
1, \hbox{ if }i+j=2l \\ 
0, \hbox{ otherwise}
\end{array}
\right.
\end{equation}
and $(wf_i(e) + xg_i(e),\, wf_j(e) + xg_j(e)) = (wh_i(e) + xk_i(e),\, wh_j(e) + xk_j(e))=0$ for all $i,j$. 

For example, when $l=1$, $u = W_\a(3)$, the basis $w_{-2},w_0,w_2,x_{-2},x_0,x_{2}$ is
\[
w,\;we+\a xe^2,\;we^2,\;x,\;\a we^2+x(e+e^2),\;xe^2,
\]
and so we can identify the polynomials:
\[
\begin{array}{|l|l|l|l|l|} \hline 
i & f_i(t)& g_i(t) & h_i(t) & k_i(t) \\
\hline
0 & 1&0&0&1 \\
1 & t&\a t^2 & \a t^2 & t+t^2 \\
2 & t^2 & 0&0&t^2 \\
\hline
\end{array}
\]
For general $l$, we construct the relevant polynomials $f_i,g_i,h_i,k_i$ by machine.

Our strategy is the usual one. Suppose $g \in G$ is conjugate to $u$. We aim to compute $y \in G$ such that $g^y = u$. We let 
$e_0 = 1+g$ and seek $w',x'\in V$ with the property that 
\begin{equation}\label{gb1}
w'f_i(e_0) + x'g_i(e_0),\;\;w'h_i(e_0) + x'k_i(e_0)\;\;(0\le i \le 2l)
\end{equation}
is a basis of $V$ with symplectic form values 
\begin{equation}\label{valbadal1}
(w'f_i(e_0) + x'g_i(e_0),\, w'h_j(e_0) + x'k_j(e_0)) = \left \{ 
\begin{array}{l} 
1, \hbox{ if }i+j=2l \\ 
0, \hbox{ otherwise}
\end{array}
\right.
\end{equation}
and $$(w'f_i(e_0) + x'g_i(e_0),\, w'f_j(e_0) + x'g_j(e_0)) 
= (w'h_i(e_0) + x'k_i(e_0),\, w'h_j(e_0) + x'k_j(e_0))=0$$ 
for all $i,j$.  Then the element $y$ sending each vector in (\ref{gb}) to the corresponding vector in (\ref{gb1}) will lie in $G$ and conjugate $g$ to $u$, as required. 

Here is our algorithm to find such vectors $w',x'$. 
First, find $z,t \in V$ such that $\la Ve_0,\,z,\,t\ra = V$, and  
compute the scalars
\begin{equation}\label{abcbad1}
\begin{array}{l}
\a_i = (z,\,zf_i(e_0) + tg_i(e_0))\;\;\;(1\le i\le 2l), \\
\b_i = (t,\,zh_i(e_0) + tk_i(e_0))\;\;\;(1\le i\le 2l), \\
\g_i = (z,\,zh_i(e_0) + tk_i(e_0))\;\;\;(0\le i\le 2l). \\
\end{array}
\end{equation}
Now we aim to solve the following equations for 
$w',x' \in V$ and $b_i,c_i,d_i \in \F_q$:
\begin{equation}\label{wxeqbad2-A}
\begin{array}{l}
z = w' + \sum_{i=0}^{l-1}b_i(w'h_{2i}(e_0)+x'k_{2i}(e_0)), \\
t = \sum_{i=0}^{l-1}c_i(w'f_{2i}(e_0)+x'g_{2i}(e_0)) + \sum_{i=0}^{2l}d_i(w'h_{i}(e_0)+x'k_{i}(e_0))
\end{array}
\end{equation}
such that $w',x'$ satisfy (\ref{valbadal1}). The values $\a_i,\b_i,\g_i$, combined with (\ref{wxeqbad2-A}) and  (\ref{valbadal1}), give $6l+1$ quadratic equations in $b_i,c_i,d_i$. Generating these quadratics is more complicated than 
in previous cases, as we do not know the polynomials 
$f_i,g_i,h_i,k_i$ explicitly. 
To generate the quadratic equations, we proceed as follows. 

Let $b_i\,(0\le i\le l-1)$, $c_i\,(0\le i \le l-1)$, and $d_i\,(0\le i\le 2l)$ 
be indeterminates, and let 
$F$ be the field of 
rational functions over $\F_q$ in the $4l + 1$ indeterminates.
Let $W,X,E$ be indeterminates over $F$, and define the ring
\[
R = \frac{F[W,X,E]}{(E^{2l+1})}.
\]
Define the following elements of $R$:
\begin{equation}\label{videf}
V_i = Wf_i(E)+Xg_i(E),\;\;V_{2l+1+i} = Wh_i(E)+Xk_i(E) \;\;\hbox{ for }0\le i \le 2l.
\end{equation}
(Strictly speaking, these should be defined as cosets of the ideal $(E^{2l+1})$, but since all the polynomials $f_i,g_i,h_i,k_i$ have degree less than $2l+1$, there is no ambiguity.) 
For $r = \sum_{i=0}^{4l+1} r_iV_i$, and $s = \sum_{i=0}^{4l+1} s_iV_i \in R$, where each $r_i,s_i \in F$, set
\begin{equation}\label{ipdef}
IP(r,s) = \sum_{i=0}^{4l+1} r_is_{4l+1-i} \in F.
\end{equation}
Note that $IP(r,s)$ is a quadratic expression in the coefficients $r_i,s_i$. 
If we specialize down to 
$W\mapsto w$, $X \mapsto x$, $E \mapsto e$, 
then  $V_0,\ldots, V_{4l+1}$ become the original basis vectors $w_{-2l},\ldots, x_{2l}$, and for $r_i,s_i \in \F_q$, we see that $r,s \in V$ 
and $IP(r,s)$ is the value of the symplectic form $(r,s)$.

Now define the following elements of $R$:
\[
\begin{array}{l}
Z = W +\sum_{i=0}^{l-1} b_i\left(Wh_{2i}(E)+Xk_{2i}(E)\right),\\
T = \sum_{i=0}^{l-1} c_i\left(Wf_{2i}(E)+Xg_{2i}(E)\right) + \sum_{i=0}^{2l} d_i\left(Wh_{i}(E)+Xk_{i}(E)\right).
\end{array}
\]
Let $A_i,B_i,C_i$ be the following elements of $F$:
\begin{equation}\label{abceqs}
\begin{array}{l}
A_i(b,c,d) = IP(Z,\; Zf_i(E)+Tg_i(E))\;\;\;(1\le i\le 2l), \\
B_i(b,c,d) = IP(T,\; Zh_i(E)+Tk_i(E))\;\;\;(1\le i\le 2l), \\
C_i(b,c,d) = IP(Z,\; Zh_i(E)+Tk_i(E))\;\;\;(0\le i\le 2l).
\end{array}
\end{equation}
As noted above, each of $A_i,B_i,C_i \in F$ is a quadratic expression in the variables $b_i,c_i,d_i$. The quadratic equations in $b_i,c_i,d_i$ arising from the equations (\ref{wxeqbad2-A}) combined with (\ref{valbadal1}) and (\ref{abcbad1}) are
\begin{equation}\label{quadbad}
\begin{array}{l}
A_i(b,c,d) = \a_i\;\;(1\le i\le 2l), \\
B_i(b,c,d) = \b_i\;\;(1\le i\le 2l), \\
C_i(b,c,d) = \g_i\;\;(0\le i\le 2l).
\end{array}
\end{equation}

As an illustration, we compute the expression $A_1(b,c,d)$ in the case where $l=1$ (so $u = W_\a(3)$). Here 
\[
Z = W+b_0X,\;\;T = c_0W+d_0X+d_1(\a WE^2+X(E+E^2))+d_2XE^2.
\]
To compute $A_1 := IP(Z,\; Zf_1(E)+Tg_1(E))$, we need to express 
$Z$ and $Zf_1(E)+Tg_1(E)$ in terms of  $V_0,\ldots, V_5$: these expressions are
\[
\begin{array}{rl}
Z & = W+b_0X = V_0+b_0V_3,\\
\\
Zf_1(E)+Tg_1(E) & = WE+b_0XE+\a c_0WE^2+ \a d_0XE^2 \\
                         & = V_1+\a (b_0+c_0)V_2 +b_0V_4+(b_0+\a d_0+\a )V_5.
\end{array}
\]
Hence 
\[
A_1 = IP(Z,\; Zf_1(E)+Tg_1(E)) = \a b_0(b_0+c_0)+b_0+\a d_0+\a.
\]

Returning to the general case, we use a Gr\"obner basis algorithm to find 
a solution for $b_i,c_i,d_i \in \F_q$ of the quadratic equations 
(\ref{quadbad}). It may fail to produce a solution -- in which 
case we re-choose the vectors $z$ and $t$ and repeat the above steps. 
(In practice, 
we obtain a solution after a small number of attempts.)

Given a solution $b_i,c_i,d_i$ to (\ref{quadbad}), define
\[
\begin{array}{l}
f(e_0) = 1 +\sum_{i=0}^{l-1} b_ih_{2i}(e_0),\\
g(e_0) = \sum_{i=0}^{l-1} b_ik_{2i}(e_0), \\
r(e_0) = \sum_{i=0}^{l-1} c_if_{2i}(e_0) + \sum_{i=0}^{2l} d_ih_{i}(e_0), \\
s(e_0) = \sum_{i=0}^{l-1} c_ig_{2i}(e_0) + \sum_{i=0}^{2l} d_ik_{i}(e_0).
\end{array}
\]
Now, by (\ref{wxeqbad2}), 
\begin{equation}\label{zweq}
z = w'f(e_0)+x'g(e_0),\;\;t = w'r(e_0)+x's(e_0).
\end{equation}
Note that $f(e_0)$ is invertible. Let $M = g(e_0)r(e_0)+f(e_0)s(e_0)$. 
If $M$ is invertible, or equivalently, $b_0c_0\ne d_0$, then
the solution of (\ref{zweq}) is
\begin{equation}\label{exa}
x' = (zr(e_0)+tf(e_0))M^{-1}, \;w' = (z+x'g(e_0))f(e_0)^{-1}.
\end{equation}
(In practice, we observe that the inequality $b_0c_0\ne d_0$ always 
holds, but we could not prove this.) 
These vectors $w',x'$ satisfy (\ref{gb1}) and (\ref{valbadal1}), as required. 

\subsubsection{General case}
Suppose $g \in G = \Sp(V)$ is conjugate to a unipotent class representative
\begin{equation}\label{canongen}
u=\bigoplus_i W(m_i)^{[c_i]}  \oplus  \bigoplus_j V(2k_j)^{[d_j]}  \oplus  \bigoplus_r W_\a(m_r')  \oplus  \bigoplus_s V_\a(2k_s'),
\end{equation}
satisfying conditions (i)--(v) listed after (\ref{canon}). 
We aim to compute an element of $G$ that conjugates $g$ to $u$ 
``block-by-block", using the work in Cases 1--4 above. However, the procedure 
is more subtle than in previous cases.

Let $2k$ be an even block size in (\ref{canongen}), and let the contribution to $u$ from blocks of this size be
\begin{equation}\label{2kcont}
W(2k)^{[a]}  \oplus  V_\a(2k)^{[b]} \oplus V(2k)^{[c]},
\end{equation}
where $b\le 1 $ and $b+c \le 2$. By Lemma \ref{conpa}, this element 
of $\Sp(V_{4ak+2bk+2ck})$ is conjugate to 
\begin{equation}\label{2knew}
\left\{\begin{array}{l}
W(2k)^{[a]}, \hbox{ if }b+c=0, \\
V(2k)^{[2a+c]} \oplus V_\a(2k)^{[b]}, \hbox{ if }b+c>0.
\end{array} \right.
\end{equation}
In (\ref{canongen}), we replace each summand (\ref{2kcont}) by a summand 
(\ref{2knew}) to get a conjugate element $u'$. 
We do this because in a space of the form (\ref{2kcont}) with $b+c>0$
it is difficult 
to construct a $W$-block, 
so instead we convert the decomposition of the 
space into $V$-blocks, which can be easily found.

\subsubsection{The algorithm}
We now describe our algorithm that performs the following task: given $g,h \in G$ that are both conjugate to $u'$ (hence also to $u$), find $x \in G$ such that $g^x=h$. 

The first step is to  compute  orthogonal decompositions of $V\downarrow g$ and $V\downarrow h$ into $V$-blocks and $W$-blocks, as follows. 

Assume first that the largest block size in $u'$ is even, say $2k$. 
Suppose there is a block $V(2k)$ or $V_\a(2k)$ present. We find $v \in V$ such that $(v,\,v(g-1)^{2k-1}) \ne 0$, and let $X = \la vg^i : 0\le i\le 2k-1\ra$. Then $X\downarrow g$ is a $V$-block of size $2k$, and we can determine whether it is $V(2k)$ or $V_\a(2k)$ using Lemma \ref{j4test} (as in Step 3 of the algorithm in Section \ref{conjprobbad}). Now we repeat this computation in $X^\perp$. Note that we need to match the $V$-blocks thus constructed with those in $u'$, which may involve rechoosing the vectors $v$ in the above process.
 
Suppose there is no  block $V(2k)$ or $V_\a(2k)$ present in $u'$. 
We find $w,x \in V$ such that $(w,\,x(g-1)^{2k-1}) \ne 0$, and let $X = \la wg^i,xg^i : 0\le i\le 2k-1\ra$. Now $X\downarrow g$ is $W(2k)$, and we repeat with $X^\perp$. 

To conclude this step, assume now that the largest block size is odd, 
say $2l+1$. We find $w,x \in V$ such that $(w,\,x(g-1)^{2l}) \ne 0$, and 
let $X = \la wg^i,xg^i : 0\le i\le 2l\ra$. Then $X\downarrow g$ is $W(2l+1)$ or $W_\a(2l+1)$; using Lemma \ref{w3test}, we determine which 
(as in Step 2 of Section \ref{conjprobbad}). 
Now we repeat this computation in $X^\perp$. As above, we need to  match the $W$-blocks thus constructed with those in $u'$, which may involve rechoosing the vectors $w,x$ in this process.

We carry out the same procedure for the element $h$.

Now we have orthogonal decompositions
\[
V\downarrow g = \bigoplus_{i=1}^r X_i,\;\;V\downarrow h = \bigoplus_{i=1}^r Y_i,
\]
where, for each $i$, the actions $g^{X_i}$ and $h^{Y_i}$ are conjugate to a $V$-block or a $W$-block, and are conjugate to each other. Moreover, both decompositions match that of the element $u'$. 

Let $g_i = g^{X_i}$ and  $h_i = h^{Y_i}$. 
Hence $g = {\rm diag}(g_1,\ldots,g_r)$ and $h = {\rm diag}(h_1,\ldots,h_r)$
with respect to the union of standard bases of the spaces $X_i,Y_i$.
Let $\dim X_i = 2n_i$. Using the work for individual blocks in 
Cases 1--4 above, 
we compute $x_i \in \Sp_{2n_i}(q)$ such that $g_i^{x_i} = h_i$ 
for $1\le i\le r$. Now $x =  {\rm diag}(x_1,\ldots,x_r)$ is an 
element of $G$ that conjugates $g$ to $h$. 

\subsection{Orthogonal groups}  \label{obadconjel}
Let $G = \Or(V)$ in characteristic 2. As for symplectic groups, 
the major part of the work is to find conjugating elements for individual 
$V$- and $W$-blocks. 
We first consider together the blocks $V(2k)$ and $V_\a(2k)$.

\subsubsection{Case 1:  $u = V_\b(2k)$}
Let $u = V_\b(2k) \in G = \Or(V)$, where $\b \in \{0,\a\}$ and $V = V_{2k}$. 
Recall that $G \cong \Or^+_{2k}(q)$ for $\b=0$, and  
$G \cong \Or^-_{2k}(q)$ for $\b=\a$. The element $u$ is defined as in Section \ref{badsec} with respect to the basis 
\[
v_{-(2k-1)},v_{-(2k-3)},\ldots, v_{2k-1},
\]
with symplectic form given by $(v_i,v_{-i}) =1$ for all $i$, and quadratic form $Q_\b(v_{-1}) = \b$, $Q_\b(v_1) = 1$, all other values $(v_i,v_j)$ and $Q_\b(v_i)$ being 0. 

We proceed as in Case 2 of Section \ref{spbadconjel}, with a few tweaks to 
take account of the quadratic form. Let $f = 1+u$, $e=f+f^2+\cdots +f^{2k-1}$, and compute $h' \in {\rm End}(V_{2k})$, a polynomial in $f$, such that $v_{-1} = v_{-3}(e+h')$. If we set $h = e^{k-2}h'$, then the basis $v_{-(2k-1)},\ldots,v_{2k-1}$ can be expressed as
\begin{equation}\label{basbadal5}
v,\,ve,\ldots,ve^{k-2},\,v(e^{k-1}+h),\,v(e^{k-1}+h)f,\ldots,v(e^{k-1}+h)f^k.
\end{equation}
Note that $h=0$ if $\b=0$; and if $\b=\a$, then $h = p(f)$, where $p(f)$ is the same polynomial as in Case 2 of Section \ref{spbadconjel}.

Suppose $g \in \Or(V)$ is conjugate to $u$. We aim to compute $y \in \Or(V)$ such that $g^y=u$. As in  Section \ref{spbadconjel}, we let $f_0=1+g$, $e_0=f_0+f_0^2+\cdots +f_0^{2k-1}$, $h_0=p(f_0)$, and seek $w \in V$ such that 
\begin{equation}\label{basbadal3}
w,\,we_0,\ldots,we_0^{k-2},\,w(e_0^{k-1}+h_0),\,w(e_0^{k-1}+h_0)f_0,\ldots,w(e_0^{k-1}+h_0)f_0^k
\end{equation}
is a basis of $V$, with respect to which the symplectic form has 
matrix (\ref{formm}), and $$Q_\b(w(e_0^{k-1}+h_0)) = \b, \;\;\; Q_\b(w(e_0^{k-1}+h_0)f_0) = 1,$$ 
the rest of the $Q_\b$-values on the basis being 0.

Denote the sequence of operators $1,e_0,\ldots, (e_0^{k-1}+h_0)f_0^k$ in (\ref{basbadal3}) by the symbols 
$m_0$, $\ldots$, $m_{2k-1}$. Then the requirements on $w$ are that the equations (\ref{req}) are satisfied, together with the equations 
\begin{equation}\label{qbw}
Q_\b(wm_k) = \b, \;Q_\b(wm_{k+1})=1, \hbox{ and } Q(wm_i)=0 \hbox{ for } i \ne k,k+1.
\end{equation}
Choose $z \in V\setminus V(g-1)$. We aim to solve the following equation 
for $w \in V$  and $a_i \in \F_q$:
\begin{equation}\label{zwa1}
z = w\left(a_0+ \sum_{i\;\text{odd},\,1\le i \le 2k-1} a_im_i \right),
\end{equation}
such that $w$ satisfies (\ref{req}) and (\ref{qbw}). There is a difference with the corresponding equation (\ref{zwa0}) in  Section \ref{spbadconjel}: 
the sum on the right hand side includes a term for $i=2k-1$; the variable $a_{2k-1}$ is needed to allow enough freedom to solve the extra equations (\ref{qbw}).

To find a solution to (\ref{zwa1}), we compute the values
\[
\begin{array}{l}
\a_i = (z,\,zm_i)\;\;(i\hbox{ odd}, 1\le i\le 2k-1), \hbox{ and} \\
\g = Q_\b(z).
\end{array}
\]
The requirements (\ref{req}) and (\ref{qbw}) give $k$ quadratic equations in the $a_i$, together with the equation 
\[
\g = \left\{\begin{array}{l} a_{k-1}^2\b+a_0a_{2k-1}, \hbox{ if } k \hbox{ is even}, \\
            a_{k}^2+a_0a_{2k-1}, \hbox{ if } k \hbox{ is odd.}
\end{array}
\right.
\]
As in Section \ref{spbadconjel}, these $k+1$ equations 
can be solved for the $a_i$.  We let 
\[
M = a_0I+ \sum_{i\;\text{odd},\,1\le i \le 2k-1} a_im_i, 
\]
and set $w = zM^{-1}$. Then $w$ is a vector satisfying (\ref{req}) and (\ref{qbw}), as required. 

\subsubsection{\bf Case 2:  $u = W(m)$}
This case is similar to 
Case 3 of Section  \ref{spbadconjel}, with some tweaks to 
take account of the quadratic form. 
Let $u = W(m) \in G = \Or(V) \cong \Or_{2m}^+(q)$, 
defined as in Section \ref{badsec} with respect to a basis 
$w_{-(m-1)},x_{-(m-1)},\ldots, w_{m-1},x_{m-1}$. 
Let $e=1+u$ and $f=e+e^2+\cdots+e^{m-1}$.
The basis $w_{-(m-1)},x_{-(m-1)},\ldots, w_{m-1},x_{m-1}$ is
\[
w,\,we,\,\ldots,we^{m-1},\,x,\,xf,\ldots ,xf^{m-1},
\]
the symplectic form takes values as in (\ref{syval}), 
and the quadratic form $Q$ preserved by $G$ takes 
the value 0 on every basis vector.

Suppose $g \in G$ is conjugate to $u$. 
We aim to compute $y \in G$ such that $g^y = u$. 
We let $e_0 = 1+g$, $f_0=e_0+e_0^2+\cdots +e_0^{m-1}$, 
and seek $w',x'\in V$ with the property that 
\begin{equation}\label{basbad5}
w',\,w'e_0,\,\ldots,w'e_0^{m-1},\,x',\,x'f_0,\ldots ,x'f_0^{m-1},
\end{equation}
is a basis of $V$, with symplectic form values as in (\ref{valbad}), 
and $Q$-values 0 on every basis vector.

Here is our algorithm to find such vectors $w',x'$. First, find $z,t \in V$ such that $\la Vf,\,z,\,t\ra = V$ and $(z,\,tf_0^{m-1}) \ne 0$. Now compute 
\begin{equation}\label{abcbad3}
\begin{array}{llll}
\a_i & = & (z,\,ze_0^i) & (1\le i\le \lfloor (m-1) / 2\rfloor), \\
\b_i & = & (t,\,tf_0^i) & (1\le i\le \lfloor (m-1) / 2\rfloor), \\
\g_i & = & (z,\,tf_0^j) & (0\le j\le m-1), \\
\d_1 & =  & Q(z), & \\
\d_2 & = & Q(t).& 
\end{array}
\end{equation}

\no {\bf (2a).}
Let $m=2k$. We aim to solve the following equations for 
$w',x' \in V$ and $b_i,c_i,d_i \in \F_q$:
\begin{equation}\label{wxeqbad1}
\begin{array}{l}
z = w' + x'\left(\sum_{i=0}^{k-2}b_if_0^i + b_{2k-1}f_0^{2k-1}\right), \\
t = w'\sum_{i=1}^{k} c_{2i-1}e_0^{2i-1} + x'\sum_{i=0}^{2k-1}d_if_0^i,
\end{array}
\end{equation}
such that $w',x'$ give a basis as in (\ref{basbad5}) with symplectic form values as in  (\ref{valbad}) and $Q$-values 0. 
Note the extra terms in these equations in comparison to (\ref{wxeqbad}): these are the terms $b_{2k-1}f_0^{2k-1}$ in the first equation and $c_{2k-1}e_0^{2k-1}$ in the second. (The extra variables $b_{2k-1},c_{2k-1}$ are needed to take care of the $Q$-values.)
The values $\a_i,\b_i,\g_i,\d_1,\d_2$, combined with (\ref{valbad}) and the $Q$-values condition, give $4k$ quadratic equations in the variables 
$b_i,c_i,d_i$, and it turns out that these have a unique solution for $b_i,c_i,d_i$. 
For this solution, define 
\[
B = \sum_{i=0}^{k-2}b_if_0^i + b_{2k-1}f_0^{2k-1},\;C = \sum_{i=1}^{k} c_{2i-1} e_0^{2i-1},\; D=\sum_{i=0}^{2k-1}d_if_0^i.
\]
By (\ref{wxeqbad1}), $z=w'+x'B$ and $t=w'C+x'D$, so letting $M = BC+D$, 
\[
x' = (zC+t)M^{-1},\;\;w' = z+x'B.
\]
These are the required vectors $w',x'$.

\vspace*{2mm} 
\no 
{\bf (2b)}.
Let $m=2k+1$. We aim to solve the following 
slightly different equations for $w',x' \in V$ and $b_i,c_i,d_i \in \F_q$:
\begin{equation}\label{wxeqbad3-first}
\begin{array}{l}
z = w' + x'\sum_{i=0}^{k}b_{2i}f_0^{2i}, \\
t = w'\sum_{i=0}^{k} c_{2i}e_0^{2i} + x'\sum_{i=0}^{2k}d_if_0^i,
\end{array}
\end{equation}
such that $w',x'$ give a basis as in (\ref{basbad5}) with symplectic form values as in  (\ref{valbad}) and $Q$-values 0. (Note the extra terms in comparison to (\ref{wxeqbad2}), namely $b_{2k}f_0^{2k}$ and $c_{2k}e_0^{2k}$.) This time the values $\a_i,\b_i,\g_i,\d_1,\d_2$, combined with (\ref{valbad}) and the $Q$-values condition, give $4k+3$ quadratic equations in 
the variables $b_i,c_i,d_i$. These can be solved using a Gr\"obner 
basis algorithm. For this solution, define 
\[
B = \sum_{i=0}^{k}b_{2i}f_0^{2i},\;C = \sum_{i=0}^{k} c_{2i}e_0^{2i},\; D=\sum_{i=0}^{2k}d_if_0^i.
\]
Letting $M = BC+D$, we deduce that 
\[
x' = (zC+t)M^{-1},\;\;w' = z+x'B.
\]
These are the required vectors $w',x'$.

\subsubsection{\bf Case 3:  $u = W_\a(2l+1)$}
Let $u = W_\a(2l+1) \in G = \Or(V) \cong \Or^-_{4l+2}(q)$, 
defined as in Section \ref{badsec} with respect to the basis 
$w_{-2l},x_{-2l},\ldots, w_{2l},x_{2l}$. Let $e = 1+u$. 
We let $w=w_{-2l}$, $x = x_{-2l}$, and compute polynomials $f_i(t)$, $g_i(t)$, $h_i(t)$, $k_i(t) \in \F_q[t]$ for $0\le i\le 2l$, each of degree at most $2l$, such that 
\begin{equation}\label{gb2}
w_{-2l+2i} = wf_i(e) + xg_i(e),\;\;x_{-2l+2i} = wh_i(e) + xk_i(e)\;\;\hbox{ for }0\le i \le 2l.
\end{equation}
 The symplectic form takes values as in (\ref{valbadal}), and the quadratic form $Q$ preserved by $G$ takes values 0 on the basis, apart from $Q(w_{-2}) = Q(w_0) = Q(x_0) = \a$, which means that 
\begin{equation}\label{qvalbad}
Q(wf_{l-1}(e) + xg_{l-1}(e)) = Q(wf_l(e) + xg_l(e)) = Q(wh_l(e) + xk_l(e)) = \a.
\end{equation}

Our method is as in Case 4 of Section \ref{spbadconjel}. Suppose $g \in G$ is conjugate to $u$. We aim to compute $y \in G$ such that $g^y = u$. We let $e_0 = 1+g$ and seek $w',x'\in V$ with the property that 
\begin{equation}\label{gb3}
w'f_i(e_0) + x'g_i(e_0),\;\;w'h_i(e_0) + x'k_i(e_0)\;\;(0\le i \le 2l)
\end{equation}
is a basis of $V$ with symplectic form values as in (\ref{valbadal1}), and $Q$-values 0 apart from 
\begin{equation}\label{qvalbad3}
Q(w'f_{l-1}(e_0) + x'g_{l-1}(e_0)) = Q(w'f_l(e_0) + x'g_l(e_0)) = Q(w'h_l(e_0) + x'k_l(e_0)) = \a.
\end{equation}

Here is our algorithm to find such vectors $w',x'$. First, find $z,t \in V$ such that $\la Ve_0,\,z,\,t\ra = V$, and  compute the scalars
\begin{equation}\label{abcbad3-second}
\begin{array}{l}
\a_i = (z,\,zf_i(e_0) + tg_i(e_0))\;\;\;(1\le i\le 2l), \\
\b_i = (t,\,zh_i(e_0) + tk_i(e_0))\;\;\;(1\le i\le 2l), \\
\g_i = (z,\,zh_i(e_0) + tk_i(e_0))\;\;\;(0\le i\le 2l), \\
\d_1 = Q(z), \\
\d_2 = Q(t).
\end{array}
\end{equation}
Now we aim to solve the following equations for $w',x' \in V$ and $b_i,c_i,d_i \in \F_q$:
\begin{equation}\label{wxeqbad3}
\begin{array}{l}
z = w' + \sum_{i=0}^{l}b_i(w'h_{2i}(e_0)+x'k_{2i}(e_0)), \\
t = \sum_{i=0}^{l}c_i(w'f_{2i}(e_0)+x'g_{2i}(e_0)) + \sum_{i=0}^{2l}d_i(w'h_{i}(e_0)+x'k_{i}(e_0))
\end{array}
\end{equation}
such that $w',x'$ satisfy (\ref{gb3}), (\ref{valbadal1}) and (\ref{qvalbad3}). (Note the extra terms for $i=l$ in both equations, in comparison to (\ref{wxeqbad2-A}).) The values $\a_i,\b_i,\g_i,\d_1,\d_2$, combined with (\ref{wxeqbad3}),  (\ref{valbadal1}) and (\ref{qvalbad3}), give $6l+3$ quadratic equations in $b_i,c_i,d_i$. 

To generate these equations, we follow the same procedure as in Case 4 of Section  \ref{spbadconjel}, with some small adjustments. We work in the same ring $R = F[W,X,E]/(E^{2l+1})$, where $F$ is the field of 
rational functions over $\F_q$ in indeterminates $b_i\,(0\le i\le l)$, $c_i \,(0\le i\le l)$, 
and $d_i\,(0\le i \le 2l)$. 
Define $V_i \in R$ as in (\ref{videf}), and for $r = \sum_{i=0}^{4l+1} r_iV_i$ and $s = \sum_{i=0}^{4l+1} s_iV_i \in R$, where each $r_i,s_i \in F$, define $IP(r,s)$ as in (\ref{ipdef}); also define 
\[
Quad(r) = \sum_{i=0}^{2l}r_ir_{4l+1-i} + \a (r_{l-1}^2+r_l^2+r_{3l+1}^2).
\]
If we specialize down to $W\mapsto w$, $X \mapsto x$, $E \mapsto e$, 
then  $V_0,\ldots, V_{4l+1}$ become the original basis vectors $w_{-2l},\ldots, x_{2l}$, and, for $r_i \in \F_q$, we see that $r \in V$ and $Quad(r)$ is the 
value of the quadratic form $Q(r)$.

Now define the following elements of $R$:
\[
\begin{array}{l}
Z = W +\sum_{i=0}^{l} b_i\left(Wh_{2i}(E)+Xk_{2i}(E)\right),\\
T = \sum_{i=0}^{l} c_i\left(Wf_{2i}(E)+Xg_{2i}(E)\right) + \sum_{i=0}^{2l} d_i\left(Wh_{i}(E)+Xk_{i}(E)\right).
\end{array}
\]
(As above, there are extra terms for $i=l$ in both elements, in comparison to $Z,T$ as defined in Section  \ref{spbadconjel}.) Let $A_i,B_i,C_i \in F$ be as in (\ref{abceqs}), and define additional elements of $F$ as follows:
\[
Q_1(b,c,d) = Quad(Z),\;Q_2(b,c,d) = Quad(T).
\]
It is easy to give explicit expressions for both.
\[
\begin{array}{lll}
Quad(Z) & = & \left\{\begin{array}{ll} b_l & \hbox{ for odd }l >1 \\
      b_l+\a b_{l/2}^2  & \hbox{ for }l \hbox{ even} \\
      b_1+\a & \hbox{ for }l=1.
\end{array} \right. \\
& & \\
Quad(T) & = & \sum_{i=0}^l c_id_{2l-2i} + \a d_l^2 + \a c_{\lfloor l/2 \rfloor}^2.
\end{array}
\]
The quadratic equations in $b_i,c_i,d_i$ arising from the equations (\ref{wxeqbad3}) combined with (\ref{valbadal1}), (\ref{qvalbad3}) and (\ref{abcbad3-second}) are
\begin{equation}\label{quadbad3}
\begin{array}{l}
A_i(b,c,d) = \a_i\;\;(1\le i\le 2l), \\
B_i(b,c,d) = \b_i\;\;(1\le i\le 2l), \\
C_i(b,c,d) = \g_i\;\;(0\le i\le 2l), \\
Q_1(b,c,d) = \d_1,\\
Q_2(b,c,d) = \d_2.
\end{array}
\end{equation}
We use a Gr\"obner basis algorithm to find solutions for $b_i,c_i,d_i$ of these equations. 
Given such a solution, define
\[
\begin{array}{l}
f(e_0) = 1 +\sum_{i=0}^{l} b_ih_{2i}(e_0),\\
g(e_0) = \sum_{i=0}^{l} b_ik_{2i}(e_0), \\
r(e_0) = \sum_{i=0}^{l} c_if_{2i}(e_0) + \sum_{i=0}^{2l} d_ih_{i}(e_0), \\
s(e_0) = \sum_{i=0}^{l} c_ig_{2i}(e_0) + \sum_{i=0}^{2l} d_ik_{i}(e_0).
\end{array}
\]
By (\ref{wxeqbad3}), 
\[
z = w'f(e_0)+x'g(e_0),\;\;t = w'r(e_0)+x's(e_0),
\]
and we solve for $w',x'$ exactly as in (\ref{exa}). This gives the required vectors $w',x'$, completing this case 
 $u = W_\a(2l+1)$. 

\subsubsection{General case} 
Suppose $g,h \in \Or(V)$ are conjugate to a unipotent class representative $u$ 
as in (\ref{canon}). Our algorithm 
mirrors that used in the general case of Section \ref{spbadconjel}: 
we first replace summands of $u$ of the form (\ref{2kcont}) by summands 
(\ref{2knew}) to obtain a conjugate element $u'$, and then 
compute $x \in \Or(V)$ such that $g^x=h$. 

If the class $u^{\Or(V)}$ does not split in $\O(V)$, then $g$ and $h$ 
are $\O(V)$-conjugate, and we can obtain an element of $\O(V)$ conjugating $g$ to $h$ by adjusting $x$ (if necessary) by an element of $\Or(V)\setminus \O(V)$ that centralizes $g$. Finally, suppose the class $u^{\Or(V)}$ splits in $\O(V)$ 
(so that $u = \bigoplus W(2k_i)^{[a_i]}$);
if $g$ and $h$ are $\O(V)$-conjugate,  then the conjugating element 
$x$ automatically lies in $\O(V)$. 

This completes our analysis of constructing conjugating elements for the 
orthogonal groups in characteristic 2.

\section{Some examples}
We close this chapter by listing the unipotent 
class representatives and centralizer orders
for 8-dimensional symplectic and orthogonal groups 
defined over fields of even size $q$.
The structures of the centralizers are given by Theorem \ref{uni}. 
Notation is as in earlier sections. 
Recall that $\a$ is a fixed element  of $\F_q$
such that $x^2+x+\a$ is irreducible in $\F_q[x]$.

\begin{table}[htbp] 
\caption{Unipotent class representatives in $G = \Sp_8(q)$, $q$ even} \label{sp8}
\[
\begin{array}{|l|l|}
\hline 
\hbox{Representative } u & |C_G(u)| \\
\hline
W(1)^{[4]} & |G| \\
W(4) & q^7|\Sp_2(q)| \\

W(3) \oplus W(1) & q^{10}|\Sp_2(q)|\,|\Or_2^+(q)| \\
W_\a(3) \oplus W(1) & q^{10}|\Sp_2(q)|\,|\Or_2^-(q)| \\
W(3) \oplus V(2) & q^9|\Sp_2(q)| \\

W(2)^{[2]} & q^{10}|\Sp_4(q)| \\
W(2) \oplus W(1) \oplus V(2) &  q^{12}|\Sp_2(q)|^2 \\
W(2) \oplus V(2)^{[2]} & q^{13}|\Sp_2(q)| \\
W(2) \oplus V_\b(4)\,(\b\in\{1,\a\}) &  2q^9|\Sp_2(q)| \\

W(1)^{[3]} \oplus V(2) &  q^7|\Sp_6(q)| \\
W(1)^{[2]} \oplus W(2) &  q^{11}|\Sp_2(q)|\,|\Sp_4(q)| \\
W(1)^{[2]} \oplus V(2)^{[2]} &  q^{12}|\Sp_4(q)| \\
W(1)^{[2]} \oplus V_\b(4)\,(\b\in\{1,\a\})  &  2q^6|\Sp_4(q)| \\
W(1) \oplus V(4) \oplus V(2) &  q^9|\Sp_2(q)| \\
W(1)\oplus V_\b(6)\,(\b\in\{1,\a\})  &  2q^5|\Sp_2(q)| \\

V(4) \oplus V(2)^{[2]} &  q^{10} \\
V(4) \oplus V_\b(4)\,(\b\in\{1,\a\})  & 2q^8 \\
V_\b(6) \oplus V(2)\,(\b\in\{1,\a\})  &  2q^6 \\
V_\b(8)\,(\b\in\{1,\a\})  &  2q^4 \\
\hline
\end{array}
\]
\end{table}

\begin{table}[htbp] 
\caption{Unipotent class representatives in $G = \O_8^+ (q)$, $q$ even} \label{o8plus}
\[
\begin{array}{|l|l|}
\hline 
\hbox{Representative } u & |C_G(u)| \\
\hline
W(1)^{[4]} & |G| \\
W(4),\,W(4)' & q^5|\Sp_2(q)| \\
W(2)^{[2]},\,W(2) \oplus W(2)' & q^6|\Sp_4(q)| \\
W(3) \oplus W(1) & 2q^8(q-1)^2 \\
W_\a(3) \oplus W_\a(1) & 2q^8(q+1)^2 \\
W(2) \oplus W(1)^{[2]} &  q^9|\Sp_2(q)|\,|\O^+_4(q)| \\
W(2) \oplus V(2)^{[2]} &  q^9|\Sp_2(q)| \\
W(1)^{[2]} \oplus V(2)^{[2]} &  q^6|\Sp_4(q)| \\
W(1) \oplus V(4) \oplus V(2) &  q^5|\Sp_2(q)| \\
V(4)^{[2]} & q^6 \\
V(6) \oplus V(2) &  2q^4 \\
V_\a(6) \oplus V_\a(2) &  2q^4 \\
\hline
\end{array}
\]
\end{table}

\begin{table}[htbp] 
\caption{Unipotent class representatives in $G = \O_8^- (q)$, $q$ even} \label{o8minus}
\[
\begin{array}{|l|l|}
\hline 
\hbox{Representative } u & |C_G(u)| \\
\hline
W(1)^{[4]} & |G| \\
W_\a(3) \oplus W(1) & 2q^8(q^2-1) \\
W(3) \oplus W_\a(1) & 2q^8(q^2-1) \\
W(2) \oplus W_\a(1) \oplus W(1) &  q^9|\Sp_2(q)|\,|\O^-_4(q)| \\
W(2) \oplus V_\a(2)  \oplus V(2) &  q^9|\Sp_2(q)| \\
W(1)^{[2]} \oplus V_\a(2) \oplus V(2) &  q^6|\Sp_4(q)| \\
W(1) \oplus V(4) \oplus V_\a(2) &  q^5|\Sp_2(q)| \\
V_\a(4) \oplus V(4) & q^6 \\
V_\a(6) \oplus V(2) &  2q^4 \\
V(6) \oplus V_\a(2) &  2q^4 \\
\hline
\end{array}
\]
\end{table}

\chapter{Semisimple classes}   \label{semisimpleChap}

Having handled unipotent conjugacy classes of the classical groups 
in Chapters \ref{unigoodchap} and \ref{unibadchap}, we now consider 
semisimple classes. In the next chapter we combine the semisimple 
and unipotent analysis with the Jordan decomposition 
to deal with general conjugacy classes.

\section{Semisimple classes and centralizers} \label{semisimpleconjclasses}

Here we use the notation for classical groups introduced in Chapter \ref{LinearChapter}. Let $F = \F_{q^u}$ with $u=1$ or 2, 
and let $\l \mapsto \bar \l$ be the field automorphism of order $u$. 
Let $V$ be an $n$-dimensional vector space over $F$ and let $\beta$ be a 
non-degenerate alternating, symmetric or hermitian form on $V$, 
and $Q$ a non-degenerate quadratic form. 
Let $\C$  be the isometry group $\C(\beta)$ or $\C(Q)$.

Let $X \in \GL(V)$ be semisimple. Theorem \ref{elementsofCG} gives 
necessary and sufficient conditions for $X$ to be similar to an element 
of $\C$  -- that is, for the intersection $X^{\GL(V)} \cap \C$ to be non-empty, 
where $X^{\GL(V)}$ denotes the conjugacy class of $X$ in $\GL(V)$. 
Our next theorem determines precisely how this intersection splits 
into $\C$-classes.

Recall from Section \ref{fincla} that for $q$ odd, there are, up to congruence, two non-degenerate symmetric bilinear forms
on a given vector space $U$ over $\F_q$, and these are distinguished by their discriminants. We refer to these as the 
two {\it types}\index{form!type}
 of such forms on $U$. If $U$ has even dimension, then we also refer to the two types as plus and minus, in accordance with the description of orthogonal groups in Section \ref{fincla}.

\begin{theorem}\label{MainConjinC}
Let $\C = \C(\b)$ or $\C(Q)$ be an isometry group as above. 
Let $X_1, X_2 \in \C$ be semisimple. 
\begin{itemize}
\item[{\rm (i)}] If $\C$ is a symplectic group, a unitary group, or an orthogonal group in characteristic $2$, then $X_1$ and $X_2$ are conjugate in $\C$ if and only if they are similar. 
\item[{\rm (ii)}] If $\C$ is an orthogonal group in odd characteristic, then $X_1$ and $X_2$ are conjugate in $\C$ if and only if they are similar and the forms induced by $\beta$ on the eigenspaces for the eigenvalues $+1$ and $-1$ have the same type.
\end{itemize}
\end{theorem}

Simultaneously, we prove the following result describing the 
centralizers of semisimple elements in $\C$.
Recall, from Chapter \ref{LinearChapter}, that for $X \in \C$, 
the generalized elementary divisors of $X$ lie in the set 
$\Phi = \Phi_1\cup \Phi_2\cup\Phi_3$ of polynomials listed in 
Definition \ref{ThreeCases}.

\begin{theorem} \label{CentrSem}
Let $\C = \C(\b)$ or $\C(Q)$ as above, and let $X \in \C$ be semisimple.
 For every polynomial $f$ in $\Phi$, let $m_f$ be the multiplicity of 
$f$ as a \ged of $X$, let $d_f = u\deg(f)/2$ 
$($where $u=2$ if $\C$ is unitary and $u=1$ otherwise$)$, and let 
 $\g_f$ denote the restriction of the form $\beta$ or $Q$ to 
the generalized eigenspace $\ker(f(X))$. Then
\begin{eqnarray}\label{CardofCentralizer}
C_{\C}(X) \cong \prod_{f \in \Phi_1} \!\! 
\C(\g_f) \times \prod_{f \in \Phi_2}\!\! \GL_{m_f}(q^{d_f}) 
\times \prod_{f \in \Phi_3}\!\! \GU_{m_f}(q^{d_f}),
\end{eqnarray}
where the products run over all \geds of $X$. 
\end{theorem}

Wall \cite{SCCSS} proved Theorem \ref{MainConjinC} for symplectic groups 
in odd characteristic. We extend his approach to all sesquilinear and 
quadratic forms in arbitrary characteristic.  In Section \ref{semisimplesesquilinear} 
we prove both theorems for sesquilinear forms,  
and complete their proof in Section \ref{sectQuad}
for orthogonal groups in characteristic 2. 

\subsection{Sesquilinear forms}\label{semisimplesesquilinear}
Let $\b$ be a non-degenerate sesquilinear form as above, and for $g \in \GL(V)$, define a form $\b g$ on $V$ by 
\[
\b g\,(u,v) = \b(ug,\,vg) \;\;\hbox{ for all }u,v \in V.
\]
Define $\mathscr{L} = \{\b g : g \in \GL(V)\}$. 
For $X \in \C(\b)$, let $\mathscr{L}^X = \{\g \in \mathscr{L} : \g X = \g\}$, the set of forms in $\mathscr{L}$ fixed by $X$. Let 
$X^{\GL(V)}$ be the conjugacy class of $X$ in $\GL(V)$, and denote by
$\mathscr{M}_X$ the set of $\C(\b)$-conjugacy classes in $X^{\GL(V)} \cap \C(\b)$. In the statement of the next lemma, $^g\!X$ denotes the conjugate $gXg^{-1}$.

\begin{lemma}\label{lxmx} Let $X \in \C(\b)$. There is a bijection from $\mathscr{M}_X$ to the set of orbits of $C_{\GL(V)}(X)$ on $\mathscr{L}^X$. For $^g\!X \in \C(\b)$ \textup{(}where $g \in \GL(V)$\textup{)}, this bijection sends the conjugacy class $(^g\!X)^{\C(\b)}$ to the orbit of $\b g$ under the action of $C_{\GL(V)}(X)$.
In particular, if $X_1,X_2 \in \C(\b)$ are similar to $X$, 
with $X_i =\, ^{g_i}\!X$ $($for $g_i \in \GL(V)$, $i=1,2)$, then $X_1,X_2$ are $\C(\b)$-conjugate if and only if the forms $\b g_1,\,\b g_2$ are in the same $C_{\GL(V)}(X)$-orbit.
\end{lemma}

\begin{proof} 
The given map is a well-defined bijection from $\mathscr{M}_X$ to the set of orbits of $C_{\GL(V)}(X)$ on $\mathscr{L}^X$.
\end{proof}

We next identify some reductions we employ in 
the proof of Theorems \ref{MainConjinC} and \ref{CentrSem}. Let $X \in \C(\b)$, and let $f_1,\ldots,f_h \in \Phi$ be the distinct generalized elementary divisors of $X$. Let the minimal polynomial of $X$ be 
$f(t) = f_1(t)^{e_1} \cdots f_h(t)^{e_h}$, with generalized eigenspaces $V_i = \ker(f_i(X)^{e_i})$. Relative to a basis of $V$ which is a union of bases of the $V_i$, the matrix of $X$ is 
$$
\left( \begin{array}{llll} X_1 & & & \\ & X_2 & & \\ & & \ddots & \\ & & & X_h \end{array}\right),
$$
where $X_i$ is the matrix of the restriction of $X$ to $V_i$. Moreover, the matrix of every form in $\mathscr{L}^X$ has shape 
$$
\left(\begin{array}{lll} B_{11} & \cdots & B_{1h} \\ \vdots & \ddots & \vdots \\ B_{h1} & \cdots & B_{hh} \end{array}\right),
$$
where $X_iB_{ij}X_j^* = B_{ij}$, and $B_{ij} = \varepsilon B_{ji}^*$ for every $1 \leq i,j \leq h$ (where $\varepsilon= -1$ in the symplectic case and $\varepsilon=1$ in the other cases). 
More generally, for every $g \in F[t]$,
\begin{eqnarray}\label{7w}
g(X_i)B_{ij} = B_{ij}g(X_j^{*-1}).
\end{eqnarray}

\begin{lemma}\label{L26w}
With the above notation, if $f_i \neq f_j^*$, then $B_{ij}=0$.
\end{lemma}

\begin{proof}
The polynomial $f_i^{e_i}$ is the minimal polynomial of $X_i$, so, 
taking $g=f_i$ in (\ref{7w}), we get 
$0 = f_i(X_i)^{e_i}B_{ij} = B_{ij}f_i(X_j^{*-1})^{e_i}$. 
Since $f_i$ and $f_j^*$ are 
coprime, $f_i(X_j^{*-1})^{e_i}$ is non-singular. 
Hence $B_{ij}f_i(X_j^{*-1})^{e_i}=0$ implies $B_{ij}=0$.
\end{proof}

\no 
This lemma implies that the matrix of every form in $\mathscr{L}^X$ has 
block diagonal shape
$$
\begin{pmatrix} B_1 & & \\ & \ddots & \\ & & B_h \end{pmatrix}.
$$
Thus the heart of the proof of Theorems \ref{MainConjinC} and \ref{CentrSem} is the case where $h=1$: namely, $X$ has precisely one generalized elementary 
divisor. We handle this case in the next result. 
Since $X$ is semisimple, $e_i=1$ for all $i$.

\begin{proposition}\label{h1cas} Let $X \in \C(\b)$ be semisimple, 
and assume that $X$ has minimal polynomial $f \in \Phi$ which occurs with  
multiplicity $m_f$. Let $d_f = u\deg(f)/2$.
\begin{itemize}
\item[{\rm (i)}] If $f \in \Phi_1$, then $X$ is a scalar and $C_{\C(\b)}(X) = \C(\b)$.
\item[{\rm (ii)}] If $f \in \Phi_2$, then $X^{\GL(V)}\cap \C(\b) = X^{\C(\b)}$, a single $\C(\b)$-conjugacy class, and 
\[
C_{\C(\b)}(X)  = \GL_{m_f}(q^{d_f}).
\]
\item[{\rm (iii)}] If $f \in \Phi_3$, then $X^{\GL(V)}\cap \C(\b) = X^{\C(\b)}$, a single $\C(\b)$-conjugacy class, and 
\[
C_{\C(\b)}(X)  = \GU_{m_f}(q^{d_f}).
\]
\end{itemize}
\end{proposition}

\begin{proof}
(i) This is clear. 

\vspace{1mm}
\no (ii) Suppose $f \in \Phi_2$, so that $f=gg^*$ where $g$ is irreducible and $g \ne g^*$. 
Then $V = \ker g(X) \oplus \ker g^*(X)$. Let $\g \in \mathscr{L}^X$. 
Relative to a suitable basis of $V$, the matrices of $X$ and $\g$ are 
\begin{eqnarray} \label{case2}
X = \left( \begin{array}{cc} X_1 & 0 \\ 0 & X_1^{*-1} \end{array}\right), \quad B = \left( \begin{array}{cc} B_{11} & B_{12}\\ \varepsilon B_{12}^* & B_{22} \end{array}\right).
\end{eqnarray}
By Lemma \ref{L26w}, $B_{11} = B_{22}=0$. The identity $XBX^*=B$ implies $X_1B_{12}X_1^{-1}=B_{12}$, so $X_1$ commutes with $B_{12}$. If we take
$$
Y = \left( \begin{array}{ll} B_{12} & 0 \\ 0 & \mathbb{I} \end{array}\right), \quad J = \left(\begin{array}{ll} 0 & \mathbb{I} \\ \varepsilon \mathbb{I} & 0 \end{array} \right),
$$
where $\mathbb{I}$ is the identity matrix of the same dimension as $X_1$, then $Y$ commutes with $X$, and $YJY^* = B$. The matrix $J$ does not depend on $B$, so given forms in $\mathscr{L}^X$ with matrices $B, B'$ and $Y,Y'$ constructed as above, 
$B = YJY^* = (YY'^{-1})B'(YY'^{-1})^*$ and $YY'^{-1} \in C_{\GL(V)}(X)$. This proves that all forms in 
$\mathscr{L}^X$ are in the same $C_{\GL(V)}(X)$-orbit, and so $X^{\GL(V)}\cap \C(\b) = X^{\C(\b)}$ by Lemma \ref{lxmx}.

It remains to compute the centralizer of $X$ in $\C(\b)$. For this we may suppose that the form $\b$ has matrix $J$,  and that $X_1$ is a block diagonal matrix
$$
\left( \begin{array}{ccc} C & & \\ & \ddots & \\ & & C\end{array} \right),
$$
where $C$ is the companion matrix of $g$. Every $Y \in C_{\GL(V)}(X)$ has the form
$$
\left( \begin{array}{cc} Y_1 & 0 \\ 0 & Y_2 \end{array}\right),
$$
where $Y_1$ commutes with $X_1$ and $Y_2$ with $X_1^{*-1}$, and the centralizer of $X_1$ is isomorphic to $\GL_{m_f}(q^{d_f})$ by Proposition \ref{largerfield}. 
If also $Y \in \C(\b)$, then $YJY^*=J$, giving $Y_2=(Y_1^*)^{-1}$. 
Hence $Y_1$ can be an arbitrary element of the centralizer of $X_1$ and it determines $Y_2$ uniquely. Thus 
$C_{\C(\b)}(X)  = \GL_{m_f}(q^{d_f})$, completing the proof of (ii).

\vspace{2mm}
\no (iii) This case requires more effort. 
Suppose $f \in \Phi_3$, so that $f$ is irreducible of degree $d \ge 2$ and $f = f^*$. For convenience write $h=m_f$, the multiplicity of the generalized elementary divisor $f$ of $X$. Observe that $X$, relative to a suitable basis of $V$, has a block diagonal matrix
$$
\left( \begin{array}{ccc} R & & \\ & \ddots & \\ & & R\end{array}  \right),
$$
where $R$ has minimal polynomial $f$ and there are $h$ diagonal blocks.  By Proposition \ref{largerfield}, 
$C_{\GL(V)}(X) \cong \GL_{h}(q^{2d_f})$, and $Y$ commutes with $X$ if and only if $Y$ is non-singular and has the block matrix shape
\begin{equation}\label{shap}
\left( \begin{array}{ccc} & \vdots & \\ \cdots & f_{ij}(R) & \cdots \\ & \vdots & \end{array} \right)
\end{equation}
for some $f_{ij} \in F[t]$; in the unitary case, where $F = \F_{q^2}$, we can take $f_{ij} \in \F_q[t]$. 
Every form in $\mathscr{L}^X$ has block matrix
$$
B = \left( \begin{array}{ccc} B_{11} & \cdots & B_{1h} \\ \vdots & \ddots & \vdots \\ B_{h1} & \cdots & B_{hh} \end{array} \right),
$$
and the equation $XBX^*=B$ is equivalent to
\begin{eqnarray}\label{15w}
RB_{ij}R^* = B_{ij} \mbox{ for } 1 \leq i,j \leq h.
\end{eqnarray}
Since $f = f^*$, it follows that $R^*$ is similar to $R^{-1}$, so there exists $T \in \GL_d(F)$ such that
\begin{eqnarray}\label{16w}
R^* = T^{-1}R^{-1}T.
\end{eqnarray}
Thus (\ref{15w}) may be rewritten as
$$
R(B_{ij}T^{-1}) = (B_{ij}T^{-1})R.
$$
This shows that $B_{ij}T^{-1}$ belongs to the centralizer of $R$, and so $B_{ij} = f_{ij}(R)T$ for a certain polynomial $f_{ij} \in F[t]$; as above, in the unitary case we can take $f_{ij} \in \F_q[t]$. We obtain the equation
\begin{eqnarray}\label{18w}
B = H \mathcal{T},
\end{eqnarray}
where $H = (f_{ij}(R)) \in C_{\GL(V)}(X)$ and $\mathcal{T} = T \oplus \cdots \oplus T$.

\vspace{2mm} \no {\bf Claim 1.} The matrix $T$ in (\ref{16w}) can be chosen such that $T = \varepsilon T^*$.

\vspace{2mm} 
\no 
{\it Proof of Claim $1$.} 
If $R^* = T^{-1}R^{-1}T$, then $T$ can be replaced by 
$g(R)T$ for some $g \in F[t]$ such that $g(R)$ is invertible. 
Note that $R$ has irreducible minimal polynomial $f(t)$, 
so every non-zero $g(R)$ is in fact invertible. 
 So the aim is to prove that there exists $g(t)$ such that $g(R)T = \varepsilon(g(R)T)^*$ and $g(R)\ne 0$. From (\ref{16w}) we see that  $R^* = T^{*-1}R^{-1}T^*$. Hence $T^*T^{-1}$ commutes with $R$, so $T^* = \phi(R)T$ for some 
$\phi \in F[t]$.
	
	If $T = -\varepsilon T^*$ and $RT = -\varepsilon (RT)^*$, then
	$$
	RT = -\varepsilon (RT)^* = (-\varepsilon T^*)R^* = TR^* = R^{-1}T
	$$
	by (\ref{16w}). But this implies $R=R^{-1}$ and so $R^2=1$, contradicting the assumption that $f\in \Phi_3$. Thus at least one of $T \neq -\varepsilon T^*$ and $RT \neq -\varepsilon (RT)^*$ holds.
	
	If $T \neq -\varepsilon T^*$, then choose $g(t) = 1+\varepsilon\phi(t)$ and deduce that
	$$
	g(R)T = (1+\varepsilon\phi(R))T = T+\varepsilon T^*.
	$$
Then $g(R)T= \varepsilon (g(R)T)^*$, as required. 

If $RT \neq -\varepsilon (RT)^*$, then take $\psi \in F[t]$ such that 
$\psi(R) = R^{-1}$ and let $g(t) = t+\varepsilon\psi(t)\phi(t)$. Then
	$$
	g(R)T = RT+\varepsilon R^{-1}T^* = RT+\varepsilon T^*R^* = RT+\varepsilon (RT)^*,
	$$
and again $g(R)T= \varepsilon (g(R)T)^*$. Claim 1 is now proved.

\vspace{2mm}
We observed in (\ref{shap}) that if $Y \in C_{\GL(V)}(X)$, then $Y$ is a block matrix $(\phi_{ij}(R))$ (where $\phi_{ij} \in \F_q[t]$), so it can be identified with a matrix in $\GL_h(E)$, where $E$ is the field $\F_q[t]/(f)$ in the symplectic and orthogonal cases, and $E$ is the field $\F_q[t]/(f\bar f)$ in the unitary case. The mapping $\phi_{ij}(R) \mapsto \phi_{ij}(R^{-1})$ is a field automorphism of $E$ of order 2. For $Y = (\phi_{ij}(R)) \in \GL_h(E)$, define
\begin{eqnarray}        \label{dag-symbol}\index{$Y^{\dag}$}
Y^{\dag}:=(\phi_{ji}(R^{-1})).
\end{eqnarray}
The map $Y \mapsto Y^{\dag}$ is a ``conjugate-transpose" involutory automorphism of $\GL_h(E)$, where the transpose sends $(\phi_{ij}(R)) \mapsto (\phi_{ji}(R))$ and the conjugate is the automorphism $(\phi_{ij}(R)) \mapsto (\phi_{ij}(R^{-1}))$.

\vspace{2mm} \no {\bf Claim 2.}	Let $B = H\mathcal{T}$ as in $(\ref{18w})$. Then $H = H^{\dag}$. Moreover, if $Y \in C_{\GL(V)}(X) \cong \GL_h(E)$, then $YBY^* = YHY^{\dag}\mathcal{T}$.

\vspace{2mm} 
\no 
{\it Proof of Claim $2$.}	
Since $B$ is the matrix of a sesquilinear form in $\mathscr{L}$, 
we conclude that $B = \varepsilon B^*$. 
By (\ref{18w}), $B^* = \mathcal{T}^* H^*$, and by Claim 1, 
we can suppose that 
$\mathcal{T} = \varepsilon \mathcal{T}^*$. Using this and (\ref{16w}), we deduce that
	\begin{eqnarray*}
		\varepsilon B^* & = & \varepsilon \mathcal{T}^*H^*\\
		& = & (\varepsilon T^* f_{ji}(R^*)) \\
		& = & (Tf_{ji}(R^*)) \\
		& = & (f_{ji}(R^{-1})T)\\
		& = & H^{\dag}\mathcal{T}.
	\end{eqnarray*}
Since $B = \varepsilon B^*$, this implies that 
$H \mathcal{T} = H^{\dag}\mathcal{T}$. 
Since $\mathcal{T}$ is invertible, $H = H^{\dag}$.
	
	Now consider the second assertion in Claim 2. Recall that $H = (f_{ij}(R)) \in C_{\GL(V)}(X)$. Write $Y = (\phi_{ij}(R))$ for $\phi_{ij} \in \F_q[t]$. Now
	\begin{empheq}{equation*}
	\begin{split}
	YBY^* = YH\mathcal{T}Y^* & =  (\phi_{ij}(R))\,(f_{ij}(R)T)\,(\phi_{ij}(R))^*\\
	& =  (\sum_{\lambda,\mu} \phi_{i \lambda}(R) f_{\lambda \mu}(R)T\phi_{j\mu}(R^*))\\
	& =  (\sum_{\lambda,\mu} \phi_{i \lambda}(R)f_{\lambda \mu}(R)\phi_{j\mu}(R^{-1})T)\\
	& =  YHY^{\dag}\mathcal{T}. 
	\end{split}
	\end{empheq}
This completes the proof of Claim 2.

\vspace{2mm}
We now finish the proof of part (iii) of the proposition. 
If $B_1$ and $B_2$ are the matrices of forms in $\mathscr{L}^X$, with $B_1=H_1\mathcal{T}$ and $B_2 = H_2\mathcal{T}$, then $H_1$ and $H_2$, viewed 
as matrices in $\GL_h(E)$, are hermitian, by Claim 2. 
Thus they are congruent: there exists $Y \in \GL_h(E) \cong C_{\GL(V)}(X)$ 
such that $H_1 = YH_2Y^{\dag}$. By Claim 2, 
this implies $$B_1 = H_1\mathcal{T} = YH_2Y^{\dag}\mathcal{T} = YB_2Y^*.$$ 
There is only one orbit of forms in $\mathscr{L}^X$ under the 
action of $C_G(X)$, and 
hence $X^{\GL(V)} \cap \C(\b)= X^{\C(\b)}$ by Lemma \ref{lxmx}. 

For the final assertion, we must identify $C_{\C(\b)}(X)$. 
Let $B$ be the matrix of the form $\b$, and let $B = H\mathcal{T}$ 
as in $(\ref{18w})$. By Claim 2, $H$ is hermitian (with respect to 
the automorphism~$\dagger$). Consider $W = E^h$ as an $h$-dimensional unitary space 
with respect to the form with matrix $H$. 
If $Y \in C_{\GL(V)}(X)$, then 
\[
\begin{array}{lll}
Y \in \C(\b) & \Leftrightarrow & YBY^* = B\\
    & \Leftrightarrow & YHY^{\dag}\mathcal{T} = H\mathcal{T} \\
    & \Leftrightarrow & YHY^{\dag} = H.
\end{array} 
\]
Hence the centralizer of $X$ in $\C(\b)$ is isomorphic to $\GU(W)$. This is 
the group $\GU_{m_f}(q^{d_f})$ in the conclusion of (iii).  \qedhere
\end{proof}

\no {\bf Proof of Theorems $\ref{MainConjinC}$ and $\ref{CentrSem}$ 
for $\C = \C(\b)$. }

\vspace{2mm}
Let $X \in \C(\b)$ be semisimple, and let $f(t) = f_1(t)\cdots f_h(t)$ be the minimal polynomial of $X$, where each $f_i \in \Phi$ (and $f_1,\ldots,f_h$ are distinct). Let $V_i = \ker f_i(X)$. 
Now $X$, relative to a basis of $V$ which is a union of bases of the $V_i$, 
has block diagonal matrix $\bigoplus_{i=1}^h X_i$, 
where $X_i$ is the matrix of the restriction of $X$ to $V_i$. 
Moreover, by Lemma \ref{L26w}, the matrix of $\b$ has the shape $B = \bigoplus_{i=1}^h B_i$, and 
every $Y \in C_{\GL(V)}(X)$ has the shape $\bigoplus_{i=1}^h Y_i$, where $[Y_i,X_i] = 1$ for all $i$.
Also $Y \in \C(\b)$ if and only if each $Y_i$ fixes the form $\b_i$ with matrix $B_i$. Hence 
\[
C_{\C(\b)}(X) = \prod_i C_{\C(\b_i)}(X_i).
\]
Each factor $C_{\C(\b_i)}(X_i)$ is given by Proposition \ref{h1cas}, and hence $C_{\C(\b)}(X) $ is as in the conclusion of Theorem \ref{CentrSem}, completing the proof of that theorem.

We now prove Theorem \ref{MainConjinC}. Let $X' \in \C(\b)$ be similar to $X$, 
so $X' = X^g$ for some $g \in \GL(V)$. Then $X \in \C(\b g)$, and, writing matrices with respect to the above basis, the matrix of $\b g$ is $B' = gBg^*$. By Lemma \ref{lxmx}, $X$ and $X'$ are conjugate in $\C(\b)$ if and only if $\b$ and $\b g$ are in the same $C_{\GL(V)}(X)$-orbit: namely, $B$ and $B'$ are congruent in $C_{\GL(V)}(X)$.
As above, since $X$ fixes $\b g$, the matrix $B'$ has block diagonal form $B' = \bigoplus_{i=1}^h B_i'$. 

Assume now that, for all $i$, the forms $B_i$ and $B_i'$ have the 
same type: both are symplectic, or both are unitary, or both are 
orthogonal of the same discriminant. This assumption holds if 
$\b$ is symplectic or unitary. By Proposition \ref{h1cas}, $B_i$ and $B_i'$ are congruent in $C_{\GL(V_i)}(X_i)$, so there exists $Y_i \in C_{\GL(V_i)}(X_i)$ such that $Y_iB_iY_i^* = B_i'$. Then $Y = \bigoplus_{i=1}^h Y_i$ is in 
$C_{\GL(V)}(X)$ and satisfies $YBY^*=B'$. Hence $X$ and $X'$ are conjugate in $\C(\b)$. This completes the proof of part (i) of Theorem \ref{MainConjinC} 
for sesquilinear forms.

It remains to consider the case where $\b$ is a symmetric form in odd characteristic. By Theorem \ref{elementsofCG}, if $f_i \in \Phi_2 \cup \Phi_3$, then $B_{i}$ and $B_i'$ have the same type. 
If $f_i = t \pm 1 \in \Phi_1$,  then $B_{i}$ can be either of the two possible types. 
If at most one of $t+1$ and $t-1$ is an elementary divisor for $X$, 
then the types of $B_{i}$ and $B_{i}'$ are uniquely determined, and 
they coincide, so 
$B$ and $B'$ are congruent in $C_{\GL(V)}(X)$. However, if both $t+1$ and $t-1$ are elementary divisors for $X$, say $f_1 = t+1$
and $f_2 = t-1$, then $B_1',B_2'$ can be either the same types as $B_1,B_2$, or the opposite types. Thus in this case the conjugacy class $X^{\GL(V)} \cap \C(\b)$ splits into two $\C(\b)$-classes. 

This completes the proof of Theorems $\ref{MainConjinC}$ and $\ref{CentrSem}$ for $\C = \C(\b)$.

\subsection{Orthogonal groups in characteristic 2} \label{sectQuad}
We now prove Theorems \ref{MainConjinC} and \ref{CentrSem} for the 
orthogonal groups in characteristic 2. 
Let $\C = \C(Q) \cong \Or_{2m}^\e (q)$, where $Q$ is a non-degenerate 
quadratic form on a $2m$-dimensional vector space $V$ over $F = \F_q$, 
with $q = 2^k$. Let $\b_Q$ be the associated symplectic form on $V$. 
As in Section \ref{semisimplesesquilinear}, the heart of the proof 
is the case where 
the elements in question ($X_1,X_2$ in Theorem \ref{MainConjinC}, and 
$X$ in Theorem \ref{CentrSem}) have precisely one elementary divisor. 
The next result is the analogue of Proposition \ref{h1cas}. 

\begin{proposition}\label{h1caso2} Let $\C(Q) \cong \Or_{2m}^\e (q)$, and let $\b_Q$ be the associated symplectic form.  Let $X \in \C(Q)$ be semisimple, 
and assume that $X$ has minimal polynomial $f \in \Phi$.
\begin{itemize}
\item[{\rm (i)}] If $f \in \Phi_1$, then $X = I$ and $C_{\C(Q)}(X) = \C(Q)$.
\item[{\rm (ii)}] If $f \in \Phi_2\cup \Phi_3$, then $X^{\GL(V)}\cap \C(Q) = X^{\C(Q)}$, a single $\C(Q)$-conjugacy class, and 
\[
C_{\C(Q)}(X)  = C_{\C(\b_Q)}(X).
\]
Hence $C_{\Or_{2m}^\e (q)}(X) = C_{\Sp_{2m}(q)}(X)$ $($which is given by Proposition $\ref{h1cas})$.
\end{itemize}
\end{proposition}

\begin{proof} 
Part (i) is trivial. 
So assume $f\in \Phi_2\cup \Phi_3$.  Recall that 
$Q$ can be represented 
by a matrix $A = (a_{ij})$ such that $Q(v) = vAv^{\tr}$ for all $v \in V$. 
	
 Let $\Or_{2m+1}(q)$ be the group of isometries for the quadratic form with matrix
$$
\widehat{A} =\begin{pmatrix} 1&0 \\0&A \end{pmatrix}.
$$
 Every matrix in $\Or_{2m+1}(q)$ has shape
\begin{equation}\label{matr}
\widehat{Y}= \begin{pmatrix} 1&0 \\v& Y \end{pmatrix},
\end{equation}
where $v$ is a $2m$-dimensional column vector. Note that if $v=0$, then $Y$ is an isometry in 
$\C(A) = \Or^{\epsilon}_{2m}(q)$. It is well known (see \cite[14.1]{CGSA}) that the map $\widehat{Y} \mapsto Y$ is an isomorphism between $\Or_{2m+1}(q)$ and $\C(A+A^{\tr}) \cong \Sp_{2m}(q)$, and the inverse of this isomorphism maps $Y$ to a matrix of the above form, where $v=0$ if and only if $Y \in \Or^{\epsilon}_{2m}(q)$. We now use this isomorphism to describe the centralizer of $X$ in $\Or^{\epsilon}_{2m}(q)$. The centralizer of $X$ in $\Sp_{2m}(q)$ is isomorphic to the centralizer of $\widehat{X}$ in $\Or_{2m+1}(q)$, where
	$$
	\widehat{X} = \begin{pmatrix} 1&0 \\0& X \end{pmatrix}.
	$$
	Since $t+1$ is not an elementary divisor of $X$, every element of the centralizer of $\widehat{X}$ in $\Or_{2m+1}(q)$ has the form $\widehat{Y}$ as in (\ref{matr}) with $v=0$ for some $Y$, and by the above observations, $Y \in \Or^\e_{2m}(q)$.  
Conversely, it is clear that for every $Y \in C_{\Or^\e_{2m}(q)}(X)$, the corresponding $\widehat{Y}$ is in the centralizer of $\widehat{X}$ in $\Or_{2m+1}(q)$. 
This proves that $|C_{\Or^\e_{2m}(q)}(X)| = |C_{\Sp_{2m}(q)}(X)|$, so these centralizers are equal, since the first is obviously contained in the second.
	
	Finally, if $X_1,X_2 \in \Or^\e_{2m}(q)$ are similar, both having minimal polynomial $f\in \Phi_2\cup \Phi_3$, then they are conjugate in $\Sp_{2m}(q)$ by Proposition \ref{h1cas}, and so $\widehat{X}_1$ and $\widehat{X}_2$ are conjugate in $\Or_{2m+1}(q)$. Hence there exists
	$$
	\widehat{Z} = \begin{pmatrix} 1&0 \\v_Z& Z \end{pmatrix} \in \Or_{2m+1}(q)
	$$
	such that $\widehat{Z}^{-1}\widehat{X}_1\widehat{Z} = \widehat{X}_2$. 
Since $t+1$ is not an elementary divisor of $X_1$,  we deduce that
$\widehat{Z}$ preserves the diagonal block structure of $\widehat{X}_1$ and $\widehat{X}_2$, so $v_Z=0$. This shows that 
$Z \in \Or^\e_{2m}(q)$, and so $X_1$ and $X_2$ are conjugate in $\Or^\e_{2m}(q)$, completing the proof.
\end{proof}
\no
As in Section \ref{semisimplesesquilinear},
the conclusions 
of Theorems \ref{MainConjinC} and \ref{CentrSem} now follow. 

\subsection{Special and Omega groups}  \label{sectSpecial}
We now determine the splitting of  
those semisimple classes in the isometry 
groups $\C = \GU_n(q)$ or $\Or_n^\e(q)$ that lie 
in the special groups $\Spec = \SU_n(q)$ or $\SO_n^\e(q)$, 
or in the Omega group\index{Omega group} $\O = \O_n^\e(q)$. 
This is a matter of computing the 
indices $|C_\C (X): C_\Spec (X)|$ and $|C_\C (X): C_\O (X)|$ 
for semisimple elements $X$ of $\Spec$ and $\O$. 
At the end of the section we show how to determine
whether a given semisimple element of $\Spec$ lies in $\O$.

\vspace*{1cm}
\begin{proposition}\label{speco}
\mbox{}
\begin{itemize}
\item[{\rm (i)}] Let $\C = \GU_n(q)$ and $\Spec = \SU_n(q)$. If $X \in \Spec$ is  semisimple, then $|C_\C(X): C_\Spec (X)| = q+1$ and $X^\C = X^\Spec$.
\item[{\rm (ii)}] Let $\C = \Or_n^\e(q)$ and $\Spec = \SO_n^\e(q)$ with $q$ odd. If $X \in \Spec$ is semisimple, then
$|C_\C(X): C_\Spec(X)| = c_S$, where 
\[
c_S = \left\{\begin{array}{l} 2, \hbox{ if $X$ has an eigenvalue }\pm 1, \\
                                          1, \hbox{ otherwise.}
\end{array}
\right.
\]
The class $X^\C$ splits into $2/c_s$ 
conjugacy classes 
in $\Spec$ of equal size; if $c_S = 1$, then representatives 
are $X$ and $X^s$ for $s \in \C\setminus \Spec$. 
If $n$ is odd, then $c_S = 2$ always.
\item[{\rm (iii)}] Let $\Spec = \SO_n^\e(q)$ and $\O = \O_n^\e(q)$. If $X \in \O$ is semisimple, then
$|C_\Spec(X): C_\O(X)| = c_\O$, where 
\[
c_\O = \left\{\begin{array}{l} 2, \hbox{ if $X$ has an eigenvalue }\pm 1, \hbox{ or if $q$ is odd,} \\
                                          1, \hbox{ otherwise.}
\end{array}
\right.
\]
The class $X^\Spec$ splits into $2/c_\O$
conjugacy classes in $\O$ of equal size; if $c_\O=1$, then representatives are $X$ and $X^v$ for $v \in \Spec \setminus \O$.
If $n$ is odd, then $c_\O = 2$ always.
\end{itemize}
\end{proposition}

\begin{proof} 
(i) 
We show that the centralizer $C_\C(X)$, given by (\ref{CardofCentralizer}), contains elements of all determinants in $D:=\{\l \in \F_{q^2}: \l^{q+1}=1\}$. This is clearly the case if $X$ has a generalized elementary divisor $f \in \Phi_1$, since $C_\C(X)$ has a factor $\C(\g_f) = \GU_{m_f}(q)$. Consider a generalized elementary divisor $f = gg^* \in \Phi_2$, where $g$ is irreducible and $g\ne g^*$. The corresponding factor of $C_\C(X)$ is $\GL_m(q^{2d})$, where $m=m_f$ and $d = \deg{g}$, and from the proof of Proposition \ref{h1cas} the
factor consists of block diagonal matrices $(Y,(Y^*)^{-1})$ for $Y \in \GL_m(q^{2d})$. 
We have an embedding $$\GL_m(q^{2d}) \le \GL_{md}(q^2) < \GU_{2md}(q) \le \C,$$ 
and, via this embedding, $Y \in \GL_m(q^{2d})$ can have any determinant $\l \in \F_{q^2}^*$. Hence $(Y,(Y^*)^{-1})$ can have any determinant $\l \bar \l^{-1}$ in $D$, proving the result for generalized elementary divisors $f = gg^* \in \Phi_2$.
Finally, consider a generalized elementary divisor $f \in \Phi_3$ of $X$, with corresponding factor $\GU_m(q^d)$ of $C_\C(X)$, where $d = \deg{f}$. Here $d$ is odd by Proposition \ref{elem}, and we have an embedding $\GU_m(q^d) \le \GU_{md}(q)$, for which \cite[(4.3.13)]{KL} shows that $\GU_m(q^d)$ contains elements of all determinants in $D$. This proves (i). 

\vspace{2mm}
(ii) Note that $|\C:\Spec| = 2$, and the elements in the coset $\C\setminus \Spec$ are those of determinant $-1$. As in the proof of (i), we need to show that $C_\C(X)$ contains elements of determinant $-1$ if and only if $X$ has an eigenvalue $\pm 1$. The centralizer $C_\C(X)$ is given by (\ref{CardofCentralizer}). Let $f$ be a generalized elementary divisor of $X$. If $f \in \Phi_1$, then $X$ has an eigenvalue $\pm 1$ and $C_\C(X)$ has a factor $\Or_{m_f}(q)$, which has elements of determinant $\pm 1$. On the other hand, if $f \in \Phi_2 \cup \Phi_3$, then the corresponding factors $\GL_{m}(q^d)$ and $\GU_{m}(q^{d})$ (where $m=m_f, d=d_f$) lie in the special orthogonal group $\Spec$ via the embeddings 
$$\GL_m(q^d) \le \SO_{2m}^+(q^d) < \SO^+_{2md}(q)$$ and 
$$\GU_m(q^d) \le \SO_{2m}^\e(q^d) < \SO^\e_{2md}(q),$$ where $\e = (-)^{m}$ (see \cite[4.2.7, 4.3.18]{KL}). This proves (ii).

\vspace{2mm}
(iii)  This follows as in (ii) once we establish that $\GL_m(q^d) \le \O^+_{2md}(q)$ and $\GU_m(q^d) \le \O^\e_{2md}(q)$ if and only if $q$ is even, where $m=m_f,d=d_f$ and $f \in \Phi_2 \cup \Phi_3$ as before. 
For even $q$ these containments are clear, since neither 
$\GL_m(q^d)$ nor $\GU_m(q^d)$ has a subgroup of index 2. 

Assume $q$ is odd. Consider the embedding 
$$\GL_m(q^d) \le \GL_{md}(q) < \SO^+_{2md}(q) = \SO(V),$$ 
where $\GL_{md}(q)$ stabilizes two maximal totally isotropic subspaces $W = \la e_1,\ldots,e_{md} \ra$ and $W' = \la f_1,\ldots, f_{md}\ra$ 
(and $e_i,f_i$ are in
a standard hyperbolic basis). 
Consider $g \in \GL_{md}(q)$ sending $e_1 \mapsto \l e_1$, $f_1 \mapsto \l^{-1}f_1$ and fixing all other $e_i,f_i$, where $\l \in \F_q$ is a non-square. Now $g$ can be expressed as the product $r_{e_1-\l f_1}r_{e_1-f_1}$ of the reflections in the vectors $e_1-\l f_1$ and $e_1-f_1$, and so  
has non-identity spinor norm. It follows that $g \not \in \O(V)$, and hence no element in the subgroup $\GL_{md}(q)$ of determinant $\l$ is in $\O(V)$, proving that $\GL_m(q^d) \not \le \O^+_{2md}(q)$, as required. 
Now consider the embedding 
$$\GU_m(q^d) \le \SO_{2m}^\e(q^d) < \SO^\e_{2md}(q) = \SO(V)$$ 
where $d = d_f$ and $f \in \Phi_3$.
Choose a natural subgroup $H = \GU_1(q^d)$ of the left-hand group, 
so $H \le \SO_{2}^-(q^d) < \SO^-_{2d}(q) = \SO(W)$, fixing $W^\perp$ pointwise. From \cite[4.3.15]{KL} we see that $N_{\O(W)}(H)$ has order $d(q^d+1)/2$, whereas $N_{\SO(W)}(H) = \GU_1(q^d).d$ has order $d(q^d+1)$. Hence $H \not \le \O(W)$, and it follows that 
$\GU_m(q^d) \not \le \O^\e_{2md}(q)$, as required. 
\end{proof}

Finally, we record a consequence of the above proof that determines which semisimple elements of $\SO(V)$ are in $\O(V)$.

\begin{lemma}   \label{oddOrderinOmega}
Let $V = V_n(q)$ be an orthogonal space with quadratic form $Q$, let $\Spec = \SO(V)$, and $\O = \O(V)$.
\begin{itemize}
\item[{\rm (i)}] Every element of odd order in $\Spec$ is also in $\Omega$.
\item[{\rm (ii)}] Let $q$ be odd and let $x \in \Spec$ be semisimple. 
Assume that $x$ has a unique \ged $f \in \Phi$ of degree $d$ and 
$f$ occurs with multiplicity $m$.
\begin{itemize}
\item[{\rm (a)}] If $f = t-1$, then $x =I \in \Omega$.
\item[{\rm (b)}] If $f = t+1$, then $x = -I \in \Omega$ if and only if $Q$ has square discriminant.
\item[{\rm (c)}] Let $f \in \Phi_2$ with $f=gg^*$, where $g \ne g^*$ is irreducible, and let $C$ be the companion matrix of $g$. Then $x \in \Omega$ if and only if either $m$ is even, or $\det(C)$ is a square in $\mathbb{F}_q^*$.
\item[{\rm (d)}] Let $f \in \Phi_3$ have degree $d$ and companion matrix $C$. 
Then $x \in \Omega$ if and only if either $m$ is even,
 or the order of $C$ divides $(q^{d/2}+1)/2$.
\end{itemize}
\end{itemize}
\end{lemma}

\begin{proof}
Part (i) is obvious, as $|\Spec:\O| = 2$. Part (ii)(b) follows from 
\cite[Prop.\ 2.5.13]{KL}. Finally, (ii)(c) and (ii)(d) follow from 
the argument for part (iii) of the previous proof.
\end{proof}

If $x \in \Spec$ is semisimple and has \geds $f_1, \dots, f_k$, then 
we can decide its membership in $\Omega$ by applying the lemma to each $f_i$. 

\section{Representatives for semisimple classes} \label{formsandelements}
We show how to write down an explicit set of representatives for the  
semisimple conjugacy classes of classical groups. We give 
the matrices of both the representatives and the corresponding sesquilinear 
or quadratic forms with respect to an appropriate basis. 
We thank Donald Taylor for his assistance with this task.

Let $\beta$ or $Q$ be a sesquilinear or quadratic form on $V$, 
and let $\C = \C(\beta)$ or $\C(Q)$. 
By Lemma \ref{L26w}, every semisimple class in $\C$ has a block diagonal representative $X = \bigoplus_{i=1}^h X_i$, where each $X_i$ has a single \ged $f_i$, and the form $\b$ or $\b_Q$ has corresponding matrix $B = \bigoplus_{i=1}^h B_i$. By Theorem \ref{MainConjinC}, it is sufficient to write down,
for each $f \in \Phi$, a matrix $X_f$, 
with \ged $f$ of multiplicity 1 
when $f \in \Phi_2\cup \Phi_3$, 
and a form $\g_f = \b_f$ or $Q_f$ with matrix $B_f$
such that $X_f \in \C(\g_f)$. We do this separately for 
the cases where $f \in \Phi_1$, $\Phi_2$ or $\Phi_3$.

Let $C(h)$ denote the companion matrix of $h(x)$, a monic polynomial 
of degree $d$.

\subsubsection{Case $f \in \Phi_1$} 
Here $X_f$ is a scalar matrix $\pm I$ in the symplectic and orthogonal cases, and is $\l I$ with $\l^{q+1}=1$ in the unitary case, and $B_f$ can be taken with respect to any basis (but we must allow both types for $Q_f$ in the orthogonal case).

\subsubsection{Case $f \in \Phi_2$} Let $f = gg^*$, where $g\ne g^*$ and $g$ is irreducible. With respect to a suitable basis we take 
$$
X_f = \begin{pmatrix} C(g) & \\ & C(g)^{*-1} \end{pmatrix}, 
\quad B_f = \begin{pmatrix} \mathbb{O} & \mathbb{I} \\ \varepsilon \mathbb{I} & \mathbb{O}
\end{pmatrix},
$$
where $\varepsilon$ is $-1$ (for $\g_f = \b_f$ symplectic), $1$ (for $\g_f = \b_f$ symmetric or unitary), or $0$ (for $\g_f = Q_f$ quadratic). 

\subsubsection{Case $f \in \Phi_3$} 
\noindent
\textit{Symplectic case}. Here $\deg{f}$ is even by Proposition \ref{elem}, and $f$ has constant term 1. Since $f=f^*$, we can write
\begin{equation}\label{f_t}
f(t) = 1+a_1t+ a_2t^2 + \cdots + a_dt^d+a_{d-1}t^{d+1}+ \cdots + a_1t^{2d-1}+t^{2d}.
\end{equation}
Define
\[
X_f = C(f)^{\tr},\;\;\;B_f = \left( \begin{array}{cc} \mathbb{O} & -P^{\tr} \\ P & \mathbb{O} \end{array}\right),
\]
where $C(f)$ is the companion matrix of $f$, and $P$ is the $d \times d$ upper triangular matrix with constant upper diagonals:

$$
P = \left( \begin{array}{llllll} 1 & a_1 & a_2 & \cdots & a_{d-2} & a_{d-1} \\ & 1 & a_1 & \ddots & \ddots & a_{d-2} \\ && \ddots & \ddots & \ddots & \vdots \\ &&  & \ddots & \ddots & a_2 \\ &&&& 1 & a_1 \\ &&&&& 1 \end{array}\right).
$$
A computation shows that $X_f \in \C(B_f)$.\\
\\
\textit{Quadratic case}. The degree of $f \in \Phi_3$ is even. 
If $\deg{f}=2$, so $f=t^2+at+1$, then we take $X_f = C(f)$, and $Q_f$ to be the quadratic form with matrix 
$$
B_f = \left(\begin{array}{ll} 1 & -a \\ 0 & 1 \end{array}\right).
$$
Then $X_f \in \C(Q_f)$.

Now suppose $\deg{f}>2$, where $f$ is defined in (\ref{f_t}), and $d \geq 2$.
Let  $X_f = C(f)$, and define $Q_f$ to be the quadratic form with $2d\times 2d$ matrix 
\[
B_f = \begin{pmatrix} 0 & A \\ 0 & 0 \end{pmatrix},
\]
where the 0's are zero matrices of the appropriate sizes, $A$ is the $(d+1) \times (d+1)$ matrix
$$
\left( \begin{array}{lllll}  1 & b_0 & b_1 & \cdots & b_{d-1} \\ & 1 & \ddots & \ddots & \vdots \\ && \ddots & \ddots & b_1 \\ &&& \ddots & b_0 \\ &&&& 1 \end{array}\right),
$$
and the coefficients $b_i$ are defined as follows. If $p$ is odd, then the vector of the $b_i$'s satisfies the linear system
$$
\left(\begin{array}{ccc} b_0 & \cdots & b_{d-1} \end{array}\right) \left( \begin{array}{lllll} 2 & a_1 & a_2 & \cdots & a_{d-1} \\ & 1 & a_1 & \ddots & \vdots \\ && \ddots & \ddots & a_2 \\ &&& \ddots & a_1 \\ &&&& 1 \end{array}\right) = \left(\begin{array}{cccccc} 2a_1 & a_2-1 & a_3 & a_4 & \cdots & a_d \end{array}\right).
$$
If $p=2$, then the vector of the $b_i$'s satisfies the linear system
$$
\left(\begin{array}{ccc} b_0 & \cdots & b_{d-1} \end{array}\right) \left( \begin{array}{lllll} c_1 & a_1 & a_2 & \cdots & a_{d-1} \\ c_2 & 1 & a_1 & \ddots & \vdots \\ \vdots && \ddots & \ddots & a_2 \\ \vdots &&& \ddots & a_1 \\ c_d &&&& 1 \end{array}\right) = \left(\begin{array}{cccccc} \delta & a_2-1 & a_3 & a_4 & \cdots & a_d \end{array}\right),
$$
where
\begin{eqnarray*}
	c_i & = & \sum_{j=1}^{d+1-i} a_{j-1}a_{d+1-j} \quad \hbox{ for } i=1, \dots, d;\\
	\delta & = & \sum_{j=0}^d a_ja_{d+j-1}
\end{eqnarray*}
and we set $a_0=1$ and $a_{d+j}=a_{d-j}$ for $j=1, \dots, d$. A computation shows that $X_f \in \C(Q_f)$.\\
\\
\textit{Unitary case}. Here $\deg{f}$ is odd by Proposition \ref{elem}. 
Since $f=f^*$ we can write
\[
f(t) = a_0+a_1t+ \cdots +a_dt^d+a_0t^{d+1}(\overline{a}_d+\overline{a}_{d-1}t+ \cdots +\overline{a}_0t^d)
\]
(an irreducible polynomial in $\mathbb{F}_{q^2}[t]$ with $a_0\overline{a}_0 =1$). Choose $b_0 \in \mathbb{F}_{q^2}$ such that $b_0^{q-1} = a_0^{-1}$, and define
\[
b_i :=  \overline{b}_0 \sum_{j=0}^i a_j \: \mbox{ for } 1 \leq i \leq d.
\]
Choose also $c \in \mathbb{F}_{q^2}$ such that $c^{q+1} = b_d+\overline{b}_d$. 
Now define the $(2d+1)\times (2d+1)$ matrix
$$
X_f = \left(\begin{array}{llll|c|llll} & 1 &&&&&&& \\ && \ddots &&&&&& \\ &&& 1 &&&&& \\ &&&&&&&& 1/\overline{b}_0 \\ \hline &&&& 1 &&&& -\overline{c}/\overline{b}_0 \\ \hline b_0 & b_1 & \cdots & b_{d-1} & c &&&& - \overline{b}_d/\overline{b}_0 \\ &&&&& 1 &&& -\overline{b}_{d-1}/\overline{b}_0 \\ &&&&&& \ddots && \vdots \\ &&&&&&& 1 & -\overline{b}_1/\overline{b}_0 \end{array}\right).
$$
Observe that 
$X_f$ has characteristic and minimal polynomial $f$ and preserves the hermitian form with matrix
$$
B_f = \begin{pmatrix}
   && 1 \\ & \iddots & \\ 1 &&
\end{pmatrix}.
$$
In Section \ref{SectListGen}, we need to compute the matrix of the form preserved by the companion matrix $C(f)$. We  can compute this by a change of basis: if $X_f$ and $B_f$ are the matrices described above, and $P \in \GL(V)$ satisfies $PX_f P^{-1}=C(f)$, then $C(f)$ preserves  the form $PB_fP^*$.

\vspace*{2mm}
This completes the description of semisimple class representatives in 
the isometry groups $\C = \C(\b)$ or $\C(Q)$. Representatives in 
$\Spec$ and $\O$ 
can be written down using Proposition \ref{speco}.

\section{Generators for the centralizer of a semisimple element}\label{centgens}
Having identified the structure of the centralizer of a semisimple element 
in a classical isometry group in Theorem \ref{CentrSem}, we now show how 
to write down a generating set for the centralizer.

We assume that the following algorithms are available.
\begin{itemize}
\item (\textit{Algorithm} 1)\index{Algorithm 1} Given similar matrices $X,Y \in \GL_n(q)$, we determine explicitly $Z \in \GL_n(q)$ such that $Z^{-1}XZ=Y$. This follows from the Jordan basis algorithm: if $J$ is the Jordan form of $X$ and $Y$ and $J=P_XXP_X^{-1}=P_YYP_Y^{-1}$ for $P_X, P_Y \in \GL(V)$, then $Z = P_X^{-1}P_Y$.
\item (\textit{Algorithm} 2)\index{Algorithm 2} Given matrices $B_1,B_2$ of two non-degenerate sesquilinear or quadratic forms on $V=V_n(q^u)$ of the same type (where $u=2$ for unitary forms, $u=1$ otherwise), we determine explicitly $T \in \GL_n(q^u)$ such that $TB_1T^*=B_2$ (or $TB_1T^*-B_2$ is alternating in the case of quadratic forms).
\end{itemize}
Descriptions of such algorithms appear in \cite{JF} and \cite{GSM} respectively. Algorithm 2 allows us to write a generating set for a classical group in any basis. Suppose a sesquilinear or quadratic form $\beta$ has matrix $B$ in a certain basis, and let $B_0$ be the matrix of $\beta$ in the basis with respect to which the standard generators for $\C(\b)$ are defined. 
If $\{Y_1, \dots, Y_r\}$ generates $\C(B_0)$, 
then $\{ TY_1T^{-1}, \dots, TY_rT^{-1} \}$ generates $\C(B)$, where $T$ is the matrix returned by Algorithm 2 such that $B = TB_0T^*$.

Let $B$ be the matrix of a non-degenerate sesquilinear form on $V$. Let $X \in \C(B)$ be semisimple and let $f_1, \dots, f_h \in \Phi$ be the \geds of $X$. 
With respect to a suitable basis,
\begin{eqnarray} \label{Xdiagform}
X = \begin{pmatrix}
X_1 && \\ & \ddots & \\ && X_h
\end{pmatrix},
\end{eqnarray}
where each $X_i$ is the matrix of the restriction of $X$ to the generalized eigenspace $\ker(f_i(X))$. By Lemma \ref{L26w}, with respect to this basis, 
\begin{eqnarray} \label{Bdiagform}
B = \begin{pmatrix}
B_1 && \\ & \ddots & \\ && B_h
\end{pmatrix}
\end{eqnarray}
and each $X_i$ is an isometry for $B_i$. If $Y \in C_{\C(B)}(X)$, 
then $Y$ is the block diagonal sum of $Y_i$ for $Y_i \in C_{\C(B_i)}(X_i)$, 
so we can write down generators for $C_{\C(B)}(X)$ 
once we have done this 
when $X$ has a single generalized elementary divisor.

Suppose that $X$ has a single \ged $f$, and let $m$ be the multiplicity 
of $f$ and $d = \deg{f}$. \\

\no 
\textbf{Case 1:} $f \in \Phi_1$. Here $X$ is a scalar matrix, so its centralizer is $\C(B)$. A generating set for $\C(B)$ is referenced in Section \ref{fincla}. \\

\no 
\textbf{Case 2:} $f \in \Phi_2$. Here $f=gg^*$, where $g$ is irreducible 
and $g\ne g^*$. Let $d' = d/2$ be the degree of $g$. 
As in Section \ref{formsandelements}, we choose 
a basis with respect to which 
$$
X = \left(\begin{array}{ll} \widehat{X} & \mathbb{O} \\ \mathbb{O} & \widehat{X}^{*-1} \end{array}\right),
$$
where $\widehat{X}$ is a block diagonal sum of $m$ copies of the companion matrix $C(g)$. By Lemma \ref{L26w}, the form preserved by $X$ is
\begin{eqnarray}    \label{formadiB}
B = \begin{pmatrix*}[l]
\mathbb{O} & A \\ \varepsilon A^* & \mathbb{O}
\end{pmatrix*},
\end{eqnarray}
where $\varepsilon= -1$ in the symplectic case and $1$ otherwise; and the centralizer of $X$ in $\C(B)$ is the group of matrices of the form
$$
\left(\begin{array}{ll} Y & \mathbb{O} \\ \mathbb{O} & A^*Y^{*-1}A^{*-1} \end{array}\right),
$$
where $Y \in C_{\GL_{d'm}(q)}(\widehat{X}) \cong \GL_m(q^{d'})$. 
So, if $Y_1,Y_2$ are standard generators for 
$C_{\GL_{d'm}(q)}(\widehat{X})$, 
then the centralizer of $X$ in $\C(B)$ is generated by
$$
\left(\begin{array}{ll} Y_i & \\ & A^*Y_i^{*-1}A^{*-1} \end{array}\right) \quad (i=1,2).
$$
\\
\textbf{Case 3:} $f \in \Phi_3$. Here $f=f^*$ is irreducible. Let $E$ be the field $\F_q[t] /(f)$ in the symplectic and orthogonal cases, and let $E = \F_q[t]/(f\bar f)$ in the unitary case. We follow the argument in the proof of part (iii) of Proposition \ref{h1cas}. Let $R$ be the companion matrix of $f$ and let $\varepsilon = -1$ if $B$ is alternating, $\varepsilon=1$ otherwise. Let $X$ be the block diagonal sum of $m$ copies of $R$. Using Algorithm 1 we find $T$ such that $R^*=T^{-1}R^{-1}T$, and by Claim 1 in the proof of Proposition \ref{h1cas}(iii), we can choose $T$ such that $T = \varepsilon T^*$. Let $\mathcal{T}$ be the block diagonal sum of $m$ copies of $T$. The matrix $H = B\mathcal{T}^{-1}$ lies in the centralizer of $X$, so it is the embedding into $\GL_{md}(F)$ of $\widetilde{H} \in \GL_m(E)$. By Claim 2 in the proof of Proposition \ref{h1cas}(iii), $\widetilde{H}$ is hermitian and, if $Y \in C_{\C(B)}(X)$, then $Y$ is the embedding into $\GL_{md}(F)$ of some $\widetilde{Y} \in \C(\widetilde{H})$. 
So, if $\widetilde{Y}_1,\widetilde{Y}_2$ are standard 
generators for the unitary group $\C(\widetilde{H})$ 
and $Y_1,Y_2$ are their embeddings into $\GL_{md}(F)$, 
then $C_{\C(B)}(X) = \langle Y_1, Y_2 \rangle$.

If $Q$ is a quadratic form in even characteristic and $X \in \C(Q)$ is semisimple, then the construction of the generators of $C_{\C(Q)}(X)$ is similar to the sesquilinear case: if $f \in \Phi_1$, then we take standard generators 
for $\C(Q)$; 
if $f \in \Phi_2 \cup \Phi_3$, then Proposition \ref{h1caso2} shows that 
we can repeat the argument of the sesquilinear case by replacing the quadratic form $Q$ with the corresponding bilinear form $\b_Q$.

\subsection{Special and Omega groups}\label{specocent}
Let $\b$ be a sesquilinear or quadratic form on $V$, and 
let $\C = \C(\b)$ be the isometry group, with corresponding special 
group $\Spec$ and Omega group $\O$ (in the orthogonal case). 
For a semisimple element $X$ of $\Spec$ or $\O$, 
we construct a generating set for $C_{\C}(X)$
using the above algorithm.
As in the unipotent case discussed in Section \ref{gencentgood}, we then
use Schreier's algorithm to construct 
generating sets for $C_\Spec(X)$ and $C_\O(X)$. 

\section{Constructing a conjugating element} \label{ConjElSemi}
Let $V$ be a vector space over $F=\F_{q^u}$,  and let $\C = \C(\b)$ or $\C(Q)$ 
be the isometry group on $V$. Given conjugate semisimple elements 
$X$ and $Y$ of $\C$, we give an algorithm to compute $Z \in \C$ such that $X^Z= Z^{-1}XZ=Y$. It uses Algorithms 1 and 2 from Section \ref{centgens}.
Then we show how to extend the algorithm to address conjugation in 
$\Spec$ and $\Omega$.

We first handle the sesquilinear form. Let $B$ be the matrix of a non-degenerate sesquilinear form on $V$ and let $\C=\C(B)$. Let $X$ and $Y$ 
be conjugate semisimple elements of $\C$. Let $f_1, \dots, f_h \in \Phi$ be the \geds of $X$ and $Y$, let $m_i$ be their multiplicities and let $d_i = \deg{f_i}$. The algorithm to compute $Z \in \C$ such that $Z^{-1}XZ=Y$ is the following. 
\begin{enumerate}

\item Use Algorithm 1 to compute matrices $P_X$ and $P_Y$ in $\GL(V)$ such 
that $P_XXP_X^{-1} = P_YYP_Y^{-1} =J$, where
\begin{eqnarray}          \label{Jstdform}
J = \begin{pmatrix}
J_1 && \\ & \ddots & \\ && J_h
\end{pmatrix}
\end{eqnarray}
and $J_i$ is the matrix of the restriction of $J$ to the generalized eigenspace $\ker(f_i(J))$ for every $i$. Note that $J$ is not necessarily the Jordan form of $X$: we can choose each $J_i$ completely freely. Hence, $J$ is an isometry for the sesquilinear forms with matrices $B_X = P_XBP_X^*$ and $B_Y = P_YBP_Y^*$.
\item Compute $W \in C_{\GL(V)}(J)$ such that $WB_YW^* = B_X$. 
\item Let $Z = P_X^{-1}WP_Y$. Now $Z \in \C(B)$ and $Z^{-1}XZ=Y$, so $Z$ is the desired conjugating element.
\end{enumerate}
We now explain Step 2. Algorithm 2 does not return an element of $C_{\GL(V)}(J)$ in general, so more work is needed. By Lemma \ref{L26w}, the forms $B_X$ and $B_Y$ have block diagonal shape
\begin{eqnarray*}
B_X = \begin{pmatrix}
B_{X,1} && \\ & \ddots & \\ && B_{X,h}
\end{pmatrix}, \quad B_Y = \begin{pmatrix}
B_{Y,1} && \\ & \ddots & \\ && B_{Y,h}
\end{pmatrix},
\end{eqnarray*}
and each $J_i$ is an isometry for $B_{X,i}$ and $B_{Y,i}$. So we need to find $W_i$ in the centralizer of $J_i$ such that $W_iB_{Y,i}W_i^*=B_{X,i}$ for every $i$, and take $W = W_1 \oplus \cdots \oplus W_h$. We distinguish the three cases.
\begin{itemize}
\item $f_i \in \Phi_1$. Here $J_i$ is a scalar matrix, so Algorithm 2 returns an element of the centralizer of $J_i$.
\item $f_i \in \Phi_2$, $f_i = g_i{g_i}^*$. If we put
$$
J_i = \begin{pmatrix}
\widehat{J}_i & \\ & \widehat{J}_i^*
\end{pmatrix},
$$
where $\widehat{J}_i$ is the restriction of $J_i$ to $\ker{g_i(J)}$, then, by Lemma \ref{L26w},
$$
B_{X,i} = \begin{pmatrix*}[l]
\mathbb{O} & A_{X,i} \\ \varepsilon A_{X,i}^* & \mathbb{O}
\end{pmatrix*}, \quad B_{Y,i} = \begin{pmatrix*}[l]
\mathbb{O} & A_{Y,i} \\ \varepsilon A_{Y,i}^* & \mathbb{O}
\end{pmatrix*},
$$
for some $A_{X,i},A_{Y,i}$ in the centralizer of $\widehat{J}_i$, where $\varepsilon= -1$ if $B$ is alternating, $1$ otherwise. We take
$$
W_i = \begin{pmatrix*}[l]
A_{X,i}A_{Y,i}^{-1} & \mathbb{O} \\ \mathbb{O} & \mathbb{I}
\end{pmatrix*}.
$$
\item $f_i \in \Phi_3$. Let $E = F[t] /(f_i)$ in the symplectic and orthogonal cases, $E = \F_q[t]/(f_i\bar f_i)$ in the unitary case. We follow the argument in the  proof of part (iii) of Proposition \ref{h1cas}.
Let $R$ be the companion matrix of $f_i$ and let $\varepsilon = -1$ 
if $B$ is alternating, $\varepsilon=1$ otherwise. 
Let $J_i$ be the block diagonal sum of $m_i$ copies of $R$. 
Using Algorithm 1, we find $T$ such that $R^*=T^{-1}R^{-1}T$, 
and, by Claim 1 in the proof of Proposition \ref{h1cas}(iii),
we can choose $T$ such 
that $T = \varepsilon T^*$. Let $\mathcal{T}$ be the block diagonal sum of $m_i$ copies of $T$. The matrices $H_{X,i} = B_{X,i}\mathcal{T}^{-1}$ and $H_{Y,i} = B_{Y,i}\mathcal{T}^{-1}$ lie in the centralizer of $J_i$, so they are the embeddings into $\GL_{m_id_i}(F)$ of matrices $\widetilde{H}_{X,i}$ and $\widetilde{H}_{Y,i}$ in $\GL_{m_i}(E)$. These two matrices are hermitian and, using Algorithm 2, we find $\widetilde{W}_i \in \GL_{m_i}(E)$ such that
\begin{eqnarray}  \label{widetildeW}
\widetilde{H}_{X,i} = \widetilde{W}_i\widetilde{H}_{Y,i}\widetilde{W}_i^{\dag},
\end{eqnarray}
where $\widetilde{W}_i \mapsto \widetilde{W}_i^{\dag}$ is defined as in (\ref{dag-symbol}). By Claim 2 in the proof of 
Proposition \ref{h1cas}(iii), this implies $B_{X,i} = W_iB_{Y,i}W_i^*$, where $W_i$ is the embedding of $\widetilde{W_i}$ into $\GL_{m_id_i}(F)$.
\end{itemize}
This completes the algorithm in the case where $\C = \C(\b)$. 

Now suppose that $\C = \C(Q)$, where $Q$ is a non-degenerate quadratic form. 
If $f_i \in \Phi_1$, then we apply directly Algorithm 2; 
if $f_i \in \Phi_2 \cup \Phi_3$, then 
Proposition \ref{h1caso2} shows that 
we can replace 
$Q$ by the corresponding sesquilinear form $\beta_Q$ and 
apply the above algorithm.

\subsection{Special and Omega groups} \label{spogps}
Let $\Spec$ and $\O$ be the special and Omega groups corresponding to 
the isometry group $\C$. We show how to find conjugating elements 
in $\Spec$ and $\O$. 

In the unitary case, let $X,Y \in \Spec$ be conjugate in $\C$; 
by Proposition \ref{speco}, they are also conjugate in $\Spec$. 
Using the above algorithm we find $Z \in \C$ such that $X^Z = Y$. 
Let $\l = \det(Z)$. As shown in the proof of Proposition \ref{speco}(i), 
we can find an element $Z_0$ of determinant $\l^{-1}$ in one of the factors 
of $C_\C(X)$. Then $ZZ_0$ conjugates $X$ to $Y$ and lies in $\Spec$, as required.

In the orthogonal case, consider first the special group $\Spec$ with $q$ odd. 
Let $X,Y \in \Spec$ be $\Spec$-conjugate. We compute $Z \in \C$ such that $X^Z=Y$. If $c_S = 2$ (using the notation of Proposition \ref{speco}(ii)), 
then necessarily $Z \in \Spec$; if $c_S=1$, then we adjust $Z$ by an 
element $Z_0$ of a factor of $C_\C(X)$ lying in the appropriate coset 
of $\C \setminus \Spec$ such that $ZZ_0 \in \Spec$. 

The case of $\O$ is handled similarly. Let $X,Y \in \O$ 
be $\O$-conjugate. Compute $Z \in \Spec$ such that $X^Z=Y$. 
If $c_\O = 2$ (using the notation of Proposition \ref{speco}(iii)), 
then necessarily $Z \in \O$; if $c_\O=1$, then we adjust $Z$ by an 
element $Z_0$ of a factor of $C_\Spec(X)$ lying in the appropriate coset 
of $\Spec \setminus \O$ such that $ZZ_0 \in \O$. 

\section{Some examples}
We close this chapter by giving tables listing the semisimple class representatives 
and centralizers for the groups $\Sp_6(3)$, $\SU_6(2)$, $\O_8^+(2)$ and $\O_8^-(2)$.
The representatives are identified using the notation
of Section \ref{formsandelements}.
The structures of the centralizers are given by
Theorem \ref{CentrSem}.

To keep the notation concise, we adopt the following conventions in the tables. In the first column, the superscripts in brackets denote the multiplicities of the elementary divisors; also the elementary divisors in the set $\Phi_2$ are of the form $f = gg^*$, and we give the polynomial $g$ rather than $f$. In the second column, $f_i \in \Phi_i$ refers to the corresponding polynomial from the first column, and $X_{f_i}$ denotes the representative defined in Section \ref{formsandelements}. Again, the superscripts in brackets denote the multiplicities of the diagonal blocks $X_{f_i}$.

Note that in Table \ref{O+82}, there are several classes in $\Or_8^+(2)$ that split into two classes in $\O_8^+(2)$ (see 
Proposition \ref{speco}(iii)); we have indicated these classes in the first column of the table, and give one representative 
in the second column.
The same applies to Table \ref{O-82}. 

Finally, in Table \ref{SU62} we denote by $\o$ a primitive element of $\F_4$.
In the third column, for convenience, we give the centralizer in $\GU_6(2)$ 
rather than $\SU_6(2)$; the latter centralizer has index 3 in the given centralizer, by 
Proposition \ref{speco}(i). 

\begin{table}[htbp]
\begin{center}
\begin{tabular}{|l|l|l|}
\hline
\textrm{Elementary divisors} & Representative  & Centralizer \\ \hline 
$(t + 1)^{[6]}$ &  $X_{f_1}^{[6]}$  & $\Sp_6(3)$ \\ \hline 

$(t + 2)^{[6]}$ &  $X_{f_1}^{[6]}$  & $\Sp_6(3)$ \\ \hline

$(t + 1)^{[4]}, (t + 2)^{[2]}$ & $X_{f_1}^{[4]} \oplus X_{f_1}^{[2]}$ & 
   $\Sp_4(3) \times \Sp_2(3)$\\ \hline 

$(t + 2)^{[4]}, (t + 1)^{[2]}$ & $X_{f_1}^{[4]} \oplus X_{f_1}^{[2]}$ &  
   $\Sp_4(3) \times \Sp_2(3)$ \\ \hline

$(t + 2)^{[4]}, (t^2 + 1)$  & $X_{f_1}^{[4]} \oplus X_{f_3}$ & 
   $\Sp_4(3) \times \GU_1(3)$ \\ \hline

$(t + 1)^{[4]}, (t^2 + 1)$ & $X_{f_1}^{[4]} \oplus X_{f_3}$ 
   & $\Sp_4(3) \times \GU_1(3)$ \\ \hline

$(t + 1)^{[2]}, (t^2 + 1)^{[2]}$  & $X_{f_1}^{[2]} \oplus X_{f_3}^{[2]}$  & 
   $\Sp_2(3) \times \GU_2(3)$  \\ \hline  

$(t + 2)^{[2]}, (t^2 + 1)^{[2]}$ & $X_{f_1}^{[2]} \oplus X_{f_3}^{[2]}$ & 
   $\Sp_2(3) \times \GU_2(3)$ \\ \hline

$(t + 1)^{[2]}, (t^2 + 2t+ 2)$ & $X_{f_1}^{[2]} \oplus X_{f_2}$  & 
   $\Sp_2(3) \times \GL_1(3^2)$ \\ \hline

$(t + 2)^{[2]}, (t^2 + 2t + 2)$ & $X_{f_1}^{[2]} \oplus X_{f_2}$ & 
   $\Sp_2(3) \times \GL_1(3^2)$ \\ \hline

$(t + 1)^{[2]}, (t^4 + t^3 + t^2 + t + 1)$ & $X_{f_1}^{[2]} \oplus X_{f_3}$ & 
   $\Sp_2(3) \times \GU_1(3^2)$ \\ \hline

$(t + 2)^{[2]}, (t^4 + t^3 + t^2 + t + 1)$ &  $X_{f_1}^{[2]} \oplus X_{f_3}$ & 
   $\Sp_2(3) \times \GU_1(3^2)$ \\ \hline

$(t + 1)^{[2]}, (t^4 + 2t^3 + t^2 + 2t + 1)$  & $X_{f_1}^{[2]} \oplus X_{f_3}$ & 
   $\Sp_2(3) \times \GU_1(3^2)$ \\ \hline

$(t + 2)^{[2]}, (t^4 + 2t^3 + t^2 + 2t + 1)$ & 
    $X_{f_1}^{[2]} \oplus X_{f_3}$ & $\Sp_2(3) \times \GU_1(3^2)$ \\ \hline

$(t + 2)^{[2]}, (t + 1)^{[2]}, (t^2 + 1)$ & 
   $X_{f_1}^{[2]} \oplus X_{f_1}^{[2]} \oplus X_{f_3}$ & 
    $\Sp_2(3) \times \Sp_2(3) \times \GU_1(3)$ \\ \hline

$(t^2 + 1)^{[3]}$ & $X_{f_3}^{[3]}$ & $\GU_3(3)$ \\ \hline

$(t^2 + 2t + 2), (t^2 + 1)$  & $X_{f_2} \oplus X_{f_3}$ & 
   $\GL_1(3^2) \times \GU_1(3)$ \\ \hline

$(t^2 + 1), (t^4 + t^3 + t^2 + t + 1)$ & $X_{f_3} \oplus X_{f_3}$ & 
   $\GU_1(3) \times \GU_1(3^2)$ \\ \hline

$(t^2 + 1), (t^4 + 2t^3 + t^2 + 2t + 1)$ & $X_{f_3} \oplus X_{f_3}$ & 
   $\GU_1(3) \times \GU_1(3^2)$ \\ \hline

$t^3 + 2t + 1$ & $X_{f_2}$ & $\GL_1(3^3)$ \\ \hline

$t^3 + 2t^2 + 2t + 2$ & $X_{f_2}$ & $\GL_1(3^3)$ \\ \hline

$t^3 + 2t^2 + t + 1$ & $X_{f_2}$ & $\GL_1(3^3)$ \\ \hline

$t^3 + 2t + 2$ & $X_{f_2}$ & $\GL_1(3^3)$ \\ \hline

$t^6 + t^5 + t^3 + t + 1$ & $X_{f_3}$ & $\GU_1(3^3)$ \\ \hline

$t^6 + t^5 + t^4 + t^3 + t^2 + t + 1$ & $X_{f_3}$ & 
$\GU_1(3^3)$ \\ \hline

$t^6 + 2t^5 + 2t^3 + 2t + 1$ & $X_{f_3}$ & 
$\GU_1(3^3)$ \\ \hline

$t^6 + 2t^5 + t^4 + 2t^3 + t^2 + 2t + 1$ & $X_{f_3}$  & 
$\GU_1(3^3)$ \\ \hline
\end{tabular}
\end{center}
\caption{Semisimple classes for $\Sp_6(3)$} \label{SpSemi}
\end{table}


\begin{table}[htbp]
\begin{center}
\begin{tabular}{|l|l|l|}
\hline
\textrm{Elementary divisors} & Representative  & Centralizer in $\GU_6(2)$ \\ \hline


$( t + 1)^{[ 6 ]}$  & 
$X_{f_1}^{[6]}$ & 
$\GU_6(2)$ \\ \hline 


$( t + \omega)^{[ 6 ]} $ & $X_{f_1}^{[6]}$ & $\GU_6(2)$ \\ \hline 


$( t + \omega^2)^{[ 6 ]} $ 
& $X_{f_1}^{[6]}$ & $\GU_6(2)$ \\ \hline 


$( t + 1)^{[4]} , ( t + \omega) , ( t + \omega^2) $ & 
$X_{f_1}^{[4]} \oplus X_{f_1} \oplus X_{f_1}$ & 
$\GU_4(2) \times \GU_1(2) \times \GU_1(2)$ \\ \hline 


$( t + 1) , ( t + \omega)^{[4]}, ( t + \omega^2) $ & 
$X_{f_1} \oplus X_{f_1}^{[4]} \oplus X_{f_1}$ &  
$\GU_1(2) \times \GU_4(2) \times \GU_1(2)$ \\ \hline 


$( t + 1), ( t + \omega) , ( t + \omega^2)^{[4]} $ & 
$X_{f_1} \oplus X_{f_1} \oplus X_{f_1}^{[4]}$ & 
$\GU_4(2) \times \GU_1(2) \times \GU_1(2)$ \\ \hline 


$( t + 1)^{[3]} , ( t + \omega)^{[3]} $ & 
$X_{f_1}^{[3]} \oplus X_{f_1}^{[3]}$ & 
$\GU_3(2) \times \GU_3(2)$ \\ \hline 


$( t + 1)^{[3]} , ( t + \omega^2)^{[3]} $ & 
$X_{f_1}^{[3]} \oplus X_{f_1}^{[3]}$ & 
$\GU_3(2) \times \GU_3(2)$ \\ \hline 


$( t + \omega)^{[3]} , ( t + \omega^2)^{[3]} $ & 
$X_{f_1}^{[3]} \oplus X_{f_1}^{[3]}$ & 
$\GU_3(2) \times \GU_3(2)$ \\ \hline


$( t + 1)^{[ 2 ]} , ( t + \omega)^{[ 2 ]} , ( t + \omega^2)^{[ 2 ]} $ & 
$X_{f_1}^{[2]} \oplus X_{f_1}^{[2]} \oplus X_{f_1}^{[2]} $ 
& 
$\GU_2(2) \times \GU_2(2) \times \GU_2(2)$ \\ \hline

$( t + 1)^{[ 2 ]} , ( t + \omega), ( t^3 + \omega^2) $ & 
$X_{f_1}^{[2]} \oplus X_{f_1} \oplus X_{f_3}$ & 
$\GU_2(2) \times \GU_1(2) \times \GU_1( 2^3)$ \\ \hline 


$( t + \omega)^{[ 2 ]} , ( t + \omega^2), ( t^3 + \omega^2) $ & 
$X_{f_1}^{[2]} \oplus X_{f_1} \oplus X_{f_3}$ & 
$\GU_2(2) \times \GU_1(2) \times \GU_1(2^3)$ \\ \hline 


$( t + 1) , ( t + \omega)^{[ 2 ]} , ( t^3 + \omega) $ & 
$X_{f_1} \oplus X_{f_1}^{[2]} \oplus X_{f_3} $ & 
$\GU_1(2) \times \GU_2(2) \times \GU_1(2^3)$ \\ \hline 


$( t + 1)^{[ 2 ]} , ( t + \omega^2), ( t^3 + \omega) $ & 
$X_{f_1}^{[2]} \oplus X_{f_1} \oplus X_{f_3} $ & 
$\GU_2(2) \times \GU_1(2) \times \GU_1(2^3)$ \\ \hline 


$( t + 1) , ( t + \omega^2)^{[ 2 ]} , ( t^3 + \omega^2) $ & 
$X_{f_1} \oplus X_{f_1}^{[2]} \oplus X_{f_3} $ & 
$\GU_1(2) \times \GU_2(2) \times \GU_1(2^3)$ \\ \hline 


$( t + \omega), ( t + \omega^2)^{[ 2 ]} , ( t^3 + \omega) $ & 
$X_{f_1} \oplus X_{f_1}^{[2]} \oplus X_{f_3} $ & 
$\GU_1(2) \times \GU_2(2) \times \GU_1(2^3)$ \\ \hline 


$( t + 1), ( t + \omega), ( t^2 + t + \omega) $ & 
$X_{f_1} \oplus X_{f_1} \oplus X_{f_2}$ &
$\GU_1(2) \times \GU_1(2) \times \GL_1(2^4)$ \\ \hline 


$( t + \omega)^{[ 2 ]} , ( t^2 + t + \omega^2) $ & 
$X_{f_1}^{[2]} \oplus X_{f_2} $ &
$\GU_2(2) \times \GL_1(2^4)$ \\ \hline 


$( t + \omega), ( t + \omega^2), 
( t^2 + \omega t + 1) $& 
$X_{f_1} \oplus X_{f_1} \oplus X_{f_2}  $ & 
$\GU_1(2) \times \GU_1(2) \times \GL_1(2^4)$ \\ \hline 


$( t + 1) , ( t + \omega^2), ( t^2 + t + \omega^2) $ & 
$X_{f_1} \oplus X_{f_1} \oplus X_{f_2}$  & 
$\GU_1(2) \times \GU_1(2) \times \GL_1(2^4)$ \\ \hline 


$( t + 1)^{[ 2 ]} , ( t^2 + \omega t + 1) $ & 
$X_{f_1}^{[2]} \oplus X_{f_2}  $ & 
$\GU_2(2) \times \GL_1(2^4)$ \\ \hline 


$( t + \omega^2)^{[ 2 ]} , ( t^2 + t + \omega) $ & 
$X_{f_1}^{[2]} \oplus X_{f_2}  $ & 
$\GU_2(2) \times \GL_1(2^4)$ \\ \hline 


$( t + \omega) , 
( t^5 + \omega^2 t^4 + \omega^2 t^3 + t^2 + t + \omega^2) $ & 
$X_{f_1} \oplus X_{f_3} $ & 
$\GU_1(2) \times \GU_1(2^5)$ \\ \hline 


$( t + \omega) ,
( t^5 + t^4 + \omega^2  t^3 + t^2 + \omega^2  t+ \omega^2) $ & 
$X_{f_1} \oplus X_{f_3} $ & 
$\GU_1(2) \times \GU_1(2^5)$ \\ \hline 


$( t + \omega^2), ( t^5 + \omega t^4 + \omega t^3 + t^2 + t + \omega) $ & 
$X_{f_1} \oplus X_{f_3} $ & 
$\GU_1(2) \times \GU_1(2^5)$ \\ \hline 


$( t + \omega^2), ( t^5 + t^4 + \omega t^3 + t^2 + \omega t + \omega) $ & 
$X_{f_1} \oplus X_{f_3} $ & 
$\GU_1(2) \times \GU_1(2^5)$ \\ \hline 


$( t + 1), ( t^5 + \omega^2 t^4 + t^3 + t^2 + \omega t + 1) $ & 
$X_{f_1} \oplus X_{f_3} $ & 
$\GU_1(2) \times \GU_1(2^5)$ \\ \hline 


$( t + 1), ( t^5 + \omega t^4 + t^3 + t^2 + \omega^2 t + 1) $ & 
$X_{f_1} \oplus X_{f_3} $ & 
$\GU_1(2) \times \GU_1(2^5)$ \\ \hline 


$ t^3 + \omega^2 t + 1 $ & 
$X_{f_2} $ &
$\GL_1(2^6)$ \\ \hline 


$ t^3 + \omega t+ 1 $ & 
$X_{f_2} $ &
$\GL_1(2^6)$ \\ \hline 


$ t^3 + t + 1 $& 
$X_{f_2} $ & 
$\GL_1(2^6)$ \\ \hline 


$( t^3 + \omega), ( t^3 + \omega^2) $ & 
$X_{f_3} \oplus X_{f_3}$ & 
$\GU_1(2^3) \times \GU_1(2^3)$ \\ \hline 

\end{tabular}
\end{center}
\caption{Semisimple classes for $\SU_6(2)$} \label{SU62}
\end{table}


\begin{table}[htbp]
\begin{center}
\begin{tabular}{|l|l|l|}
\hline
\textrm{Elementary divisors} & Representative  & Centralizer \\ \hline 

$(t + 1)^{[8]}$ & $X_{f_1}^{[8]} $ & 
   $\Omega^+_8(2)$ \\ \hline

$(t + 1)^{[6]}, (t^2 + t+ 1)$ & $X_{f_1}^{[6]} \oplus X_{f_3}$  & 
   $\Omega^{-}_6(2) \times \GU_1(2)$ \\ \hline

$(t + 1)^{[4]}, (t^2 + t + 1)^{[2]}$ & 
$X_{f_1}^{[4]} \oplus X_{f_3}^{[2]}$ 
   & $\Omega^+_4(2) \times \GU_2(2)$ \\ \hline

$(t + 1)^{[4]}, (t^4 + t^3 + t^2 + t + 1)$  & $X_{f_1}^{[4]} 
\oplus X_{f_3}$  & 
   $\Omega^{-}_4(2) \times \GU_1(2^2)$  \\ \hline  

$(t + 1)^{[2]}, (t^2 + t + 1), (t^4 + t^3 + t^2 + t + 1)$ & 
$X_{f_1}^{[2]} \oplus X_{f_3} \oplus X_{f_3}$ 
   & $\GU_1(2) \times \GU_1(2^2)$ \\ \hline

$(t + 1)^{[2]}, (t^2 + t + 1)^{[3]}$ & 
$X_{f_1}^{[2]} \oplus X_{f_3}^{[3]}$ 
   & $\Omega^-_2(2) \times \GU_3(2)$ \\ \hline

$(t + 1)^{[2]}, (t^3 +t + 1) $ & $X_{f_1}^{[2]} \oplus X_{f_2}$
& $\GL_1(2^3)$ \\ \hline

$(t + 1)^{[2]}, (t^6 + t^3 + 1)$ & $X_{f_1}^{[2]} \oplus X_{f_3}$ &  
   $\Omega^{-}_2(2) \times \GU_1(2^3)$ \\ \hline

$(t^2 + t + 1)^{[4]}$ (two classes) & $X_{f_3}^{[4]}$ & $\GU_4(2)$ \\ \hline

$(t^2 + t + 1), (t^6 + t^3+ 1)$ (two classes) & $X_{f_3} \oplus X_{f_3}$  & 
   $\GU_1(2) \times \GU_1(2^3)$ \\ \hline

$(t^4 + t^3 + t^2 + t + 1)^{[2]}$ (two classes) &  $X_{f_3}^{[2]}$  & $\GU_2(2^2)$ \\ \hline 

$(t^4 + t + 1)$ (two classes) & $X_{f_2} $  & 
   $\GL_1(2^4)$ \\ \hline

\end{tabular}
\end{center}
\caption{Semisimple classes for $\Omega^+_8(2)$} \label{O+82}
\end{table}

\begin{table}[htbp]
\begin{center}
\begin{tabular}{|l|l|l|}
\hline
\textrm{Elementary divisors} & Representative  & Centralizer \\ \hline 

$(t + 1)^{[8]}$ & $X_{f_1}^{[8]} $ & $\Omega^-_8(2)$ \\ \hline

$(t + 1)^{[6]}, (t^2 + t+ 1)$ & $X_{f_1}^{[6]} \oplus X_{f_3}$  & 
   $\Omega^{+}_6(2) \times \GU_1(2)$ \\ \hline

$(t + 1)^{[4]}, (t^2 + t + 1)^{[2]}$ & 
$X_{f_1}^{[4]} \oplus X_{f_3}^{[2]}$ 
   & $\Omega^-_4(2) \times \GU_2(2)$ \\ \hline

$(t + 1)^{[4]}, (t^4 + t^3 + t^2 + t + 1)$  & $X_{f_1}^{[4]} 
\oplus X_{f_3}$  & 
   $\Omega^{+}_4(2) \times \GU_1(2^2)$  \\ \hline  

$(t + 1)^{[2]}, (t^2 + t + 1)^{[3]}$ & 
$X_{f_1}^{[2]} \oplus X_{f_3}^{[3]}$ 
   & $\GU_3(2)$ \\ \hline

$(t + 1)^{[2]}, (t^3 +t + 1) $ & $X_{f_1}^{[2]} \oplus X_{f_2}$
& $\Omega^-_2(2) \times \GL_1(2^3)$ \\ \hline

$(t + 1)^{[2]}, (t^6 + t^3 + 1)$ & $X_{f_1}^{[2]} \oplus X_{f_3}$ &  
   $\GU_1(2^3)$ \\ \hline

$(t + 1)^{[2]}, (t^2 + t + 1), (t^4 + t^3 + t^2 + t + 1)$ & 
$X_{f_1}^{[2]} \oplus X_{f_3} \oplus X_{f_3}$ 
   & $\Omega^{-}_2(2) \times \GU_1(2) \times \GU_1(2^2)$ \\ \hline

$(t^2 + t + 1)^{[2]}, (t^4 + t^3 + t^2 + t + 1)$ (two classes) &  $X_{f_3}^{[2]} \oplus
X_{f_3}$ & $\GU_2(2) \times \GU_1(2^2)$ \\ \hline 

$(t^2 + t + 1), (t^3 +t + 1)$ (two classes)  & $X_{f_3} \oplus X_{f_2}
$ & $ \GU_1(2) \times \GL_1(2^3)$ \\ \hline

$t^8 +  t^5 + t^4 + t^3 + 1$ (two classes) & $X_{f_3}$ & $\GU_1(2^4)$ \\ \hline

$t^8 +  t^7 + t^6 + t^4 + t^2 + t  + 1$ (two classes) & $X_{f_3}$ & $\GU_1(2^4)$ \\ \hline

\end{tabular}
\end{center}
\caption{Semisimple classes for $\Omega^-_8(2)$} \label{O-82}
\end{table}

\chapter{General conjugacy classes}        \label{generalChap}

In this chapter we use the results in the semisimple and unipotent cases 
to solve the main conjugacy problems (1)-(3) of Section \ref{mainprob} for classical groups: 
list representatives for all conjugacy classes; describe the structure
of the centralizer of an arbitrary element and provide a generating set; 
decide conjugacy between arbitrary elements 
and construct explicit conjugating elements.

\section{Conjugacy classes and centralizers}   \label{SectListGen}

Let $F=\mathbb{F}_{q^u}$, where $u=2$  in the unitary case and $u=1$ otherwise, 
and let $V$ be an $n$-dimensional vector space over $F$. 
Let $\C$ be the isometry group $\C(\beta)$ or $\C(Q)$, where $\beta$ is a non-degenerate alternating, symmetric or hermitian 
form and $Q$ is a non-degenerate quadratic form on $V$.

Let $x \in \C$. Recall the Jordan decomposition $x=su=us$, with $s$ semisimple and $u$ unipotent. 
We know that $C_\C(x) = C_\C(s) \cap C_\C(u)$. Moreover, 
$x_1 = s_1u_1$ and $x_2=s_2u_2$  are conjugate in $\C$ 
if and only if $s_1$ and $s_2$ are conjugate in $\C$ and $z^{-1}u_1z$ and $u_2$ are conjugate in $C_\C(s_2)$, where $z \in \C$ is such that $z^{-1}s_1z = s_2$.
Thus our strategy to list all conjugacy classes of $\C$ is to list all 
semisimple classes and, for each semisimple representative, list all unipotent classes 
in its centralizer. 

For the semisimple classes of $\C$, we use variants of the representatives produced in Section \ref{formsandelements}. Recall the notation of Definition \ref{ThreeCases}:
\begin{eqnarray*}
	\Phi_1 & = & \{f: f \in F[t] \; | \; f=f^* \mbox{ monic irreducible}, \, \deg{f}=1 \}, \\
	\Phi_2 & = & \{f: f \in F[t] \; | \; f=gg^*, \, g \neq g^* \mbox{ monic irreducible}  \}, \\ 
	\Phi_3 & = & \{f: f \in F[t] \; | \; f=f^* \mbox{ monic irreducible}, \, \deg{f}>1 \},
\end{eqnarray*}
and $\Phi = \Phi_1 \cup \Phi_2 \cup \Phi_3$.
Each semisimple class of $\C$ can be identified with a pair $(S,B)$, where $S$ is an isometry for the form $B$. Let $f_1, \dots, f_h \in \Phi$ be the distinct \geds of $S$, and let $m_i$ be the multiplicity of $f_i$ for $i=1, \dots, h$. 
Thus, with respect to a suitable basis,
\begin{eqnarray}  \label{matricessemisimple}
S = \left( \begin{array}{ccc} S_1 && \\ & \ddots & \\ && S_h\end{array}\right), \quad B = \left(\begin{array}{ccc} B_1 && \\ & \ddots & \\ && B_h \end{array}\right),
\end{eqnarray}
where $S_i$ and $B_i$ are the matrices of the restrictions of $S$ and $B$ to $\ker(f_i(S))$.

We first describe representatives $S$, then give the corresponding $B$. For $f_i \in \Phi_1 \cup \Phi_3$, the matrix $S_i$ is a block diagonal sum of $m_i$ companion matrices $C(f_i)$; for $f_i = g_ig_i^* \in \Phi_2$, we take
$$
S_i = \left( \begin{array}{cc} Y_i & \\ & Y_i^{*-1} \end{array}\right),
$$
where $Y_i$ is the block diagonal sum of $m_i$ companion matrices $C(g_i)$. For the matrix $B$, we choose the following form:
\begin{itemize}
	\item If $f_i \in \Phi_1$, then $S_i$ is a scalar matrix, so $B_i$ can be chosen arbitrarily. Note that $C_{\C(B_i)}(S_i) = \C(B_i)$.
	\item If $f_i \in \Phi_2$, then choose
	$$
	B_i = \left(\begin{array}{cc} \mathbb{O} & \mathbb{I} \\ \varepsilon\mathbb{I} & \mathbb{O} \end{array}\right), \quad \varepsilon = \left\{ \begin{array}{rl} 1 & \mbox{ if $B$ is hermitian or symmetric}; \\ -1 & \mbox{ if $B$ is alternating}; \\ 0 & \mbox{ if $B$ is quadratic}. \end{array}\right.
	$$
Note that $C_{\C(B_i)}(S_i) \cong \GL_{m_i}(q^{d_i})$, where $d_i = u\deg (f_i)/2$ (see Theorem \ref{CentrSem}).
	\item If $f_i \in \Phi_3$, then choose
	\begin{eqnarray}  \label{choiceforBi}
	B_i = \left( \begin{array}{ccc} && B_{f_i} \\ & \iddots & \\ B_{f_i} \end{array}\right),
	\end{eqnarray}
	where $B_{f_i}$ is the matrix of a form preserved by the companion matrix $C(f_i)$, computed as discussed at the end of Section \ref{formsandelements}, and it appears $m_i$ times. Note that $C_{\C(B_i)}(S_i) \cong \GU_{m_i}(q^{d_i})$ where $d_i = u \deg(f_i)/2$ (see Theorem \ref{CentrSem}).
\end{itemize}
Let $x=su=us \in \C$, with $s$ semisimple and $u$ unipotent. Choose a basis such that $s$ and the form ($\beta$ or $Q$) have matrices $S$ and $B$ respectively as in (\ref{matricessemisimple}). Let $U$ be the matrix of $u$ with respect to
this basis. We know that $U$ belongs to $C_{\C}(S) = \prod_{i=1}^h C_{\C(B_i)}(S_i)$. Thus $U$ is a block diagonal matrix
$$
\left( \begin{array}{ccc} U_1 && \\ & \ddots & \\ && U_h \end{array}\right)
$$
where $U_i$ is the matrix of the restriction of $u$ to $\ker(f_i(S))$. 

With the above notation, we can now describe the conjugacy classes and centralizers in $\C$. The result follows from the above discussion.

\begin{theorem}    \label{CCgen}
	A complete set of representatives for the conjugacy classes of $\C$ is given by all pairs of matrices $(SU,B)$, defined by 
	$$
	S = \left(\begin{array}{ccc} S_1 && \\ & \ddots & \\ && S_h\end{array}\right), \quad U = \left( \begin{array}{ccc} U_{1j} && \\ & \ddots & \\ && U_{hj} \end{array}\right), \quad B = \left(\begin{array}{ccc} B_1 && \\ & \ddots & \\ && B_h \end{array}\right).
	$$
	where $(S,B)$ runs over all representatives of semisimple 
conjugacy classes of $\C$ and, for each $i$, 
the matrices $U_{ij}$ run over all representatives of unipotent conjugacy 
classes of $\C_i := C_{\C(B_i)}(S_i)$. 
For a given representative $SU$ as above,
\[
C_\C(SU) = \prod_{i=1}^h C_{\C_i}(U_{ij}).
\]
\end{theorem}

The factors $ C_{\C_i}(U_{ij})$ are described in 
Section \ref{Cue} when $\C_i= \GL_{m_i}(q^{d_i})$; 
Sections \ref{spcent} and \ref{centbad} 
when $\C_i$ is symplectic; 
Sections \ref{ocent} and \ref{centbad} 
when $\C_i$ is orthogonal; 
and Section \ref{ucent}
when $\C_i$ is unitary.

\subsection{Listing class representatives in the isometry group}
In Section  \ref{formsandelements} 
we showed how to list the semisimple class representatives in $\C$. 
It remains to list representatives for all conjugacy classes of 
$\C$ with fixed semisimple part $S$. With notation as above, 
this amounts to  listing representatives for all unipotent classes 
of $C_{\C(B_i)}(S_i)$. We distinguish the three cases.
	\begin{itemize}
		\item $f_i \in \Phi_1$. Here $S_i$ is a scalar matrix and $C_{\C(B_i)}(S_i)$ coincides with $\C(B_i)$. The representatives for unipotent classes of $\C(B_i)$ are given in Sections \ref{goodrep} and \ref{badsec}, and we are free to choose the form $B_i$.
		\item $f_i \in \Phi_2$, $f_i=g_i{g_i}^*$. Now $C_{\C(B_i)}(S_i)$ is isomorphic to $\GL_{m_i}(E)$, with $E = F[t]/(g_i)$, via the isomorphism
		$$
		y \mapsto \left(\begin{array}{cc} Y & \\ & Y^{*-1} \end{array}\right), \: \forall y \in \GL_{m_i}(E),
		$$
		where $Y$ is the embedding of $y$ into $\GL_{m_id_i}(F)$. Two elements of $\GL_{m_i}(E)$ are conjugate if and only if they have the same generalized elementary divisors, so the list of representatives of unipotent classes of $C_{\C(B_i)}(S_i)$ is just the list of the images in $C_{\C(B_i)}(S_i)$ of representatives of the unipotent classes of $\GL_{m_i}(E)$. We choose the block diagonal sum of unipotent Jordan blocks as our preferred form.
		\item $f_i \in \Phi_3$. Let $E = F[t]/(f_i)$ in the symplectic and orthogonal cases, and $E = \F_q[t]/(f_i\bar f_i)$ in the unitary case.  By Proposition \ref{h1cas} and its proof, 
$C_{\C(B_i)}(S_i)$ is the  embedding in $\GL_{m_id_i}(F)$ of the 
group $\GU_{m_i}(E)$ preserving the hermitian form with matrix
$$
\left( \begin{array}{ccc} && 1\\ & \iddots & \\ 1 && \end{array} \right).
$$
This follows from our choice of $B_i$ in (\ref{choiceforBi}). 
Representatives for all unipotent classes of $\GU_{m_i}(E)$ are 
the Jordan forms $\bigoplus_{i=1}^k V(m_i)^{[r_i]}$ given in 
Section \ref{sugood}.
\end{itemize}

\subsection{Listing class representatives in special and Omega groups}
We now determine the splitting of
those classes in the isometry
groups $\C = \GU_n(q)$ or $\Or_n^\e(q)$ that lie
in the special groups $\Spec = \SU_n(q)$ or $\SO_n^\e(q)$,
or in the Omega group $\O = \O_n^\e(q)$.
This is a matter of computing the
indices $|C_\C (X): C_\Spec (X)|$ and $|C_\C (X): C_\O (X)|$
for elements $X$ of $\Spec$ and $\O$.

\begin{proposition}\label{specogen}
\mbox{}
\begin{itemize}
\item[{\rm (i)}] Let $\C = \GU_n(q)$ and $\Spec = \SU_n(q)$. 
Let $X \in \Spec$ have generalized elementary divisors 
$f_1^{m_1},\ldots, f_k^{m_k}$, where the $f_i \in \Phi$ are not 
necessarily distinct. Then $|C_\C(X): C_\Spec (X)| = (q+1)/{r}$, 
where $r = \gcd(m_1,\ldots,m_k,q+1)$. 
The class $X^\C$ splits into $r$ classes in $\Spec$.
\item[{\rm (ii)}] Let $\C = \Or_n^\e(q)$ and 
$\Spec = \SO_n^\e(q)$ with $q$ odd, and let $X \in \Spec$. Then
$|C_\C(X): C_\Spec(X)| = c_S$, where 
\[
c_S = \left\{\begin{array}{l} 2, \hbox{ if $X$ has an elementary divisor } (t\pm 1)^m \hbox{ with $m$ odd}, \\
                                          1, \hbox{ otherwise.}
\end{array}
\right.
\]
The class $X^\C$ splits into $2/c_S$ classes in $\Spec$.
\end{itemize}
\end{proposition}

\begin{proof}
\no (i) This is similar to the proof of Theorem \ref{MainResult2}. As in that proof, we can reduce to the case where 
the \geds are $f^{m_1},\ldots, f^{m_k}$, all powers of a given polynomial $f \in \Phi$. Let $X=SU = US$ be the Jordan 
decomposition of $X$, and let $D:=\{\l \in \F_{q^2}: \l^{q+1}=1\} = \la \o\ra$. 

Suppose first that $f \in \Phi_1$, so that $f(t) = t-\l$ with $\l \in D$. Then $C_\C(S) = \C$
and $U \in \C$ is unipotent with Jordan form $\bigoplus_{i=1}^l J_{m_i}$. Relabel the $m_i$ so that this Jordan form is 
$\bigoplus_{i=1}^k J_{\l_i}^{[l_i]}$, where $\l_1,\ldots,\l_k$ are distinct. By Theorem \ref{sureps}, $C_\C(U) = Q_0R$, where $Q_0$ is a normal subgroup of $q$-power order and $R \cong \prod_{i=1}^k \GU_{l_i}(q)$.
The action of $R$ on $V = V_n(q^2)$ is as 
$\prod_{i=1}^k (\GU_{l_i}(q) \otimes I_{\l_i})$ 
(since there is such a subgroup in $C_\C(U)$). 
If $A \in \GU_{l_i}(q)$ has determinant $\o$, then 
the determinant of $A \otimes I_{\l_i}$ is $\o^{\l_i}$.
Hence the subgroup of $\F_{q^2}^*$ consisting of all determinants of 
elements of $R$ is generated by $\o^{\l_i}$ for $i=1,\ldots,k$. This equals 
$\la \o^r\ra$ where $r  = \mathrm{gcd}(\lambda_{1}, \ldots, \lambda_{k}, q+1)$. 
Hence the image of the determinant map $C_{\GU_n(q)}(X) \to D$ has order $(q+1)/r$, completing the proof for the case $f \in \Phi_1$.

Next suppose $f \in \Phi_2$, so that $f=gg^*$ with $g \ne g^*$ irreducible of degree $d$. By Theorem \ref{CentrSem}, 
$C_\C(S) \cong \GL_m(q^{2d})$, where $m = \sum m_i$. Relabel the $m_i$ so that the Jordan form of the unipotent element 
$U$ of $\GL_m(q^{2d})$ is $\bigoplus_{i=1}^k J_{\l_i}^{[l_i]}$, where $\l_1,\ldots,\l_k$ are distinct. By Theorem \ref{Card}, 
$C_\C(X) = C_{\GL_m(q^{2d})}(U) = Q_0R$, where $R \cong \prod_{i=1}^k \GL_{l_i}(q^{2d})$ acts on $V = V_m(q^{2d})$ as 
$\prod_{i=1}^k (\GL_{l_i}(q^{2d}) \otimes I_{\l_i})$. There exists $A \in \GL_{l_i}(q^{2d})$ whose determinant as an 
element of $\GU_{dl_i}(q)$ is $\o$, and so the determinant of $A \otimes I_{\l_i}$ is $\o^{\l_i}$. Now we complete the argument as before. 

The final case, where $ f \in \Phi_3$, is similar. Here $f$ is irreducible of odd degree $d>1$, and 
$C_\C(S) \cong \GU_m(q^{d})$ by Theorem \ref{CentrSem}. Again relabel the $m_i$ so that the Jordan form of  
$U \in \GU_m(q^{d})$ is $\bigoplus_{i=1}^k J_{\l_i}^{[l_i]}$, where $\l_1,\ldots,\l_k$ are distinct. Then, by Theorem \ref{sureps}, $C_\C(X) = C_{\GU_m(q^{d})}(U) = Q_0R$, where $R \cong \prod_{i=1}^k \GU_{l_i}(q^d)$ acts on $V = V_m(q^{2d})$ as 
$\prod_{i=1}^k (\GU_{l_i}(q^{d}) \otimes I_{\l_i})$. Now we complete the proof as before.

\vspace{2mm}
\no (ii) Let $X = SU$ be the Jordan decomposition of $X$. For every \ged $f_i \in \Phi_2\cup \Phi_3$, the corresponding factor of $C_\C(S)$ lies in $\Spec$ by Proposition \ref{speco}(ii). So we may suppose that each 
\ged $f_i^{m_i}$  is $(t-1)^{m_i}$ or $(t+1)^{m_i}$. 
If all the $m_i$ are even, then 
$$C_\C(X) \le C_{\GL(V)}(X) \le \SL(V)$$ 
(see Section \ref{SLconjclasses}).
On the other hand, if there is a \ged $(t\pm 1)^m$ with $m$ odd, then,
with respect to a suitable basis, $C_\C(X)$ contains a 
block diagonal matrix $-I_m\oplus I_{n-m}$, 
which has determinant $-1$. 
\end{proof}

The splitting of classes in $\O = \O_n^\e(q)$ is more 
complicated, and is discussed in the next two propositions. 
Let $X \in \Omega$ have Jordan decomposition $SU$ with $S$ semisimple and 
$U$ unipotent. Let $f_1, \dots, f_h$ be the \geds of $S$, with multiplicities 
$m_1, \dots, m_h$. Change the basis so that $X$ and the form $B$ are
$$
X= \left( \begin{array}{ccc} X_1 && \\ & \ddots & \\ && X_h \end{array}\right), \quad B= \left( \begin{array}{ccc} B_1 && \\ & \ddots & \\ && B_h \end{array}\right),
$$
where $X_i$ and $B_i$ are the matrices of the restrictions of $X$ and $B$ respectively to $\ker(f_i(S))$. Write $X_i = S_iU_i$, where $S_i$ and $U_i$ are the semisimple and unipotent parts respectively. Abbreviate $\C(B)$, $\Spec(B)$ and $\Omega(B)$ by $\C$, $\Spec$ and $\Omega$ respectively.

\begin{proposition}\label{splitomodd} Let $X \in  \O$ be as above, and assume $q$ is odd. Then $|C_\Spec(X): C_\O(X)| = c_\O$, where $c_\O \in \{1,2\}$, and $c_\O = 1$ if and only if the following conditions hold:
\begin{itemize}
\item[{\rm (i)}] for a \ged $t\pm 1$, 
the corresponding unipotent element $U_i$ 
 takes the following shape 
$($following the notation of $(\ref{vwso}))$:
\[
U_i=\bigoplus_{j=1}^r V_{\b_j}(2k_j+1) \oplus \bigoplus_{j=1}^s W(2l_j)^{[b_j]},
\]
where either $r=0$, or $r\ge 1$, the integers $k_1,\ldots,k_r$ are distinct, and the $\b_j(-1)^{k_j}$ are mutually congruent modulo $(\F_q^*)^2$;
\item[{\rm (ii)}] if both $t-1$ and $t+1$ are \geds of $X$, 
then the corresponding unipotent elements $U_i$ take the following shape:
\begin{equation}\label{u12}
U_1=\bigoplus_{j=1}^r V_{\b_j}(2k_j+1) \oplus \sum_{j=1}^s W(2l_j)^{[b_j]},\;\;
U_2=\bigoplus_{j=1}^{r'} V_{\b'_j}(2k'_j+1) \oplus \bigoplus_{j=1}^{s'} W(2l'_j)^{[b'_j]}
\end{equation}
$($as in {\rm (i)}$)$, and all of the quantities 
$\b_j(-1)^{k_j}$ and $\b'_j(-1)^{k'_j}$ are mutually congruent 
modulo $(\F_q^*)^2$;
\item[{\rm (iii)}] $X$ has no \geds $f^k$ with $f \in \Phi_2\cup \Phi_3$ and $k$ odd.
\end{itemize}
The class $X^\Spec$ splits into $2/c_\O$ conjugacy classes in $\O$.
\end{proposition}

\begin{proof}
Suppose $c_\O=1$. We show that (i)--(iii) hold. 

Consider a \ged $f_i \in \Phi_2$ of $S$. The centralizer 
$C_{\C(B_i)}(S_i)$ is isomorphic to $\GL_{m_i}(q^{d_i})$ 
(where $d_i= \deg(f_i)/2$). 
As shown in the proof of 
Proposition \ref{speco}, it intersects 
$\O(B_i)$ in the subgroup of index 2 consisting of elements 
whose determinant is a square in $\F_{q^{d_i}}^*$. 
The centralizer of the unipotent element $U_i$ in $\GL_{m_i}(q^{d_i})$ 
contains an element of non-square determinant if and only if it has 
a block of odd size $k$. As $c_\O=1$, this cannot be the case, 
so (iii) holds. (This also shows that if (iii) fails, then $c_\O=2$.) 
A similar argument applies for a \ged $f_i \in \Phi_3$.

Now consider a \ged $f_i = t\pm 1$. Here $S_i = \pm I$ and $U_i$, the 
corresponding unipotent element, is  
\[
\bigoplus_{j=1}^r (V_{\b_j}(2k_j+1) \oplus V_1(2k_j+1)^{[a_j-1]}) \oplus \bigoplus_{j=1}^s W(2l_j)^{[b_j]},
\]
in the notation of (\ref{vwso}).
As shown in the proof of \cite[Prop.\ 2.4(iii)]{GLOB}, 
$C_{\C(B_i)}(U_i)$ is contained in $\O(B_i)$ if and only if 
condition (i) holds for $U_i$. In particular, if $c_\O = 1$, then (i) holds.

Suppose that both 
$t-1$ and $t+1$ are generalized elementary divisors of $X$.
Since $c_\O = 1$, it follows that $U_1$ and $U_2$ satisfy condition (i). 
Suppose $U_1$ and $U_2$ are as in (\ref{u12}). If for some $a,\, b$ 
the quantities $\b_a(-1)^{k_a}$ and $\b'_b(-1)^{k'_b}$ are not congruent 
modulo $(\F_q^*)^2$, then 
$-I^T$, where $T = V_{2k_a+1} \oplus V_{2k_b'+1}$, centralizes $X$ and 
lies in $\C\setminus \O$ 
(see the proof of \cite[Prop.\ 2.4(iii)]{GLOB}). 
Hence, if $c_\O=1$, then condition (ii) must hold. 

We have shown that if $c_\O=1$, then conditions (i)--(iii) hold. 
The converse follows from various observations made in the proof.
\end{proof}

\begin{proposition}\label{splitomeven} Let $X \in  \O$ be as above, and assume $q$ is even. Then $|C_\C(X): C_\O(X)| = c_\O$, where $c_\O \in \{1,2\}$, 
and $c_\O = 1$ if and only if one of the following conditions holds:
\begin{itemize}
\item[{\rm (i)}] there is no \ged $f_i = t+1$;
\item[{\rm (ii)}] there is  a \ged $f_i = t+ 1$, and the corresponding 
unipotent element $U_i$ is in the class
$\bigoplus_{j=1}^s W(2l_j)^{[b_j]}$  $($following the notation of $(\ref{canon}))$.
\end{itemize}
\end{proposition}

\begin{proof}
By the proof of Proposition \ref{speco}(iii), if 
$X$ has a \ged $f_i \in \Phi_2\cup \Phi_3$, then $C_{\C(B_i)}(X_i) \le 
C_{\C(B_i)}(S_i) \le \O(B_i)$. So we only need to consider the case
where $X$ has a \ged $f_i = t+1$; now $X_i = U_i$ is a unipotent element 
as in (\ref{canon}). Moreover, 
as observed after (\ref{canon}), $C_{\C(B_i)}(X_i) \le \O(B_i)$ if 
and only 
$U_i = \bigoplus_{j=1}^s W(2l_j)^{[b_j]}$, as in part (ii). 
\end{proof}

\section{Generators for the centralizer of a general element}   \label{SectCentrGen}

Having described the centralizer of an arbitrary element of 
a classical group in Section \ref{SectListGen},
we now show how to construct a generating set for it.
We use Algorithms 1 and 2 described in Section \ref{centgens}, 
and the algorithms of Sections \ref{gencentgood} and \ref{centbad}
which solve the following problem:

\begin{itemize}
\item (\textit{Algorithm} 3)\index{Algorithm 3} Let $\C = \C(\b)$ or $\C(Q)$ be a classical isometry group. Given unipotent $X \in \C$, return a generating set for $C_{\C}(X)$. 
\end{itemize}

Let $B$ be the matrix of a non-degenerate sesquilinear form on $V$. 
Let $X \in \C(B)$, and let the \geds of $X$ be
\[
f_i^{m_{ij}} \;\;( i=1, \ldots, h, \, 1 \leq j \leq k_i),
\]
where $f_1,\ldots, f_h \in \Phi$ are distinct, and  $m_{i1} \geq m_{i2} \geq \cdots \geq m_{ik_i}$ for all $i$. 
We change basis so that 
$$
X=\begin{pmatrix}
X_1 && \\ & \ddots & \\ &&  X_h
\end{pmatrix}\,, \quad  B=\begin{pmatrix}
B_1 && \\ & \ddots & \\ && B_h
\end{pmatrix},
$$
where $X_i$ and $B_i$ are the matrices of the restrictions of $X$ and $B$ to $\ker(f_i(X)^{m_{i1}})$. Let $m_i = \sum_{j=1}^{k_i} m_{ij}$ for each $i$. 
We can suppose that $X_i$ is in Jordan form if $f_i \in \Phi_1 \cup \Phi_3$, or
\begin{eqnarray}    \label{formOfX}
X_i = \begin{pmatrix}
\widehat{X}_i & \\ & \widehat{X}_i^{*-1}
\end{pmatrix}
\end{eqnarray}
if $f_i = g_i{g_i}^* \in \Phi_2$, where $\widehat{X}_i$ is the Jordan form of the matrix of the restriction of $X$ to 
$\ker(g_i(X)^{m_{i1}})$. Let
$$
E_i = \left\{ \begin{array}{ll} F & \mbox{if } f_i \in \Phi_1 \\ F[t]/(g_i) & \mbox{if } f_i \in \Phi_2, \, f_i=g_i{g_i}^*\\ F[t]/(f_i) & \mbox{if } f_i \in \Phi_3, \mbox{symplectic and orthogonal cases} \\
\F_q[t]/(f_i\bar f_i) & \mbox{if } f_i \in \Phi_3, \mbox{unitary case.}
\end{array}\right.
$$
Finally, let 
$$
\begin{pmatrix}
X_1 && \\ & \ddots & \\ && X_h
\end{pmatrix} = \begin{pmatrix}
S_1 && \\ & \ddots & \\ && S_h
\end{pmatrix} \begin{pmatrix}
U_1 && \\ & \ddots & \\ && U_h
\end{pmatrix}
$$
be the Jordan decomposition of $X$, and
$$
\begin{pmatrix}
\widehat{X}_i & \\ & \widehat{X}_i^{*-1}
\end{pmatrix} = \begin{pmatrix}
\widehat{S}_i & \\ & \widehat{S}_i^{*-1}
\end{pmatrix}\begin{pmatrix}
\widehat{U}_i & \\ & \widehat{U}_i^{*-1}
\end{pmatrix}
$$
when $f_i \in \Phi_2$. 
A generating set for $C_{\C}(X)$ consists of the matrices
\begin{eqnarray}            \label{GenForGen}
y_{ij} = \begin{pmatrix}
\mathbb{I} && \\ & Y_{ij} & \\ && \mathbb{I}
\end{pmatrix},
\end{eqnarray}
where the $\mathbb{I}$'s are identity matrices of the appropriate sizes, and $Y_{ij}$ runs over a generating set for $C_{\C(B_i)}(X_i)$. These are obtained as follows.
\begin{itemize}
	\item $f_i \in \Phi_1$. The $Y_{ij}$ are the generators for $C_{\C(B_i)}(U_i)$ returned by Algorithm $3$.
	\item $f_i \in \Phi_2$. Let $d_i' = \deg(f_i)/2$. By Lemma \ref{L26w}, the form preserved by $X_i$ is
	$$
	B_i = \begin{pmatrix}
	\mathbb{O} & A_i \\ \varepsilon A_i^* & \mathbb{O}
	\end{pmatrix},
	$$
	where $\varepsilon= -1$ in the symplectic case and $1$ otherwise. 
Now $\widehat{U}_i$ is the embedding of a unipotent 
$\widetilde{U}_i \in \GL_{m_i}(E_i)$. We take
	\begin{eqnarray}         \label{GeneratoriPhi2}
	Y_{ij} = \begin{pmatrix}
	Z_{ij} & \\ & A_i^*Z_{ij}^{*-1}A_i^{*-1}
	\end{pmatrix},
	\end{eqnarray}
	where the $Z_{ij}$ are the embeddings into $\GL_{m_id_i'}(F)$ 
of the generators of $C_{\GL_{m_i}(E_i)}(\widetilde{U}_i)$, 
as described in Section \ref{Cue}. 
	\item $f_i \in \Phi_3$. Define $E_i$ as above. 
Let $R$ be the companion matrix $C(f_i)$ and let $\varepsilon = -1$ if $B$ is alternating, $\varepsilon=1$ otherwise. We can suppose that $S_i$ is the block diagonal sum of $m_i$ copies of $R$. Using Algorithm 1 we find $T$ such that $R^*=T^{-1}R^{-1}T$. By Claim 1 in the proof of Proposition \ref{h1cas}(iii), we can choose $T$ such that $T = \varepsilon T^*$. Let $\mathcal{T}$ be the block diagonal sum of $m_i$ copies of $T$. The matrix $H_i = B_i\mathcal{T}^{-1}$ lies in the centralizer of $S_i$, so it is the embedding into $\GL_{m_id_i}(F)$ of $\widetilde{H}_i \in \GL_{m_i}(E_i)$. 
By Claim 2 in the proof of Proposition \ref{h1cas}(iii), $\widetilde{H}_i$ is hermitian and $U_i$ is the embedding into 
$\GL_{m_id_i}(F)$ of a unipotent $\widetilde{U}_i \in \C(\widetilde{H}_i) \cong \GU_{m_i}(E_i)$. So $C_{\C(B_i)}(X_i)$ is generated by the embeddings into $\GL_{m_id_i}(F)$ of the generators of $C_{\C(\widetilde{H}_i)}(\widetilde{U}_i)$ returned by Algorithm 3.
\end{itemize}

This addresses the case where the classical group $\C = \C(B)$. 
For the case $\C=\C(Q)$, an orthogonal group in characteristic 2, note that in the analysis of the cases $f_i \in \Phi_2 \cup \Phi_3$, if $Q_i$ is a quadratic form, then it can be replaced by the associated bilinear form $\b_{Q_i}$, so we can proceed as above.

\subsection{Special and Omega groups}
To obtain generators for the centralizer in $\Spec$
of $X \in \Spec$, we apply Schreier's algorithm to $C_\C(X)$.
If $X \in \O$, then we apply the algorithm to $C_\Spec (X)$.

\section{The conjugacy problem}\label{conprobgen}

Having solved the conjugacy problem for unipotent elements 
in Sections \ref{conjtestgood} and \ref{conjprobbad}, 
and for semisimple elements in Section 
\ref{semisimpleconjclasses},
we now solve the problem in general. 

\begin{theorem} \label{ConjGen}
Let $\C$ be a classical group on $V = F^n$, where $F = \mathbb{F}_{q^u}$ \textup{(}with $u=2$  in the unitary case and $u=1$ 
otherwise\textup{)}. Let $X,Y \in \C$, let $\prod_{i=1}^h f_i^{m_i}$ be the 
minimal polynomial of $X$, with $f_1,\ldots, f_h  \in \Phi$ distinct, and  let $X_i$ be the restriction of $X$ to 
$\ker(f_i(X)^{m_i})$, with similar notation for $Y,Y_i$. 
\begin{itemize}
\item[{\rm (i)}] If $\C = \GU_n(q)$, then $X$ and $Y$ are conjugate in $\C$ if and only if $X \sim Y$.
\item[{\rm (ii)}] If $\C= \Sp_n(q)$ or $\Or^{\epsilon}_n(q)$, then $X$ and $Y$ are conjugate in $\C$ if and only if $X \sim Y$ and, for every $i$ such that $f_i(t) = t \pm 1$, the unipotent parts of $X_i$ and $Y_i$ are conjugate in the corresponding symplectic or orthogonal group $($see Sections $\ref{conjtestgood}$ and $\ref{conjprobbad})$.
\item[{\rm (iii)}] If $\Spec = \SU_n(q)$ and $X,Y \in \Spec$, then $X$ and $Y$ are conjugate in $\Spec$ if and only if $X$ and $Y$ are conjugate in $\SL(V)$: 
namely, $X$ and $Y$ are conjugate in $\C$ and every conjugating element in 
$\C$ has determinant a power of $\omega^r$, where $\o \in F$ has order $q+1$ and $r$ is the greatest common 
divisor of $q + 1$ and the dimensions of the Jordan blocks of $X$. 
\item[{\rm (iv)}] If $\Spec = \SO^{\epsilon}_n(q)$ with $q$ odd, and $X,Y \in \Spec$, then $X$ and $Y$ are conjugate in $\Spec$ if and only if $X$ and $Y$ are conjugate in $\C$, and either each has an elementary divisor $(t \pm 1)^e$ with $e$ odd, or every conjugating element in $\C$ has determinant $1$.
\item[{\rm (v)}] If $\Omega = \Omega^{\epsilon}_n(q)$, and $X,Y \in \Omega$, then $X$ and $Y$ are conjugate in $\Omega$ if and only if they are conjugate in $\Spec$, and either their class in $\Spec$ does not split into two distinct classes in $\Omega$ 
$($see Propositions $\ref{splitomodd}$ and $\ref{splitomeven})$, or every conjugating element in $\Spec$ lies in $\O$.
\end{itemize}
\end{theorem}

Note that for the criteria in parts (iii)-(v), conjugating elements in the relevant groups are constructed in the next section.

\begin{proof} (i) The left-to-right implication is obvious, so consider the converse. Let $X,Y \in \C = \GU(V)$ and suppose $X \sim Y$. Write $X=SU, \,Y=S'U'$ for the Jordan decompositions. Then $S\sim S'$, so $S$ is conjugate to $S'$ in $\C$ by Theorem \ref{MainConjinC}. Hence, replacing $Y$ by a conjugate, we can take $S=S'$, so $U$ and $U'$ are similar unipotent elements in $C_\C(S)$. By Theorem \ref{CentrSem}, $C_\C(S)$ is a direct product of general linear and general unitary groups, and so $U$ and $U'$ are conjugate in $C_\C(S)$. Thus $X$ and $Y$ are conjugate in $\C$, as required.

\vspace{2mm} \no (ii) Again, we only need to prove the right-to-left implication. Assume that $X \sim Y$ and that for $f_i = t\pm 1$, the unipotent parts of $X_i$ and $Y_i$ are conjugate in the corresponding symplectic or orthogonal group. Write $X=SU, \,Y=S'U'$ for the Jordan decompositions. Then $S$ is conjugate to $S'$ in $\C$ by Theorem \ref{MainConjinC}, so we can take $S=S'$. For $f_i = t\pm 1$ write $V_i = \ker(f_i(S))$, 
so by Theorem \ref{CentrSem}
\[
C_\C(S) = \prod_{f_i=t\pm 1} \C(V_i) \times D,
\]
where $D$ is a direct product of general linear and general unitary groups. Then $U$ and $U'$ are unipotent elements of this centralizer, and by assumption their projections to the first factor are conjugate in that factor. It follows that $U$ and $U'$ are conjugate in $C_\C(S)$, giving the conclusion.

\vspace{2mm} \no (iii) Suppose $X,Y\in \Spec$ are conjugate in $\Spec$, say $Y=X^Z$ with $Z \in \Spec$. If $Y = X^{Z'}$ with $Z' \in \C$, then $Z'Z^{-1} \in C_\C(X)$, and so $\det(Z')$ is a power of $\o^r$, as shown in the proof of Proposition \ref{specogen}(i).

For the converse, suppose that $X,Y\in \Spec$ and $Y=X^Z$, where $Z \in \C=\GU(V)$ has determinant a power of $\o^r$. As before, write $X=SU$ and $Y=S'U'$ for the Jordan decompositions. Then $S$ is conjugate to $S'$ in $\Spec$ by Proposition \ref{speco}, so we can take $S=S'$. Now the proof of Proposition \ref{specogen}(i) shows that $U$ and $U'$ are conjugate in $C_\Spec(S)$, as required.

\vspace{2mm} \no (iv) Let $X,Y \in \Spec$ be conjugate in $\C$. In the notation of Proposition \ref{specogen}(ii), if $c_S=2$, then $X^\C = X^\Spec$ and so $X$ and $Y$ are conjugate in $\Spec$; and $c_S=1$ if and only if $X$ has an elementary divisor $(t\pm 1)^e$ with $e$ odd. Part (iv) follows.

\vspace{2mm} \no (v) Let $X,Y \in \O$ be conjugate in $\Spec$. We use the notation of Propositions \ref{splitomodd} and 
\ref{splitomeven}. If $c_\O = 2$, then $X^\O = X^\Spec$ and so $X$ and $Y$ are conjugate in $\O$; if $c_\O=1$, then $C_\O(X) = C_\Spec (X)$, and so every conjugating element in $\Spec$ lies in $\O$. This completes the proof. 
\end{proof}

\section{Constructing a conjugating element}      \label{SectConjGen}
Let $\C = \C(B)$ or $\C(Q)$ be a classical isometry group. 
Given conjugate $X,Y \in \C$, we describe an algorithm to compute  
$Z \in \C$ such that $Y=X^Z$. 
We use Algorithms 1 and 2 described in Section \ref{centgens}, 
together the algorithms of Sections \ref{goodconjel} and \ref{badconj} 
which solve the following problem.

\begin{itemize}
\item (\textit{Algorithm} 4)\index{Algorithm 4} Given unipotent $X,Y \in \C$ that are conjugate in $\C$, return $Z \in \C$ such that $Y = X^Z$.
\end{itemize}

The algorithm is similar to the semisimple case, 
described in Section \ref{ConjElSemi}. Consider first $\C = \C(B)$.
The procedure begins as in Step 1 of Section \ref{ConjElSemi}, 
by computing matrices $P_X$ and $P_Y$ such that 
$P_XXP_X^{-1} = P_YYP_Y^{-1} =J$, where $J$ is as in (\ref{Jstdform})
and the $J_i$ are the restrictions of $J$ to the generalized eigenspaces. Then $J$ preserves the forms $B_X = P_XBP_X^*$ and $B_Y=P_YBP_Y^*$, where
$$
J = \begin{pmatrix}
J_1 && \\ & \ddots & \\ && J_h
\end{pmatrix}, \quad B_X = \begin{pmatrix}
B_{X,1} && \\ & \ddots & \\ && B_{X,h}
\end{pmatrix}, \quad B_Y = \begin{pmatrix}
B_{Y,1} && \\ & \ddots & \\ && B_{Y,h}
\end{pmatrix}
$$
and each $J_i$ is an isometry for $B_{X,i}$ and $B_{Y,i}$.
As in Section \ref{ConjElSemi}, Step 2 of the algorithm is to compute $W_i$ in the centralizer of $J_i$ such that $W_iB_{Y,i}W_i^* = B_{X,i}$ for each $i$. Taking $W = \bigoplus_{i=1}^h W_i$ and $Z = P_X^{-1}WP_Y$, we deduce that $Z \in \C(B)$ and $X^Z = Y$, as required.

To compute the $W_i$, we distinguish the three cases.
\begin{itemize}
\item $f_i \in \Phi_1$. 
Now $J_i$ is a product of a scalar and a unipotent element, and we can suppose that the scalar is the identity (so $J_i$ is unipotent) because it does not affect the computation. Using Algorithm 2, we compute $W_{B,i}$ such that $W_{B,i}B_{Y,i}W_{B,i}^* = B_{X,i}$. Now $J_i$ and $W_{B,i}J_iW_{B,i}^{-1}$ are unipotent elements of $\C(B_{X,i})$, so using Algorithm 4 we can compute ${W_{J,i} \in \C(B_{X,i})}$ such that $W_{J,i}W_{B,i}J_iW_{B,i}^{-1}W_{J,i}^{-1}=J_i$. Now take $W_i = W_{J,i}W_{B,i}$.
\item $f_i \in \Phi_2$. This case is identical to the $\Phi_2$ case in Step 2 of Section \ref{ConjElSemi}, where we never used the fact that $J_i$ was semisimple.
\item $f_i \in \Phi_3$. Following the $\Phi_3$ case in Step 2 of the semisimple case in Section \ref{ConjElSemi}, we work in $\GL_{m_i}(E)$, with $E = F[t]/(f_i)$ (or $\F_q[t]/(f_i\bar f_i)$ in the unitary case).
We compute $\widetilde{H}_{X,i}$ and 
$\widetilde{H}_{Y,i}$, matrices of hermitian forms preserved by 
$\widetilde{J}_i \in \GL_{m_i}(E)$ where $J_i$ is the embedding of $\widetilde{J}_i$ into $\GL_{m_id_i}(F)$. As in the case $f_i \in \Phi_1$, we can suppose that $\widetilde{J}_i$ is unipotent, since scalar factors do not affect the computation. The construction of $\widetilde{W}_i \in C_{\GL_{m_i}(E)}(\widetilde{J}_i)$ such that $\widetilde{W}_i\widetilde{H}_{Y,i}\widetilde{W}_i^{\dag} = \widetilde{H}_{X,i}$ is identical to that described above in the case $f_i \in \Phi_1$.
\end{itemize}

The procedure also works for orthogonal groups $\C = \C(Q)$ in 
characteristic 2 for the 
reasons stated just before Section \ref{spogps}.

\subsection{Special and Omega groups}
Here we assume that $X$ and $Y$ are conjugate in $\Spec$ or $\O$.
As above, we can compute a conjugating element $Z=P_X^{-1}WP_Y$ 
in the isometry group $\C$. If  
$Z$ has inappropriate determinant or spinor norm, then we adjust $Z$ to obtain a conjugating element in $\Spec$ or $\O$. 
As in the semisimple case, we aim to compute $D \in C_{\C(B_X)}(J)$ having the appropriate determinant or spinor norm, and then replace $Z$ by $P_X^{-1}DWP_Y$.

For $i=1, \dots, h$, let $m_{ij}$ be as defined at the beginning of Section \ref{SectCentrGen}.

Suppose first that $\Spec$ is the special orthogonal group. 
If $X$ and $Y$ are conjugate in $\Spec$, then by the proof of Proposition \ref{specogen}(ii) we may suppose 
that there is an elementary divisor $f_1^m = (t\pm 1)^m$ with $m$ odd, 
so we can choose 
\begin{eqnarray}  \label{formadiD}
D = 
\begin{pmatrix}
D_1 & \\ & \mathbb{I}
\end{pmatrix},
\end{eqnarray}
where $D_1 \in C_{\C(B_{X,1})}(J_1)$ has the same determinant as $Z$. 

Now let $\Spec$ be the special unitary group. Set
$$
r_i = \gcd({q+1}, m_{i1}, \dots, m_{ik_i}),
$$
and let $r = \gcd(q+1, r_1, \dots, r_h)$. Let $\o \in \F_{q^2}$ have order $q+1$.
From the proof of Proposition \ref{specogen}(i), we can find $H_i \in C_{\C(B_{X,i})}(J_i)$ such that $\det (H_i) = \o^{r_i}$.
Compute integers $a_i$ such that $\sum_{i=1}^h a_ir_i=r$ in $\Z_{q+1}$, and let
$$
H = \begin{pmatrix}
H_1^{a_1} && \\ & \ddots & \\ && H_h^{a_h}
\end{pmatrix}.
$$
Then $\det (H) = \o^r$, and 
since $X$ and $Y$ are conjugate in $\Spec$, there exists an 
integer $\ell$ such that $\det (Z) = \det(H^{\ell})$.
So we take $D = H^{-\ell}$.

Finally, we consider the case where $X,Y \in \Omega$ are $\O$-conjugate, where $\O$ is an orthogonal group. 
In even characteristic, in the notation of Proposition \ref{splitomeven} we can suppose that $c_\O=2$, so there exists 
$f_1 \in \Phi_1$ and 
$D_1 \in C_{\C(B_{X,1})}(J_1)\setminus \O(B_{X_1})$, and we choose
$$
D = \begin{pmatrix}
D_1 & \\ & \mathbb{I}
\end{pmatrix}.
$$
Now suppose we are in odd characteristic, and that $X$ and $Y$ are 
conjugate in $\Omega$. From the special group case handled above, we can compute a conjugating element $Z$ in $\Spec$.
Assume that ${C_{\Spec}(X) \neq C_{\Omega}(X)}$, and that $Z \in \Spec \setminus \O$.
Then, in the notation of Proposition \ref{splitomodd}, $c_\O=2$, and one of the following applies.
\begin{itemize}
\item There exists $f_i \in \Phi_2 \cup \Phi_3$ and $m_{ij}$ odd. 
Now $C_{\C(B_{X,i})}(J_i)$ contains elements of non-square spinor norm. 
Hence, we take 
$$
D = \begin{pmatrix}
\mathbb{I} && \\ & D_i & \\ && \mathbb{I}
\end{pmatrix},
$$
where $D_i \in C_{\C(B_{X,i})}(J_i)$ has the same spinor norm as 
$P_X^{-1}WP_Y$. 

\item There exists $f_i \in \Phi_1$ with 
$C_{\Spec(B_{X,i})}(J_i) \not\subseteq \Omega(B_{X,i})$. 
Proposition \ref{splitomodd} shows that 
 we can choose $D_i \in C_{\Spec(B_{X,i})}(J_i)$ with the same spinor norm as $Z$, so we take
$$
D = \begin{pmatrix}
\mathbb{I} && \\ & D_i & \\ && \mathbb{I}
\end{pmatrix}.
$$
\item 
There exist $f_1,f_2 \in \Phi_1$ such that  
$C_{\Spec(B_{X,i})}(J_i) \subseteq \Omega(B_{X,i})$ and ${C_{\C(B_{X,i})}(J_i)\not\subseteq \Spec(B_{X,i})}$ for $i=1,2$. 
Here we take
$$
D = \begin{pmatrix}
D_1 && \\ & D_2 & \\ && \mathbb{I}
\end{pmatrix}
$$
with $D_i \in C_{\C(B_{X,i})}(J_i) \setminus C_{\Spec(B_{X,i})}(J_i)$ for $i=1,2$. 
\end{itemize}

\section{Number of conjugacy classes in isometry groups}
We found it useful to check that the number of conjugacy class representatives 
for isometry groups agrees with the results of Macdonald \cite{mac81}
and Wall \cite[\S 2.6 and Thm.\ 3.7.3]{GEW}.

\begin{thm} \textit{} 
\begin{itemize}
\item The number of conjugacy classes of $\GL_n(q)$ is the coefficient of $t^n$ in
the formal power series
$$
\prod_{i=1}^{\infty} \frac{1-t^i}{1-qt^i}.
$$
\item The number of conjugacy classes of $\GU_n(q)$ is the coefficient of $t^n$ in
$$
\prod_{i=1}^{\infty} \frac{1+t^i}{1-qt^i}.
$$
\item The number of conjugacy classes of $\Sp_n(q)$ with $q$ odd is the coefficient of $t^n$ in 
$$
\prod_{i=1}^{\infty}\frac{(1+t^{2i})^4}{1-qt^{2i}}.
$$
\item Let $q$ be odd. Let $k_n^+$ and $k_n^-$ denote the numbers of conjugacy classes of $\Or^+_n(q)$ and $\Or^-_n(q)$ respectively \textup{(}with $\Or^+_n(q)=\Or^-_n(q) = \Or_n(q)$ if $n$ is odd\textup{)}. Then
\begin{eqnarray*}
\sum_{n=0}^{\infty}(k_n^++k_n^-)t^n & = & \prod_{i=1}^{\infty} \frac{(1+t^{2i-1})^4}{1-qt^{2i}} \\
\sum_{n=0}^{\infty}(k_n^+-k_n^-)t^n & = & \prod_{i=1}^{\infty} \frac{1-t^{4i-2}}{1-qt^{4i}}.
\end{eqnarray*}
\item Let $q$ be even. Define a sequence of polynomials $\chi_i = \chi_i(a,b,t)$ and a power series $\chi=\chi(a,b,t)$ as follows:
\begin{eqnarray*}
\chi_{-1} & = & a,\\
\chi_0 & = & b,\\
\chi_{2m+1}-\chi_{2m} & = & t^{2m+1}\chi_{2m-1},\\
\chi_{2m+2}-\chi_{2m+1} & = & t^{m+1}(1+t^{m+1})(\chi_{2m+1}+(1-t^{2m+1})\chi_{2m-1}),\\
\chi(a,b,t) & \equiv & \chi_{2m}\bmod{t^m} 
\end{eqnarray*}
for all $m\ge 0$. The number of conjugacy classes of $\Sp_{2m}(q)$ \textup{(}resp.\ $\Or^{\epsilon}_{2m}(q)$\textup{)} is the coefficient of $t^{2m}$ in the power series $s(t^2)$ \textup{(}resp.\ $\omega_{\epsilon}(t^2)$\textup{)}, where
\begin{eqnarray*}
s(t^2) & = & \chi(0,1,t^2) \prod_{i=1}^{\infty} (1-qt^{2i})^{-1},\\
\omega_+(t^2)+\omega_-(t^2) & = & \chi(1,1,t^2) \prod_{i=1}^{\infty} (1-qt^{2i})^{-1},\\
\omega_+(t^2)-\omega_-(t^2) & = & \prod_{i=1}^{\infty} \frac{1-t^{4i-2}}{1-qt^{4i}}.
\end{eqnarray*}
\end{itemize}
\end{thm}
In Table \ref{polys} we record the resultant polynomials in $q$ 
for the isometry groups of specified type and degree $n$ at most 10. 

\begin{longtable}{|l|l|l|}
\caption{Polynomials for number of classes in isometry groups} \label{polys} \\
\hline 
Type  &  $n$  &  Polynomial \\ \hline 
GL  &  2  &  $q^2 - 1 $ \\ \hline
GL  &  3  &  $q^3 - q $ \\ \hline
GL  &  4  &  $q^4 - q $ \\ \hline
GL  &  5  &  $q^5 - q^2 - q + 1 $ \\ \hline
GL  &  6  &  $q^6 - q^2 $ \\ \hline
GL  &  7  &  $q^7 - q^3 - q^2 + 1 $ \\ \hline
GL  &  8  &  $q^8 - q^3 - q^2 + q $ \\ \hline
GL  &  9  &  $q^9 - q^4 - q^3 + q $ \\ \hline
GL  &  10  &  $q^{10} - q^4 - q^3 + q $ \\ \hline
GU  &  2  &  $q^2 + 2q + 1 $ \\ \hline
GU  &  3  &  $q^3 + 2q^2 + 3q + 2 $ \\ \hline
GU  &  4  &  $q^4 + 2q^3 + 4q^2 + 5q + 2 $ \\ \hline
GU  &  5  &  $q^5 + 2q^4 + 4q^3 + 7q^2 + 7q + 3 $ \\ \hline
GU  &  6  &  $q^6 + 2q^5 + 4q^4 + 8q^3 + 11q^2 + 10q + 4 $ \\ \hline
GU  &  7  &  $q^7 + 2q^6 + 4q^5 + 8q^4 + 13q^3 + 17q^2 + 14q + 5 $ \\ \hline
GU  &  8  &  $q^8 + 2q^7 + 4q^6 + 8q^5 + 14q^4 + 21q^3 + 25q^2 + 19q + 6 $ \\ \hline
GU  &  9  &  $q^9 + 2q^8 + 4q^7 + 8q^6 + 14q^5 + 23q^4 + 33q^3 + 
36q^2 +  25q + 8 $ \\ \hline
GU  &  10  &  $q^{10} + 2q^9 + 4q^8 + 8q^7 + 14q^6 + 24q^5 + 37q^4 
+ 49q^3 +  50q^2 + 33q + 10 $ \\ \hline
Sp ($q$ even)  &  2  &  $q + 1 $ \\ \hline
Sp ($q$ even)  &  4  &  $q^2 + 2q + 3 $ \\ \hline
Sp ($q$ even)  &  6  &  $q^3 + 2q^2 + 5q + 4 $ \\ \hline
Sp ($q$ even)  &  8  &  $q^4 + 2q^3 + 6q^2 + 9q + 7 $ \\ \hline
Sp ($q$ even)  &  10  &  $q^5 + 2q^4 + 6q^3 + 11q^2 + 16q + 10 $ \\ \hline
Sp ($q$ odd)  &  2  &  $q + 4 $ \\ \hline
Sp ($q$ odd)  &  4  &  $q^2 + 5q + 10 $ \\ \hline
Sp ($q$ odd)  &  6  &  $q^3 + 5q^2 + 15q + 24 $ \\ \hline
Sp ($q$ odd)  &  8  &  $q^4 + 5q^3 + 16q^2 + 39q + 51 $ \\ \hline
Sp ($q$ odd)  &  10  &  $q^5 + 5q^4 + 16q^3 + 44q^2 + 90q + 100 $ \\ \hline
GO  &  3  &  $2q + 4 $ \\ \hline
GO  &  5  &  $2q^2 + 6q + 14 $ \\ \hline
GO  &  7  &  $2q^3 + 6q^2 + 20q + 28 $ \\ \hline
GO  &  9  &  $2q^4 + 6q^3 + 22q^2 + 48q + 62 $ \\ \hline
$\GO^+$ ($q$ even)  &  2  &  $1/2q + 1 $ \\ \hline
$\GO^+$ ($q$ even)  &  4  &  $1/2q^2 + 5/2q + 2 $ \\ \hline
$\GO^+$ ($q$ even)  &  6  &  $1/2q^3 + 2q^2 + 7/2q + 3 $ \\ \hline
$\GO^+$ ($q$ even)  &  8  &  $1/2q^4 + 2q^3 + 5q^2 + 8q + 7 $ \\ \hline
$\GO^+$ ($q$ even)  &  10  &  $1/2q^5 + 2q^4 + 9/2q^3 + 9q^2 + 13q + 9 $ \\ \hline
$\GO^+$ ($q$ odd)  &  2  &  $1/2q + 5/2 $ \\ \hline
$\GO^+$ ($q$ odd)  &  4  &  $1/2q^2 + 4q + 17/2 $ \\ \hline
$\GO^+$ ($q$ odd)  &  6  &  $1/2q^3 + 7/2q^2 + 23/2q + 37/2 $ \\ \hline
$\GO^+$ ($q$ odd)  &  8  &  $1/2q^4 + 7/2q^3 + 13q^2 + 63/2q + 85/2 $ \\ \hline
$\GO^+$ ($q$ odd)  &  10  &  $1/2q^5 + 7/2q^4 + 25/2q^3 + 34q^2 + 72q + 171/2 $ 
\\ \hline
$\GO^-$ ($q$ even)  &  2  &  $1/2q + 2 $ \\ \hline
$\GO^-$ ($q$ even)  &  4  &  $1/2q^2 + 3/2q + 2 $ \\ \hline
$\GO^-$ ($q$ even)  &  6  &  $1/2q^3 + 2q^2 + 9/2q + 4 $ \\ \hline
$\GO^-$ ($q$ even)  &  8  &  $1/2q^4 + 2q^3 + 4q^2 + 7q + 6 $ \\ \hline
$\GO^-$ ($q$ even)  &  10  & $1/2q^5 + 2q^4 + 9/2q^3 + 10q^2 + 15q + 10 $ \\ \hline
$\GO^-$ ($q$ odd)  &  2  &  $1/2q + 7/2 $ \\ \hline
$\GO^-$ ($q$ odd)  &  4  &  $1/2q^2 + 3q + 17/2 $ \\ \hline
$\GO^-$ ($q$ odd)  &  6  &  $1/2q^3 + 7/2q^2 + 25/2q + 39/2 $ \\ \hline
$\GO^-$ ($q$ odd)  &  8  &  $1/2q^4 + 7/2q^3 + 12q^2 + 61/2q + 83/2 $ \\ \hline
$\GO^-$ ($q$ odd)  &  10  &  $1/2q^5 + 7/2q^4 + 25/2q^3 + 35q^2 + 74q + 173/2 $ \\ \hline 
\end{longtable}

\newpage 
\phantomsection
\addcontentsline{toc}{chapter}{Index}  \label{notation-index}
\printindex


\begin{thebibliography}{99}
\addcontentsline{toc}{chapter}{Bibliography}

\bibitem{Adams} William W. Adams \& Philippe Loustaunau,
{\it An Introduction to Gr\"obner Bases}.  Grad.\ Stud.\ Math. {\bf 3}, 
American Math. Soc., 1994.

\bibitem{BHLO} Henrik B\"a\"arnhielm, Derek F.\ Holt, 
C.\ R.\ Leedham-Green \& E.\ A.\ O'Brien,
A practical model for computation with matrix groups.
{\it J.\ Symbolic Comput.} {\bf 68} (2015), 27--60.

\bibitem{MAGMA} Wieb Bosma, John Cannon \& Catherine Playoust, 
The {M}agma algebra system. {I}. The user language. \emph{J.\ Symbolic Comput.} {\bf 24} (1997), 235--265.

\bibitem{JBRIT} John R.\ Britnell, Cycle index methods for matrix groups over finite fields. DPhil Thesis, University of Oxford (2003).

\bibitem{PGPG} Timothy C.\ Burness \& Michael Giudici, 
{\it Classical groups, {d}erangements and {p}rimes}. Austral.\  Math.\ Soc.\ Lect.\ 
Ser.\ {\bf 25}, Cambridge University Press, Cambridge, 2016.

\bibitem{CannonHolt} John J.\ Cannon \& Derek F.\ Holt, Computing conjugacy class representatives in permutation groups. 
\emph{J.\ Algebra} {\bf 300} (2006), 213--222.

\bibitem{Carterclass} R.W. Carter, Centralizers of semisimple elements in the finite classical groups. 
\emph{Proc. London Math. Soc.} {\bf 42} (1981), 1--41.

\bibitem{CP} Matthew C.\ Clarke \& Alexander Premet, The Hesselink stratification of nullcones and base
change. {\it Invent.\ {M}ath.} {\bf 191} (2013), 631--669.


\bibitem{thesis} Giovanni De Franceschi, Centralizers and conjugacy classes in finite classical groups. PhD Thesis, University of Auckland (2018).
\href{https://researchspace.auckland.ac.nz/handle/2292/45197}
{researchspace.auckland.ac.nz/handle/2292/45197}

\bibitem{github}
Giovanni De Franceschi, Martin W. Liebeck \& E.A. O'Brien,
Conjugacy classes in finite classical groups.
\href{https://github.com/eamonnaobrien/ClassicalConjugacy}
{github.com/eamonnaobrien/ClassicalConjugacy}

\bibitem{black} Heiko Dietrich, C.\ R.\ Leedham-Green \& E.\ A.\ O'Brien. Effective black-box constructive recognition of classical groups. {\it J.\ Algebra} {\bf 421} (2015), 460--492.

\bibitem{FS} Paul Fong \& Bhama Srinivasan, The blocks of finite classical groups. \emph{J.\ {R}eine {A}ngew.\ {M}ath.} {\bf 396} (1989), 122--191.

\bibitem{vzg}
Joachim von zur Gathen \& J\"urgen Gerhard,
{\it Modern Computer Algebra}, Cambridge University Press, 2003.

\bibitem{GLOB} Samuel Gonshaw, Martin W.\ Liebeck \& E.\ A.\ O'Brien, 
Unipotent class representatives for finite classical groups. \emph{J.\ Group Theory} {\bf 20} (2017), 505--525.

\bibitem{CGSA} Larry C.\ Grove, {\it Classical Groups and Geometric Algebra}. 
Grad.\ Stud.\ Math. {\bf 39}, American Math. Soc., Providence, RI, 2002.

\bibitem{PFGG} I.\ N.\ Herstein, {\it Topics in Algebra, Second edition}. Xerox College Publishing, Lexington, Mass. - Toronto, Ont., 1975.

\bibitem{handbook} Derek F.\ Holt, Bettina Eick \& Eamonn A.\ O'Brien, {\it Handbook of Computational Group Theory}. Discrete mathematics and its applications, Chapman \& Hall / CRC Press, 2005.

\bibitem{Hulpke2000} Alexander Hulpke, Conjugacy classes in finite permutation groups via homomorphic images. 
\emph{Math.\ Comp.\ } {\bf 69} (2000), 1633--1651. 

\bibitem{Hulpke2013} Alexander Hulpke, Computing conjugacy classes of elements in matrix groups. \emph{J.\ Algebra} {\bf 387} (2013), 268--286. 

\bibitem{Huppert} Bertram Huppert, Isometrien von Vektorr\"{a}umen 1. \emph{Arch.\ Math.\ (Basel)} {\bf 35} (1980), 164--176.

\bibitem{KL} Peter Kleidman \& Martin Liebeck, {\it The Subgroup Structure of the Finite Classical Groups}. London Math.\ Soc.\ Lecture Note Ser. {\bf 129}, Cambridge University Press, Cambridge, 1990.

\bibitem{LS} Martin W.\ Liebeck \& Gary M.\ Seitz, {\it Unipotent and Nilpotent Classes in Simple Algebraic Groups and Lie Algebras}. Math.\ Surveys Monogr. {\bf 180}, American Math. Soc., Providence, RI, 2012.


\bibitem{mac81}
I.\ G.\ Macdonald, Numbers of conjugacy classes in some finite classical
groups, {\it Bull.\ Austral.\ Math.\ Soc.} {\bf 23} (1981), 23--48.

\bibitem{MD} I.\ G.\ Macdonald, {\it Symmetric functions and Hall Polynomials, Second Edition}. Oxford Math.\ Monogr., Clarendon Press, Oxford University Press, New York, 1995.

\bibitem{Milnor} John Milnor, On Isometries of Inner Product Spaces. \emph{Invent.\ Math.} {\bf 8} (1969), 83--97.



\bibitem{Murray}
Scott H.\ Murray, 
Conjugacy classes in maximal parabolic subgroups of general linear groups.
\emph{J.\ Algebra} {\bf 233} (2000), 135--155.

\bibitem{NP}
Max Neunh\"offer \& Cheryl E.\ Praeger, 
Computing minimal polynomials of matrices.
{\it LMS J. Comput. Math.} {\bf 11} (2008), 252--279.

\bibitem{RT} L.\ J.\ Rylands \& D.\ E.\ Taylor, Matrix generators for the orthogonal groups. \emph{J.\ Symbolic Computation} {\bf 25} (1998), 351--360.

\bibitem{JF} Allan Steel, A new algorithm for the computation of canonical forms of matrices over fields. 
\emph{J.\ Symbolic Comput.} {\bf 24} (1997), 409--432. 

\bibitem{PG} Donald E.\ Taylor, {\it The Geometry of the Classical Groups}. 
Sigma Ser.\ Pure Math. {\bf 9}, Heldermann Verlag, Berlin, 1992.

\bibitem{DET} Donald E.\ Taylor, Pairs of generators for matrix groups, I.
\href{https://arXiv.org/abs/2201.09155}{arXiv:2201.09155}, 2022.

\bibitem{GEW} G.\ E.\ Wall, On the conjugacy classes in the unitary, symplectic and orthogonal groups. \emph{J.\ Austral.\ Math.\ Soc.} {\bf 3} (1963), 1--62.

\bibitem{SCCSS} G.\ E.\ Wall, The semisimple conjugacy classes in the 
symplectic groups, \href{https://arxiv.org/abs/1512.04520}{arXiv:1512.04520} (2015).

\bibitem{GSM} James B.\ Wilson, Optimal algorithms of Gram-Schmidt type. \emph{Linear Algebra Appl.} {\bf 438} (2013), 4573--4583.



\end{thebibliography}
\end{document}